\documentclass[10pt]{report}
\usepackage{ubcmath-thesis}
\usepackage{amsmath,amssymb,amsfonts,amsthm,euscript,amsbsy,subfigure,psfrag}
\usepackage[dvips]{graphicx}
\usepackage{epic}

\newtheorem{THM}{Theorem}[section]
\newtheorem{COR}[THM]{Corollary}

\newtheorem{REM}[THM]{Remark}
\newtheorem{LEM}[THM]{Lemma}
\newtheorem{PRP}[THM]{Proposition}

\newcommand{\R}{{\cal R}}

\newcommand{\bestar} {\begin{eqnarray*}}
\newcommand{\eestar} {\end{eqnarray*}}
\newcommand{\be}   {\begin{eqnarray}}
\newcommand{\ee}   {\end{eqnarray}}
\newcommand{\nnb}       {\nonumber}
\newcommand{\Rbold} {{\mathbb R}}
\newcommand{\Zbold} {{\mathbb Z}}
\newcommand{\Nbold} {{\mathbb N}}

\newcommand{\Dcal} {{\cal D}}
\newcommand{\Fcal} {{\cal F}}
\newcommand{\Gcal} {{\cal G}}
\newcommand{\Kcal} {{\cal K}}
\newcommand{\Lcal} {{\cal L}}

\newcommand{\Ncal} {{\cal N}}
\newcommand{\Pcal} {{\cal P}}
\newcommand{\var}{\mathrm{Var}}
\newcommand{\cov}{\mathrm{Cov}}

\phdthesis

  \rightchapter

\singlespace
\thesistitle{
	Equilibrium States of Two Stochastic Models in Mathematical Ecology
}

\thesisauthor{Feng Yu}
\thesisauthorpreviousdegrees{}
\thesismonth{January} \thesisyear{2005}

\begin{document}

%%%%%%%%%%%%%%%%%%%%%%%%%%%%%%%%%%%%%%%%%%%%%%%%%%%%%%%%%%%%%%%%%%%%%%%%
%%%
%%% The preliminary pages of your thesis are defined here.  Some
%%% sections are optional, and many are automatically generated from
%%% data previously set.

%% The title page is mandatory, and is automatically generated here from
%% the data set in the preamble.
\thesistitlepage

%% An abstract is mandatory, and should be placed within the
%% thesisabstract environment.
\begin{thesisabstract}
This work deals with two problems arising in mathematical ecology.
The first problem is concerned with diploid branching particle models
and its behavior when rapid stirring is added to the interaction.
The particle models involve two types of particles, male and female,
and branching can only occur when both types of particles are present.
We show that if the branching rate is sufficiently large,
this particle model has a nontrivial stationary distribution,
i.e. one that does not concentrate all weight on the all-$0$ state, using
a comparison argument due to R. Durrett.
We also show extinction for small branching rates, thereby establishing
the existence of a phase transition.
We then add two different rapid stirring mechanisms to the interactions
and show that for the particle models with rapid stirring, there also
exist nontrivial stationary distribution(s); for this, we analyze the
limiting PDE and establish a condition on the PDE that guarantees existence
of nontrivial stationary distributions for sufficient fast stirring.

The second problem deals with a model of sympatric speciation,
i.e. speciation in the absence of geographical separation, originally
proposed by U. Dieckmann and M. Doebeli in 1999. We modify their original
model to obtain several constant-population particle models. We concentrate
on a continuous-time model that converges to a deterministic dynamical system
as the number of particles becomes large. We establish various results
regarding whether speciation occurs by studying the 
existence of bimodal stationary distributions for the limiting
dynamical system.

\end{thesisabstract}

%% A table of contents is mandatory.  If you require a list of tables
%% and/or a list of figures, uncomment the appropriate line(s).  All of
%% these items are automatically generated.
  \tableofcontents
% \listoftables
% \listoffigures

%% Acknowledgements and a dedication are optional.  If you desire to
%% include either or both of these, uncomment the appropriate lines.
%% Acknowledgements are not allowed to take more than one page.

\begin{thesisacknowledgments}
I would like to thank the external examiner,
Robert Adler, and the university examiners, Michael Doebeli and Jim Zidek,
for their feedback, suggestions, and generous comments.
I am also grateful for the help and guidance provided by members
of my supervisory committee, Martin Barlow and John Walsh,
and my research supervisor, Ed Perkins, through all my years at UBC.
Without them, this would never have been possible.
%Furthermore, I would also like to thank my officemates for their jokes,
%mostly either silly or completely inappropriate and sexist. And finally,
%I would like to thank all the students I have taught for their unparalleled
%ability to boost my self-confidence.
\end{thesisacknowledgments}

% \begin{thesisdedication}
% To Herbie the cat.
% \end{thesisdedication}

%%%%%%%%%%%%%%%%%%%%%%%%%%%%%%%%%%%%%%%%%%%%%%%%%%%%%%%%%%%%%%%%%%%%%%%%
%%%
%%% The main body of your thesis goes here.  Each chapter should start
%%% with a \chapter{Chapter Title} tag.  Probably you will want to place
%%% these in separate files and just \input them here.

\chapter{Introduction and Overview}
This work consists of two parts, each of which involves a class of models
arising from a problem
of mathematical ecology. In the first problem, we study diploid branching
particle system models. This class of particle systems differs from the
``usual'' models that one normally finds in the literature, in that
there are two types of particles, modelling the male and female populations,
and branching (i.e. birth of new particles) requires the presence of
both male and female particles. In the second problem, we study various
particle models that are all related to a sympatric speciation model
proposed in [Dieckmann and Doebeli 1999]. In both particle models, we
are mainly concerned with the equilibrium behaviour. More specifically,
we show that
the stationary distributions of the particle models have desirable properties,
e.g. nontriviality (i.e. does not concentrate all weight on the all-$0$
state) in the first problem and bimodality in the second problem.

\part{Existence of Nontrivial Stationary Distribution for the Diploid
        Branching Particle System with Rapid Stirring}
\chapter{The Particle Models}
In this part of our work, we consider a type of particle systems
that can be used to model sexual reproduction of a certain species.
This work was inspired in part by [Dawson and Perkins 1998]. In that
paper, the following system of stochastic partial differential equations
is the object of study:
\begin{eqnarray}
	\frac{\partial u}{\partial t} (t,x)
		& = & \frac{1}{2} \Delta u (t,x)
			+ (\gamma u(t,x) v(t,x))^{1/2} \dot W_1 (t,x) \nnb \\
	\frac{\partial v}{\partial t} (t,x)
                & = & \frac{1}{2} \Delta v (t,x)
                        + (\gamma u(t,x) v(t,x))^{1/2} \dot W_2 (t,x),
		\label{eq:spde:dp}
\end{eqnarray}
where $\Delta = \sum_i \frac{\partial^2}{\partial x_i^2}$ is the Laplacian,
$\gamma > 0$, and $\dot W_i(t,x)$ ($i=1,2$) are independent
space-time white noises on $\Rbold_+ \times \Rbold$. One can
associate $u(x,t)$ and $v(x,t)$ with the male and female populations
of ``particles'' (respectively) at spatial location $x$ and time $t$.
Loosely speaking,~(\ref{eq:spde:dp}) says that individual male or
female particles moves around according to Brownian motion, but branching
is only possible when both male and female particles are present at the
same spatial location. Notice that at spatial locations where
the female population is $0$, the branching rate for the male population
is also $0$, therefore the male population does not ``die'' and the
only effect on the male population at those spatial locations is the
diffusive effect of the heat kernel $\Delta$.
%with the initial condition
%$u(0,x)=u_0, v(0,x)=0$, the branching rate is always zero and $u(t,x)$ is
%deterministic; in particular, $\int_x u(t,x) dx$ remains constant, i.e.
%the total population of male particles remains constant in the absence of
%female particles.
This behaviour is not very realistic, since one would expect a ``natural''
death rate of male particles even without the presence of any female particles.
%This is in part caused by the model involving an infinite number of particles.
In this work, we study a model involving a finite number of male and
female particles with more (somewhat) ``realistic'' behaviour.

The model we study involves
two types of particles, male and female, residing on the integer grid
$S = \Zbold^d$ or $\epsilon \Zbold^d$. More specifically, each \emph{site}
$x\in S$ contains two \emph{nests}, one for the male particle and the other
for the female particle. Each nest can be inhabited by at most 1 particle,
either male or female. Let $E=\{0,1\}$ and $F=E \times E$ be the set
of possible states at each site in $S$. For $x \in S$, we write
\begin{eqnarray*}
        \xi(x) = (\xi^1(x),\xi^2(x)),
\end{eqnarray*}
where $\xi^1(x)$ denotes the number ($0$ or $1$) of male particles
at site $x$, and $\xi^2(x)$ denotes the number of female particles
at site $x$.
We define the interaction neighbourhood
\begin{eqnarray*}
	\Ncal = \{0,y_1,\ldots,y_N\},
\end{eqnarray*}
and the neighbourhood of $x$
\begin{eqnarray*}
        \Ncal_x = x + \Ncal.
\end{eqnarray*}
For example, $\Ncal = \{0,-1,1\}$ if we have nearest-neighbour
interaction on $\Zbold$.
Let $c_i(x,m,\xi)$ denote the rate at which nest $m$ ($m=1,2$) of site $x$
flips to state $i$ ($i=0,1$), and assume $c_i(x,m,\xi)$
depends only on the neighbourhood $\Ncal_x$, i.e.
\begin{eqnarray*}
	c_i(x,m,\xi) = h_{i,m}(\xi(x),\xi(x+y_1),\ldots,\xi(x+y_N))
\end{eqnarray*}
for some function $h_{i,m}:F^{N+1}\rightarrow \Rbold^+$.
The death rate $c_0$ is constant,
\begin{eqnarray}
	c_0(x,m,\xi) = \left\{ \begin{array}{ll}
		\delta, & \mbox{if $\xi^m(x) = 1$} \\
		0, & \mbox{otherwise}
        \end{array} \right. , \label{def:c0}
\end{eqnarray}
while the birth rate $c_1(x,m,\xi)$ is positive only if both male
and female particles can be found in $\Ncal_x$. For example,
the diploid branching particle model we consider in
Chapter~\ref{sec:model:static} a bit later has
\begin{eqnarray}
	c_1(x,m,\xi) = \left\{ \begin{array}{ll}
		\lambda n_1(x,\xi) n_2(x,\xi), & \mbox{if $\xi^m(x) = 0$} \\
		0, & \mbox{otherwise}
	\end{array} \right. , \label{def:c1}
\end{eqnarray}
where
\begin{eqnarray*}
        n_{m'}(x,\xi)=|\{z\in \Ncal_x:\xi^{m'}(x+z)=1\}|,
\end{eqnarray*}
i.e. at rate $\lambda$, each pair of male and female particles in $\Ncal_x$
give birth to a particle at nest $m$ of site $x$ if that nest is not
already occupied. A more stringent condition, as in the particle model
with rapid stirring we consider in Chapter~\ref{sec:2stir} a bit later,
is to require both parent particles to reside at the same site, i.e.
\begin{eqnarray}
        c_1(x,m,\xi) = \left\{ \begin{array}{ll}
                \lambda n_{1+2}(x,\xi), & \mbox{if $\xi^m(x)= 0$} \\
                0, & \mbox{otherwise}
        \end{array} \right. , \label{def1:c2}
\end{eqnarray}
where
\begin{eqnarray*}
        n_{1+2}(x,\xi)=|\{z\in \Ncal_x: \xi^1(x+z)=1 \mbox{ and }
		\xi^2(x+z)=1\}|.
\end{eqnarray*}
This more stringent condition
should not alter the behaviour of the particle system if one
allows a larger $\lambda$ than in~(\ref{def:c1}),
but it does help to simplify the analysis somewhat.

\section{Diploid Branching Particle Model}
\label{sec:model:static}
We first describe the model with birth and death rates as in~(\ref{def:c0})
and~(\ref{def:c1}), for which we will establish the existence of
nontrivial stationary distribution(s) and consequently a phase transition
later in Chapter~\ref{sec:static}.
First, we restate the model in words:
\begin{enumerate}
\item {\bf Birth.} For each nest $(x,m)$ and each pair
	$(z_1,z_2)\in \Ncal_x \times \Ncal_x$
	such that $\xi^1(z_1)=1$ and $\xi^2(z_2)=1$, where $z_1$ and $z_2$
	need not be distinct, with rate $\lambda$, a child of $(z_1,z_2)$
	is born into nest $m$ of site $x$ if $(x,m)$ is not already occupied.
\item {\bf Death.} Each particle dies at rate $\delta$.
\end{enumerate}

We can think of this particle system as a generalized spin system,
generalized in the sense that the phase space at each site is $\{0,1\}^2$
rather than $\{0,1\}$. One can refer to Chapter 3 of [Liggett 1985]
for a detailed introduction on classic spin systems.
We observe that the ``all-$0$'' state (i.e. $\xi^1(x)=\xi^2(x)=0$ for all $x$)
is an absorbing state, therefore the probability measure that concentrates
only on the ``all-$0$'' state is a trivial stationary distribution. We 
say a stationary distribution is \emph{nontrivial} if it does not concentrate
only on the ``all-$0$'' state. A major goal of this work is to establish
the existence of nontrivial stationary distributions for various particle
systems.

This interacting particle system involving the birth and death mechanisms
described above can be constructed using a countable
number of Poisson processes [Durrett 1995].
Without any loss of generality, we assume
$\lambda$ and $\delta$ to be $\leq 1$, since we can just slow down time
by $\max(\lambda,\delta)$ if either $\lambda>1$ or $\delta>1$. Define
\[ c^* = \sup_{\xi,m} \sum_i c_i (x,m,\xi). \]
We assume $c^*<\infty$.
Let $\{T^{x,i,m}_n: n\geq 1\}$ be
the arrival times of independent rate $c^*$ Poisson processes, and
$\{U^{x,i,m}_n: n\geq 1\}$ be independent uniform random variables on $[0,1]$.
At time $t=T^{x,i,m}_n$, nest $(x,m)$ flips to state $i$ if
$U^{x,i,m}_n \leq c_i(x,m,\xi_{t-}) / c^*$,
and stays unchanged otherwise.

Since the number of Poisson processes is
infinite, there is no first flip and the existence and uniqueness of
the process from this construction
is not completely trivial. One can, however, use Theorem 2.1
of [Durrett 1995] to find a small $t_0$ such that the spatial grid $S$
can be divided into an infinite number of components,
each of which is finite and 
no two of which interact during time $[0,t_0]$. This allows construction
of the process up to time $t_0$, and by iterating this procedure, we
can construct the process for all $t$. One can also see easily that
this construction is unique. Alternatively, one can explicitly write down
the generators $\Gcal^1$ and $\Gcal^2$ associated with the particle system
with death rates~(\ref{def:c0}) and birth rates~(\ref{def:c1})
and~(\ref{def1:c2}) respectively:
\begin{eqnarray}
	\Gcal^1 f(\xi) & = & \sum_{(x,m)\in S\times\{1,2\}} \left[
		\rule{0mm}{8mm}
		\delta \xi^m(x) (f(\xi-\delta_{x,m})-f(\xi)) \right. \nnb \\
	& & \ \ \ \ \ \left.
			+ \sum_{y,z\in \Ncal_x} \lambda \xi^1(y) \xi^2(z)
			(1-\xi^m(x)) (f(\xi+\delta_{x,m})-f(\xi)) \right]
		\label{def:gen1} \\
	\Gcal^2 f(\xi) & = & \sum_{(x,m)\in S\times\{1,2\}} \left[
		\rule{0mm}{8mm}
                \delta \xi^m(x) (f(\xi-\delta_{x,m})-f(\xi)) \right. \nnb \\
        & & \ \ \ \ \ \left.
			+ \sum_{y\in \Ncal_x} \lambda \xi^1(y) \xi^2(y)
                        (1-\xi^m(x)) (f(\xi+\delta_{x,m})-f(\xi)) \right]
		\label{def:gen2}
\end{eqnarray}
where $f$ has compact support and $\delta_{x,m}$ is a function
on $S\times\{1,2\}$ that is one at $(x,m)$ and zero everywhere else,
and apply Theorem B3 in [Liggett 1999] (Theorem I.3.9 in [Liggett 1985]
only gives the Markov property) to see that $\Gcal$ is a Markov
generator and therefore determines a unique $(\{0,1\}^2)^{\Zbold^d}$
Feller Markov process.

An important consequence of the construction using Poisson processes described
in the last paragraph is that the semigroup $T_t f(\xi_0)=E^{\xi_0} f(\xi_t)$
corresponding to the particle system is a Feller semigroup. As in
Corollary 2.3 of [Durrett 1995], one can show
that if $\xi^n_0 \rightarrow \xi_0$, then for $t\leq t_0$,
$E^{\xi^n_0} f(\xi^n_t) \rightarrow E^{\xi_0} f(\xi_t)$ since $S$ consists
of components that are finite and do not interact with each other during
$[0,t_0]$. One can then iterate this for as many times as one likes.
Summarizing results from the three previous paragraphs,
we have the following theorem:

\begin{THM}
	There exists a unique Feller process $\xi_t$ constructed as before
	with generator~(\ref{def:gen1}) or~(\ref{def:gen2}).
\label{thm:construct}
\end{THM}

One can represent this construction graphically, for which we give an
example with $S=\Zbold$ and $\Ncal = \{-1,0,1\}$ in figure~\ref{fig:construct}.
Let $m\in  \{1,2\}$, $x,y,z \in S$,
$\{R^{x,m}_n, n\geq 1\}$ be independent Poisson processes
with rate $\delta$, and $\{T^{x,m,y,z}_n, n\geq 1\}$, with $y,z \in \Ncal_x$,
be independent Poisson processes with rate $\lambda$. At space-time points
$((x,m),R^{x,m}_n)$, we draw a symbol $\delta$ to indicate that the particle
(if any) residing at $(x,m)$ is killed at time $R^{x,m}_n$.
At space-time points $((x,m),T^{x,m,y,z}_n)$, we draw arrows from
$((y,1),T^{x,m,y,z}_n)$ and $((z,2),T^{x,m,y,z}_n)$ to $((x,m),T^{x,m,y,z}_n)$
to indicate that a birth event will occur at nest $(x,m)$
if $(x,m)$ is not already occupied and nests
$(y,1)$ and $(z,2)$ are both occupied at time $T^{x,m,y,z}_n$.
{
\psfrag{d}{$\delta$}
\psfrag{-2}{$-2$}
\psfrag{-1}{$-1$}
\psfrag{0}{$0$}
\psfrag{1}{$1$}
\psfrag{2}{$2$}
\psfrag{t}{$t$}
\psfrag{t0}{$t=0$}
\begin{figure}[h!]
\centering
\includegraphics[height = 4in, width = 5.5in]{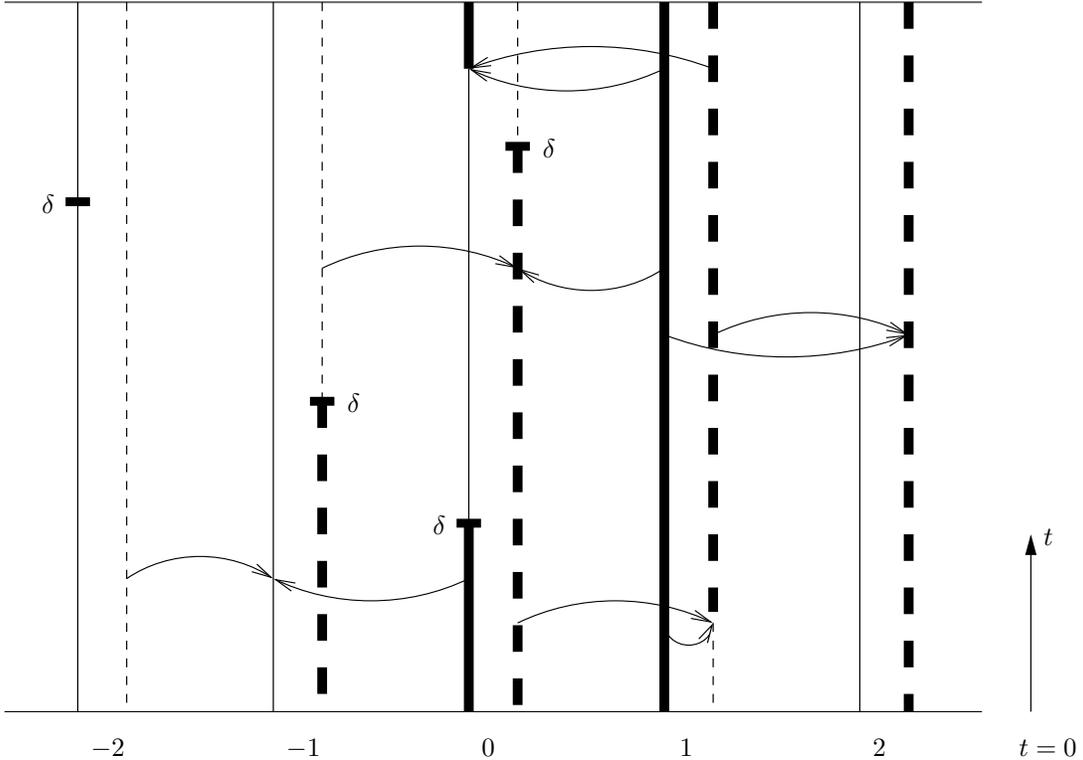}
\caption{Graphical representation: the solid lines represent nests $(x,1)$,
while the dotted lines represent nests $(x,2)$. Thick lines indicate occupied
(wet) nests, while thin lines indicate empty nests.}
\label{fig:construct}
\end{figure}
}

In figure~\ref{fig:construct}, the bottom line represent the state
(occupied or empty) of nests at $t=0$.
We use thick lines to represent occupied (wet) nests, and thin lines
to represent empty nests.
Without any birth or death event, the state of a nest remains unchanged as
$t$ increases. At a death event, i.e. at points
marked by $\delta$, a thick line is changed to a thin line, while a thin
line remains unchanged. And at a birth event, the state of nests at the origins
of the two arrows pointing at $(x,m)$ is checked -- if they are both occupied,
then a thin line at $(x,m)$ is changed to a thick line, while a thick
remains unchanged; otherwise, nothing happens.

In Chapter~\ref{sec:static}, we will use this graphical construction to
establish the existence of
a nontrivial stationary distribution for the diploid branching particle model
if $\lambda/\delta$ is sufficiently large, and
extinction if $\lambda/\delta$ is sufficiently small.

The particle system $\xi$ with generator~(\ref{def:gen1}) or~(\ref{def:gen2})
is \emph{attractive} in the sense that
$\xi$ is monotonic in initial conditions.
One can check that if $\xi_0(x)\leq \bar{\xi}_0(x)$ for all $x \in S$,
where $(0,0)< (0,1) < (1,1)$ and $(0,0)< (1,0) < (1,1)$ but
$(0,1) \not\leq (1,0)$, then $\xi_t(x)\leq \bar{\xi}_t(x)$ for all $x$ and $t$.
This is true since every birth or death event preserves the inequality $\leq$.
For example, if $\xi_{t-}(x)=(0,0)$ and $\bar{\xi}_{t-}(x) = (0,1)$, and
at time $t$ there is a male birth event at site $x$, then
$\xi_{t}(x)=(1,0)$ and $\bar{\xi}_{t}(x) = (1,1)$, so
the inequality $\xi_t(x)\leq \bar{\xi}_t(x)$ has been maintained.
Similarly, one can check that the particle system $\xi$ is increasing
in the birth rate $\lambda$ and decreasing in the death rate $\delta$,
by coupling the random variables $T^{x,i,m}_n$ and $U^{x,i,m}_n$
involved in the constructions in the obvious way. Because of this monotonicity,
along with the existence of nontrivial stationary distributions for
sufficiently large $\lambda/\delta$ and extinction for
sufficiently small $\lambda/\delta$ which we will establish a bit later
in this work, we may conclude that there is
a phase transition in the behaviour of the particle system $\xi$.

\section{Description of the Particle Model with Rapid Stirring}
\label{sec:2stir}
If we add rapid stirring to the particle system, i.e. we scale space by
$\epsilon$ and ``stir'' neighbouring particles at rate $\epsilon^{-2}$
in addition to performing the birth and death mechanisms,
then the particle system converges to the solution of a reaction-diffusion
PDE as $\epsilon \rightarrow 0$ (see Theorems 8.1 and 8.2 in [Durrett 1995]
and the beginning of Chapter~\ref{sec:conv:lily_stir}).
This PDE represents the ``mean-field'' behaviour
of the particle system and is usually easier to analyze than the
particle system itself.
As promised earlier, we will establish later in Chapter~\ref{sec:static}
that there is a phase transition for the diploid branching particle model
(i.e. without rapid stirring), but obtaining any reasonable estimates
on exactly where this transition occurs seems to be difficult. One advantage
of adding rapid stirring mechanisms is that one can get a pretty good idea
where the phase transition occurs in the rapidly stirred particle model
by analyzing the limiting PDE, or simulating this PDE on a computer.

Moreover, this convergence establishes a connection between the particle
model and PDE systems, which is of independent interest.
Since many PDE's arise out of natural systems, this connection justifies
the study of the PDE. The underlying stochastic system can also yield
information about the PDE; for example, in our case,
as we will see in Chapter~\ref{sec:conv:lily_stir},
the monotonicity of the particle system will lead to the monotonicity
of the PDE. Information about the PDE will similarly yield information
about the particle model. In Chapter~\ref{section2}, we will establish
condition ($*$) on the PDE (see page~\pageref{page:star}),
which will tell us that there exist
nontrivial stationary distributions for the particle system with sufficiently
small $\epsilon$.

For the particle models with rapid stirring, we work with
$S=\epsilon \Zbold^d$, and denote the corresponding process by $\xi^\epsilon$.
We also assume the birth and death rates in~(\ref{def:c0}) and~(\ref{def1:c2})
(i.e. generator $\Gcal^2$ in~(\ref{def:gen2}))
with $\delta = 1$, while the neighbourhood $\Ncal$ is nearest neighbour:
\begin{eqnarray*}
	\Ncal = \{y: \|y\| = 0 \mbox{ or } \epsilon \}.
\end{eqnarray*}
Here we use the $L^1$-norm: $\|y\| = \sum_{k=1}^d y_k$.
In addition to the transitions in the diploid branching model,
we introduce spatial movement
of particles between neighbouring sites called ``rapid stirring''.
We consider two rapid stirring mechanisms in this
work, one called ``lily-pad'' stirring, and the other called
``individual'' stirring:
\begin{description}
\item[$\bullet$] {\bf Lily-pad Stirring. }
        For each $x,y\in\epsilon\mathbb{Z}^d$ with $\|x-y\|_1 = \epsilon$,
        $\xi^\epsilon(x)=(\xi^{\epsilon,1}(x),\xi^{\epsilon,2}(x))$
	and $\xi^\epsilon(y)=(\xi^{\epsilon,1}(y),\xi^{\epsilon,2}(y))$
	are exchanged at rate $\epsilon^{-2}$.
\item[$\bullet$] {\bf Individual Stirring. }
        For each $i \in \{1,2\}$ and $x,y\in\epsilon\mathbb{Z}^d$
        with $\|x-y\|_1 = \epsilon$, $\xi^{\epsilon,i}(x)$ and
        $\xi^{\epsilon,i}(y)$ are exchanged at rate $\epsilon^{-2}$.
\end{description}

Just as in the particle model without rapid stirring described
in~\ref{sec:model:static}, one can construct the particle model with
either lily-pad stirring or individual stirring using a countable number
of Poisson processes. Or alternatively, one can write down the generator
explicitly and again apply Theorem B3 in [Liggett 1999] to establish:

\begin{THM}
        There exists a unique Feller process $\xi_t$ with generator
	$\Gcal^L$ for the particle model with Lily-pad stirring
	or generator $\Gcal^I$ for the particle model with individual stirring:
\begin{eqnarray}
        \Gcal^L f(\xi) & = & \Gcal^2 f(\xi)
		+ \sum_{x,y\epsilon\Zbold^d,\|x-y\|_1=\epsilon}
		\epsilon^{-2} (f(\xi^{x\leftrightarrow y})-f(\xi)) \label{def:genL}\\
	\Gcal^I f(\xi) & = & \Gcal^2 f(\xi)
		+ \sum_{m\in\{1,2\},x,y\epsilon\Zbold^d,\|x-y\|_1=\epsilon}
                \epsilon^{-2} (f(\xi^{(x,m)\leftrightarrow (y,m)})-f(\xi)), \label{def:genI}
\end{eqnarray}
where
\begin{eqnarray*}
	\xi^{x\leftrightarrow y} (z,m') = \left\{ \begin{array}{ll}
		\xi(z,m'), & \mbox{if $z\neq x,y$} \\
		\xi(x,m'), & \mbox{if $z=y$} \\
		\xi(y,m'), & \mbox{if $z=x$}
	\end{array} \right. ,
\end{eqnarray*}
and
\begin{eqnarray*}
        \xi^{(x,m)\leftrightarrow (y,m)} (z,m') = \left\{ \begin{array}{ll}
                \xi(z,m'), & \mbox{if $(z,m')\neq (x,m),(y,m)$} \\
                \xi(x,m), & \mbox{if $(z,m')=(y,m)$} \\
                \xi(y,m), & \mbox{if $(z,m')=(x,m)$}
        \end{array} \right. .
\end{eqnarray*}

\end{THM}

For lily-pad stirring, instead of thinking of a site that consists of two nests
as in the diploid branching model, we can view each site as having 4 states in
\begin{eqnarray*}
	F = \{0,1\}^2 = \{(0,0), (0,1), (1,0), (1,1)\}.
\end{eqnarray*}
We restate the dynamics of the
particle model in terms of these four states.
At any site $x\in \epsilon\mathbb{Z}^d$,
only the following transitions are possible:
$(0,0)\leftrightarrow (0,1)$, $(0,1)\leftrightarrow (1,1)$,
$(0,0)\leftrightarrow (1,0)$, and $(1,0)\leftrightarrow (1,1)$, i.e. only one particle
is born or dies at a particular time. The rates of these transitions are
as follows:
\begin{eqnarray*}
\begin{array}{ll}
	c_{(0,0)} (x,\xi^\epsilon) = 1 & \mbox{ if } \xi^\epsilon(x)=(0,1)
		\mbox{ or } \xi^\epsilon(x)=(1,0), \\
	c_{(0,1)} (x,\xi^\epsilon) = c_{(1,0)} (x,\xi^\epsilon) = 1
	    & \mbox{ if } \xi^\epsilon(x)=(1,1), \\
	c_{(0,1)} (x,\xi^\epsilon) = c_{(1,0)} (x,\xi^\epsilon) =
	    \lambda n_{1+2}(x,\xi^\epsilon)
	    & \mbox{ if } \xi^\epsilon(x)=(0,0), \\
	c_{(1,1)} (x,\xi^\epsilon) = \lambda n_{1+2}(x,\xi^\epsilon)
	    & \mbox{ if } \xi^\epsilon(x)=(0,1) \mbox{ or }
		\xi^\epsilon(x)=(1,0).
\end{array}
\end{eqnarray*}
In words, each particle, male or female, dies at rate 1. If site $x$ is
occupied by both a male and a female particle, then with rate $\lambda$,
it gives birth to a male (respectively female) child-particle
at a neighbouring site
provided that the neighbouring site is not already occupied by a male
(respectively female) particle.
%As in the diploid branching model, one can check that the particle
%system $\xi^\epsilon$ is monotone in initial conditions, i.e.
%if $\xi^\epsilon_0(x)\leq \bar{\xi}^\epsilon_0(x)$
%for all $x \in \epsilon\mathbb{Z}^d$,
%where $(0,0)< (0,1) < (1,1)$ and $(0,0)< (1,0) < (1,1)$ but
%$(0,1) \not\leq (1,0)$,
%then $\xi^\epsilon_t(x)\leq \bar{\xi}^\epsilon_t(x)$ for all $x$ and $t$.
%One can also check that $\xi^\epsilon$ is monotone in $\lambda$.

For individual stirring, we still view
the particle system with nests
$(x,m) \in \epsilon \Zbold^d \times \{1,2\}$ and each nest
assuming one of two states in $E=\{0,1\}$.

The difference between these two stirring mechanisms is that lily-pad stirring
forces male and female particles at a site to move together, but
individual stirring allows independent movement of male and female
particles. Every exchange of particles, in both lily-pad stirring
and individual stirring, is monotonicity preserving, thus neither
stirring mechanism disrupts the monotonicity property of the particle system.

\section{Convergence to a PDE for Lily-pad Stirring}
\label{sec:conv:lily_stir}
Consider the particle system with lily-pad stirring with its generator
given by~(\ref{def:genL}). For $i\in F$, define
\begin{eqnarray*}
	u^\epsilon_i (t,x) = P(\xi^\epsilon_t(x) = i)
\end{eqnarray*}
then Theorem 8.1 in [Durrett 1995] shows that if $g_i (x)$ is continuous
and $u^\epsilon_i (0,x) = g_i (x)$ for all $i$, then
\begin{eqnarray*}
	u_i (t,x) = \lim_{\epsilon\rightarrow 0} u^\epsilon_i (t,x)
\end{eqnarray*}
exists and satisfies the following system of PDE's:
\begin{eqnarray}
	\frac{\partial u_{(0,0)}}{\partial t} & = & \Delta u_{(0,0)}
	+u_{(0,1)}+u_{(1,0)}-2\lambda d u_{(0,0)} u_{(1,1)} \nnb \\
	\frac{\partial u_{(0,1)}}{\partial t} & = & \Delta u_{(0,1)}
	+u_{(1,1)}-u_{(0,1)}+\lambda d (u_{(0,0)}-u_{(0,1)}) u_{(1,1)} \nnb \\
	\frac{\partial u_{(1,0)}}{\partial t} & = & \Delta u_{(1,0)}
        +u_{(1,1)}-u_{(1,0)}+\lambda d (u_{(0,0)}-u_{(1,0)}) u_{(1,1)} \nnb \\
	\frac{\partial u_{(1,1)}}{\partial t} & = & \Delta u_{(1,1)}
	-2u_{(1,1)}+\lambda d (u_{(0,1)}+u_{(1,0)}) u_{(1,1)}.
	\label{eq:pde_init}
\end{eqnarray}
Obviously, $u_i$ must lie in $[0,1]$ for all $i$ and $t$ since it is a
limit of probabilities. We want to study the long time behaviour
of~(\ref{eq:pde_init}). The system~(\ref{eq:pde_init}) is 3-dimensional
if one takes into account the condition
$u_{(0,0)}+u_{(0,1)}+u_{(1,0)}+u_{(1,1)} = 1$.
We first do two transformations on the 3-dimensional parameter
space $(u_{(0,0)},u_{(0,1)},u_{(1,0)},u_{(1,1)})$ to obtain a
monotone 2-dimensional system, which will be easier to analyze. First,
define $u_0=u_{(0,0)}$, $u_1=u_{(0,1)}+u_{(1,0)}$, and $u_2=u_{(1,1)}$,
then $(u_0,u_1,u_2)$ satisfies:
\begin{eqnarray}
	\frac{\partial u_0}{\partial t} & = & \Delta u_0
	+ u_1 - 2\lambda d u_0 u_2 \nonumber \\
	\frac{\partial u_1}{\partial t} & = & \Delta u_1
	+ 2u_2 - u_1 + \lambda d (2u_0-u_1)u_2 \nonumber \\
	\frac{\partial u_2}{\partial t} & = & \Delta u_2	
	-2u_2 + \lambda d u_1 u_2. \label{eq:pde}
\end{eqnarray}
The above system can be written as the limiting PDE under rapid stirring
of another particle system $\zeta^\epsilon$, still on $S = \epsilon \Zbold^d$,
with state space $F=\{0,1,2\}$, and transitions
$0 \leftrightarrow 1$ and $1 \leftrightarrow 2$ at rates
\begin{eqnarray*}
\begin{array}{ll}
        c_0 (x,\zeta^\epsilon) = 1 & \mbox{ if } \zeta^\epsilon(x)=1 \\
	c_1 (x,\zeta^\epsilon) = 2 & \mbox{ if } \zeta^\epsilon(x)=2 \\
        c_1 (x,\zeta^\epsilon) = 2 \lambda n_2(x,\zeta^\epsilon)
                & \mbox{ if } \zeta^\epsilon(x)=0 \\
        c_2 (x,\zeta^\epsilon) = \lambda n_2(x,\zeta^\epsilon)
		& \mbox{ if } \zeta^\epsilon(x)=1,
\end{array}
\end{eqnarray*}
where
\begin{eqnarray*}
        n_2(x,\xi^\epsilon)=|\{z\in \Ncal:\zeta^\epsilon(x+z)=2\}|.
\end{eqnarray*}
Under this model, monotonicity still holds: if
$\zeta^\epsilon_0(x)\leq \bar{\zeta}^\epsilon_0(x)$ for all $x$ (here
the ordering of $F$ is the usual one: $0<1<2$),
then $\zeta^\epsilon_t(x)\leq \bar{\zeta}^\epsilon_t(x)$ for all $x$ and $t$,
since every transition still preserves the inequality $\leq$.
%Define $u^\epsilon_i (t,x) = P(\zeta^\epsilon_t(x) = i)$, $i=0,1,2$.
Let $(u_1^\epsilon(0,x),u_2^\epsilon(0,x))=(g_1(x),g_2(x))$ and
$(\bar{u}_1^\epsilon(0,x),\bar{u}_2^\epsilon(0,x))
=(\bar{g}_1(x),\bar{g}_2(x))$ be two sets
of initial distributions such that $g_1+g_2\leq \bar{g}_1+\bar{g}_2$
and $g_2\leq \bar{g}_2$ everywhere. Then for all $x$,
\begin{eqnarray*}
	u_2^\epsilon(0,x) & \leq & \bar{u}_2^\epsilon(0,x) \\
	u_1^\epsilon(0,x) + u_2^\epsilon(0,x) & \leq &
	\bar{u}_1^\epsilon(0,x) + \bar{u}_2^\epsilon(0,x),
\end{eqnarray*}
so it is possible to set up two initial conditions $\zeta^\epsilon_0$
and $\bar{\zeta}^\epsilon_0$, such that
$P(\zeta^\epsilon_0(x)=i)=u_i^\epsilon(0,x)$
and $P(\bar{\zeta}^\epsilon_0(x)=i)=\bar{u}_i^\epsilon(0,x)$, $i=1,2$, and
$\zeta^\epsilon_0(x)\leq \bar{\zeta}^\epsilon_0(x)$ holds for all $x$ and
$\omega$. Since $\zeta^\epsilon_t(x)\leq \bar{\zeta}^\epsilon_t(x)$
for all $t$ and $x$, the monotonicity property of $\zeta$ implies
\begin{eqnarray*}
	P(\zeta^\epsilon_t(x)\geq 1) \leq P(\bar{\zeta}^\epsilon_t(x)\geq 1)
	\mbox{ and }
	P(\zeta^\epsilon_t(x)\geq 2) \leq P(\bar{\zeta}^\epsilon_t(x)\geq 2),
\end{eqnarray*}
i.e. for all $t$ and $x$,
\begin{eqnarray*}
        u_2^\epsilon(t,x) & \leq & \bar{u}_2^\epsilon(t,x) \\
        u_1^\epsilon(t,x) + u_2^\epsilon(t,x) & \leq &
        \bar{u}_1^\epsilon(t,x) + \bar{u}_2^\epsilon(t,x).
\end{eqnarray*}
Now we transform the parameter space a second time by defining
$(\alpha,\beta)=(u_1+u_2,u_2)$ and writing $c=\lambda d$, then
$(\alpha_t,\beta_t)$ is monotone in the initial condition
since $u_i^\epsilon(t,x)\rightarrow u_i(t,x)$.
In particular, if $(1-\alpha,\alpha-\beta,\beta)$ and
$(1-\bar{\alpha},\bar{\alpha}-\bar{\beta},\bar{\beta})$ are both solutions
to~(\ref{eq:pde}) with 
$\alpha_0(x)<\bar{\alpha}_0 (x)$ and $\beta_0(x)<\bar{\beta}_0 (x)$
for all $x$, then
$\alpha_t(x)<\bar{\alpha}_t (x)$ and $\beta_t(x)<\bar{\beta}_t (x)$
for all $x$ and $t$. Straightforward calculation shows that $(\alpha,\beta)$
satisfies the following system:
\begin{eqnarray*}
        \frac{\partial \alpha}{\partial t} & = & \Delta \alpha +
	(2c-2c\alpha+1)\beta-\alpha \nonumber \\
        \frac{\partial \beta}{\partial t} & = & \Delta \beta +
	(c(\alpha-\beta)-2)\beta.
\end{eqnarray*}
Since $(u_1,u_2)\in [0,1]^2$, $(\alpha,\beta)$ lies in the triangular region
\begin{eqnarray}
	\R=\{(u,v):0\leq u,v\leq 1,u\geq v\} \label{def:R:0}
\end{eqnarray}
for all $t \geq 0$. We change variables from $(\alpha,\beta)$ to $(u,v)$ and
summarize this paragraph in the following lemma:
\begin{LEM}
The PDE
\begin{eqnarray}
        \frac{\partial u}{\partial t} & = & \Delta u +
        (2c(1-u)+1)v-u \nonumber \\
        \frac{\partial v}{\partial t} & = & \Delta v +
        (c(u-v)-2)v \label{eq:pde0}
\end{eqnarray}
is monotone in initial conditions that lie in
$\R=\{(u,v):0\leq u,v\leq 1,u\geq v\}$,
i.e.  if there are two initial conditions
$(u_0,v_0)\in \R$ and $(\bar{u}_0,\bar{v}_0)\in \R$,
with $u_0\leq \bar{u}_0$ and $v_0 \leq \bar{v}_0$ everywhere,
then $u_t\leq \bar{u}_t$ and $v_t \leq \bar{v}_t$ everywhere, for all $t$;
furthermore, both $(u_t,v_t)$ and $(\bar u_t,\bar v_t)$ lie in $\R$
for all $t$.
\label{lem:monotone}
\end{LEM}
In Chapter~\ref{section2}, we will analyze~(\ref{eq:pde0}) to establish
the following theorem:
\begin{THM}
	If $\lambda/\delta$ is sufficiently large and
	$\epsilon$ is sufficiently small, then there exists
	a nontrivial translation invariant stationary distribution
	for the diploid branching particle model with lily-pad stirring
	with generator~(\ref{def:genL}).
\label{thm:lily}
\end{THM}

\section{Convergence to a PDE for Individual Stirring}
\label{sec:conv:ind_stir}
Unlike lily-pad stirring, Theorem 8.1 in [Durrett 1995] cannot be directly
applied to get convergence to a PDE system for individual stirring.
We can, however, follow the ideas used in the proof of that theorem 8.1 to
establish a corresponding result, Theorem~\ref{thm:indp_stir}.
We consider the particle model with individual stirring described
in Chapter~\ref{sec:2stir}. For $i \in E$, define
\begin{eqnarray*}
	u^\epsilon_{i,m}(t,x) = P(\xi^\epsilon_t(x,m) = i).
\end{eqnarray*}
Then Theorem~\ref{thm:indp_stir} implies that if
$g_{i,m}: \Rbold \rightarrow [0,1]$
is continuous and $u^\epsilon_{i,m}(0,x) = g_{i,m}(x)$, then
\begin{eqnarray*}
        u_{i,m} (t,x) = \lim_{\epsilon\rightarrow 0} u^\epsilon_{i,m} (t,x)
\end{eqnarray*}
exists and satisfies the following system of PDE's:
\begin{eqnarray*}
        \frac{\partial u_{0,1}}{\partial t} & = & \Delta u_{0,1}
		+ u_{1,1} - 2 c u_{0,1} u_{1,1} u_{1,2} \nnb \\
	\frac{\partial u_{1,1}}{\partial t} & = & \Delta u_{1,1}
                - u_{1,1} + 2 c u_{0,1} u_{1,1} u_{1,2} \nnb \\
	\frac{\partial u_{0,2}}{\partial t} & = & \Delta u_{0,2}
                + u_{1,2} - 2 c u_{0,2} u_{1,1} u_{1,2} \nnb \\
        \frac{\partial u_{1,2}}{\partial t} & = & \Delta u_{1,2}
                - u_{1,2} + 2 c u_{0,2} u_{1,1} u_{1,2},
\end{eqnarray*}
where we define
\begin{eqnarray*}
	c = \lambda d^2.
\end{eqnarray*} 
Since $u_{0,1}+u_{1,1}=u_{0,2}+u_{1,2}=1$, it suffices to study the PDE
for $u_{1,1}$ and $u_{1,2}$:
\begin{eqnarray}
	\frac{\partial u_{1,1}}{\partial t} & = & \Delta u_{1,1}
                - u_{1,1} + 2 c (1 - u_{1,1}) u_{1,1} u_{1,2} \nnb \\
	\frac{\partial u_{1,2}}{\partial t} & = & \Delta u_{1,2}
                - u_{1,2} + 2 c (1 - u_{1,2}) u_{1,1} u_{1,2}.
	\label{eq:pde:indp_stir2}
\end{eqnarray}
Notice that if we start with a symmetric initial condition, i.e.
$g_{i,1}=g_{i,2}$, then the solution to~(\ref{eq:pde:indp_stir2}) is
also symmetric. And if we define $u=u_{1,1}=u_{1,2}$, then we obtain
the following PDE for $u$:
\begin{eqnarray}
	\frac{\partial u}{\partial t} & = & \Delta u + f(u),
                 \label{eq:pde:indp_stir} \\
	f(u) & = & - u + 2 c (1-u) u^2. \nnb
\end{eqnarray}
This PDE has been analyzed in [Durrett and Neuhauser 1994] as their
sexual reproduction model (example 3 on page 291). In fact, it is not
difficult to see that if $u_{1,1}=u_{1,2}$ then choosing the ``father''
from the male population is exactly the same as choosing the ``father''
from the the female population, hence it is quite natural for this reduction
to occur. Theorem 4 of [Durrett and Neuhauser 1994] states that
if $c>2.25$ then the sexual reproduction model of Durrett and Neuhauser
has nontrivial stationary distribution(s). Although this theorem does not
directly apply to our particle system $\xi^\epsilon$
with 2 types of particles because of the difference in stirring mechanisms,
one can nevertheless trace through the proof of Lemma 3.3 of
[Durrett and Neuhauser 1994]
while making obvious changes, to establish a similar result:
\begin{enumerate}
\item[$\bullet$]
	Let $0 < \rho_1 < \rho_0 < 1$ be the two nonzero roots
	of $f(u)$. Define $\beta = (\rho_0 - \rho_1)/10$ and
	$I_k = 2Lk e_1 + [-L,L)^d$.
	If $\epsilon$ is small, $L$ is large, and $\xi^\epsilon(0)$ has density
	at least $\rho_1+\beta$ of both male particles and female particles
	in $I_0$, then for sufficiently large $T$,
	with high probability $\xi^\epsilon(T)$
	will have density of at least $\rho_0-\beta$ in $I_1$ and $I_{-1}$.
\end{enumerate}
This result can then be fed into a comparison argument, comparing
the particle system with oriented percolation, as on page 312 of
[Durrett and Neuhauser 1994] or in the proof of Theorem~\ref{thm:use:comp}
later on, to establish the existence of nontrivial stationary
distribution(s) for the particle system $\xi^\epsilon$ under individual
stirring with sufficiently small $\epsilon$. Then we have the following
theorem:
\begin{THM}
        If $\lambda/\delta$ is sufficiently large and
	$\epsilon$ is sufficiently small, then there exists
        a nontrivial translation invariant stationary distribution
        for the diploid branching particle model with individual stirring
	with generator~(\ref{def:genI}).
\end{THM}

We will not explicitly write down the details of the proof, but instead
refer the interested reader to [Durrett and Neuhauser 1994] for details.

\chapter{Results on the Diploid Branching Particle Model}
\label{sec:static}
In this chapter, we assume the model with generator~(\ref{def:gen1})
described in Chapter~\ref{sec:model:static}, i.e. the
particle system with birth and death mechanisms, but no stirring.
We briefly restate the model to remind the reader: the rate
at which nest $m$ of site $x$ flips to state $i$, $c_i(x,m,\xi)$, is
\begin{eqnarray*}
        c_0(x,m,\xi) = \left\{ \begin{array}{ll}
                \delta, & \mbox{if $\xi^m(x) = 1$} \\
                0, & \mbox{otherwise}
        \end{array} \right., \ 
        c_1(x,m,\xi) = \left\{ \begin{array}{ll}
                \lambda n_1(x,\xi) n_2(x,\xi), & \mbox{if $\xi^m(x) = 0$} \\
                0, & \mbox{otherwise}
        \end{array} \right. ,
\end{eqnarray*}
where
\begin{eqnarray*}
        n_{m'}(x,\xi)=|\{z\in \Ncal_x:\xi^{m'}(x+z)=1\}|,
\end{eqnarray*}
and the $\Ncal_x$ contains the site $x$ and its $2d$ nearest neighbours.
The goal is to establish the existence of a phase transition.

\section{Existence of Stationary Distributions}
We first establish that stationary distributions exist. Define
\begin{eqnarray*}
	\overline\xi_0 (x) = (1,1) \mbox{ for all $x$}.
\end{eqnarray*}
Let $T_t f(\xi_0)=E^{\xi_0} f(\xi_t)$ be the semigroup corresponding to
the particle system, then $T_t$ is a
Feller semigroup by Theorem~\ref{thm:construct}. We begin with a lemma.

\begin{LEM}
	For any $A,B\subset S=\mathbb{Z}^d$, the function
\begin{eqnarray}
	t \mapsto P\left(\overline\xi^1_t(x)=0 \ \forall x\in A, \
		\overline\xi^2_t(y)=0 \ \forall y\in B\right)
	\label{eq:xi:increase}
\end{eqnarray}
	is increasing.
\label{lem:xi:increase}
\end{LEM}
\proof Let $\alpha_0=\overline\xi^1_s$ and $\beta_0=\overline\xi^2_s$
for an arbitrary fixed $s$. Then
$\overline\xi^1_0(x)\geq\alpha_0(x)$ and $\overline\xi^2_0(x)\geq\beta_0(x)$.
Let $(\alpha_t,\beta_t)$ be the state at time $t$ of the particle system
that started with initial condition $(\alpha_0,\beta_0)$.
Then by the fact that the particle system is
monotone in initial conditions, we have
\begin{eqnarray*}
	\overline\xi^1_t(x)\geq\alpha_t(x) \mbox{ and }
	\overline\xi^2_t(x)\geq\beta_t(x)
\end{eqnarray*}
for all $t$ and $x$.
Thus by the Markov property of $\xi$,
\begin{eqnarray}
        \lefteqn{ P\left(\overline\xi^1_t(x)=0 \ \forall x\in A, \
		\overline\xi^2_t(y)=0 \ \forall y\in B\right) } \nnb \\
	& \leq & P\left(\alpha_t(x)=0 \ \forall x\in A, \
		\beta_t(y)=0 \ \forall y\in B\right) \nonumber \\
        & = & P\left(\overline\xi^1_{s+t}(x)=0 \ \forall x\in A, \
		\overline\xi^2_{s+t}(y)=0 \ \forall y\in B\right). \nonumber
\end{eqnarray}
This implies that the function in~(\ref{eq:xi:increase}) is increasing in $t$.
\qed

\vspace{.3cm}
\begin{THM}
	As $t\rightarrow\infty$,
	$\overline\xi_t \Rightarrow \overline\xi_\infty$. The limit
	is a stationary distribution that stochastically dominates all
	other stationary distributions and called the
	\emph{upper invariant measure}.
\label{thm:upper_inv}
\end{THM}
\proof
For arbitrary subsets $A$, $B$, $C=\{x_1,\ldots,x_m\}$,
and $D=\{y_1,\ldots,y_n\}$ of $S$, we write
\begin{eqnarray}
        & & {P\left(\overline\xi^1_t(z)=0 \ \forall z\in A, \
		\overline\xi^2_t(w)=0 \ \forall w\in B, \
        	\overline\xi^1_t(x)=1 \ \forall x\in C, \ 
		\overline\xi^2_t(y)=1 \ \forall y\in D\right)} \nonumber \\
        & = & P\left(\overline\xi^1_t(z)=0 \ \forall z\in A, \
		\overline\xi^2_t(w)=0 \ \forall w\in B\right)
		- P\left(\bigcup_{i=1}^{m+n} E_i\right), \nonumber
\end{eqnarray}
where
\begin{eqnarray*}
	E_i= \{\overline\xi^1_t(z)=0 \ \forall z\in A\cup \{x_i\}, \
		\overline\xi^2_t(w)=0 \ \forall w\in B \} \
			\mbox{if $i=1,\ldots,m$}
\end{eqnarray*}
and
\begin{eqnarray*}
	E_i=	\{\overline\xi^1_t(z)=0 \ \forall z\in A, \
		\overline\xi^2_t(w)=0 \ \forall w\in B\cup\{y_{i-m}\}\} \
			\mbox{if $i=m+1,\ldots,m+n$}.
\end{eqnarray*}
We can use the inclusion-exclusion formula on
$P(\cup_{i=1}^{m+n} E_i)$, i.e.
\begin{eqnarray*}
        P\left(\bigcup_{i=1}^{m+n} E_i\right) = \sum_{i=1}^{m+n} P(E_i)
        	-\sum_{i<j} P(E_i\cap E_j)
		+\ldots+(-1)^{m+n+1} P(E_i \cap\ldots\cap E_j).
\end{eqnarray*}
Every term in the above expansion is in the form of
\begin{eqnarray*}
	P(\overline\xi^1_t(z)=0 \ \forall z\in \cdot, \
		\overline\xi^2_t(w)=0 \ \forall w\in \cdot\cdot),
\end{eqnarray*}
which is increasing in $t$ by Lemma~\ref{lem:xi:increase}. Therefore
\begin{eqnarray*}
	P(\overline\xi^1_t(z)=0 \ \forall z\in A, \
		\overline\xi^2_t(w)=0 \ \forall w\in B, \
		\overline\xi^1_t(x)=1 \ \forall x\in C, \ 
		\overline\xi^2_t(y)=1 \ \forall y\in D)
\end{eqnarray*} 
converges for all $A$, $B$, $C$, and $D$, i.e.
all finite dimensional distributions converge.
Thus a weak limit $(\overline\xi^1_\infty,\overline\xi^2_\infty)$ exists
and it follows from a standard result that since $T_t$ is a Feller semigroup,
$(\overline\xi^1_\infty,\overline\xi^2_\infty)$ is a stationary distribution.
We can also easily see that $(\overline\xi^1_\infty,\overline\xi^2_\infty)$
dominates all other stationary distributions: let
$(\tilde\xi^1_0,\tilde\xi^2_0)$ be another
stationary distribution and $(\tilde\xi^1,\tilde\xi^2)$
be the process with initial
condition $(\tilde\xi^1_0,\tilde\xi^2_0)$,
then $(\tilde\xi^1_t,\tilde\xi^2_t)$ has the same distribution
as $(\tilde\xi^1_0,\tilde\xi^2_0)$ for all $t\geq 0$ and
$(\bar\xi^1,\bar\xi^2)$ dominates $(\tilde\xi^1,\tilde\xi^2)$
because of monotonicity, therefore $(\bar\xi^1_\infty,\bar\xi^2_\infty)$
dominates $(\tilde\xi^1_0,\tilde\xi^2_0)$.
\qed

Theorem III.2.3 in [Liggett 1985] establishes the previous theorem
for spin systems with state space $\{0,1\}^S$, therefore it does not
directly apply to our case. One can, however, easily adapt the proof
of that theorem to this case, and obtain a slightly different proof.

\section{Extinction for Sufficiently Small $\lambda/\delta$}
\begin{THM}
	If $\lambda |\Ncal|^2 < \delta$, then the particle system $\xi$
	has no nontrivial stationary distribution.
\end{THM}
\proof
We compare a modification of the particle system $\xi$ with the contact
process. We recall that the contact process $\zeta$ on $\Zbold^d$ has two
states $0$ and $1$ at every site $x$, and has the following dynamics:
\begin{eqnarray*}
        c_0(x,\zeta) = \left\{ \begin{array}{ll}
                \delta, & \mbox{if $\zeta(x) = 1$} \\
                0, & \mbox{otherwise}
        \end{array} \right., \ \ \
        c_1(x,\zeta) = \left\{ \begin{array}{ll}
                \alpha n(x,\zeta), & \mbox{if $\zeta(x) = 0$} \\
                0, & \mbox{otherwise}
        \end{array} \right. ,
\end{eqnarray*}
where $n(x,\zeta)=|\{z\in \Ncal_x:\zeta(z)=1\}|$.
Theorem 2.6 of [Durrett 1995] states that if $\alpha |\Ncal|< \delta$,
then the contact process has no nontrivial stationary distribution.

We modify the mechanism of the particle model $\xi$ as follows:
for the males, when $(x,1) \in S\times \{1\}$
seeks out a pair of parents, say $(x_1,1)$ and $(x_2,2)$, it is no
longer required that $\xi^2(x_2)=1$, but $\xi^1(x_1)$ must still be 1.
Correspondingly, when a female nest $(x,2) \in S\times \{2\}$
seeks out a pair of parents $(x_1,1)$ and $(x_2,2)$, it is only
required that $\xi^2(x_2)=1$.
We denote this modified process $\tilde\xi=(\tilde\xi^1,\tilde\xi^2)$.
The result of the modification is that $\tilde\xi^1$ and $\tilde\xi^2$ are
now decoupled, and $\tilde\xi^i$ behaves exactly the same as
the contact process with birth rate $\alpha=\lambda |\Ncal|$.
Furthermore, by Theorem III.1.5 in [Liggett 1985],
the modified process $(\tilde\xi^1,\tilde\xi^2)$ stochastically
dominates the original process $(\xi^1,\xi^2)$.
If $\alpha |\Ncal| < \delta$, then $\tilde\xi$ has no nontrivial
stationary distribution and $(\tilde\xi^1_t,\tilde\xi^2_t)$ converges weakly
to the all-$0$ state as $t\rightarrow\infty$ for any initial condition.
Thus $(\xi^1_t,\xi^2_t)$ also converges to the all-$0$ state for
any initial condition if 
$\alpha |\Ncal| = \lambda |\Ncal|^2 < \delta$, as required.
\qed

\vspace{.3cm}
The proof above also shows that if $\sum_x \xi^1_0(x) + \xi^2_0(x)$
is finite, then the population dies out in finite time a.s.
if $\lambda |\Ncal|^2 < \delta$, since the contact process $\tilde\xi^1$
or $\tilde\xi^2$ has this property.

\section{Survival for Sufficiently Large $\lambda/\delta$}
We use the idea of Chapter 4 of [Durrett 1995], i.e. we compare the particle
system to an oriented percolation process. 
First, we define the oriented percolation process $W$. Let
\begin{eqnarray*}
	\Lcal_0 = \{ (x,n)\in \mathbb{Z}^2: x+n \mbox{\ is even}, n\geq 0 \}
\end{eqnarray*}
and make $\Lcal_0$ into a graph by drawing oriented edges from $(x,n)$ to
$(x+1,n+1)$ and from $(x,n)$ to $(x-1,n+1)$. Site $(x,n)$ is said to
be a parent of sites $(x+1,n+1)$ and $(x-1,n+1)$. Notice that
any site $(x,n)$ with $n\neq 0$ has two parents. We think of $n$
as the time variable. Given random
variables $\omega(x,n)$ that indicate whether site $(x,n)$ is open ($1$)
or closed ($0$), we say that $(y,l)$ can be reached from $(x,m)$ if
there is sequence of points $x=z_m,\ldots,z_l=y$ such that
$|z_k-z_{k-1}|=1$ for $m < k \leq l$ and $\omega(z_k,k)=1$ for
$m \leq k \leq l$; in this case, we write
\begin{eqnarray*}
	(x,m) \rightarrow (y,l).
\end{eqnarray*}
We say that $\omega(x,n)$ ($n\geq 1$) is
``\emph{$M$-dependent with density at least $1-\gamma$}'' if whenever
$(x_i,n_i)$, $1\leq i \leq I$, is a sequence with
$\|(x_i,n_i)-(x_j,n_j)\|_\infty > M$ for $i\neq j$, we have
\begin{eqnarray}
	P(\omega(x_i,n_i) = 0 \mbox{ for $1 \leq i \leq I$}) \leq \gamma^I.
	\label{def:Mdep}
\end{eqnarray}
Given an initial condition $W_0 \subset 2\Zbold = \{x: (x,0)\in \Lcal_0\}$,
the process
\begin{eqnarray*}
	W_n = \{y: (x,0) \rightarrow (y,n) \mbox{ for some $x\in W_0$}\}.
\end{eqnarray*}
gives all sites that can be reached from a site in $W_0$ at time $n$.
We say sites in $W_n$ are wet. Theorem 4.2 of [Durrett 1995] states:
\begin{THM}
	Let $W_n^p$ be an $M$-dependent oriented percolation
	with density at least $1-\gamma$ starting from the initial
	configuration $W_0^p$ in which the events $\{x\in W_0^p\}$,
	$x \in 2\Zbold$, are independent and have probability $p$.
	If $p>0$ and $\gamma\leq 6^{-4(2M+1)^2}$, then
\begin{eqnarray}
	\liminf_{n\rightarrow\infty} P(0\in W_{2n}^p) \geq 19/20.
	\label{est:perc:rd}
\end{eqnarray}
\label{thm:perc:rd}
\end{THM}
This theorem shows that if the density of open sites $1-\gamma$ is
sufficiently close to $1$ and we start with a Bernoulli initial condition for
$W_0$, then the probability that $0$ is wet at time $t$ does not go to 0
as $t\rightarrow\infty$. Notice that the right hand side of the
estimate~(\ref{est:perc:rd}) is a constant that does not depend on $p$.
We will construct
an oriented percolation process that is stochastically dominated by
the particle system $\xi$, such that existence of nontrivial
stationary distribution for the oriented percolation process
implies existence of nontrivial stationary distribution for
the particle system.

\begin{THM}
\label{thm:use:comp}
	If $\lambda/\delta$ is sufficiently large, then the particle system
	$\xi$ with generator~(\ref{def:gen1}) has
	a nontrivial stationary distribution.
\end{THM}
\proof
We follow the method of proof as in Chapter 4 of [Durrett 1995].
Since scaling time by an factor of $1/\delta$ does not change the behaviour
with respect to stationary distributions, we may assume without any loss of
generality that $\delta=1$. We will select an event $G_{\xi_0}$
measurable with respect to the graphical
representation in $[-1,1] \times \Zbold^{d-1} \times [0,T)$, i.e.
measurable with respect to the filtration generated by
all the Poisson arrivals $\{R^{x,m}_n\}$ and $\{T^{x,m,y,z}_n\}$ used to
construct the particle model in Chapter~\ref{sec:model:static} that
arrive at any sites in $[-1,1] \times \Zbold^{d-1}$ during the
time interval $[0,t)$.
For any $\gamma>0$ no matter how small, there is $\lambda$
and $T$ and an event $G_{\xi_0}$ with
\begin{eqnarray*}
	P(G_{\xi_0})>1-\gamma,
\end{eqnarray*}
so that on $G_{\xi_0}$, if $\xi_0(0,0,\ldots,0)=(1,1)$, then
$\xi_T(1,0,\ldots,0)=\xi_T(-1,0,\ldots,0)=(1,1)$. One can achieve this
by choosing $T$ so small that the probability of any death occurring at
any nests of sites $(-1,0,\ldots,0)$, $(0,0,\ldots,0)$, and $(1,0,\ldots,0)$
is less than $\gamma/2$; then one can
choose $\lambda$ large enough so that the probability of having
birth events from $((0,0,\ldots,0),1)$ and $((0,0,\ldots,0),2)$
to each of the four nests at sites $(-1,0,\ldots,0)$ and $(1,0,\ldots,0)$
during $[0,T)$ is larger than $1-\frac{\gamma}{2}$. In other words,
if we define the event
\begin{eqnarray*}
	G_{\xi_0} & = & \{\mbox{There are no death event during $[0,T)$ at
		sites $(-1,0,\ldots,0)$, $(0,0,\ldots,0)$, } \\
	& &	\mbox{or $(1,0,\ldots,0)$;
		and there are birth events from $((0,0,\ldots,0),1)$ } \\
	& &	\mbox{and $((0,0,\ldots,0),2)$ to each of the four nests at
		sites $(-1,0,\ldots,0)$ } \\
	& &	\mbox{and $(1,0,\ldots,0)$ during $[0,T)$}\},
\end{eqnarray*}
then $G_{\xi_0}$ satisfies the requirement and $P(G_{\xi_0}) > 1-\gamma$ for
some $\lambda$ and $T$.
$G_{\xi_0}$ is the ``good event'' that will ensure male and female
particles get born at sites $x-1$ and $x+1$ provided site $x$ is inhabited
by both a male and a female particle.
See figure~\ref{fig:goodevent} for an illustration of this event.
{
\psfrag{-1}{$-1$}
\psfrag{0}{$0$}
\psfrag{1}{$1$}
\psfrag{t1}{$t=0$}
\psfrag{t2}{$t=T$}
\begin{figure}[h!]
\centering
\includegraphics[height = 3in, width = 3in]{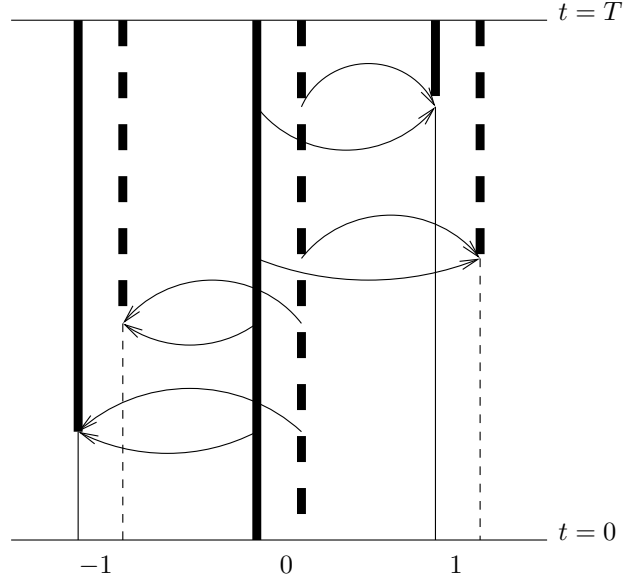}
\caption{Graphical representation of the event $G_{\xi_0}$: at least
4 birth events and no death events.}
\label{fig:goodevent}
\end{figure}
}

We start with a configuration with events $\{\xi^1_0(x)=1\}$ and
$\{\xi^2_0(x)=1\}$ all independent and having probability $p_1$ and $p_2$,
respectively. Let
\begin{eqnarray*}
	X_n = \{x:(x,n)\in \Lcal_0, \xi_{nT}(x,0,\ldots,0)=(1,1)\}.
\end{eqnarray*}
Before defining the oriented percolation process $W_n$,
we first define $V_n$ that will turn out to be slightly larger
than $W_n$ but nevertheless dominated by $X_n$.

We define $V_n$ inductively. First, we set $V_0=X_0$ and leave
$\omega(\cdot,0)$ undefined. Now assume that
$V_0$, $V_1$, $\ldots$, $V_n$, and $\omega(x,l)$ with $l \leq n-1$ have
been defined such that $V_0 \subset X_0$, $\ldots$, $V_n \subset X_n$.
If $x\in V_n$ then we set
\begin{eqnarray}
	\omega(x,n) = \left\{ \begin{array}{ll}
		1, & \mbox{if $G_{\sigma_{-x e_1}(\xi_{nT})}$ occurs in
			the graphical representation } \\
		& \mbox{($G_{\sigma_{-x e_1}(\xi_{nT})}$ is $G_{\xi_0}$
			translated by $-x e_1$ in space and} \\
		& \mbox{$-nT$ in time)} \\
		0, & \mbox{otherwise}
	\end{array} \right. , \label{def:omega}
\end{eqnarray}
which ensures that $\omega(x,n)=1$ with probability more than $1-\gamma$
if $x\in V_n$. For completeness,
if $x\notin V_n$, then we set $\omega(x,n)$ equal to an independent random
variable that is 1 with probability $1-\gamma$ and 0 with probability
$\gamma$. Now we define $V_{n+1}$ to consist of all sites $(x,n+1)$ with
either
\begin{eqnarray}
	x-1\in V_n \ \mbox{ and } \ \omega(x-1,n)=1 \label{cond:perc1}
\end{eqnarray}
or
\begin{eqnarray}
	x+1\in V_n \ \mbox{ and } \ \omega(x+1,n)=1. \label{cond:perc2}
\end{eqnarray}
For $(x',n) \in V_n \subset X_n$, the definition of $X_n$ means that
$\xi_{nT}(x',0,\ldots,0)=(1,1)$; if $\omega(x',n)=1$, then the ``good event''
occurs in the space-time rectangle
\begin{eqnarray*}
	[x'-1,x'+1] \times \Zbold^{d-1} \times [nT,(n+1)T),
\end{eqnarray*}
which, since $\xi_{nT}(x',0,\ldots,0)=(1,1)$, implies that
\begin{eqnarray*}
	\xi_{(n+1)T}(x'-1,0,\ldots,0) = (1,1) \mbox{ and }
	\xi_{(n+1)T}(x'+1,0,\ldots,0) = (1,1).
\end{eqnarray*} 
Therefore the conditions~(\ref{cond:perc1}) and~(\ref{cond:perc2})
ensure that any $(x,n+1)\in V_{n+1}$ is a member of $X_{n+1}$,
hence $V_{n+1} \subset X_{n+1}$. By induction, for all $n\in \Zbold^+$, we have
\begin{eqnarray*}
	V_n \subset X_n.
\end{eqnarray*}
{
\psfrag{n9}{$n=8$}
\psfrag{n1}{$n=1$}
\psfrag{n0}{$n=0$}
\psfrag{x1}{$(-4,0)$}
\psfrag{x2}{$(-2,0)$}
\psfrag{x3}{$(0,0)$}
\psfrag{x4}{$(2,0)$}
\psfrag{x5}{$(4,0)$}
\psfrag{text1}{$G$ occurs}
\psfrag{text2}{$G$ does not}
\psfrag{text3}{occur}
\begin{figure}[h!]
\includegraphics[height = 3in, width = 5in]{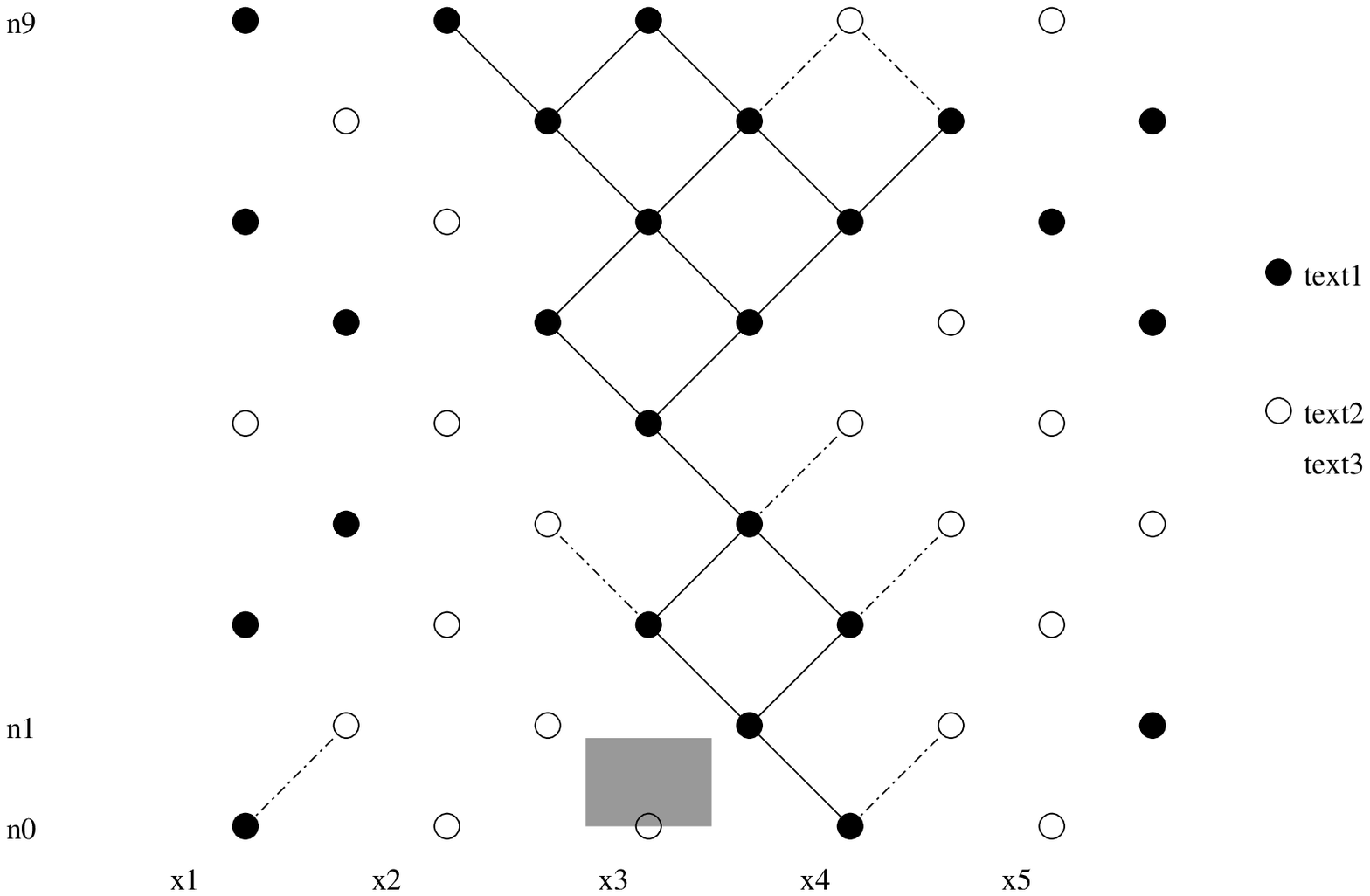}
\caption{Illustration of $V_n$ and $W_n$: all sites in $W_n$ are connected
to $(2,0)$ via a sequence of solid lines, while sites in $V_n\backslash W_n$
are connected to some site in $W_n$ via a dotted line. The shaded rectangle
indicates the space-time region that affects whether $G_{\xi_0}$ occurs}
\label{fig:perc}
\end{figure}
}   

The $\omega(x,n)$ thus defined is a 2-dependent oriented percolation
on $\Lcal_0$ with density at least $1-\gamma$;
notice that $\{\omega(0,0)=0\}$ and $\{\omega(2,0)=0\}$ are dependent events
since both require that no death occur during $[0,T)$ at nest
$((1,0,\ldots,0),1)$ or $((1,0,\ldots,0),2)$,
but $\{\omega(0,0)=0\}$ and $\{\omega(4,0)=0\}$ are clearly
independent. Also, $\{\omega(x_1,n_1)=0\}$ and $\{\omega(x_2,n_2)=0\}$ are
independent for any $x_1$ and $x_2$ provided $n_1 \neq n_2$.
Thus if $(x_i,n_i)$, $1\leq i \leq I$, is a sequence with
$\|(x_i,n_i)-(x_j,n_j)\|_\infty > 2$ for $i\neq j$, then
the events $\{\omega(x_i,n_i) = 0\}$ are all independent, which implies
\begin{eqnarray*}
        P(\omega(x_i,n_i) = 0 \mbox{ for $1 \leq i \leq I$}) \leq \gamma^I,
\end{eqnarray*}
as required by~(\ref{def:Mdep}). Notice that even though $V_{n+1}$ clearly
depends on $V_n$, $\omega(x_1,n+1)$ and $\omega(x_2,n)$ as defined
by~(\ref{def:omega}) are indeed independent since they relate to independent
Poisson arrivals in disjoint space-time rectangles.

Now we define
\begin{eqnarray*}
	W_n = \{(x,n): (x,n)\in V_n \mbox{ and } \omega(x,n)=1\},
\end{eqnarray*}
then $W_n \subset V_n$ and $W_n$ is a 2-dependent oriented percolation with
density at least $1-\gamma$. If $\gamma$ is sufficiently small, then
by Theorem~\ref{thm:perc:rd},
\begin{eqnarray*}
        \liminf_{n\rightarrow \infty} P(0\in W_{2n}) \geq 19/20.
\end{eqnarray*}
Since the particle system $\xi$ dominates $V_n$, which in turn dominates
$W_n$, we have
\begin{eqnarray}
        \liminf_{n\rightarrow \infty} P(\xi_{2nT}(0,0,\ldots,0)=(1,1))
		\geq 19/20. \label{ineq:perc3}
\end{eqnarray}
Thus if we start with initial configuration $\overline\xi_0 (x) = (1,1)$
for all $x$ (in this case, $p=1$ for Theorem~\ref{thm:perc:rd}), then 
by Theorem~\ref{thm:upper_inv},
$\overline\xi_t \Rightarrow \overline\xi_\infty$.
And~(\ref{ineq:perc3}) implies that
\begin{eqnarray*}
	P(\overline \xi_\infty(x)=(1,1)) \geq 19/20 \ \ \
	\mbox{for any $x\in \Zbold^d$},
\end{eqnarray*}
i.e. the upper invariant measure $\overline\xi_\infty$ is nontrivial.
\qed

\vspace{.3cm}
\begin{REM}
Using the same idea of comparing the particle system $\xi$ to an
oriented percolation process, one can show that
if the initial condition is finite, e.g. $\xi_0(0,0,\ldots,0)=(1,1)$ and 0
everywhere else, then for sufficiently large $\lambda$,
\begin{eqnarray*}
        \liminf_{n\rightarrow \infty} P(\xi_{nT}(x,0,\ldots,0)=(1,1)
		\mbox{\ for some \ } x)\geq 19/20
\end{eqnarray*}
for some $T$.
Here we use the following fact (see Theorem 4.1 of [Durrett 1995]) about
$C_0 = \{(y,n): (0,0)\rightarrow (y,n)\}$:
\begin{eqnarray}
	\mbox{If $\gamma\leq 6^{-4(2M+1)^2}$, then }
		P(|C_0| < \infty) \leq 1/20.
\label{ineq:orient_perc}
\end{eqnarray}
\end{REM}

\vspace{.3cm}
\begin{REM}
In fact, the proof of Theorem~\ref{thm:perc:rd} also establishes the following:
if the initial population has a positive density of male and female particles
everywhere and $\lambda/\delta$ is sufficiently large,
then at sufficiently large times, the density of $(1,1)$-sites
(i.e. sites with both a male and female particle) is larger than $9/10$.
Similarly, using~(\ref{ineq:orient_perc}), the following is also true:
if the initial population is nonzero and $\lambda/\delta$
is sufficiently large, then at sufficiently large times, the
probability of survival (i.e. existence of a site with both a male
and a female particle) is at least $9/10$.
\end{REM}

\chapter{Convergence Theorem for Individual Stirring}
In this chapter, we establish the convergence result for the individual
stirring model as promised in
Chapter~\ref{sec:conv:ind_stir}. We work in a slightly more general setting
and consider random processes
\begin{eqnarray*}
	\xi_t^\epsilon: \epsilon \mathbb{Z}^d \times \{1,2,\ldots,M\}
		\rightarrow \{0,1,\ldots,\kappa-1\}.
\end{eqnarray*}
We call each $x \in \epsilon \mathbb{Z}^d$ a site, and each
$(x,m) \in \epsilon \mathbb{Z}^d \times \{1,2,\ldots,M\}$ a nest. There
are $M$ nests at each site. We think of the set of spatial locations
$\mathbb{Z}^d \times \{1,2,\ldots,M\}$ as consisting of $M$ ``floors''
of $\mathbb{Z}^d$. Let
\begin{eqnarray*}
	\Ncal=\{0,\epsilon y_1,\ldots,\epsilon y_N\}
\end{eqnarray*}
be the interaction neighbourhood of site 0 and
$r_\Ncal^\epsilon = \max_{x\in\Ncal} \|x\|_1$ be its radius.
The process $\xi_t^\epsilon$ evolves as follows:

\begin{enumerate}
\item {\bf Birth and Death.} The state of nest $(x,m)$ flips to $i$,
$i=0,\ldots,\kappa-1$, at rate
\begin{eqnarray*}
    	c_i(x,m,\xi)= h_{i,m}(\xi(x,m),\xi(x+\epsilon z_1,m_1),
		\ldots,\xi(x+\epsilon z_L,m_L)),
\end{eqnarray*}
where $L$ is a positive integer, $z_1,\ldots,z_L \in \Ncal$,
$m_1,\ldots,m_L \in \{1,2,\ldots,M\}$,
and \[h_{i,m}:\{0,1,\ldots,\kappa-1\}^{L+1} \rightarrow K \subset \Rbold^+\] with
$h_{i,m}(i,\ldots)=0$ and $K$ compact.

\item {\bf Rapid Stirring.}
For each $m\in \{1,2,\ldots,M\}$ and
$x,y\in \epsilon \mathbb{Z}^d$ with $\|x-y\|_1 = \epsilon$,
$\xi^\epsilon(x,m)$ and $\xi^\epsilon(y,m)$ are exchanged
at rate $\epsilon^{-2}$.
\end{enumerate}

This individual stirring model differs from the Lily-pad stirring model
described in Chapter~\ref{sec:2stir} in that the stirring action between
corresponding nests at neighbouring sites are now independent. More
specifically, exchanges are allowed between neighbouring nests on the same
floor only, i.e. between $(0,1)$ and
$(\epsilon,1)$ but not between $(0,1)$ and $(\epsilon,2)$.
We will show, in the theorem below, that
this individual rapid stirring action between corresponding nests
``decouples'' (in the limit $\epsilon \rightarrow 0$)
the dependence between all nests, at neighbouring sites and
even at the same site.

As an example, for $d=1$, in the particle model with individual stirring
with generator~(\ref{def:genI}), we have $\kappa = 2$, $M=2$, $L=4$,
$\Ncal = \{0,-\epsilon,\epsilon\}$,
\begin{eqnarray*}
	c_0(x,m,\xi) = \left\{ \begin{array}{ll} 
		\delta, & \mbox{if $\xi(x,m)=1$} \\
		0,	& \mbox{otherwise}
	\end{array} \right. ,
\end{eqnarray*}
and
\begin{eqnarray*}
	c_1(x,m,\xi)= \left\{ \begin{array}{ll} 
		\lambda (\xi(x-\epsilon,1)\xi(x-\epsilon,2)
			+ \xi(x+\epsilon,2)\xi(x+\epsilon,1)),
			& \mbox{if $\xi(x,m)=0$} \\
		0,      & \mbox{otherwise}
        \end{array} \right. .
\end{eqnarray*}
In particular, we should define
\begin{eqnarray*}
        (z_1,m_1) = (-1,1), \ (z_2,m_2) = (-1,2), \
        (z_3,m_3) = (1,1), \ (z_4,m_4) = (1,2),
\end{eqnarray*} 
and $h_i=h_{i,m}$ as
\begin{eqnarray*}
        h_0(\alpha_0,\alpha_1,\alpha_2,\alpha_3,\alpha_4)
                = \delta \alpha_0,
\end{eqnarray*}
and
\begin{eqnarray*}
        h_1(\alpha_0,\alpha_1,\alpha_2,\alpha_3,\alpha_4)
                = \lambda (\alpha_1\alpha_2+\alpha_3\alpha_4) (1-\alpha_0),
\end{eqnarray*}
where $\alpha_0 = \xi(x,m)$,
$\alpha_1=\xi(x+\epsilon z_1,m_1)=\xi(x-\epsilon,1)$, \emph{etc}.

\begin{THM}
\label{thm:indp_stir}
	Suppose $\left\{\xi_0^\epsilon (x,m), (x,m)\in
		\epsilon \Zbold^d \times \{1,2,\ldots,M\} \right\}$
	are independent and let
	$u^\epsilon_{i,m} (t,x) = P(\xi^\epsilon_t (x,m) = i)$. If 
	$u_{i,m}^\epsilon (0,x) = g_{i,m} (x) = g_i (x,m)$ and
	$g_i: \Rbold^d \times \{1,2,\ldots,M\} \rightarrow [0,1]$
	is continuous, then for any smooth function $\phi$ with compact
	support, as $\epsilon\rightarrow 0$,
\begin{eqnarray}
	\epsilon^d \sum_{y\in \epsilon\Zbold^d}
		\phi(y)1_{\{\xi^\epsilon_t(y,m)=i\} }
	\rightarrow \int \phi(y) u_{i,m}(t,y) \ dy, \label{conv:part_sys}
\end{eqnarray}
where $u_{i,m} (t,x)$ is the bounded solution of
\begin{eqnarray*}
	\frac{\partial u_{i,m}}{\partial t}=\Delta u_{i,m} + f_{i,m} (u),
    	\mbox{        } u_{i,m} (0,x)=g_{i} (x,m),
\end{eqnarray*}
\begin{eqnarray*}
    	f_{i,m} (u) = \left< c_i (0,m,\xi) 1(\xi(0,m)\neq i)
		\right>_{u} - \sum_{j\neq i}
    	\left< c_j (0,m,\xi) 1(\xi(0,m)=i) \right>_{u},
\end{eqnarray*}
and $\left< \phi(\xi) \right>_u$ denotes the expected value of $\phi(\xi)$
under the product measure in which state $j$ at nest $m$ has density
$u_{j,m}$, i.e. $\xi(x,m)$, with $x\in \epsilon\Zbold^d$
and $1\leq m\leq M$, are independent with $P(\xi(x,m)=j)=u_{j,m}$.
\end{THM}

\proof
We follow the program used in the proof of Theorem 8.1 in [Durrett 1995].
We first define a dual process for the particle system in part (a), then
in part (b) we show that the dual process is almost a branching random
walk. We will not explicitly write down parts (c) and (d) of the proof,
which almost exactly resemble parts (c) and (d) of the proof of
Theorem 8.1 in [Durrett 1995]. But to summarize these two parts, part
(c) establishes that the dual process converges to a branching Brownian
motion as $\epsilon \rightarrow 0$ and defines a candidate limit
$u_{i,m}(t,x)$ of $u^\epsilon_{i,m}(t,x)$, and part (d) shows that this candidate limit
satisfies the PDE in the statement of the theorem. Finally, part (e), which we write
out explicitly, shows that convergence described in~(\ref{conv:part_sys}) does occur,
in addition to the convergence of the mean $u^\epsilon_{i,m}(t,x)$ to $u_{i,m}(t,x)$.

{\bf a. Defining the dual process.}
The dual process associated with nest $(\hat x,\hat m)$
and a fixed time $t$ is a random process
$I_\epsilon^{\hat{x},\hat{m},t}(s)$, $s \in [0,t]$, where
\begin{eqnarray}
	I_\epsilon^{\hat{x},\hat{m},t}(s) \in
	\bigcup_{k\in\Zbold^+} \{(x,m): x\in \epsilon \Zbold^d,
		m\in \{1,2,\ldots,M\} \}^k \label{def:IX}
\end{eqnarray}
consists of a finite number of particles residing at nests
$(x,m) \in \epsilon \Zbold^d \times  \{1,2,\ldots,M\}$.
The events that influence the behaviour of the dual process
$I_\epsilon^{\hat{x},\hat{m},t}$ at time $s$ are those Poisson arrivals
in the original process $\xi^\epsilon$ that occur at time $t-s$.
We start with
\begin{eqnarray*}
	I_\epsilon^{\hat{x},\hat{m},t}(0)=\{(\hat x,\hat m)\}
\end{eqnarray*}
and evolve $I_\epsilon^{\hat{x},\hat{m},t}(s)$ until $s=t$, which corresponds
to rolling back the clock in the original process $\xi^\epsilon$ from time $t$
back to time $0$.
The process $I_\epsilon^{\hat{x},\hat{m},t}$ is constructed such  
that: in order to know the value of $\xi^\epsilon(\hat{x},\hat{m})$
at time $t$, it suffices to do a computation using values of $\xi^\epsilon$
at time $0$ and nests $(x,m) \in I_\epsilon^{\hat{x},\hat{m},t}(t)$.
We call $I_\epsilon^{\hat{x},\hat{m},t}(s)$ the \emph{influence set} of
$I_\epsilon^{\hat{x},\hat{m},t}(0)$.

{
\psfrag{m=0}{$m=0$}
\psfrag{m=1}{$m=1$}
\psfrag{x}{$x$}
\psfrag{t}{$t$}
\begin{figure}[h!]
\centering       
\includegraphics[height = 6.6in, width = 5in]{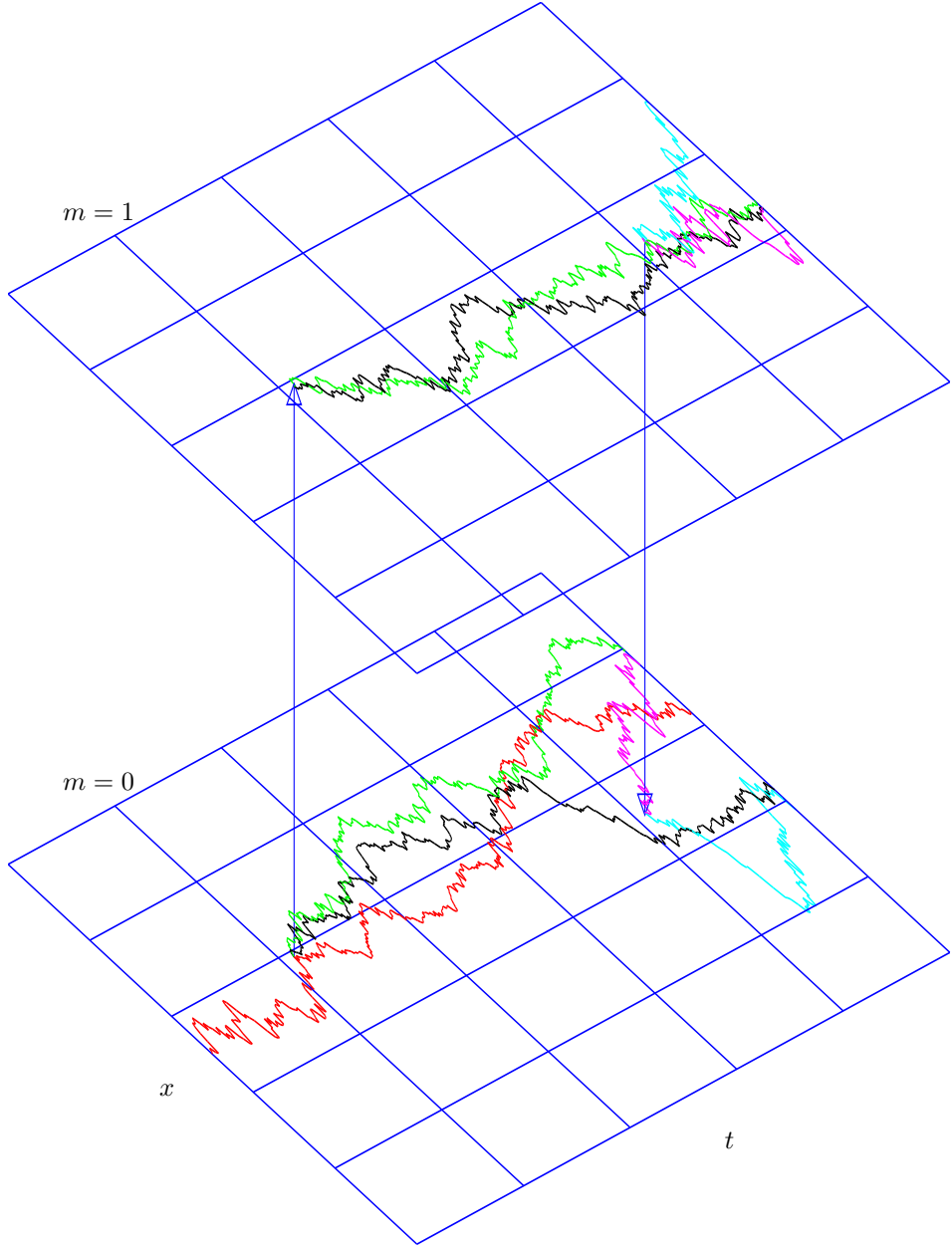}
\caption{Illustration of the dual process: there are two birth/death events
here, each giving birth to four additional particles, two on each floor;
the first birth event is from a male ($m=0$) nest, while the second
birth event is from a female ($m=1$) nest.}
\label{fig:dual}
\end{figure}                                       
}

In order to define the dual process, we need to give a graphical
construction for the particle system $\xi^\epsilon$
similar to the one given in Chapter~\ref{sec:model:static}.
We define a family of uniform $[0,1]$ random variables
$\{U_n^{x,m,i}:n\geq 1 \}$ and two families of independent Poisson processes
that respectively correspond to the birth/death flips and the rapid stirring
mechanism:
\begin{eqnarray*}
	& & \{T_n^{x,m,i}:n\geq 1 \} \ \mbox{at rate
		$c^* = \sup_{\xi,m} \sum_i c_i (x,m,\xi)$} \\
	& \mbox{ and } & \{S_n^{x,y,m}:n\geq 1\} \ \mbox{at rate
		$\epsilon^{-2}$},
\end{eqnarray*}
with $m \in \{1,2,\ldots,M\}$ and $\|x-y\|_1 = \epsilon$.
Notice that $c^* < \infty$ since $\max_{i,m} \|h_{i,m}\|_\infty < \infty$.
At $T_n^{x,m,i}$, we check all nests $(x,m)$,
$(x+\epsilon z_1,m_1),\ldots,(x+\epsilon z_L,m_L)$ and use $U_n^{x,m,i}$ to
decide whether nest $(x,m)$ should flip to state $i$ based on the
function $h_{i,m}$. And at $S_n^{x,y,m}$, we exchange
values of $\xi^\epsilon(x,m)$ and $\xi^\epsilon(y,m)$,
the corresponding nests at two neighbouring sites.
The evolution of the dual process $I_\epsilon^{\hat{x},\hat{m},t}(s)$
depends on the Poisson arrival times $\{T_n^{x,m,i}:n\geq 1 \}$
and $\{S_n^{x,y,m}:n\geq 1\}$ in the following way:
\begin{enumerate}
\item If $(x,m) \in I_\epsilon^{\hat{x},\hat{m},t}(s-)$ and
	$T_n^{x,m,i}=t-s$, then we set
\begin{eqnarray}
	I_\epsilon^{\hat{x},\hat{m},t}(s) = I_\epsilon^{\hat{x},\hat{m},t}(s-)
		\cup \{(x+\epsilon z_1,m_1),
			\ldots,(x+\epsilon z_L,m_L)\}. \label{eq:dual_birth}
\end{eqnarray}
	Therefore each particle
	$(x,m) \in I_\epsilon^{\hat{x},\hat{m},t}(s)$
	gives birth to $L$ new particles at rate $c^*$.

\item If $(x,m)\in I_\epsilon^{\hat{x},\hat{m},t}(s-)$ and either
	$S_n^{x,y,m}=t-s$ or $S_n^{y,x,m}=t-s$, then we set
\begin{eqnarray*}
	I_\epsilon^{\hat{x},\hat{m},t}(s) = \left(
		I_\epsilon^{\hat{x},\hat{m},t}(s-) \backslash
			\{(x,m)\} \right) \cup \{(y,m)\}.
\end{eqnarray*}
\end{enumerate}
Thus two dual processes $I_\epsilon^{\hat{x},\hat{m},t}(s)$             
and $I_\epsilon^{\hat{x}',\hat{m}',t}(s)$
evolve independently from each other except at time
\begin{eqnarray*}
	s=\mbox{ some } T^{x,m,i}_n
\end{eqnarray*}
where $(x,m)\in I_\epsilon^{\hat{x},\hat{m},t}(s)
\cap I_\epsilon^{\hat{x}',\hat{m}',t}(s)$ or at time
\begin{eqnarray*}
        s=\mbox{ some } S^{x,y,m}_n
\end{eqnarray*}
where $(x,m)\in I_\epsilon^{\hat{x},\hat{m},t}(s)$ and
$(y,m)\in I_\epsilon^{\hat{x}',\hat{m}',t}(s)$, or
$(y,m)\in I_\epsilon^{\hat{x},\hat{m},t}(s)$ and
$(x,m)\in I_\epsilon^{\hat{x}',\hat{m}',t}(s)$ .

An equivalent way of describing the dual process is to define the pair
\\
$(X_\epsilon^{\hat{x},\hat{m}}(s), \Kcal^{\hat{x},\hat{m}}_\epsilon(s))$
for $0\leq s \leq t$, where $X_\epsilon^{\hat{x},\hat{m}}(s)$
is the ordered set
\begin{eqnarray*}
        X_\epsilon^{\hat{x},\hat{m}}(s)=(X_\epsilon^{\hat{x},\hat{m},0}(s),
            \ldots,X_\epsilon^{\hat{x},\hat{m},
			\Kcal^{\hat{x},\hat{m}}_\epsilon(0)-1}(s))
\end{eqnarray*}
and each $X_\epsilon^{\hat{x},\hat{m},j}(s) \in \epsilon \Zbold^d
\times \{1,2,\ldots,M\}$.
For $s=0$, define the number of particles
\begin{eqnarray*}
	\Kcal^{\hat{x},\hat{m}}_\epsilon(0)=n
\end{eqnarray*}
for some $n\in \Zbold^+$ and the location of particles
\begin{eqnarray*}
	X_\epsilon^{\hat{x},\hat{m},j}(0)=(x_j,m_j)
\end{eqnarray*} 
for $0\leq j <n$, such that
\begin{eqnarray}
	I_\epsilon^{\hat{x},\hat{m},t}(0)
		= \bigcup_{j=0}^{\Kcal^{\hat{x},\hat{m}}_\epsilon(0)-1}
			\{ X_\epsilon^{\hat{x},\hat{m},j}(0) \},
	\label{eq:IX}
\end{eqnarray}
where $I_\epsilon^{\hat{x},\hat{m},t}(0)$ is defined in~(\ref{def:IX}).
Typically, $\Kcal^{\hat{x},\hat{m}}_\epsilon(0)=1$ and
$X_\epsilon^{\hat{x},\hat{m},0}(0)=\{ (\hat{x},\hat m) \}$, but
we use the more general initial condition since we will define
$X_\epsilon^{\hat{x},\hat{m},j}$ inductively. Each particle
in $X_\epsilon^{\hat{x},\hat{m},j}$,
$0 \leq j < \Kcal^{\hat{x},\hat{m}}_\epsilon$, jumps to a neighbouring nest
on the same ``floor'' of $\epsilon \Zbold^d \times \{1,2,\ldots,M\}$
%performs random walk on a ``floor''
as dictated by the stirring mechanism (i.e. $S$-arrivals),
until a $T$-arrival of type $T^{x,m,j}$ with
$0 \leq j < \Kcal^{\hat{x},\hat{m}}_\epsilon(s-)$ at some
$(x,m) \in X_\epsilon^{\hat{x},\hat{m}}(s-)$.
At time $s$ of this $T$-arrival, we set
\begin{eqnarray*}
	\Kcal^{\hat{x},\hat{m}}_\epsilon(s)
		=L+\Kcal^{\hat{x},\hat{m}}_\epsilon(s-)
\end{eqnarray*}
and
\begin{eqnarray*}
	X_\epsilon^{\hat{x},\hat{m},k+\Kcal^{\hat{x},\hat{m}}_\epsilon(s-)}(s)
		= (x+\epsilon z_k,m_k)
\end{eqnarray*}
for $0 \leq k < L$, while leaving all
existing $X_\epsilon^{\hat{x},\hat{m},j}$'s unchanged. Observe that
\begin{eqnarray*}
	I_\epsilon^{\hat{x},\hat{m},t}(s)
	=\bigcup_{j=0}^{\Kcal^{\hat{x},\hat{m}}_\epsilon(s)-1}
		\{ X_\epsilon^{\hat{x},\hat{m},j}(s) \},
\end{eqnarray*}
therefore the relation~(\ref{eq:IX}) is maintained for all $s\in [0,t]$.
Afterwards, the $\Kcal^{\hat{x},\hat{m}}_\epsilon(s)$
particles jump according to the $S$-arrivals until the next $T$-arrival
at some $(x,m) \in \{X_\epsilon^{\hat{x},\hat{m},j}(s'),
0 \leq j < \Kcal^{\hat{x},\hat{m}}_\epsilon(s'-)\}$.

A new particle may be born at a nest where a
particle already resides; if this happens, we say that a
\emph{collision} occurs, and call the new particle \emph{fictitious}.
We prescribe this fictitious particle to give birth to $L$ particles
at rate $c^*$ and jump to a
neighbouring nest on the same ``floor'' at rate $\epsilon^{-2}$,
independently from all other particles, fictitious or not; furthermore, all
offsprings of a fictitious particle are defined to be fictitious as well.
If the number of collisions is 0, then there are no fictitious particles.
As will be seen in the next paragraph , this is in fact the case in the limit
$\epsilon \rightarrow 0$.
Note that the stirring mechanism does not cause particles at different nests
to ``collide'', i.e. end up at the same nest, because only exchanges
between nests occur under the stirring mechanism.

{\bf b. Characterizing the dual process.}
Having finished defining the dual process, we proceed to show that
with high probability,
the dual process can be coupled to a branching random walk and the movement
of one dual process is independent from the movement of another dual.
This part of the proof is quite similar to
part (b) of the proof of Theorem 8.1 of [Durrett 1995].

First, we couple $X_\epsilon^{\hat x,\hat m,j}$ to independent random
walks $Y_\epsilon^{\hat x,\hat m,j}$ that start at the same location
at the time of birth of $X_\epsilon^{\hat x,\hat m,j}$ and jump to
a randomly chosen neighbour at rate $2d \epsilon^{-2}$. We define
the distance between two particles on two different ``floors'' to
only depend on the site location:
\begin{eqnarray*}
        \| (x_1,\hat m_1) - (x_2,\hat m_2) \|_1
                = \|x_1 - x_2 \|_1,           
\end{eqnarray*}
and say $X_\epsilon^{\hat x,\hat m,j}$ is \emph{crowded} if for
some $k\neq j$,
$\|X_\epsilon^{\hat x,\hat m,j}-X_\epsilon^{\hat x,\hat m,k}\|_1 \leq
r_\Ncal^\epsilon$. Recall that $r_\Ncal^\epsilon$ is the radius of
the interaction neighbourhood $\Ncal$.
When $X_\epsilon^{\hat x,\hat m,j}$ is not crowded, we define the
displacements of $Y_\epsilon^{\hat x,\hat m,j}$ to be equal to those of
$X_\epsilon^{\hat x,\hat m,j}$. But when $X_\epsilon^{\hat x,\hat m,j}$
is crowded, we use independent Poisson processes to determine the jumps of
$Y_\epsilon^{\hat x,\hat m,j}$. To estimate the difference between
$X_\epsilon^{\hat x,\hat m,j}$ and $Y_\epsilon^{\hat x,\hat m,j}$, we
need to estimate the amount of time $X_\epsilon^{\hat x,\hat m,j}$ is crowded.
If $j\neq k$, then
$V^\epsilon_s=X_\epsilon^{\hat x,\hat m,j}(s)-X_\epsilon^{\hat x,\hat m,k}(s)$
is stochastically larger than $W^\epsilon_s$, a random walk that jumps
to a randomly chosen neighbour in $\epsilon \Zbold^d \times \{1\}$
at rate $4d\epsilon^{-2}$ (see
page 175 of [Durrett 1995] for details on how to couple these two processes
so that one is stochastically larger than the other). Strictly speaking,
$V^\epsilon_s$ is only defined after the $j^\mathrm{th}$ and $k^\mathrm{th}$
particles of $X^{\hat x,\hat m}_\epsilon$ are created, but as we shall see
below, we are only interested in an upper bound on the occupation time
of $V^\epsilon_s$ and $W^\epsilon_s$ in a ball, so
the fact that these particles may only start
to exist at some positive time only helps matters.
Thus for any integer $M\geq 1$,
$v^{M\epsilon}_t = |\{s\leq t: \|V^\epsilon_s\|_1 \leq M \epsilon\}|$
is stochastically smaller than
$w^{M\epsilon}_t = |\{s\leq t: \|W^\epsilon_s\|_1 \leq M \epsilon\}|$.
Well known asymptotic results [Durrett 1995] for random walks imply that
if $t\epsilon^{-2}\geq 2$ then
\begin{eqnarray}
	E w^{M\epsilon}_t \leq \left\{ \begin{array}{ll}
		C M \epsilon t^{1/2} & d=1 \\
		CM^2 \epsilon^2 \log(t\epsilon^{-2}) & d=2 \\
		CM^d \epsilon^2 & d\geq 3
	\end{array} \right. . \label{ineq:Wocc_time}
\end{eqnarray}
In the estimates above, $d=1$ is the worst case, so we use
$E w^{M\epsilon}_t \leq C M \epsilon t^{1/2}$ if $t\epsilon^{-2}\geq 2$
from now on.

The amount of time $X_\epsilon^{\hat x,\hat m,j}$ is crowded
in $[0,t]$, denoted $\chi_\epsilon^j(t)$, can be estimated as follows:
\begin{eqnarray*}
	E[\chi_\epsilon^j(t)]
	& = & \sum_{K=0}^\infty E[\chi_\epsilon^j(t)|
		\Kcal^{\hat x,\hat m}_\epsilon(t) = K]
		P(\Kcal^{\hat x,\hat m}_\epsilon(t) = K) \\
	& \leq & \sum_{K=0}^\infty K E[w^{M\epsilon}_t]
                P(\Kcal^{\hat x,\hat m}_\epsilon(t) = K),
\end{eqnarray*} 
where we pick $M$ large enough so that $M \epsilon \geq r_\Ncal^\epsilon$,
e.g. $M = \frac{1}{\epsilon} \max_{0<i \leq N} \|y_i\|_1$.
In what follows, $C_{\ldots}$ is a constant whose value may change from
line to line. Then
\begin{eqnarray*}
	E[\chi_\epsilon^j(t)]
	\leq E[w^{M\epsilon}_t] E[\Kcal^{\hat x,\hat m}_\epsilon(t)]
	= e^{c^* L t} E[w^{M\epsilon}_t]
	\leq C_M e^{c^* L t} \epsilon t^{1/2}
\end{eqnarray*}
if $t\epsilon^{-2}\geq 2$. To see the middle equality above, we observe
that the branching mechanism of the dual process described
in~(\ref{eq:dual_birth}) occurs at rate $c^*$ for both fictitious and
nonfictitious particles, and every time a branching event occurs, one particle
is replaced by $L$ particles; therefore the mean number of branches at time $t$
is $e^{c^* L t}$. It follows that
\begin{eqnarray}
	E[\chi_\epsilon^j(t)]
	\leq  C_M e^{c^* L t} \epsilon t^{1/2} + 2\epsilon^2
	\leq C_M e^{c^* L t} \epsilon (1+t^{1/2}). \label{ineq:chi}
\end{eqnarray}
This means that the expected number of births from
$X_\epsilon^{\hat x,\hat m,j}$ while there is some other
$X_\epsilon^{\hat x,\hat m,k}$ in $\Ncal+X_\epsilon^{\hat x,\hat m,j}$
is smaller than
\begin{eqnarray}
	C_{L,M} e^{c^* L t} \epsilon (1+t^{1/2}). \label{ineq:chi2}
\end{eqnarray}
Thus
\begin{eqnarray}
	\lefteqn{ P(\mbox{at least one collision during $[0,t]$}) } \nnb \\
	& \leq & P(\mbox{at least one collision during $[0,t]$} |
		\Kcal^{\hat x,\hat m}_\epsilon(t) \leq \epsilon^{-0.5})
		P(\Kcal^{\hat x,\hat m}_\epsilon(t) \leq \epsilon^{-0.5})
		\nnb \\
	& & \ \ \ \ \
		+ P(\Kcal^{\hat x,\hat m}_\epsilon(t) > \epsilon^{-0.5}) \nnb \\
	& \leq & E[\mbox{\# of collisions during $[0,t]$} |
                \Kcal^{\hat x,\hat m}_\epsilon(t) \leq \epsilon^{-0.5}]
		P(\Kcal^{\hat x,\hat m}_\epsilon(t) \leq \epsilon^{-0.5})
		\nnb \\
	& & \ \ \ \ \	+ \epsilon^{0.5} E[\Kcal^{\hat x,\hat m}_\epsilon(t)]
		\label{ineq:1coll} \\
	& \leq & C_{L,M} e^{c^* L t} \epsilon^{0.5} (1+t^{1/2})
		+ \epsilon^{0.5} e^{c^* L t} \nnb \\
	& \leq & C_{L,M} e^{c^* L t} \epsilon^{0.5} (1+t^{1/2}), \nnb
\end{eqnarray}
which $\rightarrow 0$ as $\epsilon\rightarrow 0$.
We use~(\ref{ineq:chi2}) and the
condition $\Kcal^{\hat x,\hat m}_\epsilon(t) \leq \epsilon^{-0.5}$
to bound the first expectation in~(\ref{ineq:1coll}),
and we bound $P(\Kcal^{\hat x,\hat m}_\epsilon(t) \leq \epsilon^{-0.5})$
in~(\ref{ineq:1coll}) above by $1$.

This shows that the probability of at least one collision within a single dual
during time $[0,t]$ tends to $0$ as $\epsilon\rightarrow 0$.
Furthermore, the same argument shows that the probability of at least one
collision between
two different duals for nests $(\hat x,\hat m) \neq (\hat x',\hat m')$
during time $[0,t]$ also tends to $0$ as $\epsilon\rightarrow 0$. For that,
we observe that the estimate~(\ref{ineq:Wocc_time}) is independent
of the initial condition $W^\epsilon_0$; so in particular,
this estimate still holds even if $\|W^\epsilon_0\|_1 = 0$,
which is the case when for example one
considers two duals for two nests $(\hat{x},\hat{m})$
and $(\hat{x},\hat{m}')$ at the same site $x$.
Hence two different duals are asymptotically independent in the limit
$\epsilon \rightarrow 0$. 

The estimate~(\ref{ineq:chi}) also leads to the following estimate
on the difference between
$X_\epsilon^{\hat x,\hat m,j}$ and $Y_\epsilon^{\hat x,\hat m,j}$ (see
page 176 in [Durrett 1995] for details):
\begin{eqnarray}
	P\left( \max_{0\leq s \leq t} \|X_\epsilon^{\hat x,\hat m,j}(s)
		- Y_\epsilon^{\hat x,\hat m,j}(s)|_\infty
		\geq 2 \epsilon^{0.3} \right)
	\leq C \epsilon^{0.4} (1+t^{1/2}) e^{c^* L t}.
	\label{ineq:XY}
\end{eqnarray}
This shows that with high probability, the movements of all the particles
in a dual can be coupled to independent random walks, in addition to being
independent from the movement of any other dual.

{\bf c. and d. Defining a candidate limit and showing the limit satisfies the
PDE.}
We will not write down the details of these two parts of the argument, since
they are almost exactly the same as parts (c) and (d) of the proof of
Theorem 8.1 in [Durrett 1995]. From estimate~(\ref{ineq:XY}),
it is not too difficult
to see that the dual process converges to the branching Brownian motion
as $\epsilon\rightarrow 0$. We can then define the candidate limit
$u_{i,m}(t,x)$ (of $u^\epsilon_{i,m}(t,x)=P(\xi^\epsilon_t(x,m)=i)$)
using the limiting branching Brownian motion
as the dual process. Part (d) then establishes that the candidate limit
$u_{i,m}(t,x)$ satisfies the integral from of the PDE in the statement of
the theorem.
%$u_{i,m}(t,x) = P(\zeta_t(0)=i)$ for a random process $\zeta$ with,
%we use the limiting branching Brownian motion as the dual process and run
%it back in time until $t=0$ to find out which sites

{\bf e. The particle systems converge.}
So far we have established that
\begin{eqnarray*}
	u^\epsilon_{i,m}(t,x) \rightarrow u_{i,m}(t,x),
\end{eqnarray*}
i.e. the expected value converges. It remains to
establish~(\ref{conv:part_sys}). For this, we define for a bounded function
$\phi$ with compact support $\mathrm{supp}(\phi)$ of diameter $D$,
\begin{eqnarray*}
	\left<\xi^\epsilon_t, \phi \right>
	= \epsilon^d \sum_{y\in \epsilon\Zbold^d}
                \phi(y)1_{\{\xi^\epsilon_t(y,m)=i\} }.
\end{eqnarray*}
Then
\begin{eqnarray}
	E[ \left<\xi^\epsilon_t, \phi \right>]
	& = & \epsilon^d \sum_{y\in \epsilon\Zbold^d}
                \phi(y) P(\xi^\epsilon_t(y,m)=i)
	= \epsilon^d \sum_{y\in \epsilon\Zbold^d} \phi(y) u^\epsilon_{i,m}(t,y)
		\nnb \\
        & \rightarrow & \int \phi(y) u_{i,m}(t,y) \ dy
		\label{conv:part_sys:mean}
\end{eqnarray}
by bounded convergence. Now we compute the variance of
$\left<\xi^\epsilon_t, \phi \right>$:
\begin{eqnarray*}
	\lefteqn{\var[\left<\xi^\epsilon_t, \phi \right>]} \\
	& = & E\left[\epsilon^{2d} \left( \sum_{y\in \epsilon\Zbold^d} \phi(y)
		(1_{\{ \xi^\epsilon_t(y,m)=i \} } - P(\xi^\epsilon_t(y,m)=i))
		\right)^2 \right] \\
	& = & E\left[\epsilon^{2d} \sum_{y\in \epsilon\Zbold^d} \phi(y)^2
                (1_{\{ \xi^\epsilon_t(y,m)=i \} } - P(\xi^\epsilon_t(y,m)=i))^2
		\right]
		+ E\left[ \epsilon^{2d} \sum_{y,z\in \epsilon\Zbold^d, y\neq z}
		\phi(y)\phi(z) \right. \\
	& & \ \ \ \ \ \times \left.
		(1_{\{ \xi^\epsilon_t(y,m)=i \} } - P(\xi^\epsilon_t(y,m)=i))
		(1_{\{ \xi^\epsilon_t(z,m)=i \} } - P(\xi^\epsilon_t(z,m)=i))
                \right] \\
	& \leq & D^{2d} \|\phi\|_\infty^2 \sup_{
			y,z \in \epsilon\Zbold^d \cap \mathrm{supp}(\phi),
			 y\neq z } \cov[
		1_{\{ \xi^\epsilon_t(y,m)=i \} },
		1_{\{ \xi^\epsilon_t(z,m)=i \} } ] \\
	& & \ \ \ \ \
		+ \epsilon^{2d} \|\phi\|_\infty^2 \sum_{y\in \epsilon\Zbold^d
			\cap \mathrm{supp}(\phi)}
		\var[1_{\{ \xi^\epsilon_t(y,m)=i \} }].
\end{eqnarray*}
We observe that $1_{\{ \xi^\epsilon_t(y,m)=i \} }$ is a random variable
taking values in $\{0,1\}$ and therefore has variance $\leq 1/4$. Also,
part (b) of the proof shows that
\begin{eqnarray*}
	\mathrm{cov}_\epsilon = \sup_{y,z \in \epsilon\Zbold^d, y\neq z}
	\cov[ 1_{\{ \xi^\epsilon_t(y,m)=i \} },
	1_{\{ \xi^\epsilon_t(z,m)=i \} } ] \rightarrow 0
\end{eqnarray*} 
as $\epsilon\rightarrow 0$. The argument leading to the asymptotic independence
of two duals in part (b) works for any two nests, so
$\cov[ 1_{\{ \xi^\epsilon_t(y,m)=i \} }, 1_{\{ \xi^\epsilon_t(z,m)=i \} } ]$
goes to zero uniformly for all $y\neq z$.
Now we have the following estimate on
$\var[\left<\xi^\epsilon_t, \phi \right>]$:
\begin{eqnarray*}
        \var[\left<\xi^\epsilon_t, \phi \right>]
	\leq D^{2d} \|\phi\|_\infty^2 \mathrm{cov}_\epsilon
		+ \frac{D^d}{4} \epsilon^d \|\phi\|_\infty^2,
\end{eqnarray*}
which $\rightarrow 0$ as $\epsilon\rightarrow 0$.
Thus by~(\ref{conv:part_sys:mean}) and Chebyshev's inequality, we have
\begin{eqnarray*}
	\lefteqn{ P\left(\left|\left<\xi^\epsilon_t, \phi \right>
		- \int \phi(y) u_{i,m}(t,y) \ dy\right|
		>\delta\right) } \\
	& \leq & P\left(\left|\left<\xi^\epsilon_t, \phi \right>
		- \epsilon^d \sum_{y\in \epsilon\Zbold^d} \phi(y)
			u^\epsilon_{i,m}(t,y) \right| \right. \\
	& & + \left. \left| \epsilon^d \sum_{y\in \epsilon\Zbold^d} \phi(y)
			u^\epsilon_{i,m}(t,y)
		- \int \phi(y) u_{i,m}(t,y) \ dy\right|
                >\delta\right) \\
	& \leq & P\left(\left|\left<\xi^\epsilon_t, \phi \right>
                - \epsilon^d \sum_{y\in \epsilon\Zbold^d} \phi(y) u^\epsilon_{i,m}(t,y) \right|
		> \frac{\delta}{2} \right) \\
        & & \ \ \ \ \
		+ P\left( \left| \epsilon^d \sum_{y\in \epsilon\Zbold^d} \phi(y) u^\epsilon_{i,m}(t,y)
                - \int \phi(y) u_{i,m}(t,y) \ dy\right|
                > \frac{\delta}{2} \right) \\
	& \leq & \frac{4 \var [\left<\xi^\epsilon_t, \phi \right>]}{\delta^2}
	+ P\left( \left| \epsilon^d \sum_{y\in \epsilon\Zbold^d} \phi(y) u^\epsilon_{i,m}(t,y)
                - \int \phi(y) u_{i,m}(t,y) \ dy\right|
                > \frac{\delta}{2} \right)
	\rightarrow 0
\end{eqnarray*}            
as $\epsilon\rightarrow 0$, and the theorem follows.
\qed

\chapter{Existence of Invariant Stationary Distribution For Lily-pad Stirring}
\label{section2}
In this chapter,
we establish the existence of nontrivial stationary distribution of the
particle system with lily-pad stirring as promised by Theorem~\ref{thm:lily}.
First, we rewrite~(\ref{eq:pde0}) in the statement of Lemma~\ref{lem:monotone}:
\begin{eqnarray}
        \frac{\partial u}{\partial t} & = & \Delta u +
        (2c(1-u)+1)v-u \nonumber \\
        \frac{\partial v}{\partial t} & = & \Delta v +
        (c(u-v)-2)v. \label{eq:pde:real}
\end{eqnarray}
We show that for sufficiently large $c$, the solution to~(\ref{eq:pde:real})
with initial condition $u_0=f, v_0=g$, $f\geq g$
satisfies the following condition:

\vspace{.3cm} \hspace{-.7cm}
\label{page:star}
($*$) There are constants $0<D_1<d_1<d_2<D_2<1$, $L$, and $T$ so that if
$v_0(x)\in (D_1,D_2)$ for $x\in [-L,L]$ then $v_T(x) \in (d_1,d_2)$ for
$x\in [-3L,3L]$.

\vspace{.3cm} \hspace{-.7cm}
According to Chapter 9 of [Durrett 1995], this is a sufficient condition
for the existence of nontrivial invariant stationary distribution
for the particle system with sufficiently fast stirring, so
Theorem~\ref{thm:lily} will follow once condition ($*$) is established.
Recall that
Theorem~\ref{thm:use:comp} establishes that the diploid particle model
without rapid stirring has a nontrivial stationary distribution
if the birth rate $\lambda$ is sufficiently large. If one traces through
the proof, however, one will find that ``sufficiently large'' in that
argument means that $\lambda$ is larger than a number on the order of
$6^{100}$, which is not too informative on where exactly the critical
$\lambda$ for the phase transition is. On the other hand, one can get a far
better idea of exactly for what $\lambda$ condition ($*$) holds.

For this proof, we also establish condition~($*$) for sufficiently large
$c$ (recall that $c=\lambda d$), but here ``sufficiently large''
means that $c$ is ``only'' larger than a number on the order of $100$.
We assume dimension $d=1$; extension to $d>1$ is straightforward. The proof
consists of two parts: the first part, Chapter~\ref{sec:lb},
establishes the existence of constants
$d_1$ and $D_1$, and the second part, Chapter~\ref{sec:ub}, establishes the
existence of constants $d_2$ and $D_2$; the second part will be easy once
the first part has been established.

Theorem 9.2 in [Durrett 1995] establishes condition ($*$) for a specific
predator-prey system with phase space $\{0,1,2\}$ at each site. The critical
fact used in the proof is that the associated ODE system (i.e. the dynamical
system that results when one has constant initial conditions) has only one
interior equilibrium point and has a \emph{global} Lyapunov function.
The phase portrait of the ODE
associated with~(\ref{eq:pde:real}), however, shows that it has
two interior equilibrium points, one of which is always a saddle point.
See figure~\ref{fig:ode} for two examples. Thus it does not look likely
that the ODE system associated with~(\ref{eq:pde:real}) has a
readily identifiable global Lyapunov function.
%Furthermore, because the ODE system is 2-dimensional,
%the ``traveling wave'' method used in the Appendix of
%[Durrett and Neuhauser 1994] may also be difficult to apply in this case.
The method we use to establish condition~($*$) for~(\ref{eq:pde:real})
is \emph{ad hoc}, but does seem to apply to a wide variety of
reaction-diffusion systems where the reaction part of
the system is 2-dimensional (or even 3-dimensional), i.e.
$\frac{\partial u}{\partial t} = \Delta u + f(u)$ where
$f:\Rbold^2 \rightarrow \Rbold^2$.

As established in Lemma~\ref{lem:monotone}, the PDE~(\ref{eq:pde:real})
is monotone in initial conditions  that lie in
$\R=\{(u,v):0\leq u,v\leq 1,u\geq v\}$.
This fact is critical for the proof of existence of constants $d_1$ and $D_1$
in condition ($*$) above -- it enables us
to bound the initial condition $v_0 (x)$ below by a function,
say $\underbar{$v$}_0 (x)$, both $v_0$ and $\underbar{$v$}_0$ having
values $< D_1$ for $x\in [-L,L]$, such that if the solution
$\underbar{$v$}_t (x)$ to~(\ref{eq:pde:real}) with
initial condition $\underbar{$v$}_0 (x)$ satisfies
the condition
\begin{eqnarray*}
	\underbar{$v$}_T (x) > d_1 \ \forall x\in [-3L,3L],
\end{eqnarray*}
then $v_T (x) > d_1$ $\forall x\in [-3L,3L]$ as well.
We will also need results regarding the ODE associated
with~(\ref{eq:pde:real}):
\begin{eqnarray}
        \frac{du}{dt} & = & (2c(1-u)+1)v-u \nonumber \\
        \frac{dv}{dt} & = & (c(u-v)-2)v. \label{eq:ode}
\end{eqnarray}
The above ODE system is also monotone in initial conditions, since if
the initial condition for the PDE system in~(\ref{eq:pde:real}) is constant
in $x$, then the solution $(u_t,v_t)$ to~(\ref{eq:pde:real})
also remains constant in $x$ for all time and therefore satisfies the ODE
system in~(\ref{eq:ode}).

\section{Lower Bounds: Existence of $d_1$ and $D_1$ in Condition ($*$)}
\label{sec:lb}
First we recall the definition of the region $\R$ from~(\ref{def:R:0}):
\begin{eqnarray*}
	\R=\{(u,v):0\leq u,v\leq 1,u\geq v\}.
\end{eqnarray*}
If the initial condition $(u_0,v_0)$ lies in $\R$, then
so does the solution $(u_t,v_t)$.
We will establish the existence of constants $d_1$, $D_1$, $L$, and $T$
using the nonlinear Trotter product formula (Proposition 15.5.2 
from [Taylor 1996]):
\begin{eqnarray}
        (u_t,v_t) = \lim_{n\rightarrow \infty} \left(
        e^{(t/n)\Delta} {\cal F}^{t/n} \right)^n (f,g).
	\label{eq:trotter}
\end{eqnarray}
Here the convergence occurs in the space ${\cal BC}^1(\mathbb{R})$, the
space of functions whose first derivatives are bounded and continuous on
$\mathbb{R}$ and extend continuously to the compactification
$\widehat{\mathbb{R}}$ via the point at infinity; the norm used here is
$\|\cdot\|_\infty+\|\frac{\partial}{\partial x}(\cdot)\|_\infty$.
In~(\ref{eq:trotter}),
$e^{s\Delta}(f,g)$ gives the (independent) evolution of $f$ and $g$
for time $s$ according the heat equation
\begin{eqnarray*}
	\frac{\partial u}{\partial t} & = & \Delta u \\
        \frac{\partial v}{\partial t} & = & \Delta v,
\end{eqnarray*}
and $\Fcal^{s} (f,g)$ gives the pointwise evolution of $(f,g)$
for time $s$ according to the ODE in~(\ref{eq:ode}), i.e. for all $x$,
if $(u_0(x),v_0(x))=(f(x),g(x))$ then $\Fcal^{s} (f,g)(x) = (u_s(x),v_s(x))$
where $(u(x),v(x))$ evolves according to~(\ref{eq:ode}).
Note that both $e^{s\Delta}$ and $\Fcal^s$ are monotone in initial conditions,
therefore so is $e^{s\Delta} \Fcal^s$.

To establish the existence of constants $D_1$ and $d_1$ in condition~($*$),
it suffices to show that for any initial condition $(u_0,v_0)$
with $v_0$ dominating the function $D_1 I_{[-L,L]}(x)$,
for sufficiently large $T$, $(u_T,v_T)$ is such that $v_T$ dominates
the function $d_1 I_{[-3L,3L]}(x)$.

{
\psfrag{lminus}{$-l$}
\psfrag{lplus}{$l$}
\psfrag{Lminus}{$-L$}
\psfrag{Lplus}{$L$}
\psfrag{L-lminus}{$-L-l$}
\psfrag{L+lminus}{$-L+l$}
\psfrag{L-lplus}{$L-l$}
\psfrag{L+lplus}{$L+l$}
\psfrag{0}{$0$}
\psfrag{1}{$1$}
\begin{figure}[h!!]
\centering
\subfigure[The function $h$] {
\includegraphics[height = 1.7in, width = 1.5in]{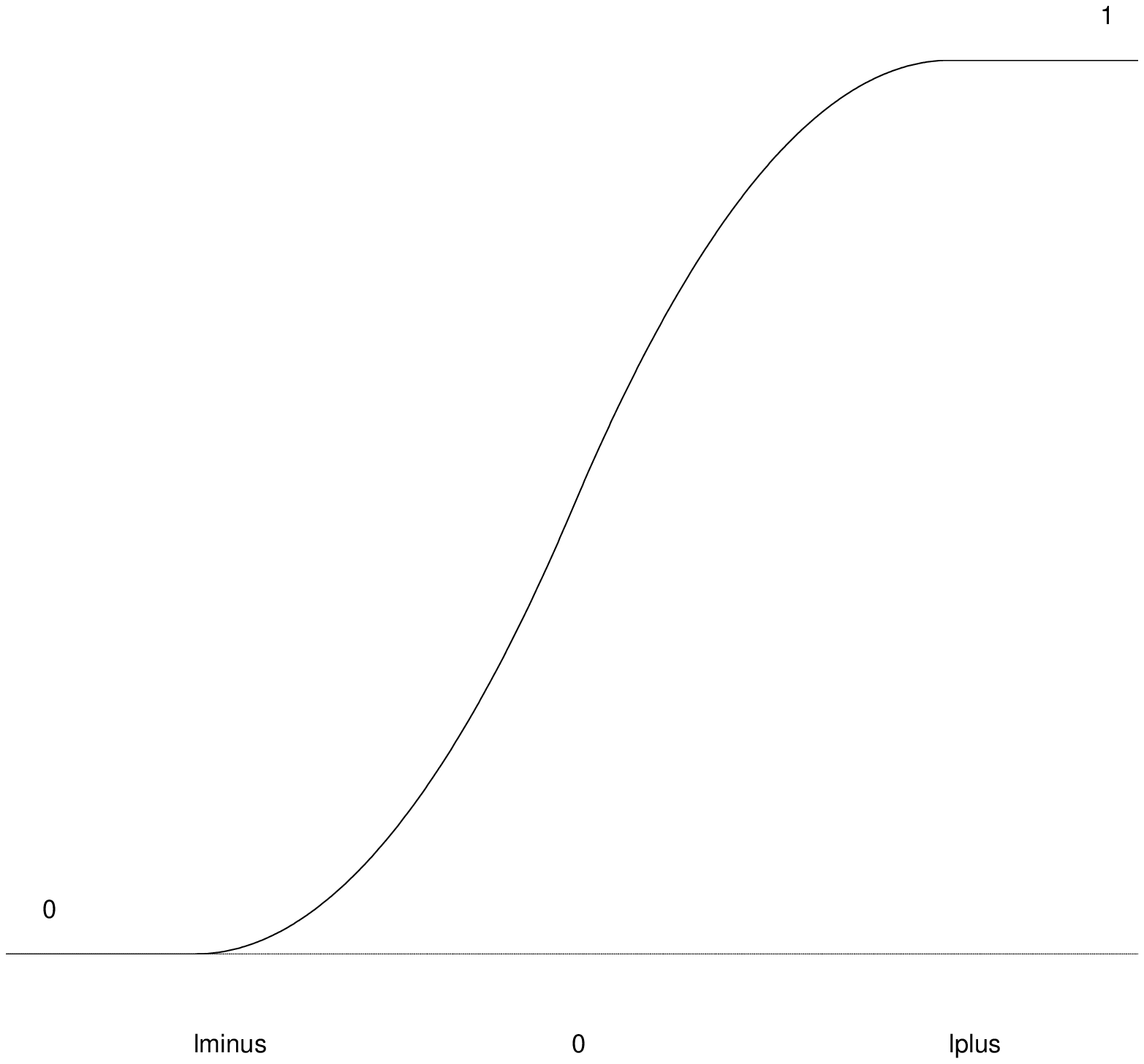}}
\subfigure[The function $f_0$] {
\includegraphics[height = 1.7in, width = 3.5in]{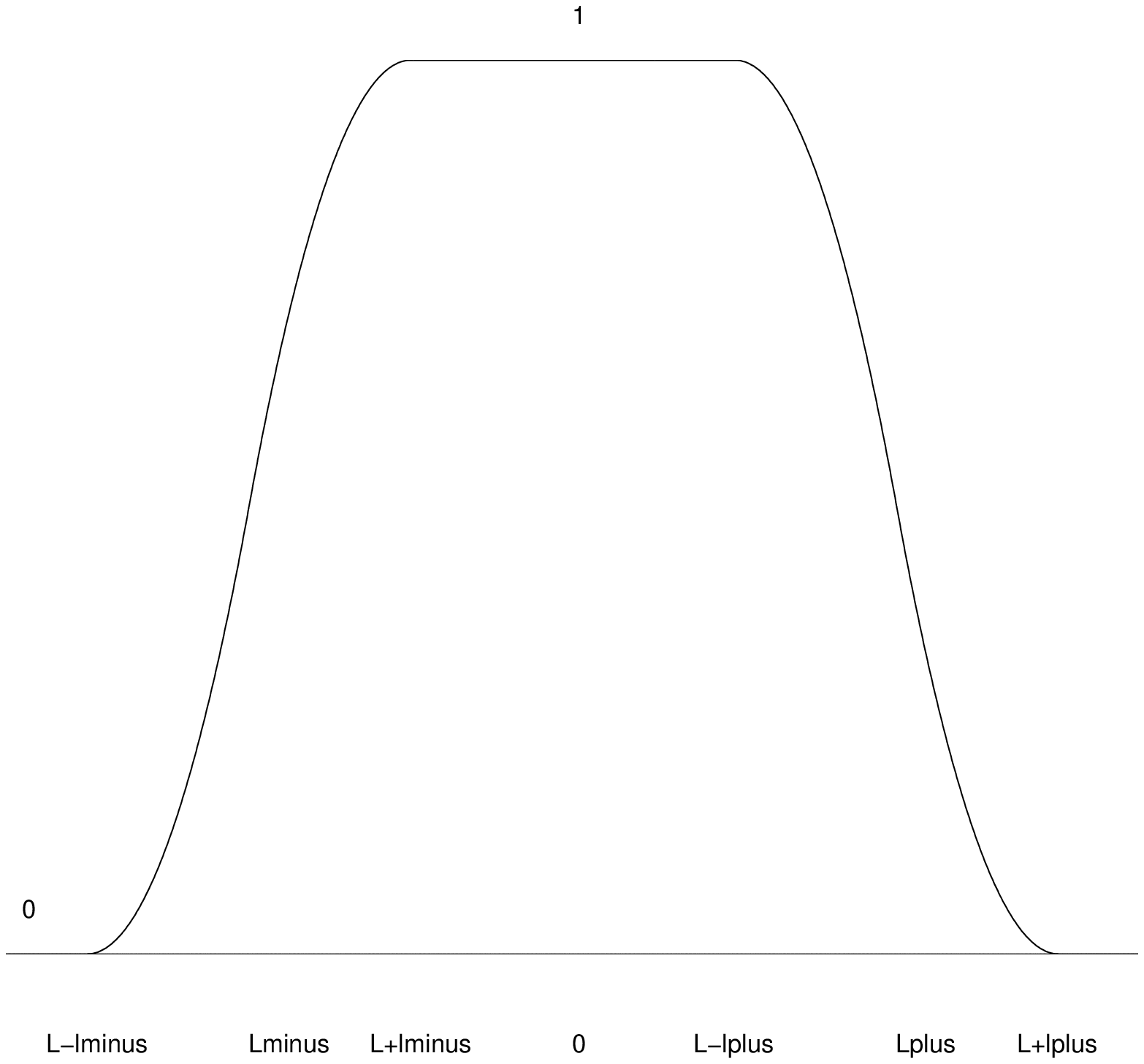}}
\caption{The functions $h$ and $f_0$.}
\end{figure}
}

Let $H^2(\Rbold)=\{f\in L^2(\Rbold): \|f\| = \int (1+\xi^2)
|\hat f(\xi)|^2 d\xi < \infty\}$ denote the Sobolev space with parameter 2,
where $\hat f$ denotes the Fourier transform of $f$. Equivalently,
$H^2(\Rbold)$ consists of $L^2$-functions with $L^2$-second derivatives.
We first define $f_0 \in H^2(\mathbb{R})$ that will be the ``shape''
of the initial conditions $(u_0,v_0)$:
\begin{enumerate}
\item[] {\bf IC 1.} $f_0(x)=1$ for $x\in [-L+l,L-l]$;
\item[] {\bf IC 2.} $f_0(x)=0$ for $x\in (\infty,-L-l]\cup [L+l,\infty)$;
\item[] {\bf IC 3.} $f_0(x)=h(x+L)$ for $x\in [-L-l,-L+l]$ and $f_0(x)=h(L-x)$
	for $x\in [L-l,L+l]$, where $h \in H^2(\mathbb{R})$ is the
	following function:
\begin{eqnarray}
        h(x) = \left \{ \begin{array}{ll}
                0,                       	& x<-l \\
                \frac{1}{2}(\frac{x+l}{l})^2, 	& -l\leq x\leq 0 \\
                1-\frac{1}{2}(\frac{l-x}{l})^2,	& 0<x\leq l \\
                1,             			& x>l
        \end{array} \right. .
\end{eqnarray}
\label{page:IC}
\end{enumerate}
In the above definition, the choice of $L$ is arbitrary as long as $L>l$,
but we will later on choose $l$ small
such that $\Delta f_0$ is large in $[-L-l,-L+l] \cup [L-l,L+l]$.
We call the intervals $[-L-l,-L+l]$ and $[L-l,L+l]$ the ``transition regions''.
We observe that $h$ is continuous at $x=0$, with
\begin{eqnarray*}
	h''(x)= \left\{ \begin{array}{ll}
			\frac{1}{l^2}	& \mbox{if $-l<x<0$} \\
			-\frac{1}{l^2}	& \mbox{if $0<x<l$}
		\end{array} \right. ,
\end{eqnarray*}
so the graph of $h$ in the plane
is symmetric about the point $(0,\frac{1}{2})$ and also,
\begin{eqnarray}
	|\Delta f_0| \leq \frac{1}{l^2} \label{eq:boundf0}
\end{eqnarray}
everywhere.

{
\psfrag{g1}{$\gamma_1$}
\psfrag{R0}{$\R_0$}
\psfrag{a0b0}{$(a_0,b_0)$}

\begin{figure}[b!]
\centering
\includegraphics[height = 3.5in, width = 3.5in]{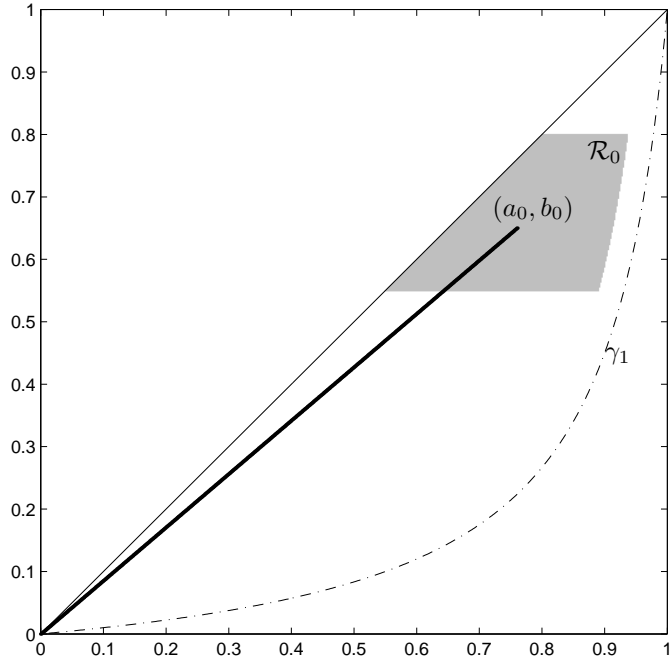}
\caption{Shaded region is $\R_0$, and $\gamma_1$ will be
defined in~(\ref{def:gamma1}).}
\label{fig:R0}
\end{figure}
}

We pick the initial condition to be $u_0=a_0 f_0,v_0=b_0 f_0$ with
$(a_0 f_0(0),b_0 f_0(0)) = (a_0,b_0) \in \R_0$, where $\R_0$ is the region
for the top tip of the line segment $\overline{O(a_0,b_0)}$ to be defined
later in~(\ref{def:R0}). See figure~\ref{fig:R0} for an illustration of
$\R_0$ and the line segment $\overline{O(a_0,b_0)}$.
We will define a family of parallel
piecewise-linear curves $\overline{ABC}{(u_0,v_0)}$
(see figure~\ref{fig:ode}) that satisfy the following requirements:
\begin{enumerate}
\item[] {\bf ABC 1.} $\overline{ABC}{(u_0,v_0)}$ lies in $\R_0$;
\item[] {\bf ABC 2.} $\overline{ABC}{(u_0,v_0)}$ passes through the points $A$, $B$, $C$,
	and $(u_0,v_0)$;
\item[] {\bf ABC 3.} $A$ lies on the line $u=v$, $B$ lies on the line $u=v+0.1$,
	and $C$ lies on the line $u=1$;
\item[] {\bf ABC 4.} $\overline{ABC}{(u_0,v_0)} = \overline{AB}(u_0,v_0)
		\cup \overline{BC}(u_0,v_0)$, where line segments
	$\overline{AB}(u_0,v_0)$ connects $A$ and $B$, and
	$\overline{BC}(u_0,v_0)$ connects $B$ and $C$;
\item[] {\bf ABC 5.} $\overline{AB}(u_0,v_0)$ makes an angle of $-\theta$
	($0 < \theta < \arctan (0.05)$) with the positive $u$-axis and
	$\overline{BC}(u_0,v_0)$ is a horizontal line segment.
\end{enumerate}
\label{page:ABC}

We will establish the following:

\begin{PRP}
	If $c$ is sufficiently large and
	$(\tilde a_0,\tilde b_0) \in \{(u,v) \in \R: 0.55 \leq v \leq 0.8 \}$,
	then for sufficiently small $s$, we have
\begin{eqnarray*}
	e^{s\Delta} \Fcal^{s} (\tilde a_0 f_0,\tilde b_0 f_0)
		\geq (\tilde a_s f_s,\tilde b_s f_s),
\end{eqnarray*}
where $f_0$ is defined in (IC 1-3) on page~\pageref{page:IC}, $\geq$ means
that $\geq$ holds in each component, and
$\tilde a_s$, $\tilde b_s$, and $f_s$ satisfy the following conditions:
\begin{enumerate}
\item $\tilde a_s$ and $\tilde b_s$ are constants depending on $s$, such that
	the curve $\overline{ABC}{(\tilde a_s,\tilde b_s)}$
        lies above the curve $\overline{ABC}{(\tilde a_0,\tilde b_0)}$,
        and the vertical distance separating them is at least $\delta_2 s$.
        In particular, since $\overline{ABC}{(\tilde a_0,\tilde b_0)}$
        lies above the horizontal line $v=0.5$, so does
        $\overline{ABC}{(\tilde a_s,\tilde b_s)}$.
\item
\begin{eqnarray}
	f_s (x) = \left\{ \begin{array}{ll}
		1, & x\in [-L+l-\delta_1 s,L-l+\delta_1 s] \\
		h(x+L+\delta_1 s), & x\in [-L-l-\delta_1 s,-L+l-\delta_1 s] \\
		h(L+\delta_1 s-x), & x\in [L-l+\delta_1 s,L+l+\delta_1 s] \\
		0, & x\in (\infty,-L-l-\delta_1 s]\cup [L+l+\delta_1 s,\infty)
		\end{array} \right. , \label{eq:Dd}
\end{eqnarray}
	i.e. $f_s$ is $f_0$ with each of the two transition regions is
	translated by $\delta_1 s$ away from the origin.
\item $\delta_1$ and $\delta_2$ are positive constants independent of $s$ and
        $(\tilde a_0,\tilde b_0)$.
\end{enumerate}
\label{prp:Dd}
\end{PRP}

The proof of the above proposition requires a few lemmas and will be deferred
until the end of Chapter~\ref{sec:pde}. Proposition~\ref{prp:Dd} states that
$e^{s\Delta} \Fcal^{s} (\tilde a_0 f_0,\tilde b_0 f_0)$ is bounded below by
$(\tilde a_s f_s,\tilde b_s f_s)$, where $(\tilde a_s,\tilde b_s)$
moves $\delta_2 s$ above $\overline{ABC}{(\tilde a_0,\tilde b_0)}$.
Furthermore, by~(\ref{eq:Dd}), the region where the value of
$(\tilde a_s f_s,\tilde b_s f_s)$
(and hence $e^{s\Delta} \Fcal^{s} (\tilde a_0 f_0,\tilde b_0 f_0)$)
is equal to or above $(\tilde a_s,\tilde b_s)$ has expanded
by $\delta_1 s$, both to the left and to the right, while the transition
regions of $(\tilde a_s f_s,\tilde b_s f_s)$ are shifted left or right
by $\delta_1 s$
but maintains exactly the same profile as in the initial condition.
Thus by the monotonicity of $e^{s\Delta} \Fcal^s$, we can iterate
$e^{s\Delta} \Fcal^s$ enough times and obtain
information about the evolution of the PDE~(\ref{eq:pde:real}) for large time.
From the construction of the piecewise linear curves $\overline{ABC}$,
requirement (5) implies that $B(0.8,0.8)$ on $\overline{AB}(0.8,0.8)$
has $v$-coordinate $\geq 0.75$. Therefore Corollary~\ref{cor:Dd} below
is an easy consequence of Proposition~\ref{prp:Dd}.
Let
\begin{eqnarray}      
        \pi_v (u_0,v_0) = v_0, \label{def:pi}
\end{eqnarray}
and $\lfloor x \rfloor = \max \{z\in \Zbold: z\leq x\}$.

\begin{COR}
If $c$ and $T$ are sufficiently large,
and $v_0 (x) > 0.55$ for $x\in [-L+l,L-l]$
then \[ \pi_v((e^{s\Delta}\Fcal^s)^{\lfloor T/s \rfloor}(u_0,v_0)(x)) > 0.7 \]
for $x\in [-3L,3L]$ and sufficiently small $s$.
\label{cor:Dd}
\end{COR}

In other words, the constants $D_1$ and $d_1$ in condition ($*$) are
picked to be $D_1 = 0.55$ and $d_1=0.7$.
Note that we restrict $v_0(x)=\tilde b_0 f_0(x)$ to be $>0.55$ for
$x\in [-L+l,L-l]$ in the above corollary because Prop~\ref{prp:Dd}
only works for
$(\tilde a_0,\tilde b_0) \in \{(u,v) \in \R: 0.55 \leq v \leq 0.8 \}$.
Also note that the ``$L$'' in condition ($*$) is picked to be
$L-l$.

\subsection{Analysis of the ODE~(\ref{eq:ode})}
\label{sec:ode}
We first characterize the phase portrait of the ODE. We carry this out for
sufficiently large $c$. See figure~\ref{fig:ode} for phase portraits
with $c=5$ and $c=25$. Define
\begin{eqnarray}
	\eta(u,v)=(\eta_1(u,v),\eta_2(u,v))=((2c(1-u)+1)v-u,(c(u-v)-2)v),
	\label{def:eta}
\end{eqnarray}
such that the solution to the ODE~(\ref{eq:ode}),
$(u_s,v_s)=\Fcal_\eta^s (u_0,v_0)$,
flows along the vector field $\eta$.
Define the curves $\gamma_1$, $\gamma_2$, and $\gamma_3$:
\begin{eqnarray}
	\gamma_1 & = & \left\{(u,v): u\in [0,1], v=\frac{u}{1+2c-2cu}
		\right\}, \label{def:gamma1} \\
	\gamma_2 & = & \left\{(u,v): u\in [0,1], v=u-\frac{2}{c}
                \right\}, \label{def:gamma2} \\
	\gamma_3 & = & \left\{(u,v): u\in [0,1], u=v \right\}.
		\label{def:gamma3}
\end{eqnarray}
We have $\eta_1=0$ on $\gamma_1$ and $\eta_2=0$ on $\gamma_2$.
An easy calculation shows that for all $c$, $(0,0)$ and
$(1,1)$ pass through $\gamma_1$. We observe that $\eta_1$ is a
linear function in $u$ for fixed $v$, so $\eta_1>0$ to the left of
$\gamma_1$ and $\eta_1<0$ to the right of $\gamma_1$. 
By similar reasoning, we also have $\eta_2<0$ to the left of $\gamma_2$,
while $\eta_2>0$ to its right. The two
intersection points of $\gamma_1$ and $\gamma_2$,
\begin{eqnarray*}
	P_+ & = & \left(
	\frac{1}{2}+\frac{1}{c}+\sqrt{\frac{1}{4}-\frac{1}{c}},
	\frac{1}{2}-\frac{1}{c}+\sqrt{\frac{1}{4}-\frac{1}{c}}\right), \\
	P_- & = & \left(
        \frac{1}{2}+\frac{1}{c}-\sqrt{\frac{1}{4}-\frac{1}{c}},
        \frac{1}{2}-\frac{1}{c}-\sqrt{\frac{1}{4}-\frac{1}{c}}\right),
\end{eqnarray*}
are the only equilibrium points of $\eta$ in the interior of $\R$, with
$O=(0,0)\in \partial \R$ being the third equilibrium point.
Elementary computation shows that
$O$ and $P_+$ are stable, but $P_-$ is a saddle point, thus one would expect
any point that lies significantly above $P_-$ to flow toward $P_+$
under $\eta$. Elementary calculations also show the following:
\begin{eqnarray}
	& P_+ \rightarrow (1,1) & \mbox{ as } c\rightarrow \infty,
		\label{eq:Pplus} \\
	& P_- \rightarrow (0,0) & \mbox{ as } c\rightarrow \infty, \nnb \\
	& \frac{P_{-,u}}{P_{-,v}} \rightarrow \infty & \mbox{ as }
		c\rightarrow \infty, \nnb
\end{eqnarray}
where $P_{-,u}$ and $P_{-,v}$  denote the $u$- and $v$- coordinates of
$P_-$, respectively.

{
\psfrag{P1}{$P_+$}
\psfrag{P2}{$P_-$}
\psfrag{P3}{$P_3$}
\psfrag{P4}{$P_4$}
\psfrag{A1}{$A_1$}
\psfrag{A2}{$A_2$}
\psfrag{B1}{$B_1$}
\psfrag{B2}{$B_2$}
\psfrag{C1}{$C_1$}
\psfrag{C2}{$C_2$}
\psfrag{Q}{$Q$}
\psfrag{E}{}
\psfrag{g1}{$\gamma_1$}
\psfrag{g2}{$\gamma_2$}
\psfrag{uu}{$u$}
\psfrag{vv}{$v$}

\begin{figure}[t!]
\subfigure[$c=5$] {
\includegraphics[height = 2.5in, width = 2.5in]{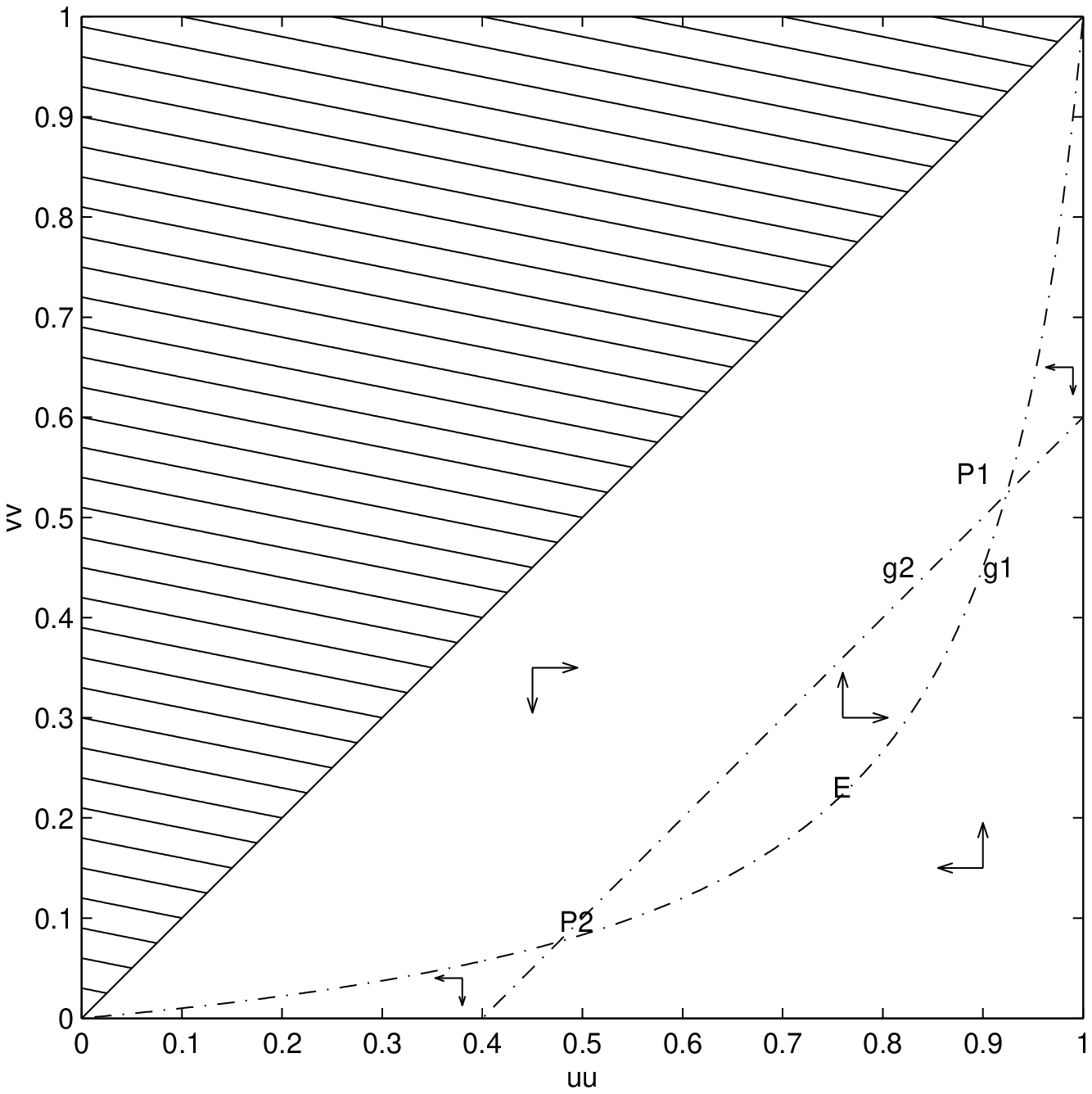}}
\subfigure[$c=25$] {
\includegraphics[height = 2.5in, width = 2.5in]{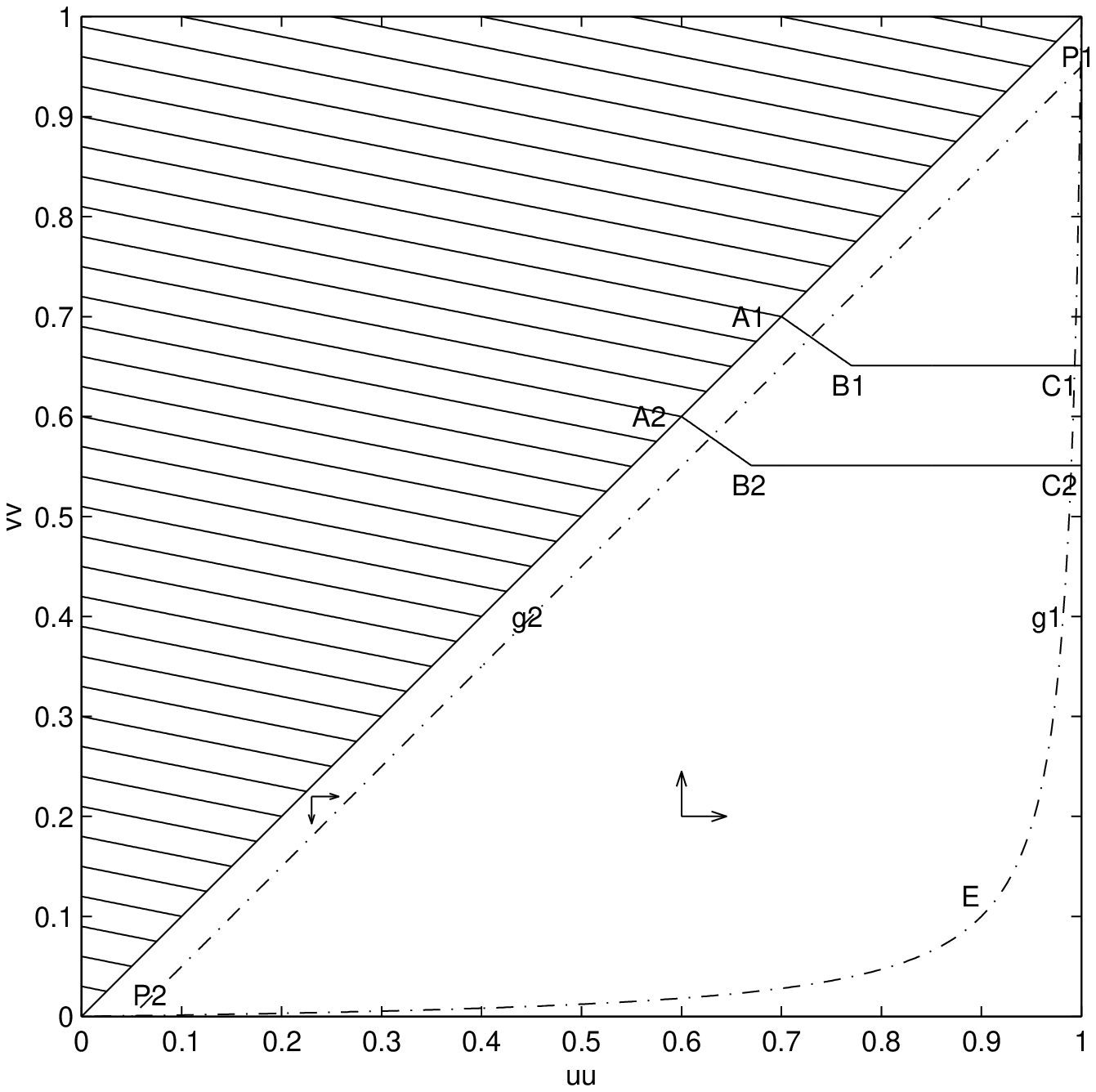}}
\caption{Phase space of the ODE}
\label{fig:ode}
\end{figure}
}

%We now define $(P_3,P_4)=\left(0,\frac{1}{2c}(c-2\mp\sqrt{c^2-4c})\right)$
%and characterize $\eta$ in the rectangle $P_+ P_- P_3 P_4$.
%First, we compute: $\frac{d}{du} \frac{\eta_1(u,v)}{\eta_2(u,v)} =
%2 \frac{1+(2c-c^2)v+c^2 v^2}{v(2-cu)^2}$. We observe that the numerator
%$1+(2c-c^2)v+c^2 v^2$ is a quadratic function of $v$ and has roots
%$\frac{1}{2c}(c-2\pm\sqrt{c^2-4c})$; since $c^2>0$, the numerator is $<0$
%for $v \in \left(\frac{1}{2c}(c-2-\sqrt{c^2-4c}),
%\frac{1}{2c}(c-2+\sqrt{c^2-4c})\right)$; notice that the endpoints of this
%interval are the $v$-coordinates of $P_+$ and $P_-$, respectively.
%Thus in the rectangle $P_+ P_- P_3 P_4$, the angle $\theta(u_0,v_0)$,
%between the horizontal line $v=v_0$ and the vector $\eta(u_0,v_0)$, is
%negative and decreasing (i.e. getting more horizontal) along horizontal lines
%as $u$ increases from 0 to $2/c$. Furthermore, along the $v$-axis,
%$\tan \theta(0,v_0) = -\frac{1}{c(1-v_0)+1}$ is negative and decreasing
%as $v_0$ increases, i.e. getting less horizontal.

We will need some crude estimates of $\eta_1$ and $\eta_2$.
First of all, since $u - v \geq 0$ everywhere in $\R$, we have
\begin{eqnarray}
        \eta_2 \geq -2v. \label{eq:eta2:easy}
\end{eqnarray}
If the point $(u,v)$ is at least $\delta$ to the right of $\gamma_2$, then
since $\eta_2(u,v)=0$ on $\gamma_2$,
\begin{eqnarray}
	\eta_2(u,v)>\delta c v, \label{est:eta2}
\end{eqnarray}
Similar reasoning shows that if the point $(u,v)$ is
at least $\delta$ to the left of $\gamma_1$, then
\begin{eqnarray}
	\eta_1(u,v)>2\delta cv. \label{est:eta1}
\end{eqnarray}
Now the horizontal distance between $\gamma_1$ and $\gamma_3$
at a fixed $v$ is
\begin{eqnarray*}
	d(v)=\frac{(1+2c)v}{1+2cv} - v.
\end{eqnarray*}
Simple calculations show that
\begin{eqnarray*}
	d'(v) = -\frac{2c(2cv^2+2v-1)}{(1+2cv)^2} \ 
	\mbox{ and } \ d''(v) = -\frac{4c(1+2c)}{(1+2cv)^3}.
\end{eqnarray*}
Notice that $d''(v) < 0$ if $v,c>0$, so for $v\in [0,1]$, $d(v)$
is a strictly concave function, and
\begin{eqnarray*}
	\hat v=\frac{1}{2c}(\sqrt{1+2c}-1)
\end{eqnarray*}
is the unique point where $d(v)$ attains its maximum for $v\in [0,1]$.
We observe that
\begin{eqnarray}
	\hat v \rightarrow 0 \mbox{ as } c\rightarrow \infty. \label{eq:hatv}
\end{eqnarray}
If $c$ is sufficiently large such that
the horizontal line $v=\epsilon$ lies above the line $v=\hat v$
but below the line $v=0.8$,
then for $v\in [\epsilon,0.8]$, the minimum of $d(v)$
occurs at $v=0.8$, and
\begin{eqnarray}
	d(0.8)= \frac{0.8(1+2c)}{1+2(0.8)c} - 0.8
	= \frac{0.8+1.6c}{1+1.6c} - 0.8
	\rightarrow 1-0.8 = 0.2 \label{ineq:d}
\end{eqnarray}
as $c\rightarrow \infty$.
This shows that for arbitrary $\epsilon>0$ and sufficiently large $c$,
the minimum horizontal distance between $\gamma_1$ and $\gamma_3$
for $v\in [\epsilon,0.8]$, is larger than 0.19. This fact will
be needed a bit later on in the proof of Lemma~\ref{lem:case1}. We also
observe that since $d(v)$ is strictly concave for $v \in [0,1]$, the curve
$\gamma_1$ written as $u = u(v)$ is also strictly concave for $v \in [0,1]$;
then~(\ref{eq:hatv}) and the fact
\begin{eqnarray*}
	u(\hat v) = \frac{1+2c-\sqrt{1+2c}}{2c} \rightarrow 1
	\mbox{ as } c \rightarrow \infty
\end{eqnarray*}
imply that
\begin{eqnarray}
	\gamma_1 \rightarrow \{(u,v): v=0, u \in [0,1]\} \cup
		\{(u,v): u=1, v \in [0,1]\} 
		\mbox{ as } c \rightarrow \infty. \label{eq:gamma1:cbig}
\end{eqnarray}
This finishes the characterization of the phase portrait of ODE~(\ref{eq:ode}).
These facts, which are admittedly tedious and not much fun to establish,
will all be required later in the proof of Lemmas~\ref{lem:case1}
and~\ref{lem:case2}.

We define the region (see figure~\ref{fig:R0})
\begin{eqnarray}
	\R_0=\{(u,v)\in \R: u < \frac{(1+2c)v}{1+2cv}-0.04,
		0.55\leq v\leq 0.8\}, \label{def:R0}
\end{eqnarray}
Recall that $\gamma_1$ defined in~(\ref{def:gamma1}) is the curve
$v=\frac{u}{1+2c-2cu}$ or $u=\frac{(1+2c)v}{1+2cv}$, and therefore
the region $\R_0$ lies at least 0.04 to the left of the curve $\gamma_1$.
By (ABC 5), the line segment $AB$ forms a small negative
angle with the positive $u$-axis, so by requiring $v\geq 0.55$
in the definition of $\R_0$, we can be sure that all of
$\overline{ABC}(a_0,b_0)$ lies above the horizontal line $v=0.5$ if
$(a_0,b_0)\in \R_0$. For any point $(u,v)\in \R_0$, we have $u<2v$
and $(u,v)$ lies at least 0.04 to the left of $\gamma_1$.
Since our initial condition has form $(a_0 f_0,b_0 f_0)$,
%where $\overline{ABC}(a_0,b_0)$ stays above the line $v=0.5$,
the set of values in each ``transition region''
\begin{eqnarray*}
	\{(a_0 f_0(x), b_0 f_0(x)): x \in [L-l,L+l] \},
\end{eqnarray*}
forms a line segment with endpoints $O$ and $(a_0,b_0)$ in the $(u,v)$-plane.
We require that the tip of this line segment $(a_0,b_0)$ lies in the region
$\R_0$. We do not need to worry about the case where the initial condition
for the PDE~(\ref{eq:pde:real}) is such that
$(a_0,b_0) \in \R\cap \{0.55\leq v\leq 0.8\}$ lies to the right of $\R_0$. 
\label{page:mono}
If we want to establish condition ($*$) for
that initial condition, then by the monotonicity of the
PDE~(\ref{eq:pde:real}), it is sufficient to pick
$a_0' < a_0$ such that $(a_0',b_0)\in \R_0$
and prove condition ($*$) for the initial condition $(a_0' f_0,b_0 f_0)$.
Therefore we only consider $(a_0,b_0)$ in $\R_0$.

Assuming $(a_0,b_0)$ lies in $\R_0$, a part of
the line segment $\overline{O(a_0,b_0)}$ still
lies below the horizontal line $v=0.55$. To study the evolution of the
whole line segment under $\Fcal_\eta$, we will a bit later consider two
cases: 1. $\epsilon \leq v \leq 0.8$, and 2. $0\leq v < \epsilon$,
where we will pick $\epsilon=0.24$ in Chapter~\ref{sec:pde}.
We will construct piecewise linear curves $\overline{ABC}{(v_0,v_0)}$,
$v_0 \in [0.55,0.8]$, with $A=(v_0,v_0)$, $B$, and $C$ satisfying the
requirements laid down in (ABC 1-5) on page~\pageref{page:ABC},
in the proof of the following two lemmas. See figures~\ref{fig:case1a}
and~\ref{fig:case2a} for an illustration of each lemma.

\begin{LEM}
	(Case 1) If $(a_0,b_0)$ lies on
	$\overline{AB}{(v_0,v_0)}=\overline{AB}{(a_0,b_0)}$ with
	$a_0 - b_0 < 0.09$ and $0.55 < b_0 \leq 0.8$,
	then for sufficiently small $s$,
	there exist $a_s$ and $b_s$ with $b_s > 0.5$,
	and a positive number $\tilde{K}$ independent of $s$
	such that $\overline{AB}{(a_s,b_s)} =\overline{AB}{(v_0,v_0)}$ and
\begin{eqnarray}
	\Fcal^s_\eta (\alpha a_0,\alpha b_0) \geq
	\left((1+\tilde{K}s)\alpha a_s,(1+\tilde{K}s)\alpha b_s\right)
	\label{eq:case1}
\end{eqnarray}
	for all $\alpha\in [0,1]$. Moreover, the constant $\tilde K$ can
	be chosen to be arbitrarily large if $c$ is also allowed to be
	arbitrarily large.
\label{lem:case1}
\end{LEM}
\begin{REM}
	Using some easy geometric considerations, one can say the following:
	if $\epsilon$ is fixed and $c$ is allowed to be
        arbitrarily large, then there exists an arbitrarily large
	positive number $K$ depending on $\epsilon$ but independent of $s$,
	such that if $\alpha \in [\frac{\epsilon}{b_s},1]$ then
\begin{eqnarray}
	((1+\tilde{K}s)\alpha a_s,(1+\tilde{K}s)\alpha b_s)
	-(\alpha a_s,\alpha b_s)>\left(\frac{a_s}{b_s}Ks,Ks\right),
	\label{eq:case1a}
\end{eqnarray}
	and if $\alpha \in [0,\frac{\epsilon}{b_s})$ then
\begin{eqnarray}
	((1+\tilde{K}s)\alpha a_s,(1+\tilde{K}s)\alpha b_s)
	-(\alpha a_s,\alpha b_s) \geq (0,0). \label{eq:case1b}
\end{eqnarray}
\label{rem:case1}
\end{REM}

\begin{LEM}
	(Case 2) If $(a_0,b_0)$ lies on $\overline{ABC}{(v_0,v_0)}=
\label{page:case2}
	\overline{ABC}{(a_0,b_0)}$ with
	$0.08 + b_0 < a_0 < \frac{(1+2c)v}{1+2cv}-0.04$ and
	$0.55 < b_0 \leq 0.8$, then for sufficiently small $s$,
	there exists a positive number
	$K$ such that, if $\alpha \in [\frac{\epsilon}{b_0},1]$, then
\begin{eqnarray}
	\Fcal^s_\eta (\alpha a_0,\alpha b_0)-(\alpha a_0,\alpha b_0)
	\geq \left(\frac{a_0}{b_0} Ks,Ks\right), \label{eq:case21}
\end{eqnarray}
	and if $\alpha \in [0,\frac{\epsilon}{b_0})$, then
\begin{eqnarray}
	\Fcal^s_\eta (\alpha a_0,\alpha b_0)-(\alpha a_0,\alpha b_0)
	\geq (-2\alpha a_0 s,-2\alpha b_0 s). \label{eq:case22}
\end{eqnarray}
	Moreover, the constant $K$ can be chosen to be arbitrarily
	large if $c$ is also allowed to be arbitrarily large.
\label{lem:case2}
\end{LEM}

In case 1 above (Lemma~\ref{lem:case1}),
$(a_s,b_s)$ changes with $s$, but in case 2 \\
(Lemma~\ref{lem:case2}), $(a_s,b_s)$ remains
fixed and equal to $(a_0,b_0)$, thus $(a_s,b_s)$ is not explicitly defined.
As will be seen later on, $l$ is picked small so that the lower part
of the ``transition region'' (of $(a_0 f_0,b_0 f_0)$) moves up at a
sufficiently large speed under the heat kernel to cancel out the downward
movement as described in~(\ref{eq:case22}).
But the heat kernel pushes down the top part of the ``transition region'',
so $K$ (and thus $c$) is picked large to cancel out
that effect. And finally $s$ is picked small so that the movement
caused by $\Fcal_\eta^s$ is small.

In case 1, we assume $(a_0,b_0)$, the top tip of the line segment formed by
the ``transition region'', lies to the
left of the line $u=v+0.09$, while in case 2, we assume that $(a_0,b_0)$
lies to the right of $u=v+0.08$. There
is a thin strip, i.e. $0.08 < u-v < 0.09$, where both cases apply, so we
can apply either case 1 or case 2 there. Let $(u_0,v_0)$ be a point on
the line segment $\overline{O(a_0,b_0)}$. Intuitively, we would like to
view ``progress'' as an increase in $v$-coordinate, which is measured by
$\pi_v(\Fcal_\eta^s (u_0,v_0) - (u_0,v_0))$.
In case 1, however, it is not always possible for the $v$-coordinate to
increase. So instead, we measure progress with respect to the family of
parallel lines $\overline{AB}$, each of
which makes a small negative angle with the positive $u$-axis,
and thus allowing the $v$-coordinate to decrease slightly while still
making ``progress'' with respect to $\overline{AB}$.
More specifically, we compare $\Fcal_\eta^s (u_0,v_0)$ not with
$(u_0,v_0)=(\frac{u_0}{a_0}a_0,\frac{v_0}{b_0}b_0)$, but with
$(u_s,v_s)=(\frac{u_0}{a_0}a_s,\frac{v_0}{b_0}b_s)$.
We show that with respect to the lines $\overline{AB}$,
the entire line segment makes progress with respect to the lines
$\overline{AB}$ when moving under $\eta$.
The $u$-coordinate actually increases very rapidly, so we move very quickly
into where case 2 applies. In case 2, we compare $\Fcal_\eta^s (u_0,v_0)$
directly with $(u_0,v_0)$, and show that the part of the
line segment with $v$-coordinate $>\epsilon$ makes progress,
but the part of the line segment with
$v$-coordinate $<\epsilon$ actually makes small negative progress. This 
negative progress will be compensated by positive progress made under
evolution according to the heat kernel (to be shown in Chapter~\ref{sec:pde}).

\vspace{.3cm}
\noindent {\bf Proof of Lemma~\ref{lem:case1}}. \hspace{2mm}
We first define
\begin{eqnarray}
        \R_1 & = & \{(u,v)\in \R:u-v\in [0,0.1], v\in [\epsilon,0.8]\},
        \label{def:R1p} \\
        \R_1' & = & \{(u,v)\in \R_1:u-v\in [0,0.09]\}. \label{def:R1pp} \\
	\R_3' & = & \{(u,v)\in \R:u-v\in [0,0.1],
		u\leq 2v,v\in [0,\epsilon)\}. \label{def:R3}
\end{eqnarray}
Later in the proof of Lemma~\ref{lem:case2}, we will also define
the follow three regions, which we include here for easy reference
(See figure~\ref{fig:regions}).
\begin{eqnarray*}
        \R_2 & = & \{(u,v)\in \R: v+0.02<u<\frac{(1+2c)v}{1+2cv}-0.04,
                v\in [\epsilon,0.8]\}, \\
        \R_2' & = & \{(u,v)\in \R_2: v+0.08<u<\frac{(1+2c)v}{1+2cv}-0.04,
                v > 0.55 \}, \\
        \R_3 & = & \{(u,v)\in \R: u\leq 2v,v\in [0,\epsilon)\}.
\end{eqnarray*}
Notice that $\R_1' \subset \R_1$, $\R_2' \subset \R_2$, and
$\R_3' \subset \R_3$.

{
\psfrag{y=8}{$v=0.8$}
\psfrag{y=e}{$v=\epsilon$}
\psfrag{u=2v}{$u=2v$}
\psfrag{R1}{$\R_1$}
\psfrag{R2}{$\R_2$}
\psfrag{R3}{$\R_3$}
\psfrag{uu}{$u$}
\psfrag{vv}{$v$}

\begin{figure}{t}
\centering
\subfigure[$\R_1$ and $\R_3$] {
\includegraphics[height = 2.6in, width = 2.6in]{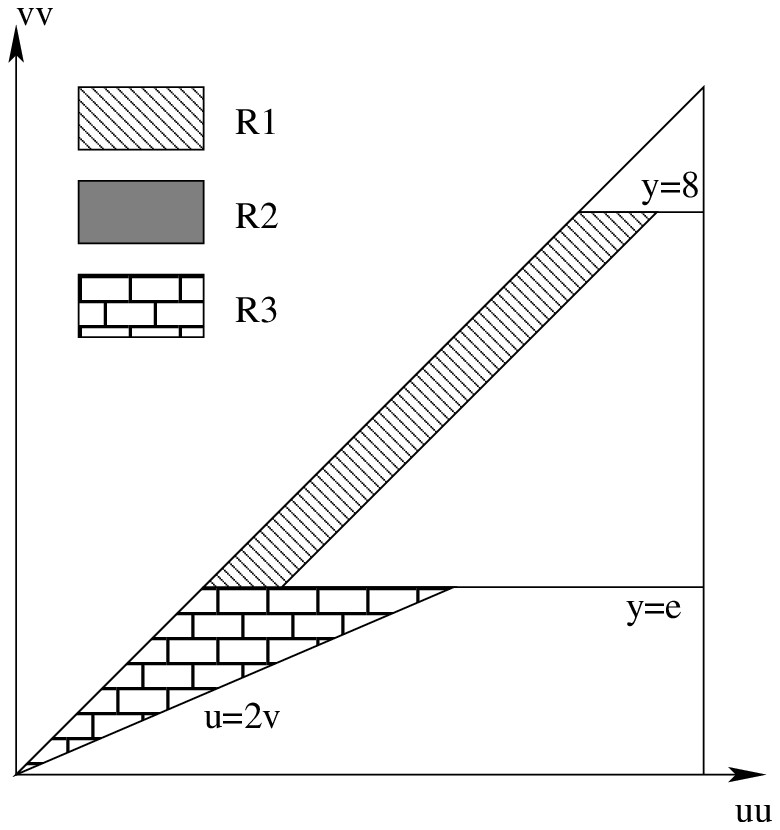}}
\subfigure[$\R_2$ and $\R_3$] {
\includegraphics[height = 2.6in, width = 2.6in]{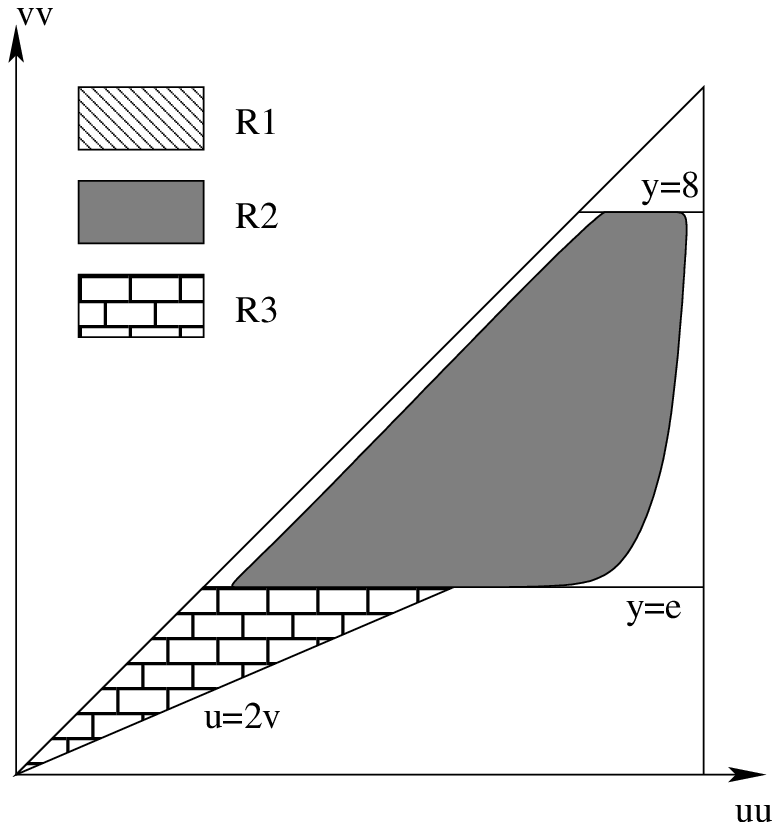}}
\caption{Various regions.}
\label{fig:regions}
\end{figure}
}

To study the evolution of $\Fcal_\eta$
of a line segment $\overline{OP}$ that passes
through the origin $O$ and has its top tip $P$ in the region
$\R_1'\cap \R_0$, we will define a new vector field
$\xi=(\xi_1,\xi_2)$ for $(u,v)\in \R_1 \cup \R_3'$, such that
$\xi_1 \leq \eta_1$ and $\xi_2 \leq \eta_2$ everywhere in $\R_1 \cup \R_3'$,
which means that
\begin{eqnarray}
        \Fcal^s_\xi(u,v) \leq \Fcal^s_\eta (u,v) \label{ineq:xieta}
\end{eqnarray}
in both $u$- and $v$- coordinates, for small $s$ and all
$(u,v) \in \R_1' \cup \R_3'$. Notice that $\R_1'$ is required to stay
a finite distance left of the line $u-v=0.1$, the
right edge of the parallelogram $\R_1$; this is such that we can still
control how much $(u,v) \in \R_1'$ moves under $\Fcal_\xi$,
even if it leaves $\R_1'$ and enters the strip $\R_1\backslash \R_1'$.
Also, we define $\xi$ in $\R_1 \cup \R_3'$ because this is the region where
the line segments $\overline{OP}$ lie.
We only consider the region $\R_3'$ for $(u,v)$ close to the origin,
rather than the region
\begin{eqnarray*}
        \{(u,v)\in \R:u-v\in [0,0.1],v\in [0,\epsilon)\},
\end{eqnarray*}
since top tip of the line segment $b_s$ will be $>0.5$,
which implies that $u \leq 2v$ if $(u,v)\in \overline{OP}$.

{
\psfrag{A}{$A$} 
\psfrag{B}{$B$}
\psfrag{Ap}{$A'$}
\psfrag{Bp}{$B'$}
\psfrag{e}{$\epsilon$}
\psfrag{R1}{$\R_1$}
\psfrag{R3}{$\R_3'$}
\psfrag{u}{$u$}
\psfrag{v}{$v$}
\psfrag{u=v}{$u=v$}
\psfrag{u=v+.1}{$u=v+0.1$}
\psfrag{u=2v}{$u=2v$}
\psfrag{F}{$v F_1$}
\psfrag{-2v}{$-2v$}
\psfrag{P1}{$\Fcal^s_\xi(a_0,b_0) = (1+\tilde{K}s)(a_s,b_s)$}
\psfrag{P2}{$(a_s,b_s)$}
\psfrag{P6}{$(a_0,b_0)$}
\begin{figure}[h!]
\centering
\includegraphics[height = 5in, width = 5in]{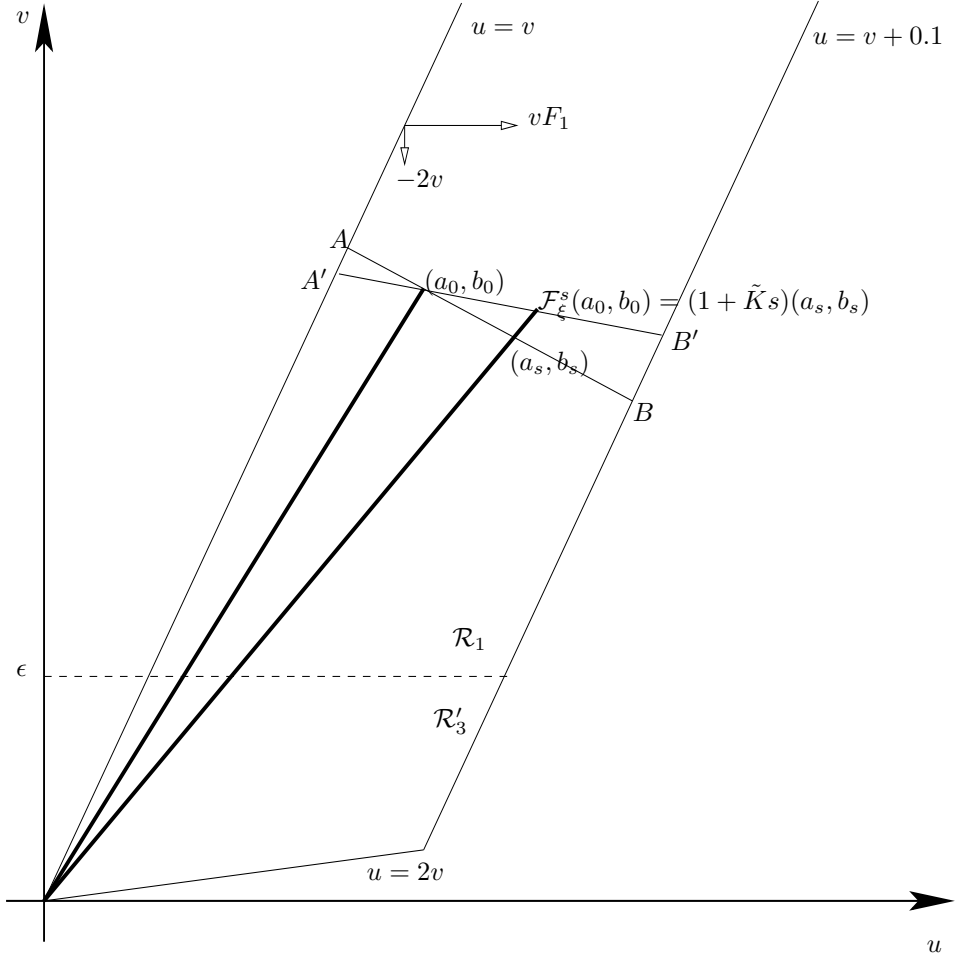}
\caption{Illustration of Lemma~\ref{lem:case1} (Case 1)}
\label{fig:case1a}                  
\end{figure}    
}               

{\bf Step 1: Defining $\xi$.}
We first define $\xi$ on the diagonal line
\[ \gamma_3 \cap \{(u,v): v\leq 0.8\} = \{(u,v): u=v, 0 \leq v\leq 0.8\}, \]
where $\gamma_3$ is defined in~(\ref{def:gamma3}).
Let $\theta_0 \leq \frac{1}{2} \arctan (0.05)$ be a small angle
such that a line passing through $(v_0,v_0)$ with $v_0 > 0.55$ and making
an angle of $-2\theta_0$ with the positive $u$-axis
intersects the vertical line $u=1$
above the horizontal line $v=0.5$. Define
\begin{eqnarray}
	\delta = 0.1, \label{def1:delta}
\end{eqnarray}
and pick $F_1$ to be large but                      
\begin{eqnarray}                                                          
        F_1<\delta c. \label{eq:pick1}
\end{eqnarray}                       
We also define
\begin{eqnarray}
        \xi(0.8,0.8)=(0.8F_1,0.8(-2)), \label{eq:pick1a}
\end{eqnarray}                               
where $F_1$ is large enough such that $\xi(0.8,0.8)$ makes an
angle of $-\theta_1$ with the positive $u$-axis, and $0<\theta_1<\theta_0$.
This can be done for sufficiently large $c$. We define
\begin{eqnarray*}                            
        \xi(v,v) = \frac{v}{0.8}\xi(0.8,0.8) = (v F_1, -2 v),
\end{eqnarray*}                             
for $v\in [0,0.8]$. This defines $\xi$ on the line segment
$\gamma_3 \cap (\R_1 \cup \R_3')$.

Finally, along straight lines that make angles of $-\theta_1$             
with the positive $u$-axis, denoted $\overline{A'B'}(v,v)$, we define $\xi$ to
be equal to
$\xi(v,v)$, i.e. for all $(u',v') \in \overline{A'B'}(v,v)$, we define
\begin{eqnarray*}
        \xi(u',v')=\xi(v,v).
\end{eqnarray*}
Here $A'=(v,v)$ is the point where $\overline{A'B'}(v,v)$ intersects
the line $u=v$ and $\xi$ has already been defined, and $B'$ is the
intersection point of $\overline{A'B'}(v,v)$ and the right/bottom edge
of the 4-gon $\R_1\cup \R_3'$, i.e. either the line $u=v+0.1$ or the line
$u=2v$. Thus we have defined $\xi$ on all points in $\R_1 \cup \R_3'$.
To summarize, $\xi$ in $\R_1$ is defined in a way such that:
\begin{enumerate}
\item On the line $u=v$ where $0 \leq v \leq 0.8$,
\begin{eqnarray}
        \xi=(vF_1,v(-2)), \label{def:xi}
\end{eqnarray}
        such that $\xi$ makes a small angle of $-\theta_1$ with the positive
	$u$-axis.
\item $\xi$ is constant along lines that start at a point on the line $u=v$,
        and make angles of $-\theta_1$ with the positive $u$-axis, where
\begin{eqnarray}
	\theta_1 < \theta_0 < \frac{1}{2} \arctan (0.05). \label{def:theta}
\end{eqnarray}
\end{enumerate}

{\bf Step 2: Verifying $\xi_1 \leq \eta_1$.}
We divide into two sub-cases:
\begin{enumerate}
\item $\R_1$: $\epsilon \leq v \leq 0.8$;
\item $\R_3'$: $0 \leq v < \epsilon$.
\end{enumerate}

We first deal with sub-case 1.
By the discussion following~(\ref{ineq:d}), the minimum horizontal distance
between $\gamma_1$ and $\gamma_3$, defined in~(\ref{def:gamma1})
and~(\ref{def:gamma3}), is at least $0.19$ for $v\in [\epsilon,0.8]$ and
sufficiently large $c$. Therefore the region $\R_1$ is more than
$\delta = 0.1$ left of $\gamma_1 \cap \{(u,v): \epsilon \leq v \leq 0.8\}$.
Thus by~(\ref{est:eta1}),
\begin{eqnarray}
	\eta_1>2\delta cv \label{eq:case1:req1}
\end{eqnarray}
in $\R_1$. On the line segment
$\gamma_3 \cap \R_1 = \{(u,v): u=v, v\in [\epsilon,0.8]\}$,
$\xi$ is defined by~(\ref{def:xi}). Condition~(\ref{eq:pick1}) then implies
\begin{eqnarray*}
	\xi_1 = v F_1 < \delta c v.
\end{eqnarray*}
Thus~(\ref{eq:case1:req1}) shows that
$\xi_1 \leq \eta_1$ on the line segment $\gamma_3 \cap \R_1$.

We also need to
verify that $\xi_1 \leq \eta_1$ everywhere in $\R_1$. For that,
recall that $\xi$ in $\R_1$ is constant along line segments
$\overline{A'B'}(v,v)$, each of which makes an angle of
$-\theta_1$ with the positive $u$-axis, so we
estimate how much the $v$-coordinate can decrease
along the line segments $\overline{A'B'}$, to make sure that $\xi \leq \eta$
even on the line $u=v+0.1$. We observe the following:
since $\tan\theta_1<0.05$, the amount by which the $v$-coordinate decreases,
from the point $A'(v_0,v_0)$ on the line $u=v$ to the point $B'$
on the line $u=v+0.1$,
is $0.1 \tan\theta_1 \leq \frac{\epsilon}{2} \leq \frac{v_0}{2}$, where
we will pick $\epsilon=0.24$ in Chapter~\ref{sec:pde}.
Thus even for $(u_0,v_0)$ lying on the line $u=v+0.1$,
we still have $\overline{A'B'}(u_0,v_0)=\overline{A'B'}(v_1,v_1)$
(i.e. $(v_1,v_1)$ lies on the line $u=v$)
for some $v_1$ with $v_1<2 v_0$.
Because $\xi$ is constant along $\overline{A'B'}(u_0,v_0)$, we have
\begin{eqnarray*}
	\xi_1(u_0,v_0)=v_1 F_1< 2 v_0 F_1<2\delta c v_0,
\end{eqnarray*}
where we use the requirement~(\ref{eq:pick1}) in the last inequality.
Hence for all $(u_0,v_0)\in \R_1$,
\begin{eqnarray*}
	\xi_1(u_0,v_0)<\eta_1(u_0,v_0)
\end{eqnarray*}
by~(\ref{eq:case1:req1}).

Now we deal with sub-case 2. For $v < \epsilon < \frac{1}{4}$ and $u\leq 2v$,
we have
\begin{eqnarray}
        \eta_1 & = & (2c(1-u)+1)v-u \geq (2c(1-2v)+1)v-2v \nnb \\
	& = & (2c(1-2v)-1)v \geq cv
        \label{est:eta1R3}
\end{eqnarray}
if $c$ is sufficiently large.
This estimate applies to both Case 1 (this lemma) and a bit later on
Case 2 (Lemma~\ref{lem:case2}), and
shows that on $\gamma_3 \cap \R_3' = \{(u,v): u=v, v\in [0,\epsilon)\}$,
\begin{eqnarray*}
	\xi_1 = v F_1 < \delta c v = 0.1 c v < \eta_1,
\end{eqnarray*}
where we use~(\ref{def1:delta}),~(\ref{eq:pick1}), and~(\ref{est:eta1R3})
in the second, third, and fourth steps, respectively.
For the rest of $\R_3'$, we make the observation that the $v$-coordinate
of any point on $\overline{A'B'}(v_1,v_1)$ is larger than the $v$-coordinate
$v_2$ of the intersection point of $\overline{A'B'}(v_1,v_1)$
and the line $u=2v$, which we obtain by solving $v-v_1=-\tan\theta_1(2v-v_1)$,
i.e. $v_2 = \frac{1+\tan\theta_1}{1+2\tan\theta_1} v_1>\frac{v_1}{2}$.
Therefore for any $(u_0,v_0) \in \overline{A'B'}(v_1,v_1)$,
\begin{eqnarray*}
        \xi_1(u_0,v_0) = v_1 F_1 < 0.1 c v_1 < 0.2 c v_0 < \eta_1,
\end{eqnarray*}
from which we conclude that $\xi_1 < \eta_1$ in $\R_3'$.

{\bf Step 3: Verifying $\xi_2 \leq \eta_2$.}
This is considerably easier than verifying $\xi_1 \leq \eta_1$.
From~(\ref{def:xi}), $\xi_2 = -2v$ on the line segment
$\gamma_3 \cap (\R_1 \cup \R_3')$,
so $\xi_2 \leq \eta_2$ by~(\ref{eq:eta2:easy}).
Along $\overline{A'B'}(v_1,v_1)$, the $v$-coordinate decreases. So
for any $(u_0,v_0) \in \overline{A'B'}(v_1,v_1)$,
\begin{eqnarray*}
        \xi_2(u_0,v_0) = -2 v_1 < -2 v_0 < \eta_2.
\end{eqnarray*}
Thus $\xi_2 < \eta_2$ in $\R_1 \cup \R_3'$.

{\bf Step 4: Defining $a_s$, $b_s$, and $\overline{AB}$}.
The vector field $\xi$ is defined such that any point
$(u_0,v_0) \in \R_1 \cup \R_3'$
moves under $\xi$ at a constant speed (linear in $v_0$) along the line
$\overline{A'B'}(u_0,v_0)$.
Thus any line segment $\overline{OP}$
lying in $\R_1 \cup \R_3' $ remains a line segment
(i.e. does not become a curve)
under the flow $\xi$, and if $(u_1,v_1)$ and $(u_2,v_2)$ are
two points on such a line segment, then the ratio
$\frac{|\Fcal^s_\xi(u_1,v_1)|}{|\Fcal^s_\xi(u_2,v_2)|}$ remains constant.
We define $\overline{AB}{(v_1,v_1)}$, $0 \leq v_1 \leq 0.8$,
to be the line segment that makes an angle of $-2\theta_0$ with the positive
$u$-axis and connects points $A=(v_1,v_1) \in \gamma_3$ and $B$, with
$B$ lying on the right/bottom boundary of the 4-gon $\R_1\cup \R_3'$.
Recall from~(\ref{def:theta}) that $\theta_0$ is a small angle
and $\overline{A'B'}(a_0,b_0)$ makes an angle of $-\theta_1$
with the positive $u$-axis, where $0 < \theta_1 < \theta_0$. Thus
the part of $\overline{A'B'}(a_0,b_0)$ to the right of the point $(a_0,b_0)$
lies strictly above $\overline{AB}(a_0,b_0)$, and the angle between
$\overline{A'B'}(a_0,b_0)$ and $\overline{AB}(a_0,b_0)$ is at least $\theta_0$
by the choice of $\theta_0$ and $\theta_1$ in~(\ref{def:theta}).
Also, we define
\begin{eqnarray*}
	(a_s,b_s) = \overline{AB}{(a_0,b_0)} \cap
		\overline{O\Fcal^s_\xi(a_0,b_0)},
\end{eqnarray*}
where $(a_0,b_0) \in \R_0 \cap \R_1'$ is the top tip of any line segment
that we consider for this lemma.
We collect various facts for later use:
\begin{enumerate}
\item $\Fcal_\xi$ moves the point $(a_0,b_0)$ to the right.
\item The part of $\overline{A'B'}(a_0,b_0)$ to the right of the point
	$(a_0,b_0)$ lies strictly above $\overline{AB}(a_0,b_0)$.
\item $(a_s,b_s)$ lies on $\overline{AB}(a_0,b_0)$.
\item $\Fcal^s_\xi(a_0,b_0)$ lies on $\overline{A'B'}(a_0,b_0)$
\item $(a_s,b_s)$ and $\Fcal^s_\xi(a_0,b_0)$ both lie on
	$\overline{O\Fcal^s_\xi(a_0,b_0)}$.
\end{enumerate}
The above facts imply that $b_s$ is a lower bound for
$\pi_v(\Fcal^s_\xi(a_0,b_0))$.

We now estimate the speed at which $\Fcal^s_\xi(a_0,b_0)$ separates from
$(a_s,b_s)$. First of all, since $\overline{A'B'}(a_0,b_0)$ makes a negative
angle with the positive $u$-axis, the intersection point $(v_1,v_1)$ of
$\overline{A'B'}(a_0,b_0)$ and the line $u=v$ must lie above $(a_0,b_0)$,
i.e.
\begin{eqnarray*}
	v_1 \geq b_0.
\end{eqnarray*}
This means that
\begin{eqnarray*}
        |\xi(a_0,b_0)| = |\xi(v_1,v_1)| = |(v_1 F_1,-2v_1)|
	\geq b_0 \sqrt{F_1^2+2^2}.
\end{eqnarray*}
Let
\begin{eqnarray*}
	F_2(b_0) = b_0 \sqrt{F_1^2+2^2}.
\end{eqnarray*}

Since $\overline{A'B'}(a_0,b_0)$ and $\overline{AB}(a_0,b_0)$
make angles of $-\theta_1$ and $-2\theta_0$ with the positive
$u$-axis, respectively, where $0 < \theta_1 < \theta_0$,
the angle $\theta_2$ between $\overline{A'B'}(a_0,b_0)$
and $\overline{AB}(a_0,b_0)$ is larger than $\theta_0$.
Let $\alpha=s F_2(b_0)$ be the distance between
$(a_0,b_0)$ and $\Fcal^s_\xi(a_0,b_0)$,
and $\beta$ be the distance between $(a_0,b_0)$ and $(a_s,b_s)$, then
the Euclidean distance $\gamma$ between $\Fcal^s_\xi(a_0,b_0)$ and $(a_s,b_s)$
(the thick line in figure~\ref{fig:case1b}) is
$\gamma = \sqrt{\alpha^2 + \beta^2 - 2 \alpha \beta \cos\theta_2}$,
which attains the minimum $u \sin\theta_2$ when $\beta=\alpha \cos\theta_2$,
therefore $\gamma\geq s F_2(b_0)\sin\theta_0$.
Since $s$ is small and $(a_0,b_0)$ lies in $\R_0 \cap \R_1'$ (in particular,
to the left of the line $u-v=0.09$), $\Fcal^s_\xi(a_0,b_0)$ lies in
$\{(u,v)\in \R: u-v \in [0,0.1], v\geq 0.5\}$, thus
the smallest angle between $\overline{O\Fcal^s_\xi(a_0,b_0)}$
(portion of which is $\overline{(a_s,b_s)\Fcal^s_\xi(a_0,b_0)}$)
and the positive $u$-axis is greater than
$\arctan \frac{0.5}{0.6} > \frac{\pi}{6}$.
Therefore the vertical distance between
$\Fcal^s_\xi(a_0,b_0)$ and $(a_s,b_s)$ is at least
$s F_2(b_0)\sin\theta_0 \sin\frac{\pi}{6}$.
Similarly, $\Fcal^s_\xi(a_0,b_0)$ lies in $\R$, so
the largest angle between $\overline{O\Fcal^s_\xi(a_0,b_0)}$
and the positive $u$-axis
is less than $\frac{\pi}{4}$, hence
the horizontal distance between these two points
is at least $s F_2(b_0)\sin\theta_0 \cos\frac{\pi}{4}$, which is
larger than $s F_2(b_0)\sin\theta_0 \cos\frac{\pi}{3}$.
{
\psfrag{A}{$A$}
\psfrag{B}{$B$}
\psfrag{Ap}{$A'$}
\psfrag{Bp}{$B'$}
\psfrag{u=v}{$u=v$}
\psfrag{u=v+.1}{$u=v+0.1$}
\psfrag{P1}{$\Fcal^s_\xi(a_0,b_0)$}
\psfrag{P2}{$(a_s,b_s)$}
\psfrag{P6}{$(a_0,b_0)$}
\psfrag{tt}{$\theta_2$}
\psfrag{tt2}{$\leq \frac{\pi}{4}, \geq \frac{\pi}{6}$}
\psfrag{alpha}{$\alpha$}
\psfrag{beta}{$\beta$}
\psfrag{gamma}{$\gamma$}
\begin{figure}[h!]
\centering
\includegraphics[height = 3in, width = 5.5in]{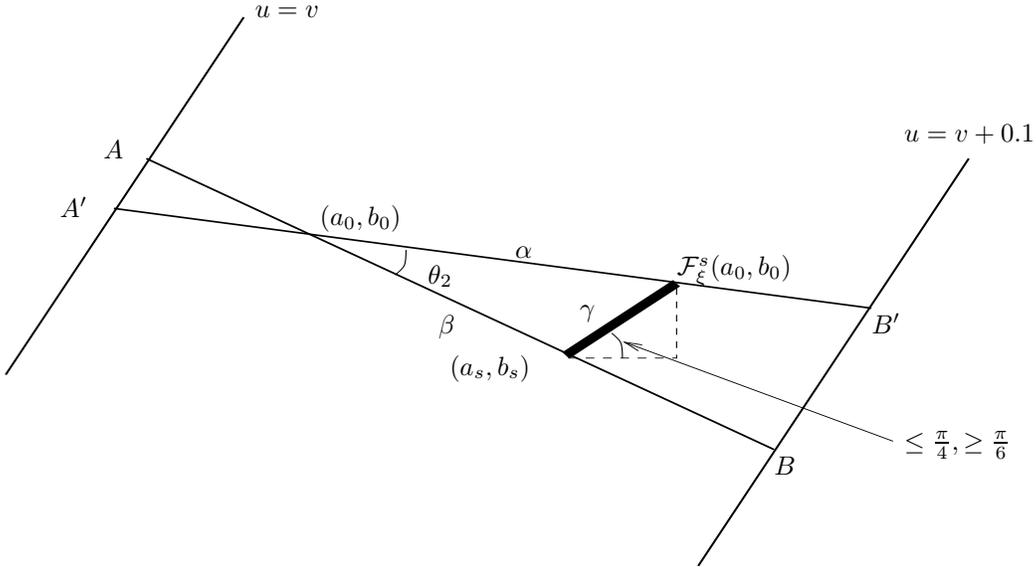}
\caption{Another illustration of Lemma~\ref{lem:case1} (Case 1)}
\label{fig:case1b}
\end{figure}
}

More precisely,
\begin{eqnarray*}
        \Fcal^s_\xi(a_0,b_0) - (a_s,b_s) & \geq &
        	\left(s F_2(b_0)\sin\theta_0 \cos\frac{\pi}{3},
        	s F_2(b_0)\sin\theta_0 \sin\frac{\pi}{6} \right) \\
	& = &  \left(\frac{s}{2} b_0 \sqrt{F_1^2+2^2} \sin\theta_0,
                \frac{s}{2} b_0 \sqrt{F_1^2+2^2} \sin\theta_0 \right)
\end{eqnarray*}
Since $b_s \leq b_0$, the above inequality implies
\begin{eqnarray}
        \Fcal^s_\xi(a_0,b_0) - (a_s,b_s) & \geq &
                \left(\frac{s}{2} b_s \sqrt{F_1^2+2^2} \sin\theta_0,
                \frac{s}{2} b_s \sqrt{F_1^2+2^2} \sin\theta_0 \right).
	\label{ineq:case1:compute}
\end{eqnarray}
By the choice of a small $\theta_0$ in~(\ref{def:theta}), the
entire line segment $\overline{AB}(a_0,b_0)$ lies above the horizontal
line $v=0.5$ if $(a_0,b_0) \in \R_0 \cap \R_1'$. Thus $(a_s,b_s)$ lies
above the horizontal line $v=0.5$, which means that $a_s \leq 2 b_s$.
Thus~(\ref{ineq:case1:compute}) implies
\begin{eqnarray}
	\Fcal^s_\xi(a_0,b_0) - (a_s,b_s) & \geq &
		\left(\frac{s}{4} a_s \sqrt{F_1^2+2^2} \sin\theta_0,
                \frac{s}{4} b_s \sqrt{F_1^2+2^2} \sin\theta_0 \right).
	\label{ineq:case1:compute2}
\end{eqnarray}
If we define
\begin{eqnarray*}
	\tilde{K} = \frac{1}{4} \sqrt{F_1^2+2^2} \sin \theta_0,
\end{eqnarray*}
then~(\ref{ineq:xieta}) and~(\ref{ineq:case1:compute2})
verifies condition~(\ref{eq:case1})
for the point $(a_0,b_0)$. We recall from~(\ref{eq:pick1}) that $F_1$
can be chosen to be arbitrarily large if we allow $c$ to be large.
Therefore $\tilde K$ can also be chosen to be arbitrarily large
if $\theta_0$ is fixed.

For any other point $(\alpha a_0,\alpha b_0)$ on the line segment
$\overline{O(a_0,b_0)}$, $\alpha \in [0,1]$, notice that from~(\ref{def:xi}),
$\xi(\alpha a_0,\alpha b_0) = \alpha \xi(a_0,b_0)$, i.e. $\xi$ is
linear on the line $\overline{O(a_0,b_0)}$. Also, the entire line segment
$\overline{O(a_0,b_0)}$ lies in $\R_1 \cup \R_3'$.
Thus by~(\ref{ineq:xieta}), condition~(\ref{eq:case1})
holds for all $\alpha \in [0,1]$.
\qed

\vspace{.3cm}
\noindent {\bf Proof of Lemma~\ref{lem:case2}}. \hspace{2mm}
We define
\begin{eqnarray}
	\R_2 & = & \{(u,v)\in \R: v+0.02<u<\frac{(1+2c)v}{1+2cv}-0.04,
		v\in [\epsilon,0.8]\}, \label{def:R2} \\
	\R_2' & = & \{(u,v)\in \R_2: v+0.08<u<\frac{(1+2c)v}{1+2cv}-0.04,
		v > 0.55 \}, \label{def:R2p} \\
	\R_3 & = & \{(u,v)\in \R: u\leq 2v,v\in [0,\epsilon)\}. \label{def:R3p}
\end{eqnarray}
For this lemma, the region $\R_2'$ is the region for $P$,
the top tip of the line segment $\overline{OP}$ that connects
the origin $O$ and the point $P$. Since we pick $\epsilon=0.24$ later
in Chapter~\ref{sec:pde}, any line segment $\overline{OP}$,
with $P \in \R_2'$, lies entirely in $\R_2 \cup \R_3$.
We follow the same steps as in the proof of
Lemma~\ref{lem:case1} (Case 1).

{
\psfrag{A}{$A$}
\psfrag{B}{$B$}
\psfrag{C}{$C$}
\psfrag{u}{$u$}
\psfrag{v}{$v$}
\psfrag{u=v}{$u=v$}
\psfrag{u=v+1}{$u=v+0.02$}
\psfrag{R2}{$\R_2$}
\psfrag{R3p}{$\R_3$}
\psfrag{F}{$\zeta(a_0,b_0)$}
\psfrag{-2v}{$\zeta(\alpha a_0,\alpha b_0)$}
\psfrag{P5}{$(a_0,b_0)$}
\psfrag{P6}{$(\alpha a_0,\alpha b_0)$}
\psfrag{eee}{$\epsilon$}
\begin{figure}[h!]
\centering
\includegraphics[height = 5in, width = 5in]{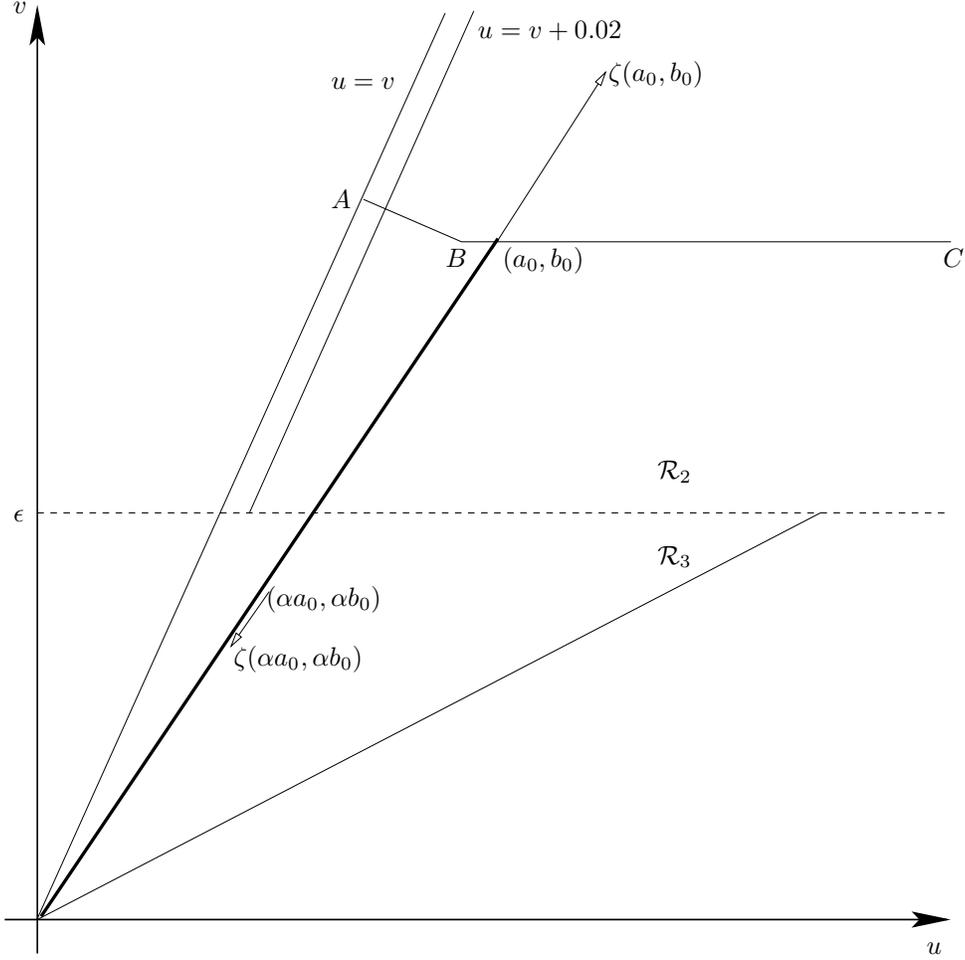}
\caption{Illustration of Lemma~\ref{lem:case2} (Case 2):
$\zeta(\alpha a_0,\alpha b_0) = (-2 \alpha a_0, -2 \alpha b_0)$ and
$\zeta(a_0,b_0) = (a_0 F_3, b_0 F_3)$.}
\label{fig:case2a}
\end{figure}
}

{\bf Step 1: Defining $\zeta$.}
As in the the proof of Case 1, we will define
a vector field $\zeta$, such that $\zeta_1 \leq \eta_1$ and
$\zeta_2 \leq \eta_2$ everywhere in $\R_2 \cup R_3$.
First, we define
\begin{eqnarray*}
	\zeta(1,1)=(F_3,F_3),
\end{eqnarray*}
where we pick $F_3$ large but with
\begin{eqnarray}
	0 < F_3 < 0.01 c. \label{eq:pick2}
\end{eqnarray}
It is convenient to define $\zeta$ at the
point $(1,1)$, even though this point is not even in $\R_2 \cup R_3$.
For $(u,v) \in \R_2$, we define
\begin{eqnarray}
	\zeta(u,v)=(u F_3,v F_3). \label{eq:pick2a}
\end{eqnarray}
And for $(u,v)\in \R_3$,
\begin{eqnarray}
	\zeta(u,v)=(u(-2),v(-2)). \label{eq:pick2b}
\end{eqnarray}
Notice that $\zeta$ is discontinuous across the horizontal line $v=\epsilon$.

{\bf Step 2: Verifying $\zeta \leq \eta$.}
The region $\R_2$ is at least 0.01 to the right of $\gamma_2$
(the line $u=v+2/c$) for sufficiently large $c$. By~(\ref{est:eta2}), we have
\begin{eqnarray}
        \eta_2>0.01 c v \label{eq:case2:req2}
\end{eqnarray}
for all points in $\R_2$. The region $\R_2$ is also at least 0.04 left
of $\gamma_1$, where $\gamma_1$ is the curve of $u=\frac{(1+2c)v}{1+2cv}$.
Thus by~(\ref{est:eta1}), we have
\begin{eqnarray}
        \eta_1>0.08 c v>0.08 c \epsilon = 0.08 c (0.24) > 0.01 c.
	\label{eq:case2:req1}
\end{eqnarray}
for all points in $\R_2$.
By~(\ref{eq:pick2}) and~(\ref{eq:pick2a}),
\begin{eqnarray*}
	\zeta(u,v) < (0.01 c u, 0.01 c v) < (0.01 c, 0.01 c v),
\end{eqnarray*}
which implies, by~(\ref{eq:case2:req2}) and~(\ref{eq:case2:req1}),
\begin{eqnarray*}
	\zeta(u,v) < \eta(u,v)
\end{eqnarray*}
for $(u,v) \in \R_2$.

Recall from~(\ref{eq:gamma1:cbig}) that the curve $\gamma_1$
approaches the degenerate curve
\[ \{(u,v): v=0, u \in [0,1]\} \cup \{(u,v): u=1, v \in [0,1]\} \]
as $c\rightarrow \infty$. Therefore $\R_3$ stays to the left of $\gamma_1$
if $c$ is sufficiently large, and by the discussion below~(\ref{def:gamma3})
regarding the sign of $\eta_1$, $\eta_1>0$ for $(u,v) \in \R_3$. The
definition of $\zeta$ in~(\ref{eq:pick2b}) says that $\zeta_1 < 0$
for $(u,v) \in \R_3$, therefore
\begin{eqnarray*}
	\zeta_1 (u,v) < \eta_1 (u,v)
\end{eqnarray*}
for $(u,v) \in \R_3$.
Furthermore,~(\ref{eq:eta2:easy}) implies that $\zeta_2 (u,v) < \eta_2 (u,v)$.
Thus for all $(u,v) \in \R_3$.
\begin{eqnarray*}
        \zeta(u,v) < \eta(u,v).
\end{eqnarray*}

{\bf Step 3: Defining $\overline{BC}$}.
We define $\overline{BC}{(v_0,v_0)}$ to be the \emph{horizontal} line segment
$v=v_1$ starting at point $B(v_0,v_0)=(v_1+0.1,v_1)$ on the line $u=v+0.1$
and ending on the vertical line $u=1$; notice that $(v_0,v_0)$ itself does
not lie on the line segment $\overline{BC}{(v_0,v_0)}$. This definition
of $\overline{BC}(v_0,v_0)$ means that $B(v_0,v_0)=(v_1+0.1,v_1)$ is the
right end point of the line segment $\overline{AB}{(v_0,v_0)}$, which
was defined in Step 4 of the proof of Lemma~\ref{lem:case1}. We also define
\begin{eqnarray}
	K = \frac{F_3}{2}. \label{def:K}
\end{eqnarray}
Notice that $K$ can be made arbitrarily large since $F_3$
(picked in~(\ref{eq:pick2})) is allowed to be arbitrarily large.

The vector field $\zeta$ is defined such that any point $(u,v)\in \R_2$
moves in the direction $\overrightarrow{O(u,v)}$, i.e. $\zeta$ is a dilation
for points in $\R_2$. But any point $(u,v) \in \R_3$ moves in the 
direction of $\overrightarrow{(u,v)O}$, i.e. $\zeta$ is a contraction
for points in $\R_3$.
Thus any line segment $\overline{OP}$ with $P\in \R_2'$ immediately
splits into two line segments under $\zeta$; the two line segments, however,
lie on the same straight line through the origin $O$.
For the top tip of the line segment $P = (a_0,b_0) \in \R_2'$, we have
\begin{eqnarray}
        \Fcal^s_\zeta(a_0,b_0) - (a_0,b_0)
		\geq \left(s a_0 F_3, s b_0 F_3 \right),
		\label{ineq:case2:upper}
\end{eqnarray}
since the fact that both $u$- and $v$- coordinates increases under $\zeta$
in $\R_2$ implies that $s a_0 F_3$ and $s b_0 F_3$ are lower bounds
for the increase in $u$- and $v$- coordinates, respectively. Since $b_0 > 0.5$
if $(a_0,b_0) \in \R_2'$,~(\ref{def:K}) and~(\ref{ineq:case2:upper})
implies that
\begin{eqnarray}
        \Fcal^s_\zeta(a_0,b_0) - (a_0,b_0)
                \geq \left(s \frac{a_0}{b_0} b_0 F_3, s b_0 F_3 \right)
		\geq \left(\frac{a_0}{b_0} K s, K s\right).
\end{eqnarray}
This verifies~(\ref{eq:case21}) for $P = (a_0,b_0) \in \R_2'$. For any
other point $(\alpha a_0, \alpha b_0)$ on $\overline{OP}$ that lies in $\R_2$
(i.e. $\alpha \in [\frac{\epsilon}{b_0},1]$), linearity in the definition
$\zeta(\alpha a_0, \alpha b_0) = \alpha \zeta(a_0,b_0)$ in~(\ref{eq:pick2a})
implies~(\ref{eq:case21}).

The verification of~(\ref{eq:case22}) for points in $\R_3$ is similar.
Recall from~(\ref{eq:pick2b}) the definition of $\zeta$ in $\R_3$:
\begin{eqnarray*}
	\zeta(u,v) = (-2u,-2v).
\end{eqnarray*}
Let $P = (a_0,b_0) \in \R_2'$ and $\alpha \in [0,\frac{\epsilon}{b_0})$.
Then both $\alpha a_0$ and $\alpha b_0$ decrease under $\zeta$, initially
at speed $2\alpha a_0$ and $2\alpha b_0$, respectively. The speed
of decrease immediately becomes smaller than $2\alpha a_0$ and $2\alpha b_0$
(respectively) after the initial movement. Thus $2\alpha a_0$
and $2\alpha b_0$ are upper bounds of the speed of decrease:
\begin{eqnarray*}
        \Fcal^s_\zeta(\alpha a_0,\alpha b_0) - (\alpha a_0,\alpha b_0)
                \geq \left(- 2 \alpha a_0 s, -2 \alpha b_0 s \right),
\end{eqnarray*}
as required by~(\ref{eq:case22}). \qed

\vspace{.3cm}
To summarize the results in Lemmas~\ref{lem:case1} (Case 1)
and~\ref{lem:case2} (Case 2),
if the top tip of the line segment $\overline{OP}$
at time 0, $P=(a_0,b_0)$ with $b_0 > 0.55$, 
lies to the left of the line $u=v+0.09$, then we use Case 1 to define
$(a_s,b_s)\in \overline{AB}(a_0,b_0)$ and $\xi \leq \eta$ where
$(a_s,b_s)$ is below but to the right of $(a_0,b_0)$, such that $\xi$ moves
$(a_0,b_0)$ at an arbitrarily large speed below and to the right of
$(a_0,b_0)$, but above $(a_s,b_s)$.
Once the top tip of the line segment has moved
to the right of the line $u=v+0.08$ but to the left of the curve
$u = \frac{(1+2c)v}{1+2cv}-0.04$, (or it lies between those two
at time 0 to start with), then using Case 2, where
$(a_s,b_s)\equiv (a_0,b_0)$, we define $\zeta \leq \eta$ such that $\zeta$
moves $(a_0,b_0)$ above and to the right of $(a_0,b_0)$, in fact,
along the same direction as $\overrightarrow{O(a_0,b_0)}$,
again at an arbitrarily large speed. Finally, if $(a_0,b_0)$ lies to
the right of the curve $u = \frac{(1+2c)v}{1+2cv}-0.04$, then we
move the initial condition to the left of this curve and apply Case 2.

\subsection{Analysis of the PDE~(\ref{eq:pde:real})}
\label{sec:pde}
Now we use the results obtained in the previous section about the evolution
of the ODE~(\ref{eq:ode}), together with some results on the heat equation,
to study the evolution of the PDE~(\ref{eq:pde:real}).
First, we need to characterize
how values in the transition region evolve according to the heat equation.
We will establish two technical lemmas to that end.

\begin{LEM}
	If $l$ is fixed and $f = f_0$ is as defined in (IC 1-3)
	on page~\pageref{page:IC}, then for
\begin{eqnarray*}
	x\in (-L-l-s,-L-\frac{l}{200}) \cup (L+\frac{l}{200},L+l+s)
\end{eqnarray*}
	and $s$ small, we have
\begin{eqnarray}
	e^{s \Delta} f(x)>f(x)+\frac{s}{5l^2}. \label{ineq:lem3}
\end{eqnarray}
\label{lem:heat}
\end{LEM}

\proof
First, we shift $f$ right by $L$ such that $f(x)=h(x)$ for
$x\in [-l,l]$ and $f(x)=h(2L-x)$ for $x\in [2L-l,2L+l]$. Thus
$\Delta f(x) = \Delta h(x)$ for $x\in (-l,l)$ and
$\Delta f(x) = \Delta h(2L-x)$ for $x\in (2L-l,2L+l)$, where
\begin{eqnarray}
        \Delta h (x) = \left \{ \begin{array}{ll}
                0,              & x<-l \\
                \frac{1}{l^2},  & -l<x<0 \\
                -\frac{1}{l^2}, & 0<x<l \\
                0               & x>l
        \end{array} \right. . \label{eq:Deltah}
\end{eqnarray}
We make the following observation to aid our computation:
if $u=e^{s\Delta} h$ gives the evolution of the heat equation
with initial condition $h$, then $\Delta u = \Delta (e^{s\Delta} h)$ gives the
evolution of the heat equation with initial condition $\Delta h$, i.e.
\begin{eqnarray}
	\Delta (e^{s\Delta} h) = e^{s\Delta} (\Delta h). \label{eq:observ:D}
\end{eqnarray}
Define
\begin{eqnarray*}
	k(x)=\Delta f(x),
\end{eqnarray*}
then
\begin{eqnarray}
	k(x) = \left \{ \begin{array}{ll}
                \Delta h(x),	& x \in (-l,l)  \\
                \Delta h(2L-x), & x \in (2L-l,2L+l) \\
                0               & \mbox{otherwise}
        \end{array} \right. . \label{eq:k}
\end{eqnarray}
By equation (5.5.10) in [Taylor 1996], the solution of the heat equation
can be expressed in terms of an integral.
\begin{eqnarray}
        e^{s\Delta} k (x) & = & \int_{-\infty}^\infty
        	\frac{1}{\sqrt{4\pi s}} e^{-\frac{y^2}{4s}} k(x-y) dy.
		\label{eq:heat1}
\end{eqnarray}
Using the above formula and the expression of $k$ in~(\ref{eq:k}), we can
estimate $e^{s\Delta} k (x)$ for $x \in (-\frac{5}{4}l,-l]$ and $s$ small:
\begin{eqnarray}
        \lefteqn{ e^{s\Delta} k (x)
	= \frac{1}{l^2} \left(
	\int_{-x-l}^{-x} \frac{e^{-\frac{y^2}{4s}}}{\sqrt{4\pi s}} dy
        - \int_{-x}^{-x+l} \frac{e^{-\frac{y^2}{4s}}}{\sqrt{4\pi s}} dy
	- \int_{-x-l+2L}^{-x+2L} \frac{e^{-\frac{y^2}{4s}}}{\sqrt{4\pi s}} dy
	\right. } \nnb \\
	& & \ \ \ \ \ \left.
	+ \int_{-x+2L}^{-x+l+2L} \frac{e^{-\frac{y^2}{4s}}}{\sqrt{4\pi s}} dy
	\right) \nonumber \\
	& \geq & \frac{1}{l^2} \left(
        \int_{-x-l}^{\infty} \frac{e^{-\frac{y^2}{4s}}}{\sqrt{4\pi s}} dy
	- \int_{-x}^{\infty} \frac{e^{-\frac{y^2}{4s}}}{\sqrt{4\pi s}} dy
        - \int_{-x}^{-x+l} \frac{e^{-\frac{y^2}{4s}}}{\sqrt{4\pi s}} dy
        - \int_{-x-l+2L}^{-x+2L} \frac{e^{-\frac{y^2}{4s}}}{\sqrt{4\pi s}} dy
	\right)	\nnb \\
        & \geq & \frac{1}{l^2} \left(
        \int_{-x-l}^{\infty} \frac{1}{\sqrt{4\pi s}} e^{-\frac{y^2}{4s}} dy
        - 3\int_{-x}^{\infty} \frac{1}{\sqrt{4\pi s}} e^{-\frac{y^2}{4s}} dy
        \right) \nonumber \\
        & \geq & \frac{1}{l^2} \left(
	\int_{|x|-l}^{\infty} \frac{1}{\sqrt{4\pi s}} e^{-\frac{y^2}{4s}} dy
        - 3 \int_{l}^{\infty} \frac{1}{\sqrt{4\pi s}} e^{-\frac{y^2}{4s}} dy
        \right), \label{ineq:heat_on_k}
\end{eqnarray}
where in the last step we use the fact that $x\in (-\frac{5}{4}l,-l]$
implies $|x|=-x \geq l$.
We can take $s$ to be sufficiently small such that
$\int_{l}^{\infty} \frac{e^{-\frac{y^2}{4s}}}{\sqrt{4\pi s}} dy < 10^{-5}/3$,
then with a substitution of variable in the first integral
in~(\ref{ineq:heat_on_k}), we obtain
\begin{eqnarray}
        e^{s\Delta} k (x) & \geq & \frac{1}{l^2} \left(
	\int_{(|x|-l)/\sqrt{s}}^\infty \frac{1}{\sqrt{4\pi}}
	e^{-\frac{y^2}{4}} dy - 10^{-5} \right), \nnb \\
        & = & \frac{1}{l^2} \left(
	\frac{1}{2}- 10^{-5} - \int_0^{(|x|-l)/\sqrt{s}}
        \frac{1}{\sqrt{4\pi}} e^{-\frac{y^2}{4}} dy \right). \label{eq:k1}
\end{eqnarray}
If $x \in (-l-s,-l]$, then $|x|-l < s <\sqrt{s}$ if $s < 1$, and~(\ref{eq:k1})
implies
\begin{eqnarray}
        e^{s\Delta} k (x) \geq \frac{1}{l^2} \left( \frac{1}{2}- 10^{-5}
	- \int_0^1 \frac{1}{\sqrt{4\pi}} e^{-\frac{y^2}{4}} dy \right)
	> \frac{1}{5l^2}. \label{ineq:k1}
\end{eqnarray}
On the other hand, for $x \in (-l,-\frac{l}{200})$ and $s$ small, we also have
\begin{eqnarray}
        e^{s\Delta} k (x) > \frac{1}{5l^2} \label{ineq:k2}
\end{eqnarray}
since for $x \in (-l,l)$, $k(x) = \Delta h(x)$ is a step function with
discontinuity at $0$, where $\Delta h$ is given in~(\ref{eq:Deltah}).

Estimates~(\ref{ineq:k1}) and~(\ref{ineq:k2}) on the behaviour of $k$
under the heat kernel implies that
for $s$ small and $x \in (-l-s,-\frac{l}{200})$,
\begin{eqnarray*}
        \frac{\partial e^{s\Delta} f(x)}{\partial s} = (\Delta(e^{s\Delta} f))
	(x) = (e^{s\Delta} (\Delta f)) (x) > \frac{1}{5l^2},
\end{eqnarray*}
where we use~(\ref{eq:observ:D}) in the second equality. This
establishes~(\ref{ineq:lem3}) for $x \in (-L-l-s,-L-\frac{l}{200})$.
Verification of~(\ref{ineq:lem3}) for $x \in (L+\frac{l}{200},L+l+s)$
is similar.
\qed

\vspace{.3cm}
\begin{LEM}
	Let $t>0$ be fixed and $f_1 = f_0$ be as defined in (IC 1-3)
	on page~\pageref{page:IC}, i.e.
\begin{eqnarray}
        f_1 (x) = \left \{ \begin{array}{ll}
                h(x+L),		& x\in (-L-l,-L+l) \\
                1,		& x\in [-L+l,L-l] \\
		h(L-x),		& x\in (L-l,L+l) \\
                0               & \mbox{otherwise}
        \end{array} \right. .
\end{eqnarray}
	Let
\begin{eqnarray}
        f_3 (x) = \left \{ \begin{array}{ll}
                f_1(x)+mt,      & -L-l-t<x<L+l+t \\
                0,              & \mbox{otherwise}
        \end{array} \right. \label{eq:lemma2}
\end{eqnarray}
where $m > 0$. Then there exist positive
constants $\delta_1$ and $\delta_2$ depending on $m$ but independent of $t$
such that if
\begin{eqnarray}
        f_2 (x) = \left \{ \begin{array}{ll}
		(1+\delta_2 t) h(x+L+\delta_1 t),
                        & x\in (-L-l-\delta_1 t,-L+l-\delta_1 t) \\
                1+\delta_2 t,	& x\in [-L+l-\delta_1 t,L-l+\delta_1 t] \\
		(1+\delta_2 t) h(L+\delta_1 t-x),
                        & x\in (L-l+\delta_1 t,L+l+\delta_1 t) \\
		0,	& \mbox{otherwise}
        \end{array} \right. , \label{eq:lemma2a}
\end{eqnarray}
then $f_2\leq f_3$.
\label{lem:space}
\end{LEM}

\proof
Without any loss of generality, assume $m<1$.
Let $M=1 \wedge \sup_{x\in\mathbb{R}}|f_1'(x)|$, then $M=1 \wedge 1/l = 1/l$
since $l$ will be picked to be $<1$ in~(\ref{eq:pickl}) a bit later. Define
\begin{eqnarray*}
	g_1(x) = \left \{ \begin{array}{ll}
	f_1(x+\frac{mt}{3M}), & x\in (-L-l-\frac{mt}{3M},-L+l-\frac{mt}{3M}) \\
	f_1(0),	& x\in [-L+l-\frac{mt}{3M},L-l+\frac{mt}{3M}] \\
	f_1(x-\frac{mt}{3M}), & x\in (L-l+\frac{mt}{3M},L+l+\frac{mt}{3M}) \\
	0, & \mbox{otherwise}
	\end{array} \right. ,
\end{eqnarray*}
then any small piece of the curve of $g_1$ is $f_1$ shifted by either
$0$, $\frac{mt}{3M}$, or $-\frac{mt}{3M}$, with $\frac{mt}{3M} < t$.
In particular, the two transition
regions in $g_1$ are the two transition regions in $f_1$ shifted by
$\frac{mt}{3M}$ or $-\frac{mt}{3M}$, and the ``middle'' region (i.e. the
region sandwiched between the two transition regions) in $g_1$ is
the middle region of $f_1$ expanded left and right by $\frac{mt}{3M}$.
Since $M=\sup_{x\in\mathbb{R}}|f_1'(x)|$, we have
\begin{eqnarray*}
	g_1(x)-f_1(x)\leq \frac{mt}{3}
\end{eqnarray*}
for all $x\in (-L-l-\frac{mt}{3M},L+l+\frac{mt}{3M})$,
therefore $f_3 \geq g_1$ everywhere and in particular,
since $f_3(x) - f_1(x) = mt$ for $x\in (-L-l-t,L+l+t)$, we have
\begin{eqnarray}
	f_3(x)-g_1(x) \geq \frac{2mt}{3} \label{ineq:f3.g1}
\end{eqnarray}
for $x\in (-L-l-\frac{mt}{3M},L+l+\frac{mt}{3M})$.

Next we define
\begin{eqnarray*}
        f_2(x)=\left(1+\frac{mt}{3}\right) g_1(x).
\end{eqnarray*}
Then
\begin{eqnarray*}
	f_2(x) - g_1(x) = \frac{mt}{3} g_1(x) \leq \frac{mt}{3}
\end{eqnarray*}
since $g_1(x) \leq 1$ everywhere.
The above inequality and~(\ref{ineq:f3.g1}) imply that
\begin{eqnarray*}
	f_2(x) < f_3(x)
\end{eqnarray*} 
for $x \in (-L-l-\frac{mt}{3M},L+l+\frac{mt}{3M})$. Then
$\delta_1=\frac{m}{3M}$ and $\delta_2=\frac{m}{3}$ satisfy
the requirement of $f_2$ in~(\ref{eq:lemma2a}), and the proof
is complete. \qed

\vspace{.3cm}
For the remainder of this section, we will establish Proposition~\ref{prp:Dd}. 
We assume that the initial condition of the PDE~(\ref{eq:pde:real})
is $(a_0 f_0,b_0 f_0)$, where $f_0$ is as defined in (IC 1-3) on
page~\pageref{page:IC} and $(a_0,b_0)=(\tilde a_0,\tilde b_0)$ lies in
the region $\R_0$ defined in~(\ref{def:R0}).
By Remark~\ref{rem:case1} and Lemma~\ref{lem:case2},
we can pick $(a_s,b_s) \in \overline{ABC}{(a_0,b_0)}$,
with $(a_s,b_s)=(a_0,b_0)$ if $a_0 - b_0 > 0.08$
(i.e. Lemma~\ref{lem:case2}/Case 2),
such that estimates~(\ref{eq:case1})-(\ref{eq:case22})
regarding the evolution of the ODE~(\ref{eq:ode})
are valid. We will use this, together
with Lemma~\ref{lem:heat} at the beginning of this section, to show that
there is a positive constant $m$ such that for $x\in (-L-l-s,L+l+s)$
and sufficiently small $s$,
\begin{eqnarray*}
	e^{s\Delta} \Fcal_\eta^s (a_0 f_0,b_0 f_0)(x) - (a_s f_0, b_s f_0)(x)
	\geq \left(a_s ms,b_s ms\right).
\end{eqnarray*}
Finally we will apply Lemma~\ref{lem:space} to complete the proof of
Proposition~\ref{prp:Dd}.

We divide this task into proving two lemmas, which
correspond to the two cases in Lemmas~\ref{lem:case1} and~\ref{lem:case2}, respectively. Before we proceed, we first pick
\begin{eqnarray}
        l=\sqrt{\frac{0.1}{3}} \label{eq:pickl}
\end{eqnarray}
and
\begin{eqnarray}
	\epsilon=0.24<0.5 h\left(-\frac{l}{100}\right) . \label{eq:picke}
\end{eqnarray}

\begin{LEM}
	(Case 1) Recall that
\begin{eqnarray*}
        \R_1' \cup \R_3' = \{(u,v)\in \R: u-v \in [0,0.1] \mbox{ and
                $u\leq 2v$ for $v\in [0,\epsilon)$} \}
\end{eqnarray*}
	and
\begin{eqnarray*}
        \R_1' \cap \R_0 = \{(u,v)\in \R: u-v \in [0,0.09] \mbox{ and
                $v\in [0.55,0.8]$} \}.
\end{eqnarray*}
	If $\{(a_0 f_0(x),b_0 f_0(x)): x\in [-L-l,-L+l]\}
	\subset \R_1' \cup \R_3'$ and $(a_0,b_0) \in R_1' \cap \R_0$,
	where $\R_0$, $\R_1'$ and $\R_3'$ are defined
	in~(\ref{def:R0}),~(\ref{def:R1p}), and~(\ref{def:R1pp}) respectively,
	and $f_0$ is defined in (IC 1-3) on page~\pageref{page:IC},
	then the conclusion of Proposition~\ref{prp:Dd} holds.
\label{lem:case1pde}
\end{LEM}

\proof
By Lemma~\ref{lem:case1} and Remark~\ref{rem:case1},
for sufficiently small $s$,
we can pick $K$ and $\tilde K$ large enough such that
\begin{eqnarray}
	K>\frac{3}{l^2}, \label{eq:pickK}
\end{eqnarray}
and a point $(a_s,b_s)\in \R$ with $b_s>0.5$ such that
$(a_s f_0(0),b_s f(0)) \in \overline{AB}{(a_0 f_0(0),b_0 f_0(0))}$
and
\begin{eqnarray}
	\Fcal^s_\eta (a_0 f_0(x),b_0 f_0(x))\geq
	((1+\tilde{K}s)a_s f_0(x),(1+\tilde{K}s)b_s f_0(x))
	\label{eq:pde:case1}
\end{eqnarray}
for all $x$; furthermore, if $b_s f_0(x)\geq \epsilon$,
\begin{eqnarray}
        (1+\tilde{K}s) b_s f_0(x) - b_s f_0(x) > K s,
		\label{eq:pde:case1c}
\end{eqnarray}
and if $b_s f_0(x)\in [0,\epsilon)$,
\begin{eqnarray}
        (1+\tilde{K}s) b_s f_0(x) - b_s f_0(x) \geq 0. \label{eq:pde:case1d}
\end{eqnarray}
Here, $(a_s f_0,b_s f_0)$ is the function to which we compare 
$\Fcal^s_\eta (a_0 f_0,b_0 f_0)$ to see how much
``progress'' we are making in increasing the $v$-coordinate.

For $x\in \left[-L-\frac{l}{200},L+\frac{l}{200}\right]
		= \left[-L-\frac{l}{200},-L+l\right] \cup (-L+l,L-l)
			\cup \left[L-l,L+\frac{l}{200}\right]$,
where the intervals $\left[-L-\frac{l}{200},-L+l\right]$ and
$\left[L-l,L+\frac{l}{200}\right]$
are in the transition region, we have
\begin{eqnarray*}
	b_s f_0 (x) \geq 0.5 f_0 \left(-L-\frac{l}{200}\right)
	> 0.5 f_0 \left(-L-\frac{l}{100}\right)
	= 0.5 h\left(-\frac{l}{100}\right),
\end{eqnarray*}
therefore by~(\ref{eq:picke}),
\begin{eqnarray}
	b_s f_0 (x) > \epsilon. \label{ineq:middle:eps}
\end{eqnarray}
Therefore by~(\ref{eq:pde:case1c}) we have,
for $x\in [-L-\frac{l}{200},L+\frac{l}{200}]$,
\begin{eqnarray*}
	(1+\tilde{K}s)b_s f_0(x)-b_s f_0(x)>Ks.
\end{eqnarray*}
For $x\in [-L-l,-L-\frac{l}{200}) \cup (L+\frac{l}{200},L+l]$
where $b_s f_0$ is possibly smaller than $\epsilon$,
by~(\ref{eq:pde:case1d}), we have
\begin{eqnarray*}
	(1+\tilde{K}s)b_s f_0(x)-b_s f_0(x) \geq 0.
\end{eqnarray*}
To summarize, combining~(\ref{eq:pde:case1}) and the two inequalities
above, we have
\begin{eqnarray}
	\pi_v (\Fcal^s_\eta (a_0 f_0,b_0 f_0))(x)
		\geq (1+\tilde{K}s)b_s f_0(x) \label{eq:pde:case12}
\end{eqnarray}
and
\begin{eqnarray}
	(1+\tilde{K}s)b_s f_0(x) - b_s f_0(x) \left\{ \begin{array}{ll}
                > Ks ,& x \in [-L-\frac{l}{200},L+\frac{l}{200}] \\
                \geq 0, & x \in [-L-l,-L-\frac{l}{200})
			\cup (L+\frac{l}{200},L+l] \\
                = 0,    & x \notin [-L-l,L+l]
        \end{array} \right. . \label{eq:pde:case11}
\end{eqnarray}

{
\psfrag{ee}{$\epsilon$}
\psfrag{0}{$0$}
\psfrag{1}{$(1+\tilde K s)b_s$}
\psfrag{lminus}{$-l$}
\psfrag{lplus}{$l$}
\psfrag{Lminus}{$-L-\frac{l}{200}$}
\psfrag{Lplus}{$L+\frac{l}{200}$}
\psfrag{L-lminus}{$-L-l$}
\psfrag{L+lminus}{$-L+l$}
\psfrag{L-lplus}{$L-l$}
\psfrag{L+lplus}{$L+l$}
\begin{figure}[h]
\centering
\includegraphics[height = 1.5in, width = 5.5in]{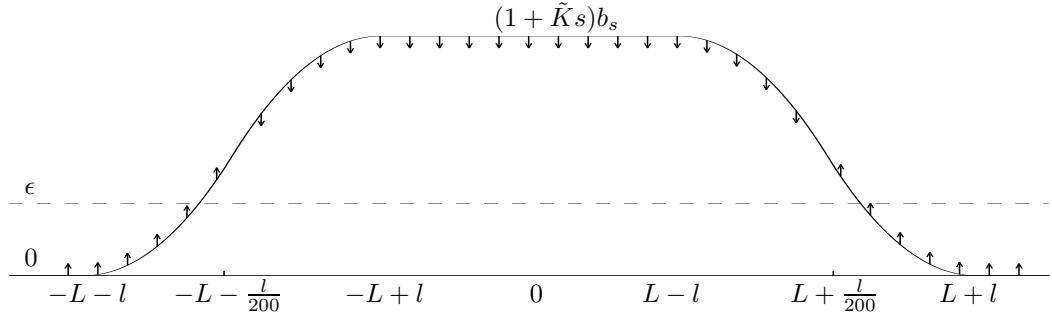}
\caption{The effect of the heat kernel on the function $(1+\tilde Ks) b_s f_0$.
The arrows indicates whether $(1+\tilde Ks)b_s f_0(x)$ increases or
decreases. The effects illustrated here are lower bounds.
In $[-L-\frac{l}{200},L+\frac{l}{200}]$, the function decreases, which is
why we need the first line~(\ref{eq:pde:case11}) to be $>K s$ to cancel
out this decrease.}
\label{fig:pde:case1}
\end{figure}
}
As stated in~(\ref{eq:boundf0}), $|\Delta f_0| \leq \frac{1}{l^2}$.
Therefore the heat operator $e^{s\Delta}$ applied to $b_s f_0$
may decrease its value by at most $\frac{b_s}{l^2}s$. More precisely,
\begin{eqnarray}
	e^{s\Delta} b_s f_0 - b_s f_0 \geq -\frac{b_s}{l^2}s.
	\label{eq:bound:heat}
\end{eqnarray}
everywhere. We can use~(\ref{eq:pde:case11}) and~(\ref{eq:bound:heat}) to
obtain estimates on $e^{s\Delta} \Fcal^s_\eta (a_0 f_0,b_0 f_0)$.
We estimate the ``progress'' made after applying the heat kernel:
by~(\ref{eq:pde:case12})
and the monotonicity of the heat kernel $e^{s\Delta}$,
\begin{eqnarray*}
	\pi_v(e^{s\Delta}\Fcal^s_\eta (a_0 f_0,b_0 f_0) - (a_s f_0,b_s f_0))(x)
	\geq e^{s\Delta}(1+\tilde{K}s)b_s f_0(x)-b_s f_0(x).
\end{eqnarray*}
For $x\in [-L-\frac{l}{200},L+\frac{l}{200}]$,
\begin{eqnarray*}
	\lefteqn{ e^{s\Delta}(1+\tilde{K}s)b_s f_0(x)-b_s f_0(x) } \\
	& = & (1+\tilde{K}s)\left( e^{s\Delta}b_s f_0(x) - b_s f_0(x) \right)
		+ \left( (1+\tilde{K}s)b_s f_0(x) - b_s f_0(x)\right) \\
	& > & -(1+\tilde{K}s)\frac{b_s}{l^2}s+\frac{3b_s}{l^2}s 
\end{eqnarray*}
by~(\ref{eq:bound:heat}), the first line of~(\ref{eq:pde:case11}),
and the fact $K>\frac{3}{l^2}>\frac{3b_s}{l^2}$. Therefore
\begin{eqnarray}
	\pi_v(e^{s\Delta}\Fcal^s_\eta (a_0 f_0,b_0 f_0) - (a_s f_0,b_s f_0))(x)
	> \frac{b_s}{l^2}s \label{eq:pde:case1a}
\end{eqnarray} 
for sufficiently small $s$.
On the other hand, for
$x\in (-L-l-s,-L-\frac{l}{200}) \cup (L+\frac{l}{200},L+l+s)$,
by Lemma~\ref{lem:heat} we have,
\begin{eqnarray*}
	e^{s\Delta}(1+\tilde{K}s)b_s f_0(x)-(1+\tilde{K}s)b_s f_0(x)
	>(1+\tilde{K}s)\frac{b_s}{5l^2}s \geq \frac{b_s}{5l^2}s.
\end{eqnarray*}
Therefore by~(\ref{eq:pde:case1}) and the above inequality,
for $x\in (-L-l-s,-L-\frac{l}{200}) \cup (L+\frac{l}{200},L+l+s)$,
\begin{eqnarray*}
	\pi_v(e^{s\Delta}\Fcal^s_\eta (a_0 f_0,b_0 f_0)-(a_s f_0,b_s f_0))(x)
	& \geq & e^{s\Delta} (1+\tilde{K}s)b_s f_0(x) - b_s f_0(x) \nnb \\
	& > & (1+\tilde{K}s)b_s f_0(x)+\frac{b_s}{5l^2}s-b_s f_0(x) \nnb \\
        & > & \frac{b_s}{5 l^2}s,
\end{eqnarray*}
where the last line is due to the second and third
lines of~(\ref{eq:pde:case11}). Hence for $x\in (-L-l-s,L+l+s)$,
\begin{eqnarray*}
	\pi_v(e^{s\Delta}\Fcal^s_\eta (a_0 f_0,b_0 f_0))(x)
	> b_s \left(f_0(x)+\frac{1}{5 l^2}s \right).
\end{eqnarray*}
Then Lemma~\ref{lem:space} implies that there exist positive constants
$\delta_1$ and $\delta_2$ independent of $s$ such that
\begin{eqnarray}
        \pi_v(e^{s\Delta}\Fcal^s_\eta (a_0 f_0,b_0 f_0))(x) > b_s f_2(x),
	\label{eq:pde:case1b}
\end{eqnarray}
where $f_2$ is defined in~(\ref{eq:lemma2a}):
\begin{eqnarray*}
        f_2 (x) = \left \{ \begin{array}{ll}
                (1+\delta_2 s) h(x+L+\delta_1 s),
                        & x\in (-L-l-\delta_1 s,-L+l-\delta_1 s) \\
                1+\delta_2 s,   & x\in [-L+l-\delta_1 s,L-l+\delta_1 s] \\
                (1+\delta_2 s) h(L+\delta_1 s-x),
                        & x\in (L-l+\delta_1 s,L+l+\delta_1 s) \\
                0,      & \mbox{otherwise}
        \end{array} \right. .
\end{eqnarray*}

Similarly, the estimates in~(\ref{eq:pde:case1a}) to~(\ref{eq:pde:case1b})
also hold for the $u$-coordinate of \\
$e^{s\Delta}\Fcal^s_\eta (a_0 f_0,b_0 f_0) - (a_s f_0,b_s f_0)$,
if $b_s$ on the right hand side of each inequality
is replaced by $a_s$. So for all $x\in (-L-l-s,L+l+s)$, we have
\begin{eqnarray*}
	e^{s\Delta} \Fcal^s_\eta (a_0 f_0,b_0 f_0) (x)
	> (a_s f_2(x),b_s f_2(x)),
\end{eqnarray*}
as required. In particular, $(\tilde a_s,\tilde b_s)$ in the statement of
Proposition~\ref{prp:Dd} should be $((1+\delta_2 s) a_s,(1+\delta_2 s) b_s)$.
\qed

\vspace{.3cm}

\begin{LEM}
        (Case 2) Recall that
\begin{eqnarray*}
        \R_2' & = & \{(u,v)\in \R_2: v+0.08<u<\frac{(1+2c)v}{1+2cv}-0.04,
                v > 0.55 \},
\end{eqnarray*}
	and
\begin{eqnarray*}
        \R_2 \cup \R_3 & = & \{(u,v)\in \R:
                v+ 0.02 < u < \frac{(1+2c)v}{1+2cv} - 0.04 \\
	& & \ \ \ \ \ \mbox{ for
                        $v\in [\epsilon,0.8]$ and
                $u \leq 2v$ for $v \in [0,\epsilon)$} \} .
\end{eqnarray*}
	If $\{(a_0 f_0(x),b_0 f_0(x)): x\in [-L-l,-L+l]\}
        \subset \R_2 \cup \R_3$ and $(a_0,b_0) \in R_2'$,
	where $\R_2$, $R_2'$ and $\R_3$ are defined
        in~(\ref{def:R2}),~(\ref{def:R2p}), and~(\ref{def:R3p})
	respectively, and $f_0$ is defined in (IC 1-3) on
	page~\pageref{page:IC}, then the conclusion of Proposition~\ref{prp:Dd}
	holds.
\label{lem:case2pde}
\end{LEM}

Before we prove this lemma, we observe that the union of all the regions
where $(a_0,b_0)$ may lie is $(\R_1'\cap \R_0) \cup \R_2'$, which
is exactly $\R_0$ as defined in~(\ref{def:R0}).
If $(a_0,b_0)\in \R_1'\cap \R_0$, then the part of the line segment
$\overline{O(a_0,b_0)}$ ($\overline{O(a_0,b_0)}$ consists of values of
$(a_0 f_0,b_0 f_0)$)
above $y=\epsilon$ lies in $\R_1'$. On the other hand,
if $(a_0,b_0)\in  \R_2'$, then the part of the line segment
$\overline{O(a_0,b_0)}$ above $y=\epsilon$ lies in $\R_2$. But in both
these cases, for the part of the line segment $\overline{O(a_0,b_0)}$
below $y=\epsilon$, if suffices to consider
$\R_3 = \{(u,v)\in \R: u\leq 2v,v\in [0,\epsilon)\}$,
because the top tip of $\overline{O(a_0,b_0)}$ in the
$(u,v)$-plane lies above the horizontal line $v=0.5$, where      
$u \leq 1 \leq 2v$. The sufficiency
of restricting to $\{(u,v)\in \R: u < \frac{(1+2c)v}{1+2cv}-0.04\}$
has been discussed below~(\ref{def:R0}) on page~\pageref{page:mono}.

\vspace{.3cm}
\noindent {\bf Proof of Lemma~\ref{lem:case2pde}}. \hspace{2mm}
Under this case, the line segment formed by \\
$\{(a_0 f_0(x),b_0 f_0(x)): x\in [-L-l,-L+l]\}$ lies in $\R_2 \cup R_3$,
i.e. the portion of the line segment above the horizontal line $v=\epsilon$
lies in $\R_2$ and right of the line $u-v=0.02$,
and the portion below $v=\epsilon$ lies in $\R_3$. Furthermore,
the top tip $(a_0,b_0)$ lies in $\R_2'$, i.e.
to the right of the line $u-v = 0.08$ and above the horizontal line $v=0.55$.
For $x\in [-L-\frac{l}{100},L+\frac{l}{100}]$, we have
\begin{eqnarray*}
	b_0 f_0 (x) > 0.55 f_0 \left(-L-\frac{l}{100}\right)
	> 0.5 h\left(-\frac{l}{100}\right) > \epsilon
\end{eqnarray*}
by~(\ref{eq:picke}). Therefore, by Lemma~\ref{lem:case2},
we can construct functions $g_2$ and $g_3$:
\begin{eqnarray}
        g_2(x) = \left \{ \begin{array}{ll}
		b_0 f_0(x)+Ks,  & x\in [-L-\frac{l}{100},L+\frac{l}{100}] \\
                b_0 f_0(x)(1-2s), & x\in 
			(-L-l,-L-\frac{l}{100})\cup (L+\frac{l}{100},L+l) \\
                0,              & \mbox{otherwise}
        \end{array} \right. , \label{eq:g2}
\end{eqnarray}
and
\begin{eqnarray*}
        g_3(x) = \left \{ \begin{array}{ll}
		b_0 f_0(x)+Ks,    & x\in [-L-\frac{l}{200},L+\frac{l}{200}] \\
                b_0 f_0(x)(1-2s), & x\in
			(-L-l,-L-\frac{l}{200})\cup (L+\frac{l}{200},L+l) \\
                0,              & \mbox{otherwise}
        \end{array} \right. ,
\end{eqnarray*}
such that both $(\frac{a_0}{b_0}g_2,g_2)$ and $(\frac{a_0}{b_0}g_3,g_3)$
are lower bounds of $\Fcal^s_\eta (a_0 f_0,b_0 f_0)$.
Notice that $g_3\leq g_2$ everywhere, and $g_2$ has discontinuities at
$-L-l/100$ and $L+l/100$, while $g_3$ has discontinuities at
$-L-l/200$ and $L+l/200$. See figure~\ref{fig:pde:case2} for graphs of
$g_2$ and $g_3$.
{
\psfrag{Lminus}{$-L$}
\psfrag{L-lminus}{$-L-l$}
\psfrag{L+lminus}{$-L+l$}
\psfrag{x0}{$-L-\frac{l}{100}$}
\psfrag{x1}{$-L-\frac{l}{200}$}
\psfrag{0}{$0$}
\psfrag{b0}{$b_0+Ks$}
\begin{figure}[h!!]
\centering
\subfigure[The function $g_2$] {
\includegraphics[height = 1.5in, width = 3in]{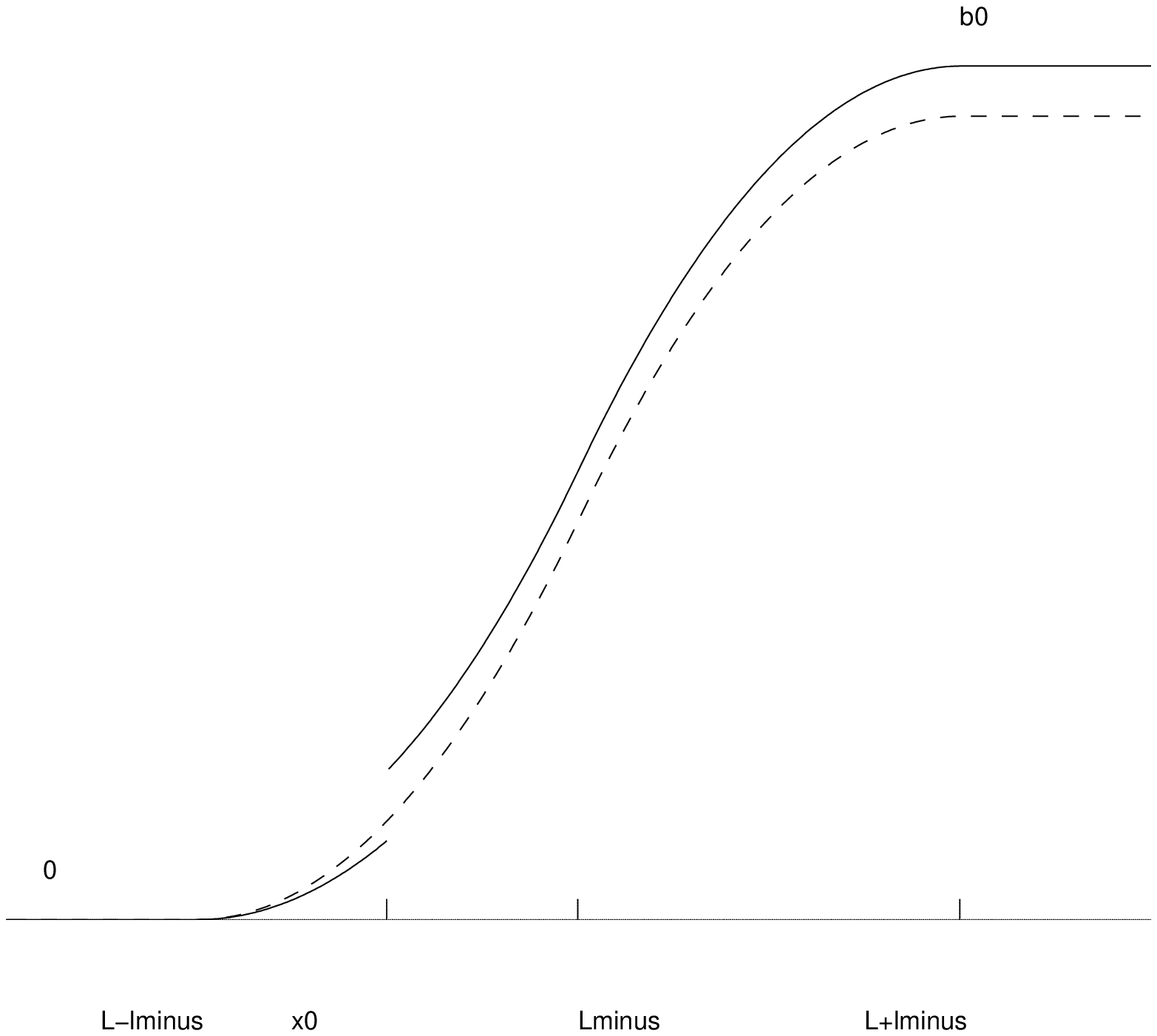}}
\subfigure[The function $g_3$] {
\includegraphics[height = 1.5in, width = 3in]{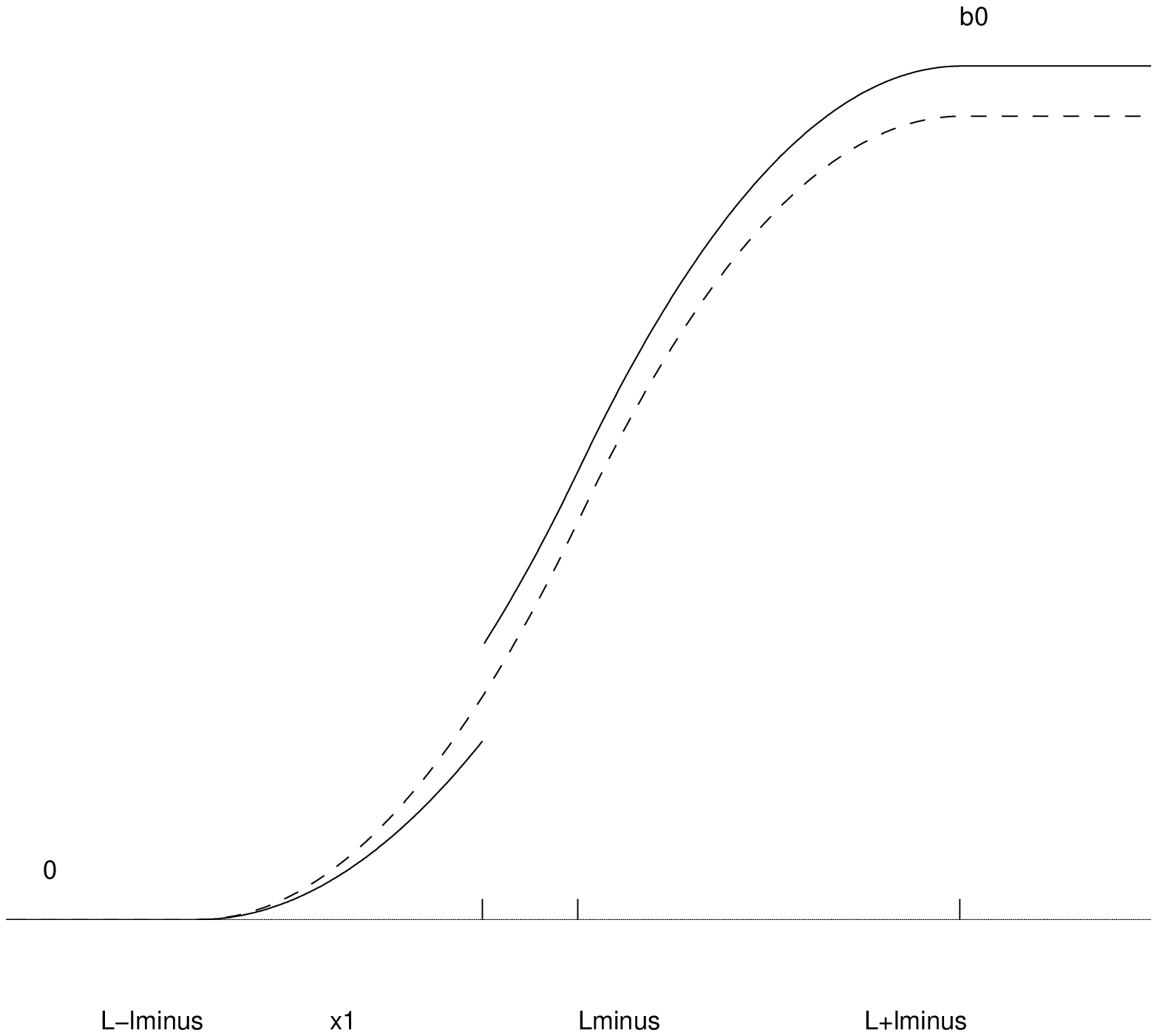}}
\subfigure[The function $g_4$] {
\includegraphics[height = 1.5in, width = 3in]{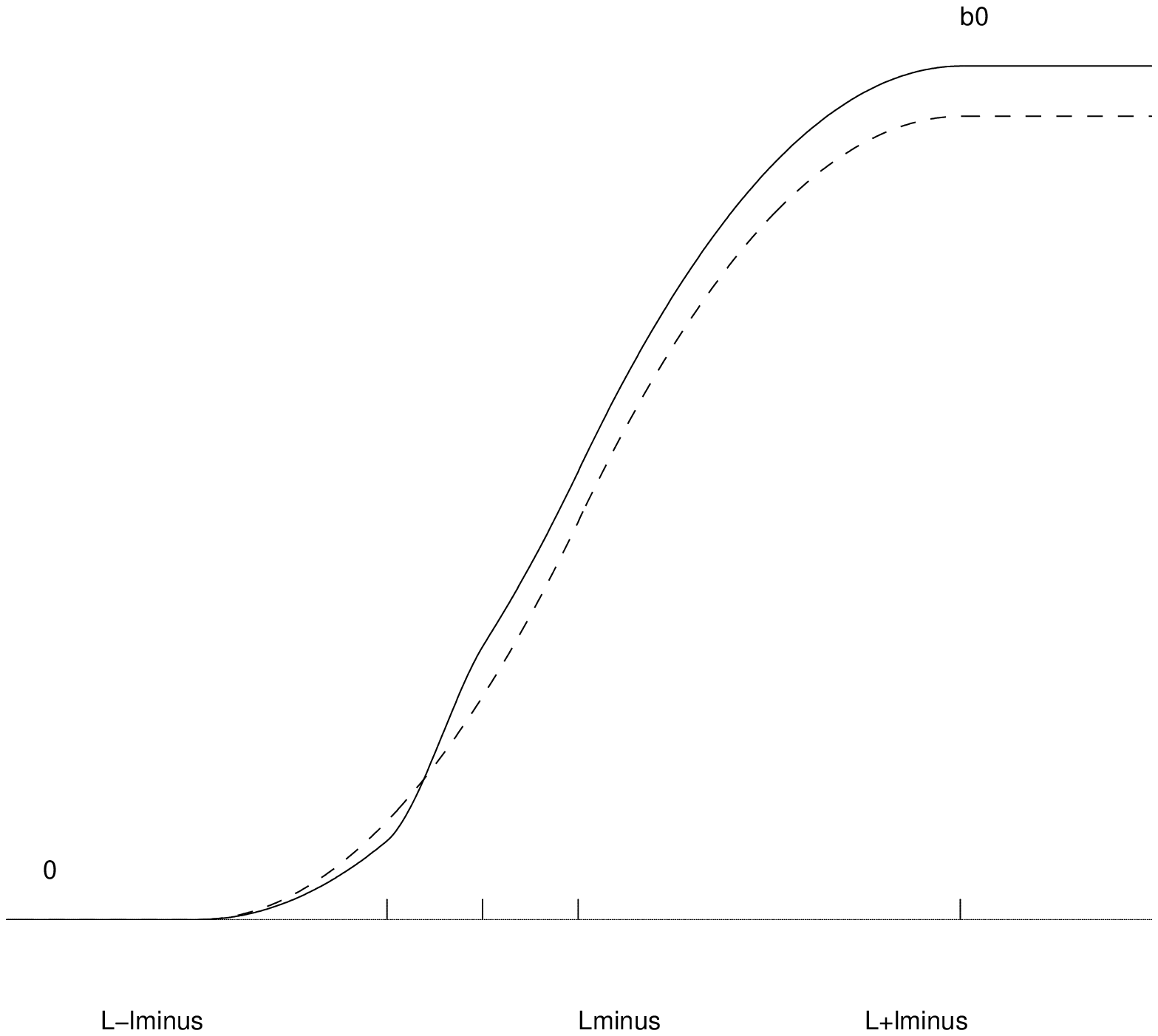}}
\caption{The functions $g_2$, $g_3$, and $g_4$; dotted lines denote the function
$b_0 f_0$.}
\label{fig:pde:case2}
\end{figure}
}

Now we construct the function $g_4$:
\begin{eqnarray}
	g_4(x) = \left\{ \begin{array}{ll}
		g_2(x),	& x\in [-L-\frac{l}{200},L+\frac{l}{200}] \\
		g_3(x), & x\in (-\infty,-L-\frac{l}{100}]
			\cup [L+\frac{l}{100},\infty)
	\end{array} \right. , \label{def:g4}
\end{eqnarray}
furthermore $g_4$ is required to be $C^\infty$, monotone
in $[-L-\frac{l}{200},-L-\frac{l}{100}] \cup
[L+\frac{l}{100},L+\frac{l}{200}]$, and lying between $g_2$ and $g_3$, with
\begin{eqnarray}
	|\Delta g_4(x)| < \frac{2}{l^2} \label{ineq:g4}
\end{eqnarray}
everywhere. For
$x \notin [-L-\frac{l}{200},-L-\frac{l}{100}]
\cup [L+\frac{l}{100},L+\frac{l}{200}]$, the last requirement
is automatic by~(\ref{eq:boundf0}); but for
$x \in [-L-\frac{l}{200},-L-\frac{l}{100}]
\cup [L+\frac{l}{100},L+\frac{l}{200}]$, it can be achieved
for sufficiently small $s$.
Notice that since $g_4 \leq g_2$, $(\frac{a_0}{b_0}g_4,g_4)$
is a lower bound of $\Fcal^s_\eta (a_0 f_0,b_0 f_0)$.
Furthermore, for $x\in [-L-\frac{l}{200},L+\frac{l}{200}]$,
\begin{eqnarray}
	g_4(x) - b_0 f_0(x) = Ks. \label{eq:g4bigger}
\end{eqnarray}
Here $(a_0 f_0,b_0 f_0)$ is the function to which we compare
$\Fcal^s_\eta (a_0 f_0,b_0 f_0)$ to see how much
``progress'' we are making increasing the $v$-coordinate.

We now turn to evolution according to the heat equation. First we deal
with $x\notin [-L-\frac{l}{200},L+\frac{l}{200}]$. For
this, we use $g_3$ as the lower bound for
$\pi_v(\Fcal^s_\eta (a_0 f_0,b_0 f_0))$. We observe that
$g_3$ dominates $(1-2s)b_0 f_0$, therefore monotonicity of the heat
kernel implies
\begin{eqnarray*}
	e^{s\Delta} g_3 \geq e^{s\Delta} ((1-2s) b_0 f_0)
	= (1-2s) b_0 e^{s\Delta} f_0.
\end{eqnarray*}
By Lemma~\ref{lem:heat}, $e^{s\Delta} f_0 (x) > f_0(x) + \frac{s}{5l^2}$
for all $x\in (-L-l-s,-L-\frac{l}{200})\cup (L+\frac{l}{200},L+l+s)$,
therefore
\begin{eqnarray*}
	e^{s\Delta} g_3 (x)> (1-2s) b_0 \left(f_0(x) + \frac{s}{5l^2}\right).
\end{eqnarray*}     
Since $b_0>0.55$, we have $(1-2s)b_0>0.5$ for sufficiently small $s$.
Also recall that we pick $l=\sqrt{\frac{3}{0.1}}$ in~(\ref{eq:pickl}) such that
$\frac{1}{l^2}=\frac{3}{0.1}$.
Thus the inequality above can be written as
\begin{eqnarray*}
        e^{s\Delta} g_3 (x) > (1-2s) b_0 f_0(x) + 0.5 \frac{3s}{0.5}
	= b_0 f_0(x) + (3-2 b_0 f_0(x))s.
\end{eqnarray*}
Finally, since $3-2 b_0 f_0 > 3-2= 1$, we have
\begin{eqnarray}
	e^{s\Delta} g_3 (x) > b_0 f_0(x) + s. \label{eq:pde:case2a}
\end{eqnarray}
for all $x\in (-L-l-s,-L-\frac{l}{200})\cup (L+\frac{l}{200},L+l+s)$.

For $x\in [-L-\frac{l}{200},L+\frac{l}{200}]$, we use
$g_4$ as the lower bound for $\pi_v(\Fcal^s_\eta (a_0 f_0,b_0 f_0))$.
By~(\ref{ineq:g4}), the heat operator $e^{s\Delta}$ may decrease values of
$g_4(x)$ by at most $\frac{2}{l^2}s$, i.e.
\begin{eqnarray}
        e^{s\Delta} g_4 (x) - g_4 (x) \geq -\frac{2}{l^2}s.
	\label{eq:bound:heat2}
\end{eqnarray}
Therefore for $x\in [-L-\frac{l}{200},L+\frac{l}{200}]$, we have
\begin{eqnarray*}
	e^{s\Delta}g_4(x)-b_0 f_0(x)
	= (e^{s\Delta}g_4(x)-g_4(x))+(g_4(x)-b_0 f_0(x))
	\geq -\frac{2}{l^2}s + Ks,
\end{eqnarray*}
where we apply~(\ref{eq:bound:heat2}) to $e^{s\Delta}g_4(x)-g_4(x)$
and~(\ref{eq:g4bigger}) to $g_4(x)-b_0 f_0(x)$.
Thus for $x\in [-L-\frac{l}{200},L+\frac{l}{200}]$, we have
\begin{eqnarray*}
	e^{s\Delta}g_4(x)-b_0 f_0(x) > \frac{1}{l^2}s
\end{eqnarray*}
since $K$ is chosen to be larger than $3/l^2$ in~(\ref{eq:pickK}).

The estimates in~(\ref{eq:pde:case2a}) and~(\ref{eq:pde:case2b}), together
with the fact that $g_3$ and $g_4$ are lower bounds of
$\pi_v(\Fcal^s_\eta (a_0 f_0,b_0 f_0))(x)$, imply
that there is a positive constant $m$, such that for $x \in (-L-l-s,L+l+s)$,
\begin{eqnarray*}
        \pi_v (e^{s\Delta} \Fcal^s_\eta (a_0 f_0,b_0 f_0)) (x)
        > b_0 f_0(x) + s \geq b_0 (f_0(x) + s).
\end{eqnarray*}
As in Lemma~\ref{lem:case1pde}, we apply Lemma~\ref{lem:space} to obtain
the estimate
\begin{eqnarray}
        \pi_v(e^{s\Delta}\Fcal^s_\eta (a_0 f_0,b_0 f_0))(x) > b_0 f_2(x),
	\label{eq:pde:case2b}
\end{eqnarray}
where $f_2$ is defined in~(\ref{eq:lemma2a}).

For the $u$-coordinate of $e^{s\Delta}\Fcal^s_\eta (a_0 f_0,b_0 f_0)$,
we can obtain estimates~(\ref{eq:pde:case2a}) to~(\ref{eq:pde:case2b})
if we replace $b_0$ on the right hand side of each inequality by $a_0$.
So we conclude that for all $x\in (-L-l-s,L+l+s)$,
\begin{eqnarray*}
        e^{s\Delta} \Fcal^s_\eta (a_0 f_0,b_0 f_0) (x)
	> (a_0 f_2(x),b_0 f_2(x)),
\end{eqnarray*}
as required. In particular, $(\tilde a_s,\tilde b_s)$ in the statement of
Proposition~\ref{prp:Dd} should be $((1+\delta_2 s) a_0,(1+\delta_2 s) b_0)$.
\qed

\vspace{.3cm}
\noindent {\bf Proof of Proposition~\ref{prp:Dd}}. \hspace{2mm}
The proposition follows from Lemmas~\ref{lem:case1pde} and~\ref{lem:case2pde},
and the discussion below~(\ref{def:R0}) on page~\pageref{page:mono}
regarding the sufficiency of restricting the region for $(a_0,b_0)$ to $\R_0$.
\qed

\section{Upper Bounds: Existence of $d_2$ and $D_2$ in Condition ($*$)}
\label{sec:ub}
We establish the following proposition, which,
together with Corollary~\ref{cor:Dd}, verifies condition ($*$) on
page~\pageref{page:star}. As in Corollary~\ref{cor:Dd},
the ``$L$'' in condition ($*$) is picked to be $L-l$.
\begin{PRP}
If $c$ is sufficiently large, then there exist constants $d_2 < D_2 < 1$
and $T$ such that if $v_0 (x) < D_2$ for $x\in [-L+l,L-l]$
then $v_t (x) < d_2$ for $x\in [-3L,3L]$, where $(u_t,v_t)$ solves
the PDE~(\ref{eq:pde:real}).
\end{PRP}

\proof
Because of the monotonicity of the PDE~(\ref{eq:pde:real}), it suffices to
pick a uniform initial condition
\begin{eqnarray*}
	u_0 & \equiv & \mbox{some } \bar{u}, \\
	v_0 & \equiv & D_2,
\end{eqnarray*}
and show that at time $T$,
\begin{eqnarray*}
	u_T & \equiv & \mbox{some } \tilde{u}, \\
	v_T & < & d_2.
\end{eqnarray*}
Therefore we need only concern ourselves with the ODE~(\ref{eq:ode}).
%Since $\gamma_2 = \left\{(u,v): u\in [0,1], v=u-\frac{2}{c} \right\}$ and the
%line $u=1$ intersect at the point $(1,1-\frac{2}{c})$.
We can bound $\eta_2(u,v)$ defined in~(\ref{def:eta})
for any $v>1-\frac{1}{c}$ as follows:
\begin{eqnarray*}
	\eta_2(u,v) = (c(u-v)-2)v
	\leq \left(c\left(1-\left(1-\frac{1}{c}\right)\right)-2\right)v = -v
	< -\left(1-\frac{1}{c}\right) < 0
\end{eqnarray*}
if $c>1$.
Thus for any two numbers $D_2$ and $d_2$ that satisfy
$1>D_2>d_2>1-\frac{1}{c}$, there exists $T$, such that
if $v_0 \equiv D_2$, then $v_T < d_2$.
\qed

\part{Stationary Distributions of A Model of Sympatric Speciation}
\chapter{A Model on Sympatric Speciation}
\section{Introduction}
Understanding Speciation is one of the great problems in the field of
evolution. According to Mayr [Mayr 1963], speciation means the splitting
of a single species into several, that is, the multiplication of species.
It is believed that many species originated through geographically isolated
populations of the same ancestral species [Dieckmann and Doebeli 1999].
This phenomenon is relatively easy to understand. In contrast, sympatric
speciation, in which new species arise without geographical isolation,
is theoretically much more difficult.

\subsection{The Dieckmann-Doebeli Model}
Dieckmann and Doebeli [Dieckmann and Doebeli 1999] proposed a general
model for sympatric speciation, for both asexual and sexual populations.
We will describe their model for the asexual population first. Each individual
in the population is assumed to have a quantitative character (phenotype)
$x\in \Rbold$ determining how effectively this individual can
make use of resources in the surrounding environment.
A typical example is the beak size of a certain bird species,
which determines the size of seeds that can be consumed by an individual bird.
The function $K: \Rbold \rightarrow \Rbold^+$ (carrying
capacity) is associated with the surrounding environment, where
$K_x$ denotes the number of individuals of phenotype $x$ that can be
supported by the environment. For example, since birds with small beak size
(say $x_1$) are more adapted to eating small seeds than birds with large
beak size (say $x_2$, $x_2>x_1$), $K_{x_1}$ will be larger than $K_{x_2}$
if the surrounding environment produces more small seeds than large seeds.
In the Dieckmann-Doebeli model, $K_x$ is taken to be
$c \exp(-\frac{(x-\hat x)^2}{2\sigma_K^2})$. Moreover, every pair
of individuals compete at an intensity determined by the phenotypical
distance of these two individuals. More specifically, an individual
of phenotype $x_1$ competes with an individual of phenotype $x_2$ at
intensity $C_{x_1-x_2}$, where
$C_x = \exp(-\frac{x^2}{2\sigma_C^2})$. Therefore each
individual in the population interacts with the environment via the
carrying capacity $K$, and interacts with the population via the competition
kernel $C$.

Let $N_x(t)$ denote the number of individuals with phenotype $x$ at time $t$.
At any time, an individual of phenotype $x$ gives birth at a constant rate,
and dies at a rate proportional to $\frac{(C*N_\cdot(t))_x}{K_x}$,
i.e. inversely proportional to the $x$-carrying capacity,
but proportional to the intensity of competition exerted by the
population on phenotype $x$, the numerator
$(C*N_\cdot(t))_x = \sum_y C_{x-y} N_y(t)$ being
how much competition (from every individual in the population)
individuals with phenotype $x$ suffer.
In addition, every time an individual gives birth,
there is a small probability that a mutation occurs and the phenotype of the
offspring is different from that of the parent; in this case,
the phenotypical distance between the offspring and the parent is
then random and assumed to have a Gaussian distribution.

Since the number of individuals of a certain phenotype
increases via the birth mechanism at a linear rate, but
decreases via the death mechanism at a quadratic rate,
extinction of  all phenotypes will occur in finite time
with probability one, i.e. $N \equiv 0$ eventually. For large initial
populations, however, extinction will happen far enough into the future that
interesting behaviour does arise before the population becomes extinct.

Monte-Carlo simulations, shown in figure~\ref{fig:model:dd}, give a fairly
good idea of the behaviour of the Dieckmann-Doebeli model for asexual
populations. If the initial population is monomorphic
($t=1$ in figure~\ref{fig:model:dd}), i.e. concentrated
near a certain phenotype $x_0$
($\frac{N_\cdot(0)}{\sum_x N_x(0)} \approx \delta_{x_0}$),
then the entire population first moves ($t=30,100,200$ in
figure~\ref{fig:model:dd}) toward $\hat x$, the
phenotype with maximum carrying capacity.
If $\sigma_C > \sigma_K$ (this
includes the case $\sigma_C = \infty$, i.e. equal competition
between all phenotypes), then the population stabilizes near
phenotype $\hat x$. But if $\sigma_C < \sigma_K$, then the monomorphic
population concentrated at phenotype $\hat x$ splits into two groups, one
group concentrating on a phenotype $< \hat x$,
while the other concentrating on a phenotype $> \hat x$
($t=330,370,400,500$ in figure~\ref{fig:model:dd}). In the latter
case, one can say that one species has evolved into two distinct species.

{
\psfrag{iter1}{$t=1$}
\psfrag{iter30}{$t=30$}
\psfrag{iter100}{$t=100$}
\psfrag{iter200}{$t=200$}
\psfrag{iter300}{$t=300$}
\psfrag{iter330}{$t=330$}
\psfrag{iter370}{$t=370$}
\psfrag{iter400}{$t=400$}
\psfrag{iter506}{$t=500$}
\psfrag{iter600}{$t=600$}
\begin{figure}[h!!]
\centering
\includegraphics[height = 2in, width = 5.5in]{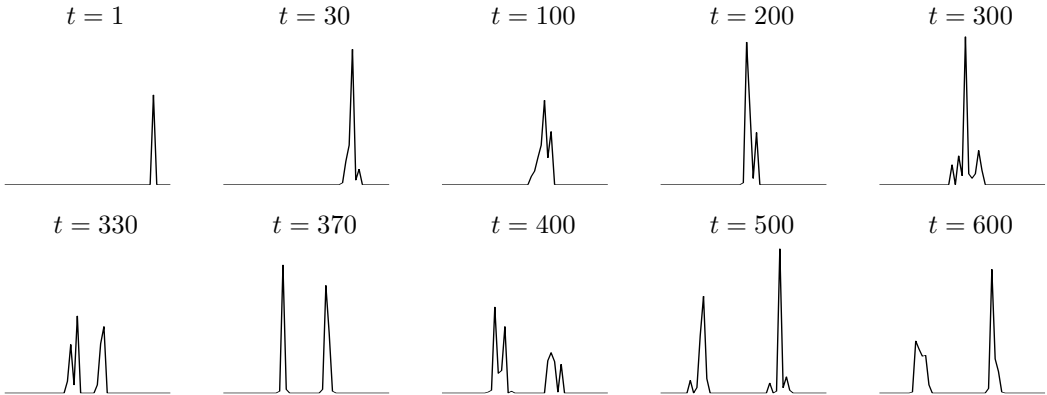}
\caption{Simulation of the Dieckmann-Doebeli model with
$E=[-50,50]\cap \Zbold$, $\sigma_K = \sqrt{1000}$,
	and $\sigma_C=\sqrt{600}$.}
\label{fig:model:dd}
\end{figure}
}

We now give a qualitative description of the Dieckmann-Doebeli
model for sexual populations. Each individual in a sexual population
is assigned three diploid genotypes
with (say) five diallelic loci each. The first set of loci determines
the ecological character $x$ (i.e. phenotype in the asexual model). The
second set of loci determines the marker trait, which is ecologically
neutral, i.e. individuals with different marker traits but the same
ecological character have exactly the same birth and death rates.
The third set of loci determines mating probabilities $m$; if $m>0$,
then such individuals prefer to mate with individuals of similar phenotypes;
if $m=0$, then such individuals have no preference;
and if $m<0$, then such individuals prefer to mate with individuals
of a distant phenotype. In addition, $|m|$ determines the strength of this
preference.

The birth rates and the death
rates are calculated the same way as in the asexual model; in particular,
only information from the first set of loci is used to calculate these
rates, as this is the only genotype that determines the
phenotype of the individual.
Dieckman and Doebeli considered two cases in their sexual model:
1. mating depends on the ecological character; and 2. mating depends on
the ecologically neutral marker trait. For example, in the second case,
individuals with $m>0$ prefers to mate with individuals of similar
marker traits.

Monte-Carlo simulations show that case 1 of the sexual model exhibits very
similar behaviour to the asexual model, i.e.
speciation if $\sigma_C < \sigma_K$ and
no speciation if $\sigma_C > \sigma_K$. A caveat: if
$\sigma_C < \sigma_K$, then only individuals, who prefer to mate with
individuals of similar phenotypes, survive after the population splits into
two groups. Hence in the end, there are two distinct groups of
individuals who refuse to mate with individuals from the other group.
For case 2 of the sexual model, Monte-Carlo simulations indicate that
$\sigma_C < \sigma_K$ is not enough for speciation to occur. In this case,
$\sigma_C < c \sigma_K$ is needed, where $c<1$ is a constant.

As the sexual model exhibits similar behaviour to the asexual model,
we will concentrate on the analysis of the simpler asexual model.
It is our hope that understanding the asexual model will give some
insights in explaining the behaviour of the sexual model as well.

\subsection{A conditioned Dieckmann-Doebeli model}
Although the Dieckmann-Doebeli model for asexual populations is considerably
less complicated than their model for sexual populations, it still seems too
complicated for rigorous analysis. Thus we will attempt to
simplify the model while preserving its key ingredients. Henceforth we
refer to the Dieckmann-Doebeli model for asexual populations simply as
the Dieckmann-Doebeli model.

Before we describe our simplified Dieckmann-Doebeli model, we first
introduce the concepts of \emph{fitness} and \emph{selection}.
Selection occurs when individuals of different genotypes leave different
numbers of offspring because their probabilities of surviving
to reproductive age are different [B\"urger 2000]. If we define fitness to be
a measure of how likely a particular individual produces offspring
that will survive to reproductive age, then individuals with higher fitness
should have higher probability of being selected for reproduction.
Along these lines, it is natural to define fitness of a phenotype
as the difference between the birth rate and the death rate of individuals
of this phenotype. It is also natural to require the fitness function to be
bounded between $0$ and $1$.

The key feature of the Dieckmann-Doebeli model is that each individual
has a fitness that depends on both the carrying capacity associated with
its phenotype and the configuration of the entire population.
More specifically, the fitness of a phenotype $x$ is an increasing function
of $K_x$, the carrying capacity, but a decreasing function
of $(C*N)_x$, the competition it suffers.
Here $N_x$ is the number of individuals of phenotype $x$.

In the Dieckmann-Doebeli model, the number of individuals can fluctuate
with time. As mentioned before, since the birth rate is linear but the
death rate is quadratic, extinction will
occur in finite time with probability one, which makes it somewhat
meaningless to analyze the equilibrium behaviour of the system.
We make the assumption that the number of individuals $N$
is fixed over time, reflecting a constant carrying capacity of the
overall population. The mechanism by which we achieve this is to require
that death of an individual and birth of its single offspring occur at
the same time, called \emph{replacement sampling} in
Moran particle models [Dawson 1993, Chapter 2.5].
This way, the number of individuals remains constant, and
analyzing the behaviour of the population is then equivalent to analyzing
the empirical distribution
\bestar
        \pi^N = \frac{1}{N} \sum_{n=1}^N \delta_{x_n},
\eestar
where $x_n$, $n=1\ldots N$, denotes the phenotype of the $n^\mathrm{th}$
individual in a population of size $N$.

Before we describe our simplified Dieckmann-Doebeli model, we say a few
words about our terminologies and notations:
we refer to individuals in a population
as ``particles'', and sometimes refer to a phenotype as a ``site''.
We consider multiple models, both discrete-time and continuous-time;
for discrete-time models, we use $V$ to denote the fitness function;
but for continuous-time models, we use $m$ instead.
Our simplified discrete-time Dieckmann-Doebeli model is as follows:
\begin{enumerate}
\item $E=[-L,L]\cap \Zbold$ is the phenotype space, and
        $\pi, \pi^N \in \Pcal(E)$ is a probability measure on $E$;
\item $K: E \rightarrow [0,1]$ is the carrying capacity, and
        $C: \Zbold \rightarrow \Rbold^+$ is the competition kernel;
\item $V_x(\pi)$ is the fitness
        of phenotype $x$ in a population with empirical distribution $\pi$
        (sometimes we notationally suppress the dependence on $\pi$); we define
        two possible fitness functions below;
\item $A$ is a Markov transition matrix associated with mutation,
        with $A(y,x)$ denoting the probability of a particle of
        phenotype $y$ mutating to a particle of phenotype $x$;
\item At every time step $t\in \Zbold^+$, the entire population is replaced by
	a new population of $N$ particles, each particle
	chosen independently, according to the distribution
	$p_\cdot(t,\pi^N)$:
\be
        p_x(t,\pi^N) = \sum_y A(y,x) \frac{\pi^N_y(t) V_y(\pi^N(t))}{
                \sum_z \pi^N_z(t) V_z(\pi^N(t))}
        \label{eq:badmodel1}
\ee
\end{enumerate}
In~(\ref{eq:badmodel1}), the denominator $\sum_z \pi^N_z(t) V_z(\pi^N(t))$
is simply the normalization factor
such that $\sum_x p_x(t,\pi^N) = 1$. In words, at every time step,
the entire population dies and is replaced by a new population, each
individual $x$ choosing an individual $x'$ from the original population as
its parent with a probability proportional to its fitness $V_{x'}$, after
which the new individual $x$ undergoes mutation according to $A$.

We consider two fitness functions:
\be
        V^{(1)}_x(\pi) & = & 0 \vee \left(
                1 - \frac{\sum_z C_{x-z} \pi_z}{K_x}\right), \nnb \\
        V^{(2)}_x(\pi) & = & \frac{K_x}{\sum_z C_{x-z} \pi_z}.
	\label{eq:badmodel:V2}
\ee
Each of the two fitness function defined above is an increasing function
of $K_x$ and a decreasing function of $(C*\pi)_x$. $V^{(1)}$ resembles
more closely the original Dieckmann-Doebeli model, but it has
the disadvantage of being in a more complicated form than $V^{(2)}$ and
it is also not differentiable.

By Theorem 1 in [Del Moral 1998], which we state below,
$\{\pi^N_t,t\in [0,T]\} \Rightarrow \{\pi_t,t\in [0,T]\}$
as $N\rightarrow \infty$, where $\Rightarrow$ denotes weak convergence
and $\pi_t$ evolves according to the following
\emph{deterministic} dynamical system:
\be
        \pi_x(t+1) = \sum_y A(y,x) \frac{\pi_y(t) V_y(\pi(t))}{
                \sum_z \pi_z(t) V_z(\pi(t))}. \label{eq:badmodel}
\ee
\begin{THM}
        Suppose $E$ is compact and $M$ is a Feller-Markov transition kernel
        on $\Pcal(E)$, i.e. $M: \Pcal(\Pcal(E)) \rightarrow \Pcal(\Pcal(E))$.
	Let $M^{(N)} = M C_N$, where $C_N$ is a Markov
	transition kernel on $\Pcal(E)$ given by
\begin{eqnarray*}
	C_N F(\pi) = \int_{E^{N}} F\left( \frac{1}{N} \sum_{j=1}^N \delta_{x_j}
		\right) \pi(dx_1) \ldots \pi(dx_N),
\end{eqnarray*}
	i.e. a probability measure $\pi$ replaced with an empirical measure
	$\pi^N$ formed by $N$ particles chosen independently according
	to $\pi$.  Then
\begin{eqnarray*}
        \left(M^{(N)}\right)^n \rightarrow M^n
\end{eqnarray*}
        as $N \rightarrow \infty$.
\end{THM}
This theorem is easy to understand: take $n=1$, then it says that the mapping
$M C_N$ converges to $M$, i.e. changing the input measure of $M$ by an
empirical measure of $N$ particles makes almost no difference if $N$ is large.

Analyzing the dynamical system~(\ref{eq:badmodel}) is not easy,
partly because it is of
a complicated form that is nonlinear in $\pi$, and we cannot
find any Lyapunov function that associates with~(\ref{eq:badmodel}).
A continuously differentiable function $V: U \rightarrow \Rbold$
is called a Lyapunov function
if $V$ is nondecreasing (or nonincreasing) along orbits. For a
discrete-time dynamical system such as~(\ref{eq:badmodel}), this means that
$V(\pi(t+1))-V(\pi(t)) \geq 0$ (or $\leq 0$) for all $t\geq 0$.
For a continuous-time dynamical system, it means that
$\partial_t V(\pi(t)) \geq 0$ (or $\leq 0$).
Simulations of~(\ref{eq:badmodel}), however, seem to display some interesting
behaviour, which we will describe after carrying out some non-rigorous
analysis of~(\ref{eq:badmodel}).

Without mutation, any site $x$ with $\pi_x = 0$ at
any time $\tau$ will stay $0$ for all $t \geq \tau$. Mutation enables
individuals of phenotype $x$ to be born in future generations
even if there are no individuals of phenotype $x$ in the present generation.
But if we start with a polymorphic initial measure, i.e. $\pi_x(0) \neq 0$
for all $x$, then adding small mutation to the system should not
cause significant changes in the behaviour of~(\ref{eq:badmodel}).
Therefore we assume that $A=I$ and $\pi(0)$ is polymorphic.
In this case,~(\ref{eq:badmodel}) can be simplified to
\bestar
        \pi_x(t+1) = \frac{\pi_x(t) V_x(\pi(t))}{
                \sum_z \pi_z(t) V_z(\pi(t))}.
\eestar
Thus if $A=I$, then $\hat\pi$ is a stationary
distribution of~(\ref{eq:badmodel}) if and only if
\be
        \hat\pi_x = \frac{1}{c} \hat\pi_x V_x(\hat\pi) \label{cond:badmodel}
\ee
for some constant $c$. Condition~(\ref{cond:badmodel}) is equivalent to
\be
        V_x(\hat\pi) = c \ \mbox{for all $x$ where $\hat\pi(x) \neq 0$},
        \label{cond:Vconst}
\ee
Let $K$ and $C$ be in the form considered by Dieckmann and Doebeli,
i.e. $K_x=\exp(-x^2/2\sigma_K^2)$ and $C_x=\exp(-x^2/2\sigma_C^2)$.
If $V=V^{(2)}$, then condition~(\ref{cond:Vconst}) means that
\bestar
        K_x = c (C*\hat\pi)(x) \ \mbox{for all $x$ where $\hat\pi(x) \neq 0$},
\eestar
which seems to indicate that if $\sigma_C < \sigma_K$, then $\hat\pi$
should be close to $\Ncal(0,\sigma_K^2-\sigma_C^2)$. On the other hand,
if $V=V^{(1)}$, then $\hat\pi$ is a stationary distribution if
$1 - \frac{\sum_z C_{x-z} \hat\pi_z}{K_x}$ is a strictly positive constant.
Notice that if $K$ and $C$ are both Gaussian-shaped with $K_0=C_0=1$
then $\hat\pi = \Ncal(0,\sigma_K^2-\sigma_C^2)$ makes
$1 - \frac{\sum_z C_{x-z} \hat\pi_z}{K_x}$ constant; furthermore,
this constant is strictly positive since $(C*\hat\pi)(0) < K_0 = 1$
if $\sigma_C < \sigma_K$.

Therefore for both $V^{(1)}$ and $V^{(2)}$, the
dynamical system~(\ref{eq:badmodel}) should have Gaussian-shaped
stationary distributions if $\sigma_C < \sigma_K$. In simulations carried out
by Dieckmann and Doebeli [Dieckmann and Doebeli 1999],
however, $\sigma_C < \sigma_K$ is the case that
leads to speciation, i.e. the stationary distribution supposedly has
two sharp well-separated peaks, which contradicts the analysis carried out
in the previous paragraph. Simulations of~(\ref{eq:badmodel}) with
$V=V^{(1)}$, shown in Figure~\ref{fig:badmodel1},
reveal that if $\pi(0) \approx \delta_0$, initially the population does
split into two groups and begins to move apart, but as $t\rightarrow \infty$,
the empirical measure converges to a Gaussian-shaped hump. This suggests
the possibility that in the original Dieckmann-Doebeli model, conditioning
on the population surviving long enough for convergence to equilibrium to
occur (recall that in the original Dieckman-Doebeli model, extinction
occurs in finite time), speciation is also a transitory phenomenon,
rather than an equilibrium phenomenon. Simulations of~(\ref{eq:badmodel})
with $V=V^{(2)}$, shown in Figure~\ref{fig:badmodel2}, does not even display
transitory speciation behaviour.
Instead, the initial spike at $0$ simply widens to a Gaussian hump centred
at $0$. Hence the particular form of the dependence on $K$ and $C*\pi$ seems
to affect whether or not speciation occurs.

{
\psfrag{iter0}{$t=0$}
\psfrag{iter10}{$t=10$}
\psfrag{iter20}{$t=20$}
\psfrag{iter30}{$t=30$}
\psfrag{iter40}{$t=40$}
\psfrag{iter80}{$t=80$}
\psfrag{iter200}{$t=200$}
\psfrag{iter400}{$t=400$}
\psfrag{iter800}{$t=800$}
\psfrag{iter1200}{$t=1200$}
\psfrag{iter1600}{$t=1600$}
\psfrag{iter2000}{$t=2000$}
\psfrag{iter3000}{$t=3000$}

\begin{figure}[h!!]
\centering
\includegraphics[height = 2in, width = 5.5in]{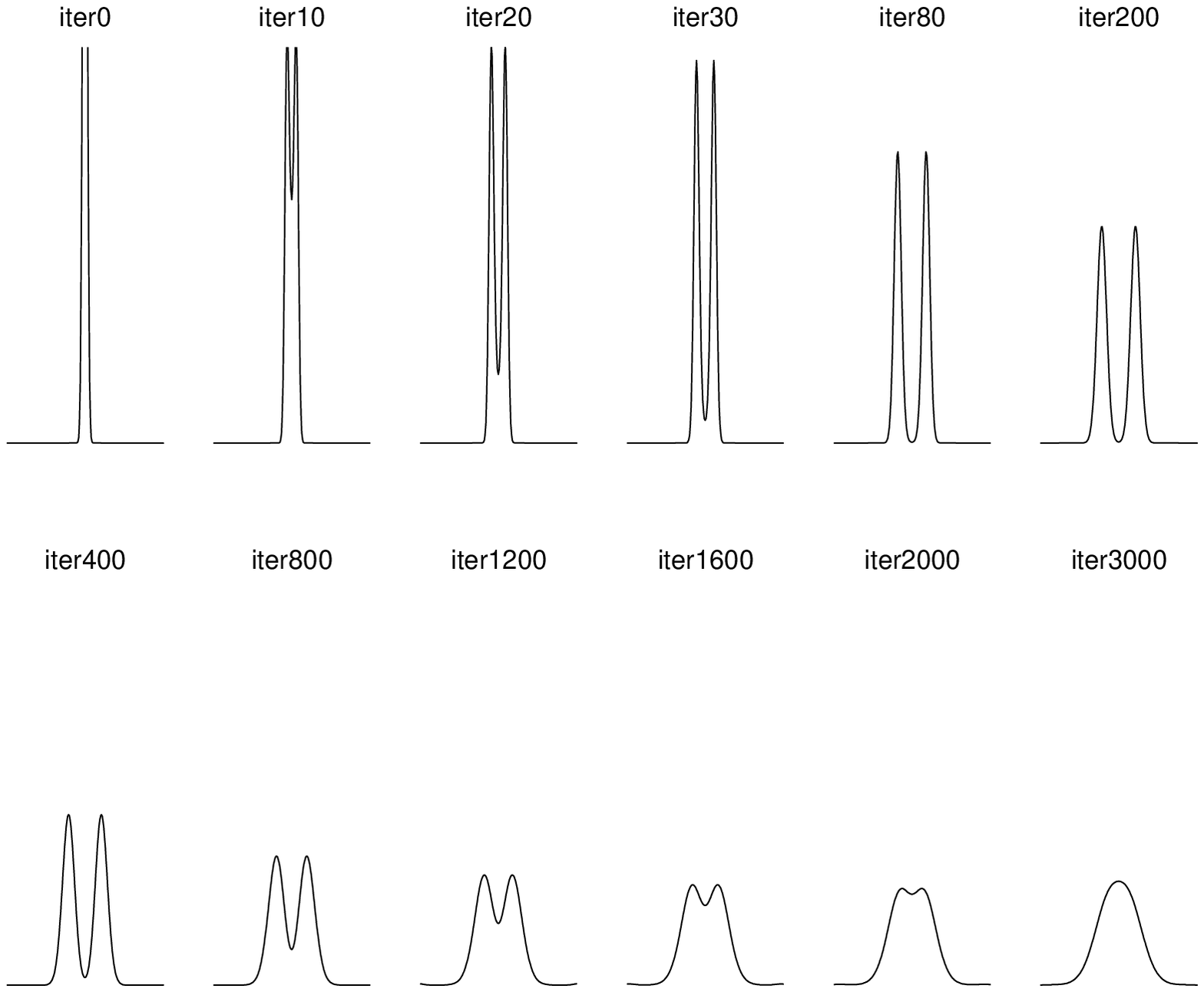}
\caption{Simulation of~(\ref{eq:badmodel}) with
$E=[-149,149]\cap \Zbold$, $\sigma_K = 60$, $\sigma_C=55$, and $V=V^{(1)}$.}
\label{fig:badmodel1}
\end{figure}

\begin{figure}[h!!]
\centering
\includegraphics[height = 1in, width = 5.5in]{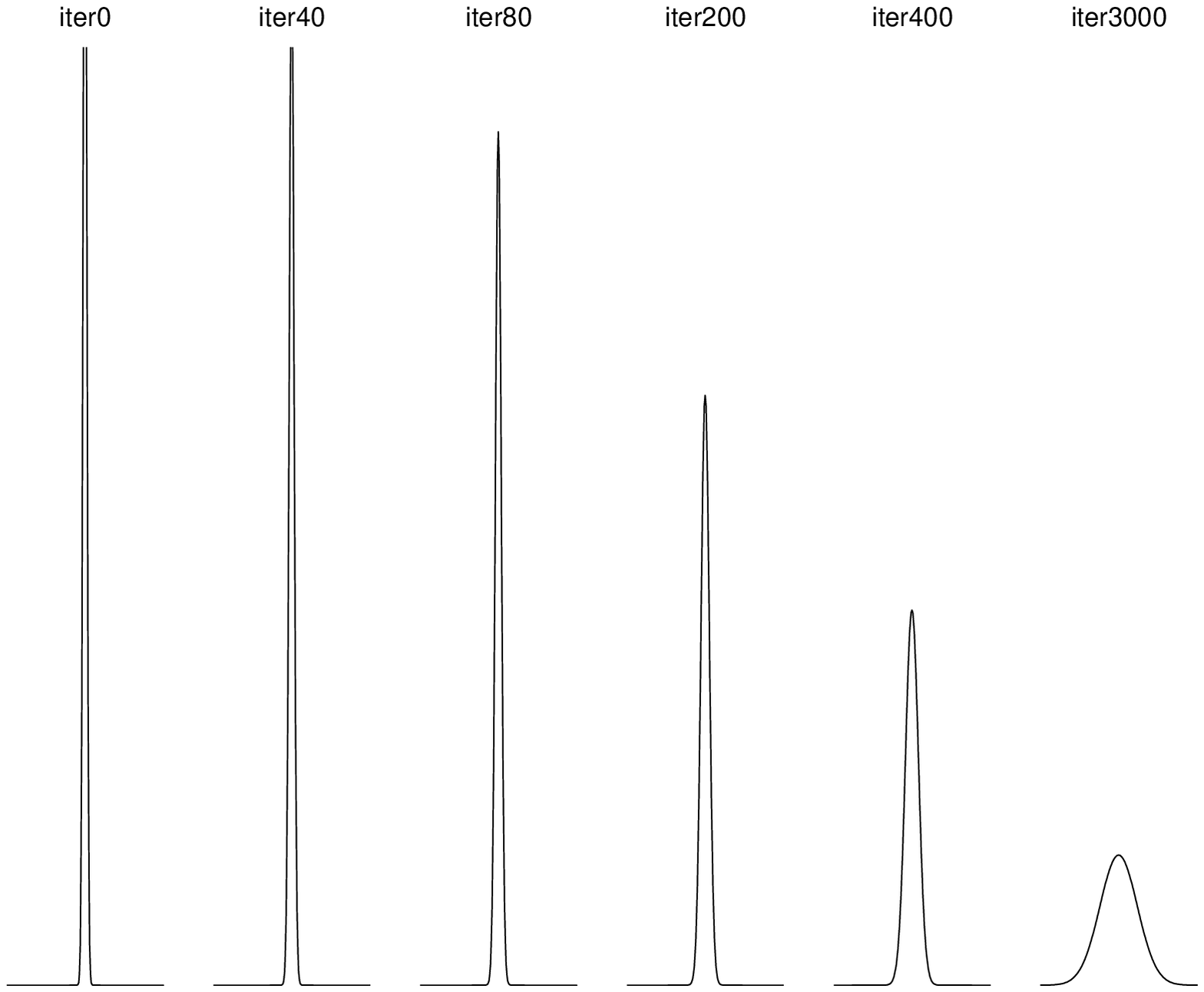}
\caption{Simulation of~(\ref{eq:badmodel}) with
$E=[-149,149]\cap \Zbold$, $\sigma_K = 60$, $\sigma_C=55$, and $V=V^{(2)}$.}
\label{fig:badmodel2}
\end{figure}
}

From the simulations and non-rigorous analysis above, it seems that
the dynamical system in~(\ref{eq:badmodel}) does not have a
bimodal stationary distribution if both $K$ and $C$ are taken to be
Gaussian-shaped. If $K$ and $C$ are taken to be rectangular (i.e.
$K_x = 1_{ \{|x|\leq L \} }$ and $C_x = 1_{ \{|x| \leq M\} }$
for some integers $L$ and $M$), however, results from
Appendix~\ref{ch:appen} shows that there exist bimodal stationary distributions.
More specifically, Theorem~\ref{thm:append} says that if
$\nu^n$ is a convergent sequence of \emph{symmetric} stationary distributions
for the conditioned Dieckmann-Doebeli model with mutation parameter $\mu^n$
with $\mu^n \rightarrow 0$,
then $\nu^n_{[-l,l]} \rightarrow 0$, where $l=M-L+1$; in words, the
mass in the middle gets very small as the mutation parameter approaches zero.

\subsection{A Moran Model with Competitive Selection}
As discussed earlier, the dynamical system~(\ref{eq:badmodel}) cannot
be easily associated with a Lyapunov function, which makes analyzing its
behaviour difficult. Keeping in mind that the essential ingredient of
the original Dieckmann-Doebeli model is that the fitness function is
an increasing function of $K_x$ and a decreasing function of $(C*\pi)_x$,
we define the fitness $m_x(\pi)$ to have the following form:
\begin{eqnarray}
        m_x(\pi) = K_x \sum_z B_{x-z} K_z \pi_z, \label{def:m1}
\end{eqnarray}
where the ``cooperation'' kernel $B$ can be taken to be $1-C$. We assume $B$
is symmetric. In the original Dieckmann-Doebeli model,
pairs of individuals with small phenotypical distance \emph{compete}
at a higher intensity than
pairs of individuals with large phenotypical distance;
in our model, pairs of individuals with small phenotypical distance
\emph{cooperate} at a lower intensity than
pairs of individuals with large phenotypical distance.
To make our formulation cleaner, we also adopt a continuous-time model.
The advantage of adopting $m_x$ in~(\ref{def:m1}) as fitness
and using a continuous time model is that the mean fitness of the
population
\bestar
        \overline m_\pi = \sum_x \pi_x m_x
        = \sum_{x,z} \pi_x K_x B_{x-z} K_z \pi_z
\eestar
is a Lyapunov function [B\"urger 2000] for the dynamical system
\be
        \partial_t \pi_x = \pi_x (m_x - \overline m_\pi). \label{eq:nomut}
\ee
This assertion can be verified by the following calculation:
\be
        \partial_t \overline m_\pi & = & 2 \sum_x m_x \partial_t \pi_x
		\nnb \\
	& = & 2 \sum_x m_x \pi_x (m_x - \overline m_\pi) \nnb \\
        & = & 2 \sum_x m_x \pi_x (m_x - \overline m_\pi)
         - 2 \sum_x \overline m_\pi \pi_x (m_x - \overline m_\pi) \nnb \\
        & = & 2 \sum_x \pi_x (m_x - \overline m_\pi)^2, \label{eq:nomut:ver}
\ee
where in the second line we use the fact
\bestar
	\sum_x \overline m_\pi \pi_x (m_x - \overline m_\pi)
	= \overline m_\pi^2 - \overline m_\pi^2 \sum_x \pi_x
	= \overline m_\pi^2 - \overline m_\pi^2 = 0. \nnb
\eestar
Since $\partial_t \overline m_\pi \geq 0$ for any $\pi$, $\overline m_\pi$
is a Lyapunov function for the dynamical system~(\ref{eq:nomut}),
and in particular, the mean fitness $\overline m_\pi$
increases at a rate proportional to the variance of the fitness.
We call~(\ref{eq:nomut}) the selection-only equation, as it does not
have a part that corresponds to mutation.
In Chapter~\ref{ch:st-model}, we will derive~(\ref{eq:nomut}) as the
deterministic limit of particle systems as the number of particles
tends to infinity.

\section{The Particle Model}
\label{ch:model}
We introduce two particle models, one with ``strong selection'' that
yields a deterministic limit, and another with ``weak selection''
that yields a stochastic limit. We work on space
$E = [-L,L] \cap \Zbold$. Let
\begin{eqnarray*}
	\Delta = \{(p_{-L},\ldots,p_0,\ldots,p_L):
                p_i > 0 \ \forall i \mbox{ and } \sum_{i=-L}^L p_i = 1 \}
\end{eqnarray*}
be the space of probability measures on $E$, i.e. $\Delta = \Pcal(E)$.
Members of $\Delta$ are usually denoted by $\pi$, $\hat\pi$, $\pi^N$, etc.
We endow $\Delta$ with the following metric:
\begin{eqnarray*}
	d(\hat\pi,\tilde\pi) = \max_x |\hat\pi(x)-\tilde\pi(x)|.
\end{eqnarray*}
Let $K: E \rightarrow [0,1]$ be the carrying capacity function,
and $B: \Zbold \rightarrow [0,1]$ be the cooperation kernel, with
$B_z=0$ meaning that sites separated by phenotypical
distance $z$ do not cooperate at all (i.e. compete at full intensity),
and $B_z=1$ meaning that they cooperate at full intensity
(i.e. do not compete at all). We assume $B$ to be symmetric. The fitness
of site $x$ in a population with distribution $\pi$ is defined as
\begin{eqnarray*}
	m_x(\pi) = K_x \sum_z B_{x-z} K_z \pi_z.
\end{eqnarray*}
If one abuses notation by writing $K$ as a diagonal matrix, $B$ as a matrix,
and $\pi$ as a vector,
then the vector formed by $m_\cdot(\pi)$ can be written as $KBK\pi$.
The mean fitness of a population with distribution $\pi$ is defined as
\begin{eqnarray*}
	\overline m_\pi = \sum_{x,z} \pi_x K_x B_{x-z} K_z \pi_z.
\end{eqnarray*}
If one abuses notation again, then $\overline m_\pi$ can be written
as a quadratic form $\pi^t K B K \pi$.

Throughout this work, we will use symmetric or house-of-cards mutation,
which means that the rate $\mu_{xy} = \mu_y$ at which phenotype $x$
mutates to phenotype $y$ depends on $y$ only.
This is a common assumption in population genetics [B\"urger 2000],
and it is precisely this assumption
that allows one to explicitly write down a Lyapunov function
for the selection-mutation equation (to be defined in~(\ref{eq:ODE})).
As a further simplification, we assume that $\mu_y = \mu$ is constant in $y$,
which makes the proofs a bit cleaner.
Let $X^N(t)=(X^N_1(t),\ldots,X^N_N(t))$, $t \in \Rbold^+$, be an $N$-particle
system, with $X^N_i (t) \in E$ for all $t$ and $i$.
Define the empirical measure
\begin{eqnarray*}
        \pi^N (t) = \frac{1}{N} \sum_{i=1}^N \delta_{X^N_i(t)}.
\end{eqnarray*}

\subsection{The Strong Selection Model}
\label{ch:st-model}
For our model with strong selection, the particle system
undergoes the following:
\begin{enumerate}
\item[$\bullet$] Selection: At rate $N\overline m_{\pi^N}$, a particle,
        say $X^N_i$,
        is chosen at random from the $N$ particles and killed; at
        the same time, a new particle is born at $x$ with probability
        $\frac{m_x(\pi^N) \pi^N_x}{\overline m_{\pi^N}}$. Since
        $\sum_x \frac{m_x(\pi^N) \pi^N_x}{\overline m_{\pi^N}} = 1$,
        $\frac{m_\cdot(\pi^N) \pi^N_\cdot}{\overline m_{\pi^N}}$
        is a probability distribution;

\item[$\bullet$] Mutation:
        At rate $N(2L+1)\mu$, a particle, say $X^N_i$, is chosen
        at random from the $N$ particles and killed; at the same time,
        a new particle is born at a site $y$ with probability
        $\frac{1}{2L+1}$.
\end{enumerate}
A particle at $x$ gets replaced via selection by a particle at $y$
at rate $N \pi^N_x m_y(\pi^N) \pi^N_y$, and gets replaced via mutation
by a particle at $y$ at a rate of $N \mu \pi^N_x$.
Let $\bar K = \sup_{x\in [-L,L]} K_x$, so that
$m_x(\pi) = K_x \sum_z B_{x-z} K_z \pi_z \leq \bar K \sum_z \bar K \pi_z
= \bar K^2$. Let $l(dx)$ denote the Lebesgue measure on $\Rbold^+$.
The process described above can be constructed using a
Poisson point process $\Lambda^N(dt,dx,dy,d\xi,de)$ on
\[ \Rbold^+ \times \{(x,y) \in E^2: x\neq y \} \times [0,1] \times \{1,2\}, \]
with intensity measure
\begin{eqnarray*}
	\lambda^N(A \times B \times C \times D)
		= l(A) (\# B) (\# C) k(D),
\end{eqnarray*}
where $\#$ denotes the counting measure, $D \subset [0,1] \times \{1,2\}$, and
$k = l \times (N \bar K^2 \delta_1 + N \mu \delta_2)$.
For all $x,y\in E^2$ with $x\neq y$, jumps of $\Lambda^N(dt,x,y,[0,1],\{1\})$
give possible times at which a particle at $y$ may be replaced by
a particle at $x$ by the selection mechanism, while
jumps of $\Lambda^N(dt,x,y,[0,1],\{2\})$
give possible times at which a particle at $y$ may be replaced by
a particle at $x$ by the mutation mechanism.
The strong selection model can be expressed in terms of the following
formula for $\pi^N_x(t)$:
\begin{eqnarray}
        \lefteqn{ \pi^N_x(t) = \pi^N_x(0) } \nnb \\
        & &	+ \frac{1}{N} \left[
                        \int_0^t \int 1\left(\xi \leq \frac{
                                \pi^N_y(s-) m_x(\pi^N(s-)) \pi^N_x(s-)}{
                                \bar K^2}\right) \Lambda^N(ds,x,dy,d\xi,1)
				\right. \nnb \\
        & & \ \ \ \ \	- \left. \int_0^t \int 1\left(\xi \leq \frac{
                                \pi^N_x(s-) m_y(\pi^N(s-)) \pi^N_y(s-)}{
                                \bar K^2}\right) \Lambda^N(ds,dy,x,d\xi,1)
				\right] \nnb \\
        & &     + \frac{1}{N} \left[
                        \int_0^t \int 1(\xi \leq \pi^N_y(s-))
				\Lambda^N(ds,x,dy,d\xi,2) \right. \nnb \\
	& & \ \ \ \ \	- \left. \int_0^t \int 1(\xi \leq \pi^N_x(s-))
				\Lambda^N(ds,dy,x,d\xi,2) \right] .
        \label{eq:stsel}
\end{eqnarray}
A solution to~(\ref{eq:stsel}) exists because the total jump rate
is finite for a fixed $N$.
The two integrals inside the first set of brackets corresponds to
selection, i.e. a particle at $x$ gets replaced by a particle at $y$
at rate $N \pi^N_y m_x(\pi^N) \pi^N_x$ due to selection. In particular,
$\Lambda^N(ds,x,dy,d\xi,1)$ in the first integral accounts for the killing
of a particle at $y$ and a new particle being born at $x$, and
$\Lambda^N(ds,dy,x,d\xi,1)$ in the second integral accounts for the killing
of a particle at $x$ and a new particle being born at $y$. The two integrals
inside the second set of brackets corresponds to mutation, i.e.
a particle at $x$ gets replaced by a particle at $y$
at rate $N \mu \pi^N_y$ due to mutation.

\begin{PRP}
\label{prp:stsel}
	As $N\rightarrow \infty$, the processes $\pi^N$ converge weakly
	to a deterministic process $\pi$
	that takes values in $\Pcal(E)$ and obeys the following system of
	ODE's:
\begin{eqnarray}
        \partial_t \pi_x = \pi_x ( m(x,\pi) - \overline m_{\pi} )
                + \mu ( 1 - (2L+1) \pi_x). \label{eq:ODE}
\end{eqnarray}
\end{PRP}
\proof First, we rewrite~(\ref{eq:stsel}) by decomposing $\Lambda^N$ into
a martingale term $\tilde \Lambda^N$ and a deterministic drift term:
\begin{eqnarray}
        \pi^N_x(t) & = & \pi^N_x(0) + M^N_x(t) + \sum_{y=-L}^L \int_0^t
			[\pi^N_y(s-) m_x(\pi^N(s-)) \pi^N_x(s-) \nnb \\
	& &	- \pi^N_x(s-) m_y(\pi^N(s-)) \pi^N_y(s-)] \ ds
		+ \mu \sum_{y=-L}^L \int_0^t [\pi_y(s-) - \pi_x(s-)] \ ds
		\nnb \\
	& = & \pi^N_x(0) + M^N_x(t)
        	+ \int_0^t \pi^N_x(s-) [m_x(\pi^N(s-))
                        - \overline m_{\pi^N(s-)}] \nnb \\
        & &	+ \mu [1 - (2L+1)\pi_x(s-)] \ ds
	\label{eq:stsel1}
\end{eqnarray}
where we define $\tilde\Lambda^N=\Lambda^N-\lambda^N$ to be the martingale
part of $\Lambda^N$ and
\begin{eqnarray*}
	M^N_x(t) & = & \frac{1}{N} \left[
                        \int_0^t \int 1\left(\xi \leq \frac{
                                \pi^N_y(s-) m_x(\pi^N(s-)) \pi^N_x(s-)}{
				\bar K^2}\right)\tilde\Lambda^N(ds,x,dy,d\xi,1)
                                \right. \nnb \\
        & &             - \left. \int_0^t \int 1\left(\xi \leq \frac{
                                \pi^N_x(s-) m_y(\pi^N(s-)) \pi^N_y(s-)}{
                                \bar K^2}\right)\tilde\Lambda^N(ds,dy,x,d\xi,1)
                                \right] \nnb \\
        & &     + \frac{1}{N} \left[
                        \int_0^t \int 1(\xi \leq \pi^N_y(s-))
                                \tilde\Lambda^N(ds,x,dy,d\xi,2) \right. \nnb \\
        & &	\left. - \int_0^t \int 1(\xi \leq \pi^N_x(s-))
                                \tilde\Lambda^N(ds,dy,x,d\xi,2) \right].
\end{eqnarray*}
We estimate the quadratic variation of the martingale term $M^N_x(t)$
and show that it converges to $0$ as $N \rightarrow \infty$:
\begin{eqnarray}
	\left< M^N_x \right>_t & = & \frac{1}{N^2} \sum_{y=-L}^L \left[
                        \int_0^t 1\left(\xi \leq \frac{
                                \pi^N_y(s-) m_x(\pi^N(s-)) \pi^N_x(s-)}{
                                \bar K^2}\right) N \bar K^2 \ ds \right.
				\nnb \\
        & &              + \left. \int_0^t 1\left(\xi \leq \frac{
                                \pi^N_x(s-) m_y(\pi^N(s-)) \pi^N_y(s-)}{
                                \bar K^2}\right) N \bar K^2 \ ds \right]
                                \nnb \\
        & &     + \frac{1}{N^2} \sum_{y=-L}^L \left[
                        \int_0^t 1(\xi \leq \pi^N_y(s-)) N \mu \ ds
                        + \int_0^t 1(\xi \leq \pi^N_x(s-)) N \mu \ ds
                                \right]. \nnb \\
	& \leq & \frac{2 \bar K^2}{N} \sum_{y=-L}^L \int_0^t \ ds
		+ \frac{2 \mu}{N} \sum_{y=-L}^L \int_0^t \ ds \nnb \\
	& \leq & \frac{2 (\bar K^2+\mu)(2L+1) t}{N},
\end{eqnarray}
which $\rightarrow 0$ as $N \rightarrow \infty$. Since the maximum
jump size in $M^N_x$ is $\frac{1}{N}$, by Burkholder's inequality
(see for example Theorem 21.1 in [Burkholder 1973]),
we have
\begin{eqnarray}
	E\left( \left(\sup_{t\leq T} M^N_x(t)\right)^2 \right)
	\leq C \left( E \left< M^N_x \right>_T + \frac{1}{N^2} \right)
	\rightarrow 0 \label{ineq:MNx}
\end{eqnarray}
as $N \rightarrow 0$. The other two terms in~(\ref{eq:stsel1}), i.e.
$\pi^N_x(0)$ and $\int_0^t \pi^N_x(s-) [m_x(\pi^N(s-))
- \overline m_{\pi^N(s-)}] + \mu [1 - (2L+1)\pi_x(s-)] \ ds$, are both
$C$-tight, $\pi^N_x(0)$ being constant and the integrand in the integral
bounded uniformly between constants. Therefore $\pi^N_x$ is $C$-tight
for each $x$.

If there exists a sequence $N_n$ such that $\pi^{N_n}$ converges to
$\pi$ weakly and $\pi$ is continuous, then by bounded convergence,
$\int_0^t \pi^{N_n}_x(s-) [m_x(\pi^{N_n}(s-))
- \overline m_{\pi^{N_n}(s-)}] + \mu [1 - (2L+1)\pi_x(s-)] \ ds$
converges to
\begin{eqnarray*}
	\int_0^t \pi_x(s) [m_x(\pi(s))
	- \overline m_{\pi(s)}] + \mu [1 - (2L+1)\pi_x(s)] \ ds,
\end{eqnarray*}
and since $\pi^N$ has representation~(\ref{eq:stsel1}), $\pi$ satisfies
the following deterministic integral equation:
\begin{eqnarray}
        \pi_x(t) = \pi_x(0) + \int_0^t \pi_x(s) [m_x(\pi(s))
                        - \overline m_{\pi(s)}]
                + \mu [1 - (2L+1)\pi_x(s)] \ ds. \label{eq:ODE:int}
\end{eqnarray}
A continuous $(\pi(t),0\leq t\leq \infty)$
solves the integral equation~(\ref{eq:ODE:int}) if and only
if it satisfies the ODE system~(\ref{eq:ODE}).
By well-known results from ODE theory
(e.g. Theorem 1.1.1 from [Wiggins 1988]), the solution to~(\ref{eq:ODE})
is unique because its right-hand-side is $C^\infty$ in $\pi$.
Therefore solution to~(\ref{eq:ODE:int}) is unique as well,
and the proof is complete. \qed

\subsection{The Weak Selection Model}
For our model with weak selection, the particle system with $N$
particles undergoes the following:
\begin{enumerate}
\item[$\bullet$] Selection: A particle at $x$ gets replaced by a particle
	at $y$ at rate $N \pi^N_x m_y(\pi^N) \pi^N_y$;

\item[$\bullet$] Mutation: A particle at $x$ gets replaced by a particle
        at $y$ at rate $N \mu \pi^N_x$;

\item[$\bullet$] Replacement sampling: A particle at $x$ gets replaced
	by a particle at $y$ at rate $\frac{N^2}{2} \pi^N_x \pi^N_y$.
\end{enumerate}
Just as in the strong selection model, this process can be constructed using
a Poisson point process $\Lambda^N(dt,dx,dy,d\xi,de)$ on 
\[ \Rbold^+ \times \{(x,y)\in E^2: x\neq y \} \times [0,1] \times \{1,2,3\}, \]
with intensity measure
\begin{eqnarray*}
        \lambda^N(A \times B \times C \times D)
                = l(A) (\# B) (\# C) k(D),
\end{eqnarray*}
where $l$ denotes the Lebesgue measure on $\Rbold^+$, $\#$ denotes the counting
measure, $D \subset [0,1] \times \{1,2,3\}$, and $k = l \times
(N \bar K^2 \delta_1 + N \mu \delta_2 + \frac{N^2}{2} \delta_3)$.
The weak selection model can then be expressed in terms of the following
formula for $\pi^N_x(t)$:
\begin{eqnarray}
        \lefteqn{ \pi^N_x(t) = \pi^N_x(0) } \nnb \\
        & &	+ \frac{1}{N} \left[
                        \int_0^t \int 1\left(\xi \leq \frac{
                                \pi^N_y(s-) m_x(\pi^N(s-)) \pi^N_x(s-)}{
                                \bar K^2}\right) \Lambda^N(ds,x,dy,d\xi,1)
                                \right. \nnb \\
        & & \ \ \ \ \	- \left. \int_0^t \int 1\left(\xi \leq \frac{
                                \pi^N_x(s-) m_y(\pi^N(s-)) \pi^N_y(s-)}{
                                \bar K^2}\right) \Lambda^N(ds,dy,x,d\xi,1)
                                \right] \nnb \\
        & &     + \frac{1}{N} \left[
                        \int_0^t \int 1(\xi \leq \pi^N_y(s-))
                                \Lambda^N(ds,x,dy,d\xi,2) \right. \nnb \\
        & & \ \ \ \ \	- \left. \int_0^t \int 1(\xi \leq \pi^N_x(s-))
                                \Lambda^N(ds,dy,x,d\xi,2) \right] \nnb \\
	& &	+ \frac{1}{N} \left[
			\int_0^t \int 1(\xi \leq \pi^N_y(s-) \pi^N_x(s-))
				\Lambda^N(ds,x,dy,d\xi,3) \right. \nnb \\
        & & \ \ \ \ \	- \left. \int_0^t \int
				1(\xi \leq \pi^N_x(s-) \pi^N_y(s-))
                                \Lambda^N(ds,dy,x,d\xi,3) \right] .
        \label{eq:wksel}
\end{eqnarray}
The two integrals inside the first set of brackets correspond to
selection, those inside the second set of brackets correspond
to mutation, and those inside the third set of brackets correspond to
replacement sampling. By carrying out computations similar to those done
in the proof of Proposition~\ref{prp:stsel} (and the tightness proof
in [Perkins 2002]), one can conclude the
processes $\{\pi^N: N \in \Nbold\}$ is $C$-tight, and that each weak limit
point $\pi$ satisfies the following martingale problem:
\begin{eqnarray}
	& & \pi_x(t) = \pi_x(0) + \int_0^t \pi_x(s) [m_x(\pi(s))
                        - \overline m_{\pi(s)}]
                +\mu [1 - (2L+1)\pi_x(s)] \ ds + M_x(t), \nnb \\
	& & \mbox{where $M_x$ is a continuous $(\Fcal^\pi_t)$-martingale
		such that $M_x(0) = 0$, and} \nnb \\
	& & \ \ \ \ \ \ \ \ \ \ \ \ \left<M_x,M_y\right>_t = \int_0^t
		\delta_{xy}\pi_x(s) - \pi_x(s) \pi_y(s) \ ds.
	\label{eq:wksel:mart}
\end{eqnarray}
Here $\delta_{xy} = 1(x=y)$.
Ito's formula shows that for $F \in C^2(E)$,
\begin{eqnarray*}
	F(\pi(t)) - F(\pi(0)) - \int_0^t \Gcal F(\pi(s)) \ ds
\end{eqnarray*}
is a bounded continuous martingale, where
\begin{eqnarray*}
	\Gcal F(\pi) & = & \sum_{x=-L}^L [ \pi_x (m_x(\pi)-\overline m_{\pi})
                	+\mu (1 - (2L+1)\pi_x) ]
			\frac{\partial F}{\partial \pi_x} \nnb \\
	& & 	+ \frac{1}{2} \sum_{x=-L}^L \sum_{y=-L}^L
		\pi_x (\delta_{xy} - \pi_y)
		\frac{\partial^2 F}{\partial \pi_x \partial \pi_y}.
\end{eqnarray*}
In particular, if $F$ takes on the
form $F(\pi) = f(\left<\pi,\phi\right>)$ where
$\left<\pi,\phi\right> = \sum_z \pi_z \phi_z$ and $\phi:E\rightarrow \Rbold$
is a function, then
\begin{eqnarray}
	\Gcal F(\pi) & = & \sum_{x=-L}^L f'(\left<\pi,\phi\right>)
		\phi_x [ \pi_x (m_x(\pi) - \overline m_{\pi})
                        +\mu (1 - (2L+1)\pi_x) ] \nnb \\
        & &     + \frac{1}{2} \sum_{x=-L}^L \sum_{y=-L}^L
		f''(\left<\pi,\phi\right>) \phi_x \phi_y
                \pi_x (\delta_{xy} - \pi_y). \label{def:Gweak}
\end{eqnarray}
This is a special (finite $E$ and symmetric mutation) case of generator
for the Fleming-Viot process with selection. The martingale problem
associated with $\Gcal$ has a unique solution
(see Chapter 10.1.1 from [Dawson 1993]), so $\pi^N$ converges weakly to $\pi$
on $D(\Rbold^+,M_1(E))$.

For the process described by the martingale problem~(\ref{eq:wksel:mart}),
Lemma 4.1 from [Ethier and Kurtz 1994] says that
\begin{eqnarray*}
	\nu(d\pi) = C \left( \prod_{x=-L}^L \pi_x \right)^{\mu-1}
		e^{\overline m_\pi} \ d\pi_{-L} \cdots d\pi_L
\end{eqnarray*}
is the unique stationary distribution. Here $C$ is the normalizing constant
such that $\nu$ is a probability measure on $\Pcal(E)$. Notice that if
$m_x(\pi) = 0$ for all $x$ and $\pi$, then $\overline m_\pi = 0$ for all
$\pi$, and in that case, $\nu(d\pi) = C
\left( \prod_{x=-L}^L \pi_x \right)^{\mu-1} \ d\pi_{-L} \cdots d\pi_L$
is the unique stationary distribution of the Fleming-Viot process on
$E$ with symmetric mutation [Ethier and Kurtz 1981]. Adding selection
to the model has the effect of putting weight $e^{\overline m_\pi}$ (where
$\overline m_\pi$ is the mean fitness of $\pi$) on the phenotypical
distribution $\pi$. But even the least fit $\pi$ has weight at least $1$,
since $\exp(\min_\pi \overline m_\pi) = e^0 = 1$, and the fittest $\pi$ has
weight at most $e$. The effect of fitness on the density at a particular
phenotypical distribution $\pi$ is thus only marginal. In contrast,
as we shall see later, fitness will have
a much more pronounced effect in the strong selection model.
For the rest of this part of the thesis, we deal with the strong selection
model outlined in Chapter~\ref{ch:st-model}.

\newpage
\chapter{The Selection-Mutation Equation}
\label{ch:selmut}
As discussed in Chapter~\ref{ch:st-model}, the particle model with strong
selection converges weakly to the selection-mutation equation as the number of
particles $N$ tends to infinity:
\begin{eqnarray}
        \partial_t \pi_x & = & \pi_x ( m_x - \overline m_{\pi} )
                + \mu ( 1 - (2L+1) \pi_x), \label{eq:ODE2}
\end{eqnarray}
where, as before, $m_x = m_x(\pi)$ is the fitness of site $x$
in population $\pi$, and $\overline m_\pi$ is the mean fitness of
the population:
\begin{eqnarray}
	m_x & = & K_x \sum_z B_{x-z} K_z \pi_z \label{eq:mx1} \\
        \overline m_\pi & = & \sum_x \pi_x m_x
                = \sum_{x,z} \pi_x K_x B_{x-z} K_z \pi_z. \label{def:mbar}
\end{eqnarray}
First, we state a few assumptions on the parameters involved:
\begin{enumerate}
\item $x \in E = [-L,L] \cap \Zbold = \{-L,\ldots,-1,0,1,\ldots,L\}$,
\item $K: E \rightarrow (0,1]$,
\item $B: \Zbold \rightarrow [0,1]$.
\end{enumerate}
We establish a few basic facts about the system~(\ref{eq:ODE2}).
First of all, $m_x$ (uniformly in $x$) and $\overline m_\pi$
both lie in $[0,1]$, therefore $m_x - \overline m_\pi \geq -1$.
Thus
\begin{eqnarray}
	\partial_t \pi_x & = & \pi_x ( m_x - \overline m_{\pi} )
                + \mu ( 1 - (2L+1) \pi_x) \nnb \\
	& \geq & - \pi_x + \mu ( 1 - (2L+1) \pi_x) \nnb \\
	& = & \mu - (2L\mu+\mu+1) \pi_x. \label{ineq:muneq0}
\end{eqnarray}
If $\pi_x < \frac{\mu}{2L\mu+\mu+1}$, then $\partial_t \pi_x > 0$ regardless
of $\pi$. This means that there can be no stationary points of
the system~(\ref{eq:ODE2}) with any $\pi_x < \frac{\mu}{2L\mu+\mu+1}$,
and furthermore, for those $x$ where $\pi_x < \frac{\mu}{2(2L\mu+\mu+1)}$
initially, $\pi_x$ will increase at a positive (bounded away from 0) speed and
eventually $\pi_x \geq \frac{\mu}{2(2L\mu+\mu+1)}$ for all $x$.
Secondly, since the sum over all $x$ of the right hand side of~(\ref{eq:ODE2})
is
\begin{eqnarray}
	\lefteqn{ \sum_x \pi_x ( m_x - \overline m_{\pi} )
		+ \mu ( 1 - (2L+1) \pi_x) } \nnb \\
	& = & \sum_x \pi_x m_x - \overline m_{\pi} \sum_x \pi_x + (2L+1) \mu
		- \mu (2L+1) \sum_x \pi_x \nnb \\
	& = & \overline m_{\pi} - \overline m_{\pi} + (2L+1) \mu - (2L+1) \mu
		\nnb \\
	& = & 0, \label{eq:sumto0}
\end{eqnarray}
we have
\begin{eqnarray*}
	\partial_t \sum_x \pi_x = 0,
\end{eqnarray*}
hence the total mass $\sum_x \pi_x$ remains constant. This, together
with the first observation we just made, imply that if
at $t=0$, $\pi(0)$ is a probability measure, i.e. $\sum_x \pi_x(0) = 1$,
then $\pi(t)$ remains a probability measure for all $t$.

Since $\pi_x$ will become instantly nonzero at $x$ where $\pi_x = 0$
initially, we restrict our attention to polymorphic initial conditions,
i.e. $\pi_x (0) > 0$ for all $x$, which is equivalent to saying that
$\pi(0)$ lies in the interior $\overset{_\circ}{\Delta}$ of the set $\Delta$.
For polymorphic initial conditions,
$\pi_x(t) > 0$ for all $x$ and $t$, and therefore~(\ref{eq:ODE2})
can be written as:
\begin{eqnarray}
	\partial_t \pi_x & = & \pi_x \left( m_x - \overline m_{\pi}
		+ \frac{\mu}{\pi_x} - \mu(2L+1) \right). \label{eq:ODE3}
\end{eqnarray}
Furthermore, there is a Lyapunov function
$V_\pi: \overset{_\circ}{\Delta} \rightarrow \Rbold$ for the
dynamical system~(\ref{eq:ODE3}):
\begin{eqnarray}
	V_\pi = \frac{1}{2} \overline m_\pi + \mu \sum_x \log \pi_x.
	\label{def:V}
\end{eqnarray}
Notice that $\partial_{\pi_x} \overline m_\pi = 2 \sum_z K_x B_{x-z} K_z \pi_z
= 2 m_x$.
The assertion that $V_\pi$ is a Lyapunov function for~(\ref{eq:ODE3})
can then be verified by the following simple calculus exercise:
\begin{eqnarray}
        \partial_t V_\pi & = & \sum_x (\partial_{\pi_x} V_\pi)
                (\partial_t \pi_x) \nonumber \\
        & = & \sum_x \left(m_x + \frac{\mu}{\pi_x}\right) \pi_x
                \left(m_x+\frac{\mu}{\pi_x}-\overline m_\pi-\mu(2L+1)\right)
                \nonumber \\
        & = & \sum_x \left(m_x + \frac{\mu}{\pi_x}\right) \pi_x
                \left(m_x+\frac{\mu}{\pi_x}-\overline m_\pi-\mu(2L+1)\right)
		\nnb \\
        & & - \left(\overline m_\pi + \mu(2L+1)\right) \sum_x \pi_x
                \left(m_x+\frac{\mu}{\pi_x}-\overline m_\pi-\mu(2L+1)\right).
		\nnb
\end{eqnarray}
Notice that
$\sum_x \pi_x \left(m_x+\frac{\mu}{\pi_x}-\overline m_\pi-\mu(2L+1)\right)=0$
by~(\ref{eq:sumto0}), therefore
\begin{eqnarray*}
        \partial_t V_\pi & = & \sum_x \pi_x \left(m_x + \frac{\mu}{\pi_x}
                - \overline m_\pi - \mu(2L+1)\right)^2 \geq 0,
\end{eqnarray*}
hence $V_\pi$ is a Lyapunov function for~(\ref{eq:ODE3}) as claimed.

In fact, according to Theorem A.9 of [B\"urger 2000],~(\ref{eq:ODE3}) is
a so-called Svirezhev-Shahshahani gradient system with potential $V_\pi$
as defined in~(\ref{def:V}), i.e. $\partial_t \pi = \tilde\nabla V(\pi)$, where
$\tilde\nabla V(\pi) = G_\pi \nabla V(\pi)$ and $G_\pi$ is 
the matrix formed by entries $g^{xy} = \pi_x (\delta_{xy} - \pi_y)$.
Any gradient system, such as~(\ref{eq:ODE3}), has the property that
all orbits, regardless of initial condition,
converge to some point in the $\omega$-limit set
\begin{eqnarray*}
	D_\omega = \{p: \mbox{$p$ is an accumulation point of $\pi_x(t)$
		as $t\rightarrow\infty$} \}.
\end{eqnarray*}
All points in $D_\omega$ are stationary points of~(\ref{eq:ODE3}).
Since $\hat\pi$ is a stationary point of~(\ref{eq:ODE3})
if and only if
\begin{eqnarray}
	\hat m_x - \overline m_{\hat\pi} + \mu \left(
		\frac{1}{\hat\pi_x} - (2L+1) \right) = 0
	\mbox{ for all } x \in E, \label{cond:ODE2}
\end{eqnarray}
where we write $\hat m_x = m_x(\hat\pi)$, all points in $D_\omega$
satisfies~(\ref{cond:ODE2}). We observe that if
$\hat\pi$ is a stationary distribution and $\hat m_x > \hat m_y$, then
condition~(\ref{cond:ODE2}) means that
$\frac{1}{\hat\pi_x} - (2L+1) < \frac{1}{\hat\pi_y} - (2L+1)$,
therefore $\hat\pi_x > \hat\pi_y$, i.e.
\begin{eqnarray}
	\hat m_x > \hat m_y \Longrightarrow \hat\pi_x > \hat\pi_y.
	\label{mfitmmass}
\end{eqnarray}
In words, fitter sites have more mass.

We will try to characterize the stationary points of the
dynamical system~(\ref{eq:ODE3}) for $K$ and $B$ satisfying the following
conditions:
\begin{eqnarray}
	& & \mbox{$K$ is symmetric and unimodal with $K_0 = 1$,
		and $B$ is of the form } \nnb \\
	& & \mbox{$B_x=b+(1-b) 1_{\{|x| \geq M\}}$
		with $b \in [0,1]$ and $L < M \leq 2L$.} \label{cond:KB}
\end{eqnarray}
Define $l=M-L$, then
for $x\in [-l+1,l-1]$, the cooperation intensity between $x$ and any other
site $z \in [-L,L]$, $B_{x-z}$, is equal to $b$, which means that
\begin{eqnarray}
	m_x =  K_x \sum_z B_{x-z} K_z \pi_z = b K_x \sum_z K_z \pi_z.
	\label{eq:mx2}
\end{eqnarray}

\section{Mild Competition: $b$ close to $1$}
\label{sec:mild}
If $b=1$, then $B_x = b$ for all $x$, hence there is equal competition
between all sites. This actually means that competition plays no part in
how fit site $x$ is and $m_x$ is proportional to $K_x$. Therefore, since
$K_x$ is unimodal (hence $K_x$ is strictly increasing in $[-L,0]$ and strictly
decreasing in $[0,L]$), the fitness should be unimodal, too.
Recall from~(\ref{mfitmmass}) that stationary
distributions of~(\ref{eq:ODE3}) has the property of fitter sites having
more mass, thus we expect the stationary distribution $\hat\pi$
to be unimodal as well. In particular, $\hat\pi$ should attain its maximum
at $x=0$.
As $\mu \rightarrow 0$, we expect the ``peak'' of $\hat\pi$
concentrated around $0$ to become sharper and sharper, approaching $\delta_0$,
the $\delta$-measure concentrated at $0$. In fact, as we shall see,
$b$ only needs to be somewhat close to $1$ for this behaviour to occur.

We now show that for any stationary distribution $\hat\pi$ of~(\ref{eq:ODE3}),
site $0$ is fitter than any other site for $b \in (\frac{1}{2},1]$
sufficiently close to 1, i.e. $\hat m_0 > \hat m_x$ if $x\neq 0$.
This will mean that as $\mu \rightarrow 0$, any $\hat\pi$ approaches the
$\delta$-measure. Recall from~(\ref{eq:mx2}) that for $x\in [-l+1,l-1]$,
we have
\begin{eqnarray*}
        m_x = b K_x \sum_z K_z \pi_z.
\end{eqnarray*}
We recall that $K_0=1$ and $K_x$ is assumed to
be strictly increasing in $[-L,0]$ and strictly decreasing in $[0,L]$,
therefore $K_x$ attains its maximum at $x=0$, and thus for
$x\in [-l+1,-1]\cup [1,l-1]$, $K_x - K_0 \leq K_1 - K_0 < 0$.
Therefore for $x\in [-l+1,-1]\cup [1,l-1]$,
\begin{eqnarray*}
        m_x - m_0 = b (K_x-K_0) \sum_z K_z \pi_z
	\leq - b (K_0 - K_1) \sum_z K_z \pi_z.
\end{eqnarray*}
We now apply the bound $\sum_z K_z \pi_z \geq
\inf_{\pi\in\Delta} \sum_z K_z \pi_z = K_L$
to the above inequality to obtain
\begin{eqnarray*}
	m_x - m_0 \leq - b (K_0 - K_1) K_L.
\end{eqnarray*}
Since $b$ is assumed to be $> \frac{1}{2}$, we have
for $x\in [-l+1,-1]\cup [1,l-1]$,
\begin{eqnarray}
        m_x - m_0 < - \frac{1}{2} (K_0 - K_1) K_L. \label{ineq:bd_fit1}
\end{eqnarray}

We also bound the fitness for sites $x$ in $[l,L]$:
\begin{eqnarray}
        m_x & = & b K_x \sum_z K_z \pi_z
                + (1-b) K_x \sum_{|x-z|\geq M} K_z \pi_z \nonumber \\
        & = & b K_x \sum_z K_z \pi_z + (1-b) K_x \sum_{x-z\geq M} K_z \pi_z
                \ \ \ \ \ \mbox{ since } x \geq l \nonumber \\
        & \leq & b K_x \sum_z K_z \pi_z + (1-b) K_x \sum_{z=-L}^{x-M} K_z
                \ \ \ \ \ \mbox{ since } \pi_z \leq 1 \nonumber \\
	& \leq & b K_x \sum_z K_z \pi_z + (1-b) K_x \sum_{z=-L}^{-l} K_z
		\nnb \\
        & \leq & b K_x \sum_z K_z \pi_z + (1-b) (L-l+1) K_x K_l,
		\label{eq:bd_fit3}
\end{eqnarray}
where in the last line, we use the following: for
$z\in [-L,x-M] \subset [-L,-l]$, $K_z \leq K_{-l} = K_l$.
Similarly, for $x\in [-L,-l]$, we have the same bound~(\ref{eq:bd_fit3}).
Therefore for $x\in [-L,-l]\cup [l,L]$, we have
\begin{eqnarray}
        m_x & \leq & b K_x \sum_z K_z \pi_z + (1-b) (L-l+1) K_x K_l \nnb \\
        & \leq & b K_x \sum_z K_z \pi_z + (1-b) (L-l+1) K_l^2,
        \label{eq:bd_fit0}
\end{eqnarray}
where again we use the bound $K_x \leq K_l$ for $x \in [-L,-l]\cup [l,L]$.
We use~(\ref{eq:mx2}) and~(\ref{eq:bd_fit0}) to estimate $m_x - m_0$ for
$x\in [-L,-l]\cup [l,L]$:
\begin{eqnarray*}
	m_x - m_0 & \leq & b K_x \sum_z K_z \pi_z + (1-b) (L-l+1) K_l^2
		- b K_0 \sum_z K_z \pi_z \\
	& = & b (K_x - K_0) \sum_z K_z \pi_z + (1-b) (L-l+1) K_l^2.
\end{eqnarray*}
Since $K_x - K_0 < 0$ for $x\in [-L,-l]\cup [l,L]$, we have
\begin{eqnarray*}
	m_x - m_0 & \leq & - b \inf_{x\in [-L,-l]\cup [l,L]}
		(K_0 - K_x) \inf_{\pi\in\Delta} \sum_z K_z \pi_z
		+ (1-b) (L-l+1) K_l^2 \nnb \\
	& = & - b (K_0 - K_l) K_L + (1-b) (L-l+1) K_l^2.
\end{eqnarray*}
Since $b$ is assumed to be $> \frac{1}{2}$, the above bound can be
simplified to:
\begin{eqnarray}
	m_x - m_0 \leq (1-b) (L-l+1) K_l^2 - \frac{1}{2} (K_0 - K_l) K_L.
		\label{ineq:bd_fit2}
\end{eqnarray}
Thus if $b$ is so close to $1$ that
\begin{eqnarray}
	1-b \leq \frac{(K_0 - K_l) K_L}{4 (L-l+1) K_l^2}, \label{cond:bd_fit2}
\end{eqnarray}
then let $\delta_1 = \frac{1}{4} (K_0 - K_l) K_L$ be a positive constant
and we have
\begin{eqnarray}
	(1-b) (L-l+1) K_l^2 \leq \frac{1}{4}(K_0 - K_l) K_L
	= \frac{1}{2} (K_0 - K_l) K_L - \delta_1. \label{ineq:bd_fit4}
\end{eqnarray}
Thus if condition~(\ref{cond:bd_fit2}) holds, then~(\ref{ineq:bd_fit2})
and~(\ref{ineq:bd_fit4}) imply that for $x \in [-L,-l]\cup [l,L]$,
\begin{eqnarray}
	m_x - m_0 \leq - \delta_1. \label{ineq:bd_fit5}
\end{eqnarray}
Define $\delta = \min\left( \delta_1, \frac{1}{2} (K_0 - K_1) K_L \right)$,
then the estimates in~(\ref{ineq:bd_fit1}) and~(\ref{ineq:bd_fit5})
mean that for all $x \in E \backslash \{ 0 \}$ and any $\pi$,
\begin{eqnarray}
	m_x - m_0 \leq -\delta. \label{eq:bd_fit6}
\end{eqnarray}
This shows that for $b$ satisfying condition~(\ref{cond:bd_fit2}) and
for any $\pi$, site $0$ is fitter than any other site.
We will use this bound to establish the following.

\begin{THM} If $K$ is symmetric and unimodal with $K_0 = 1$,
and $B_x=b+(1-b) 1_{\{|x| \geq M\}}$ with $L < M \leq 2L$, $l=M-L$, and
$b \in \left[1-\frac{(K_0 - K_l) K_L}{4 (L-l+1) K_l^2},1\right]$,
then as $\mu \rightarrow 0$,
\[\sup \{\|\hat\pi^\mu - \delta_0\|_\infty: \hat\pi^\mu
\mbox{ is an stationary distribution for mutation parameter } \mu \}
\rightarrow 0. \]
\end{THM}
\proof Recall condition~(\ref{cond:ODE2}): $\hat\pi^\mu$ is a stationary
distribution of~(\ref{eq:ODE3}) if and only if for all $x \in [-L,L]$,
\begin{eqnarray}
        \hat m^\mu_x - \overline m_{\hat \pi^\mu}
                + \mu \left(\frac{1}{\hat\pi^\mu_x} - (2L+1)\right)
        = 0, \label{eq:lesscomp}
\end{eqnarray}
where we write $\hat m^\mu_x = m_x (\hat \pi^\mu)$.
We make the following observations using condition~(\ref{eq:lesscomp}): if
$\hat m^\mu_x \geq \overline m_{\hat \pi^\mu}$, then
$- \frac{1}{\hat\pi^\mu_x} + (2L+1) \geq 0$,
which implies that $\hat\pi^\mu_x \geq 1/(2L+1)$;
similarly, if $m_x < \overline m_{\hat \pi^\mu}$
then $\hat\pi^\mu_x < 1/(2L+1)$.

Since $1/\hat\pi^\mu_x \geq 1$, the following bound holds for all $x$:
\begin{eqnarray*}
        \hat m^\mu_x - \overline m_{\hat \pi^\mu} + \mu (1 - (2L+1)) \leq 0,
\end{eqnarray*}
which implies that
\begin{eqnarray}
        \hat m^\mu_x - \overline m_{\hat \pi^\mu} \leq 2 L \mu.
		\label{eq:lesscomp2}
\end{eqnarray}
We consider $\mu$ small enough such that $2L\mu < \frac{\delta}{2}$,
where $\delta$ is defined right after~(\ref{ineq:bd_fit5}).
Then the estimate~(\ref{eq:lesscomp2}) applied to $x=0$ means that
\begin{eqnarray*}
        \hat m^\mu_0 - \overline m_{\hat \pi^\mu}
		\leq 2 L \mu < \frac{\delta}{2},
\end{eqnarray*}
which implies that
\begin{eqnarray*}
	\hat m^\mu_0 < \overline m_{\hat \pi^\mu} + \frac{\delta}{2}.
\end{eqnarray*}
Applying the above to estimate~(\ref{eq:bd_fit6}), we get,
for all $x\neq 0$,
\begin{eqnarray}
	\hat m^\mu_x \leq \hat m^\mu_0 - \delta
		< \overline m_{\hat \pi^\mu} - \frac{\delta}{2},
	\label{eq:lesscomp3}
\end{eqnarray}
In particular, the only $x$ where
$\hat m^\mu_x \geq \overline m_{\hat \pi^\mu}$ is $x=0$.

Using the bound~(\ref{eq:lesscomp3}), condition~(\ref{eq:lesscomp})
implies that for $x\neq 0$,
\begin{eqnarray*}
        \mu \left(\frac{1}{\hat\pi^\mu_x} - (2L+1)\right)
                \geq \frac{\delta}{2}.
\end{eqnarray*}
Therefore
\begin{eqnarray*}
        \frac{1}{\hat\pi^\mu_x} \geq \frac{\delta}{2\mu} + (2L+1),
\end{eqnarray*}
which $\rightarrow \infty$ as $\mu \rightarrow 0$.
Hence $\hat\pi^\mu_x \rightarrow 0$ as
$\mu \rightarrow 0$ for all $x\neq 0$.
Notice that the proof does not depend on which $\hat\pi^\mu$ we pick,
therefore we are done.
\qed

\vspace{.3cm}

\begin{REM}
The crucial estimate for the above proof is~(\ref{eq:bd_fit6}). In the case
of equal competition, i.e. $b=1$, it is very easy to
derive~(\ref{eq:bd_fit6}):
\begin{eqnarray*}
	m_x - m_0 = b (K_x - K_0) \sum_z K_z \pi_z
	& \leq & - b \inf_{x\neq 0} (K_0 - K_x)
		\inf_{\pi\in\Delta} \sum_z K_z \pi_z \\
	& = & - b (K_0-K_1) K_L,
\end{eqnarray*}
which is a positive constant independent of $\mu$.
\end{REM}

\section{Intense Competition: $b$ close to $0$}
Results from Chapter~\ref{sec:mild} show that if the competition between
pairs of sites that are far away from each other is not intense, i.e. $b$
is close to $1$, then as $\mu \rightarrow 0$, the stationary distribution(s)
converge to $\delta_0$, and therefore, there is no speciation. In this section,
we show that if there is the most intense competition between pairs of
sites that are far away from each other, i.e. $b=0$, then
we do see speciation in the stationary distribution(s).
More interesting behaviour
arises in the case of positive but small $b$. We will show that, in this
case, if $\mu$ is small enough, then
there are at least two vastly different stationary distributions, one
resembling the $\delta$-measure, the other bimodal and having almost zero
mass in the middle; on the other hand, if $\mu$ is sufficiently large, then
all stationary distributions are bimodal and have little mass in the middle.
Thus for small $\mu$, whether speciation occurs eventually in
the dynamical system depends on the initial state of the system.
But for large enough $\mu$, speciation will occur eventually if one waits
long enough.
We first illustrate this behaviour in a system with 3 phenotypes
$\{-1,0,1\}$, whose stationary points we can calculate explicitly.

\subsection{Study of A One-dimensional System}
\label{sec:one-dim}
Here we will take the simplest possible scenario, and show that the
dynamical system~(\ref{eq:ODE3}) has exactly two stable stationary points.
Let $L=1$, $E=\{-1,0,1\}$, $K_0=1$, and $K_1 = K_{-1} = \frac{1}{2}$.
Let $B_x=b+(1-b) 1_{\{|x|\geq 2\}}$, i.e. phenotype $-1$
cooperates with phenotype $0$ and itself at level $b$ (i.e. some competition),
and cooperates with phenotype $1$ at level $1$ (i.e. no competition).
We only consider symmetric distributions, i.e. $\pi_{-1}=\pi_1$. Taking into
account that $\pi_{-1}+\pi_0+\pi_1=1$, the dynamical system~(\ref{eq:ODE3})
has only one variable, say $\pi_0$. The fitness of the 3 sites in $E$ and
the mean fitness are:
\begin{eqnarray*}
m_0 & = & K_0 b (K_0 \pi_0 + 2 K_1 \pi_1) = b(\pi_0 + \pi_1) \\
m_1 & = & K_1 (b K_0 \pi_0 + (1+b) K_1 \pi_1)
	= \frac{1}{2} (b\pi_0 + \frac{1+b}{2} \pi_1) \\
\overline m_\pi & = & b\pi_0 (\pi_0 + \pi_1)
	+ \pi_1 (b\pi_0 + \frac{1+b}{2} \pi_1)
\end{eqnarray*}
Thus the dynamical system~(\ref{eq:ODE3}) with variable $\pi_0$ can
be written as a single ordinary differential equation:
\begin{eqnarray*}
\partial_t \pi_0 & = & \pi_0 (m_0 - \overline m_\pi)
	+ \mu (1-3 \pi_0) \\
& = & \pi_0 \left(b(\pi_0 + \pi_1)
	- b\pi_0 (\pi_0 + \pi_1)
	- \pi_1 (b\pi_0 - \frac{1+b}{2} \pi_1)\right)
	+ \mu (1-3\pi_0)
\end{eqnarray*}
%       & = & \pi_0 \left( b (\pi_0 + \pi_1) - b\pi_0 (\pi_0 + \pi_1)
%               - \pi_1 (b\pi_0 - \frac{b+1}{2} \pi_1)\right)
%                + \mu (1-3\pi_0) \\
%       & = & \pi_0 \left( b (\pi_0 + \pi_1) - b\pi_0^2 - 2b\pi_0\pi_1
%               + \frac{b+1}{2} \pi_1^2 \right)
%                + \mu (1-3\pi_0) \\
%       & = & \pi_0 \left( b (\pi_0 + \frac{1-\pi_0}{2}) - b\pi_0^2
%               - b\pi_0 (1-\pi_0) + \frac{b+1}{8} (1-\pi_0)^2 \right)
%                + \mu (1-3\pi_0) \\
%       & = & \pi_0 \left( b \frac{1+\pi_0}{2} - b\pi_0
%               + \frac{b+1}{8} (1-2\pi_0+\pi_0^2) \right)
%               + \mu (1-3\pi_0) \\
%       & = & \pi_0 \left( b \frac{1-\pi_0}{2}
%               - \frac{b-2b\pi_0+b\pi_0^2 + 1-2\pi_0+\pi_0^2}{8} \right)
%                + \mu (1-3\pi_0) \\
%       & = & \pi_0 \frac{-(b+1)\pi_0^2 + (2-2b)\pi_0 + (3b-1)}{8}
%               + \mu (1-3\pi_0) \\
Substituting in $\pi_1 = \frac{1}{2} (1-\pi_0)$, we get
\begin{eqnarray}
        \partial_t \pi_0 = -\frac{b+1}{8}\pi_0^3 + \frac{1-b}{4}\pi_0^2
                + \frac{3b-1-24\mu}{8}\pi_0 + \mu. \label{eq:ODE_1d}
\end{eqnarray}
Now we take $b=\frac{1}{5}$, then~(\ref{eq:ODE_1d}) can be simplified:
\begin{eqnarray*}
        \partial_t \pi_0 & = & -\frac{3}{20}\pi_0^3 + \frac{1}{5}\pi_0^2
                - \left(\frac{1}{20} + 3\mu\right)\pi_0 + \mu \\
        & = & -\frac{3}{20} \left(\pi_0^3 - \frac{4}{3}\pi_0^2
                + \left(\frac{1}{3}+20\mu\right)\pi_0 - \frac{20}{3}\mu\right).
\end{eqnarray*}
                                                                                
\begin{figure}[t]
\centering
\subfigure[$\mu=\frac{1}{150}$]{
\includegraphics[height = 2.5in, width = 2.5in]{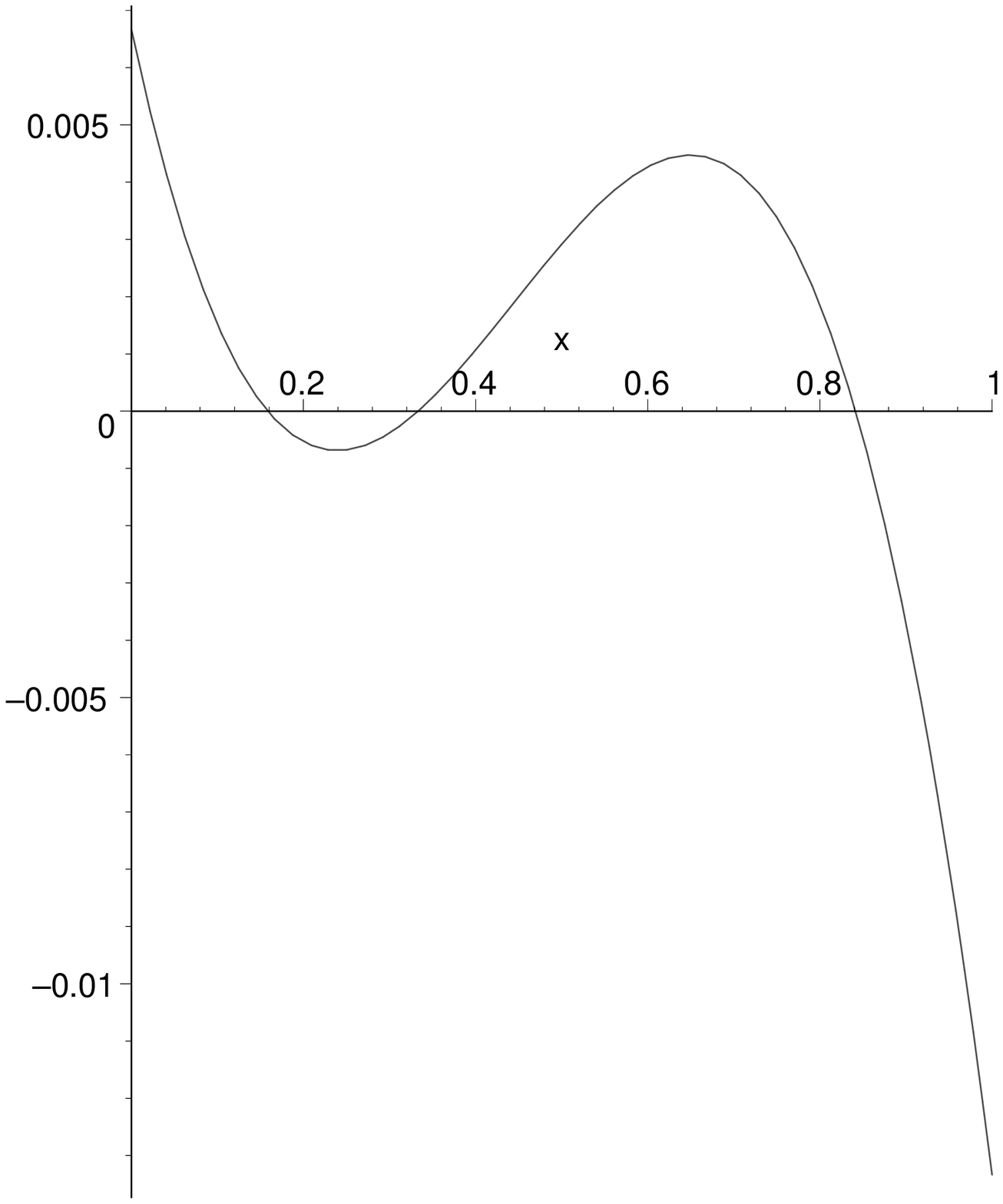}}
\subfigure[$\mu=\frac{1}{70}$]{
\includegraphics[height = 2.5in, width = 2.5in]{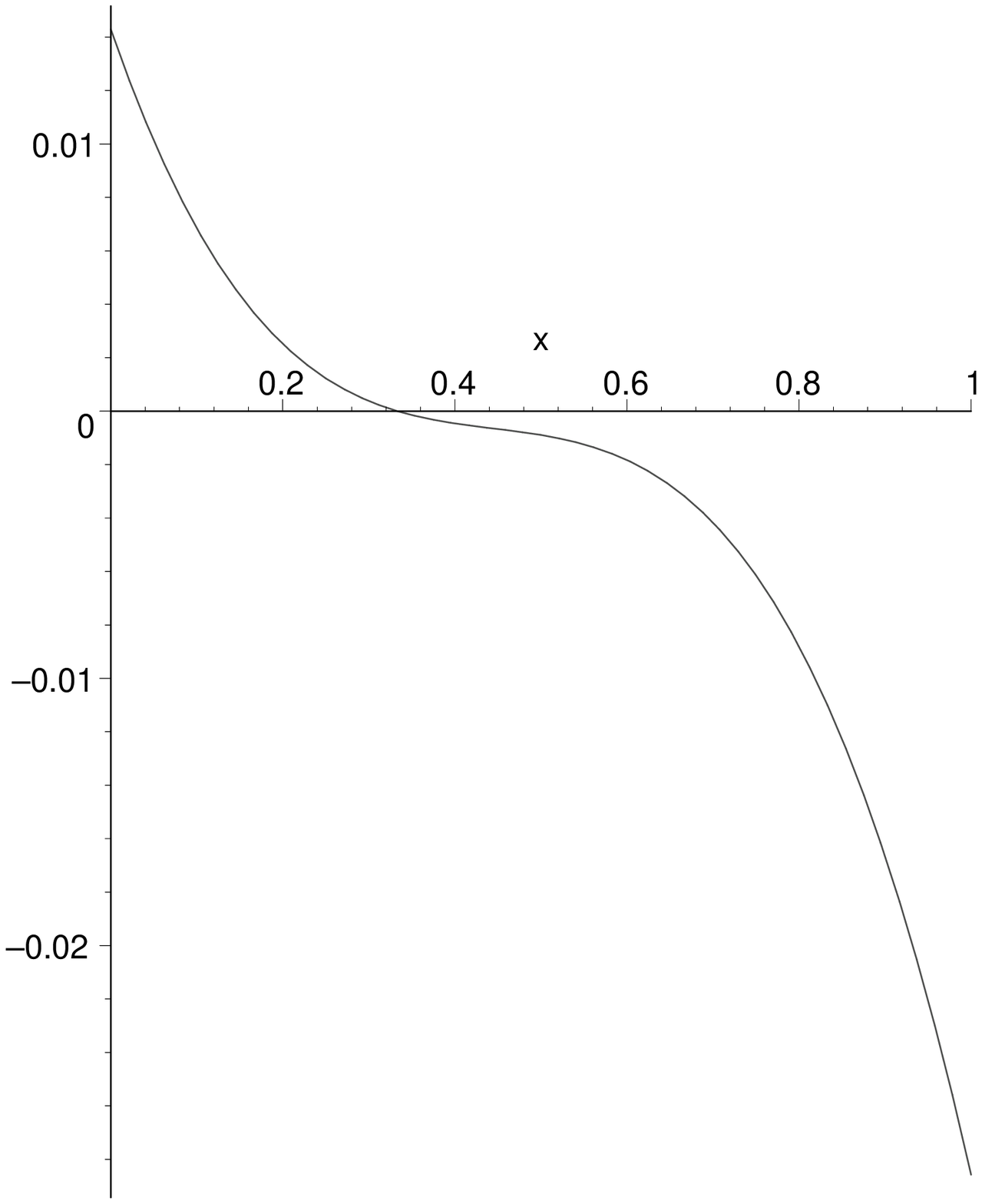}}
\caption{One dimensional case: $b=\frac{1}{5}$}
\label{fig:smallb0}
\end{figure}

We define the polynomial $p(x)=x^3 - \frac{4}{3}x^2 + (\frac{1}{3}+20\mu)x
- \frac{20}{3}\mu$, then $p$ has roots $x_1 = \frac{1}{3}$,
$x_2 = \frac{1}{2}(1+\sqrt{1-80 \mu})$,
and $x_3 = \frac{1}{2}(1-\sqrt{1-80 \mu})$. Notice that the root $x_1$
does not depend on $\mu$, although this is only true for $b=\frac{1}{5}$.
For other $b$'s, all roots depend on $\mu$. Two plots of
$-\frac{3}{20} p(x)$ are shown in figure~\ref{fig:smallb0}.
                                                                                
For $\mu<\frac{1}{80}$, there are three real roots, as in
figure~\ref{fig:smallb0}(a). In this case, if $\pi_0(0) < \frac{1}{3}$, then
$\pi_0(t) \rightarrow x_3 \approx 0.16$
as $t\rightarrow \infty$; since $\pi_{-1} = \pi_1$
and $\pi_{-1}+\pi_0+\pi_1=1$, $\pi_{-1} = \pi_1 \rightarrow 0.42$ as
$t \rightarrow \infty$. This distribution has much larger mass on sites $-1$ and
$1$ than on site $0$, thus we can say that speciation occurs eventually.
But if $\pi_0(0) > \frac{1}{3}$, then $\pi_0(t) \rightarrow x_2 \approx 0.84$
as $t\rightarrow \infty$, then $\pi_{-1} = \pi_1 \rightarrow 0.08$.
This distribution has much larger mass on site $0$ than on sites $-1$ and $1$,
and we say that speciation never occurs. On the other hand,
if $\mu>\frac{1}{80}$, there is only one real root, as in
figure~(\ref{fig:smallb0}(b). In this case, regardless of initial condition,
$\pi_0(t) \rightarrow \frac{1}{3}$ as $t\rightarrow \infty$.
                                                                                
\newpage
Two plots for $b=\frac{1}{10}$ are shown below.
For sufficiently large $\mu$, e.g. $\mu=\frac{1}{300}$
in figure~\ref{fig:smallb}(b), regardless of initial condition,
$\pi$ converges to a configuration with little mass in the middle
(approximately $0.037$), i.e.
speciation occurs. But for $\mu$ sufficiently small,
e.g. $\mu=\frac{1}{500}$ in figure~\ref{fig:smallb}(a), there may or may not
be speciation depending on the initial condition. The time evolution
of $\pi_0$ with different initial conditions for both
$\mu=\frac{1}{300}$ and $\mu=\frac{1}{500}$ are also shown.
{
\psfrag{ttt}{$t$}
\psfrag{p0}{$\pi_0$}
\begin{figure}[h!!]
\centering
\subfigure[$\mu=\frac{1}{500}$]{
\includegraphics[height = 2.3in, width = 2.5in]{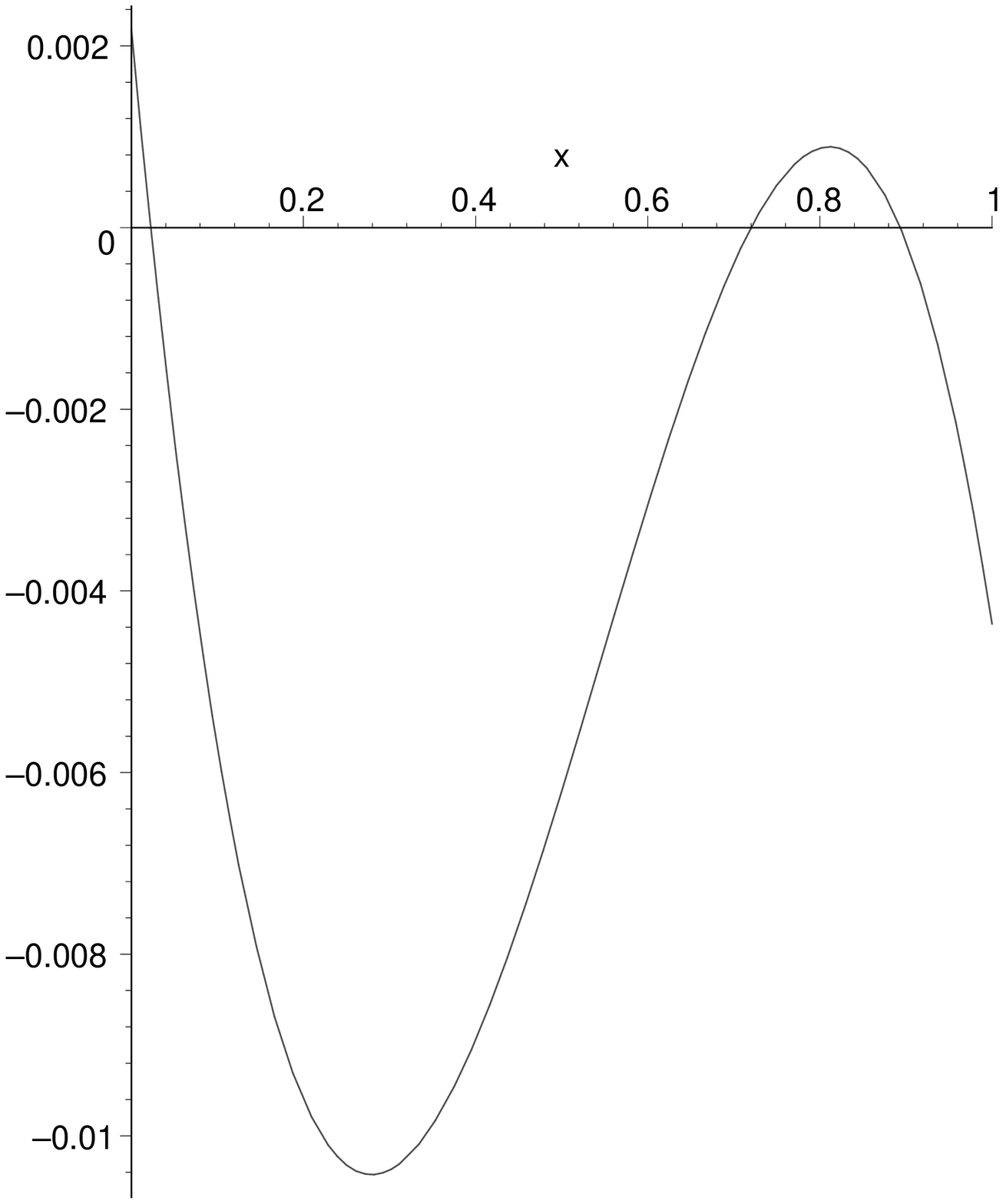}}
\subfigure[$\mu=\frac{1}{300}$]{
\includegraphics[height = 2.3in, width = 2.5in]{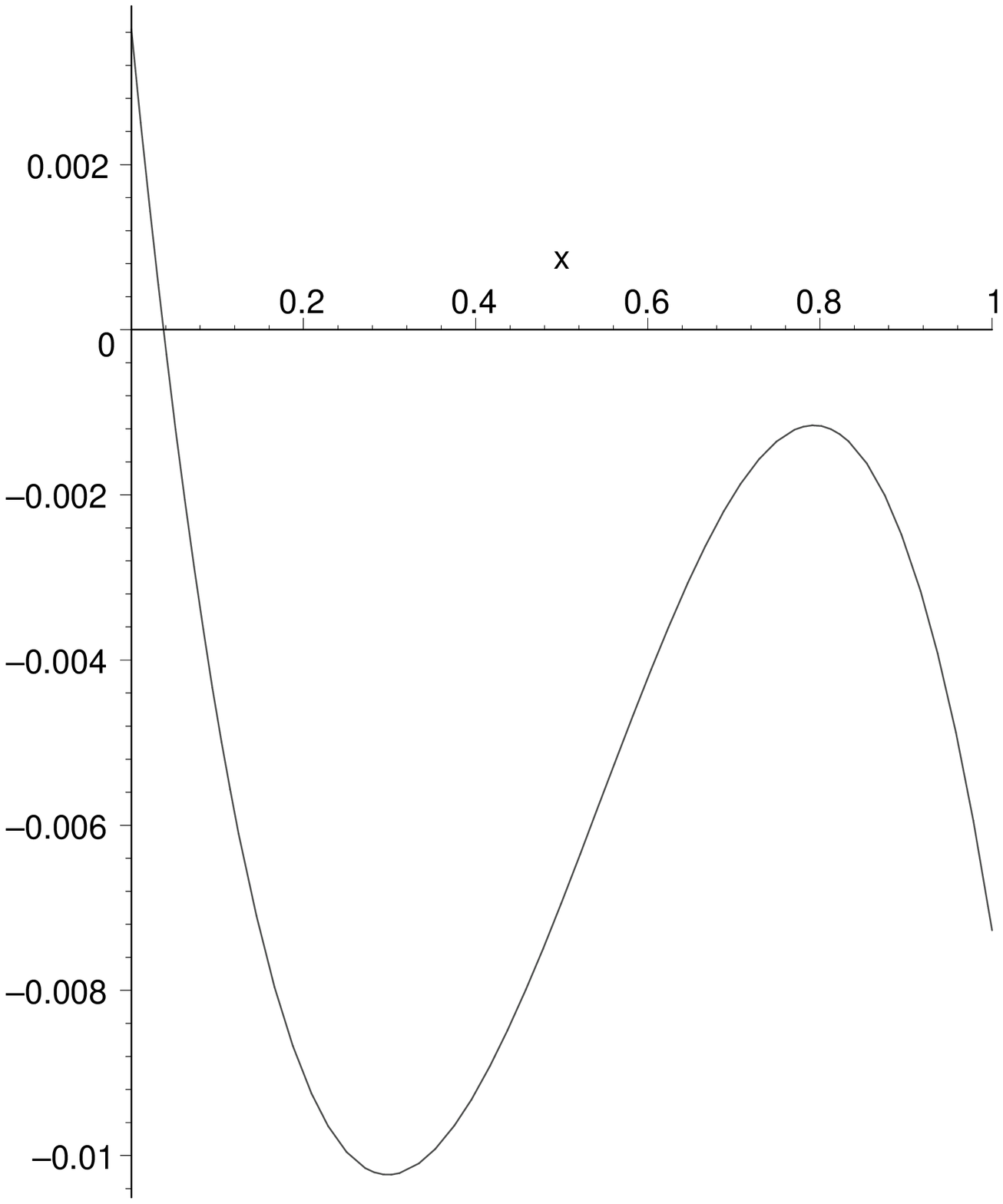}}
\subfigure[$\mu=\frac{1}{500}$]{
\includegraphics[height = 2.3in, width = 2.5in]{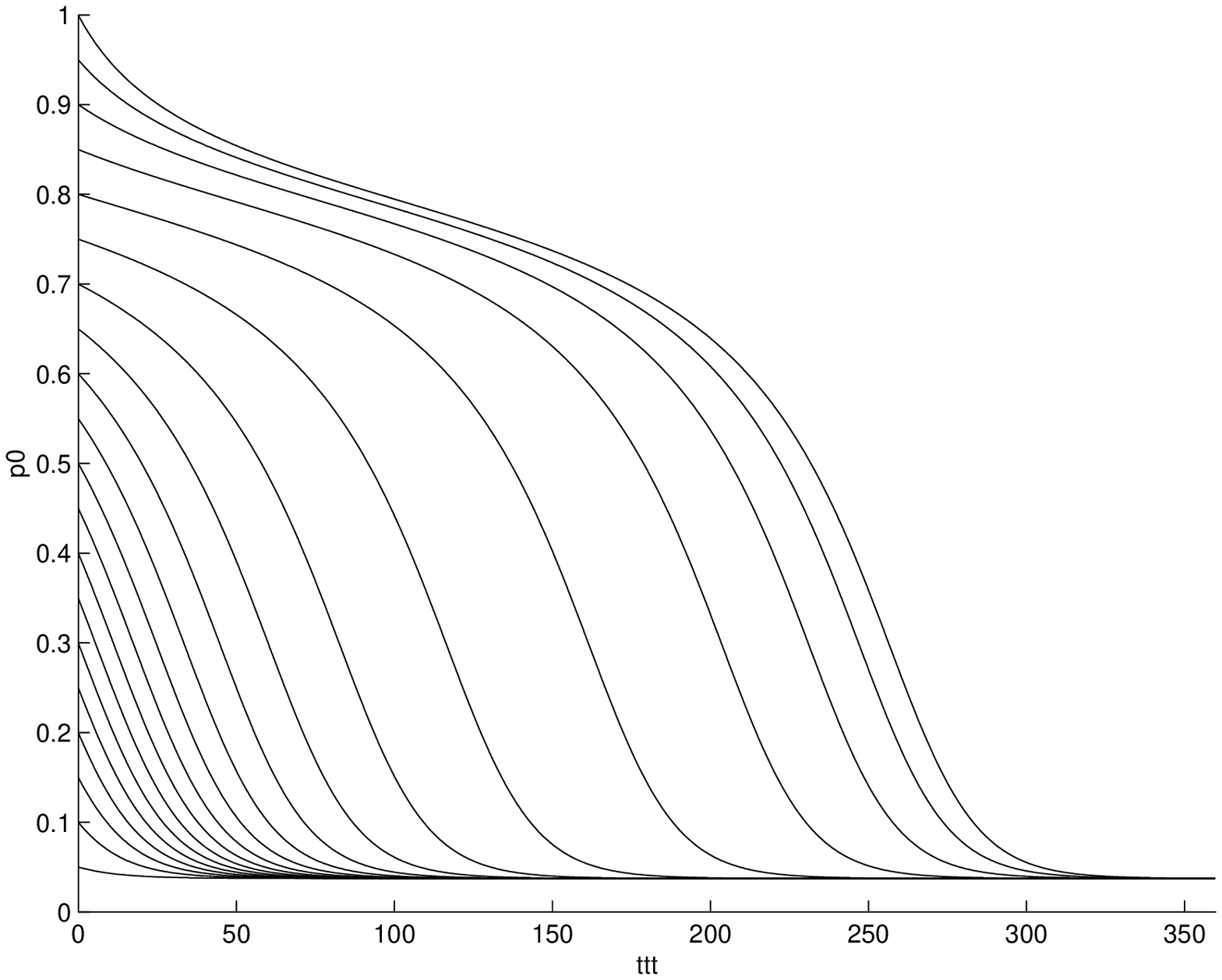}}
\subfigure[$\mu=\frac{1}{300}$]{
\includegraphics[height = 2.3in, width = 2.5in]{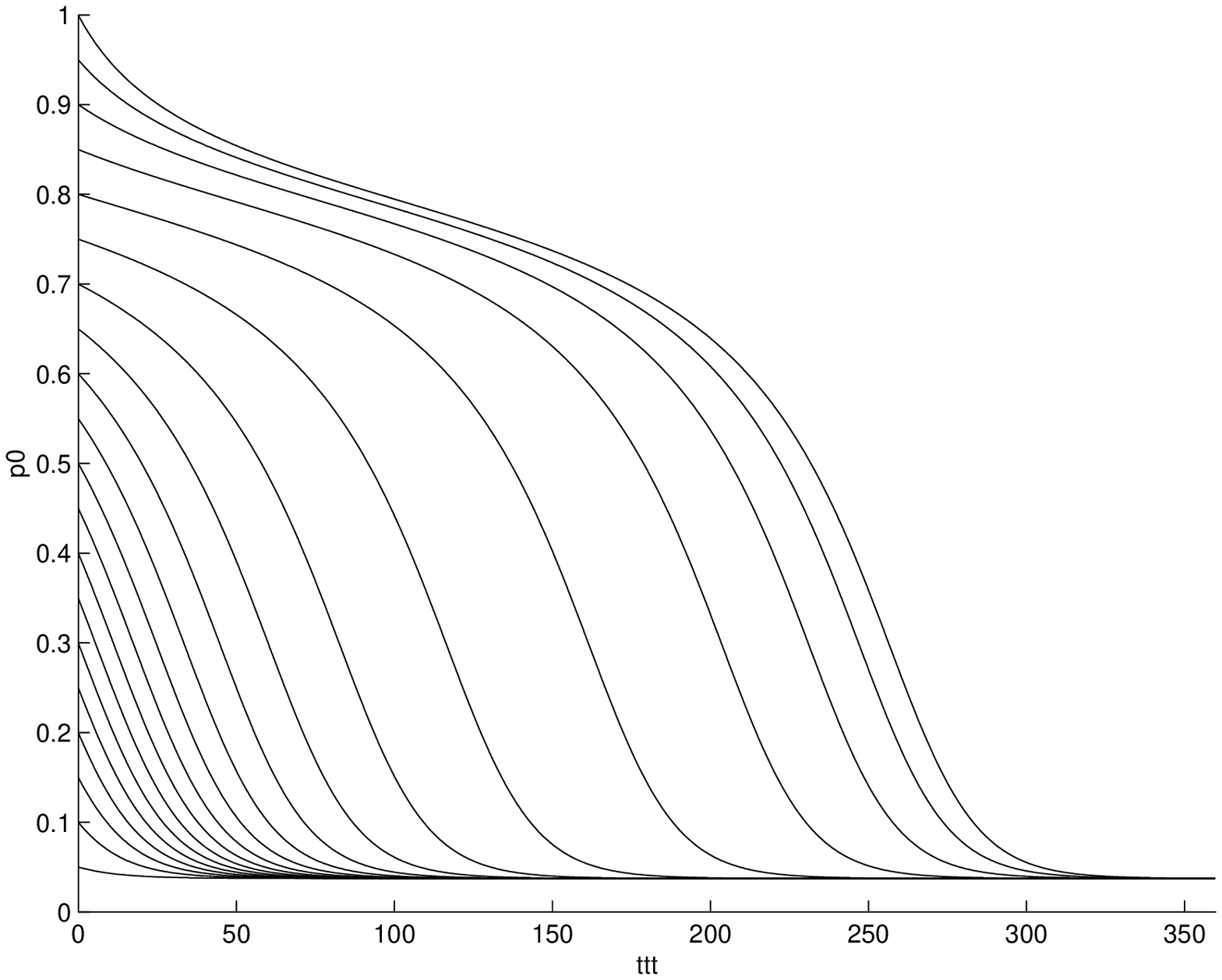}}
\caption{One dimensional case: $b=\frac{1}{10}$}
\label{fig:smallb}
\end{figure}
}
                                                                                
\clearpage
\subsection{Large enough $\mu$}
Analysis of the dynamical system~(\ref{eq:ODE3}) in its simplest form in
the last subsection shows that if $\mu$ is small, there may be two vastly
different types of stationary points for~(\ref{eq:ODE3}).
In subsections~\ref{sec:sm_mu1} and~\ref{sec:sm_mu2},
we establish this for~(\ref{eq:ODE3}) in its general form.
But in this section, we examine the behaviour of~(\ref{eq:ODE3})
when $\mu \geq \frac{b}{4 K_L^2 (L-l)}$ and establish that all stationary
points of~(\ref{eq:ODE3}) are bimodal.

We maintain the assumption in~(\ref{cond:KB}) that $K$ is symmetric and
unimodal with $K_0=1$, and $B$ is of the form $B_x=b+(1-b) 1_{\{|x| \geq M\}}$
with $b \in [0,1]$ and $L < M \leq 2L$.
We will need a uniform (in $\mu$) lower bound on
$\overline m_{\hat\pi^\mu}$. We first establish a crude lower bound
on $\overline m_{\hat\pi^\mu}$ that does depend on $\mu$.
Recall from~(\ref{cond:ODE2}) that $\hat\pi^\mu$ is a stationary distribution
if and only if for all $x \in [-L,L]$,
\begin{eqnarray}
        \hat m^\mu_x - \overline m_{\hat\pi^\mu}
        + \mu \left(\frac{1}{\hat\pi^\mu_x} - (2L+1)\right) = 0.
        \label{eq:cond3}
\end{eqnarray}
Therefore
\begin{eqnarray*}
        \mu \left(\frac{1}{\hat\pi^\mu_x} - (2L+1)\right)
        = \overline m_{\hat\pi} - \hat m^\mu_x
        \leq \overline m_{\hat\pi}
	\leq \sup_{x\in E, \pi\in\Delta} m_x(\pi).
\end{eqnarray*}
Since $K$ and $B$ lie in $[0,1]$, $m_x = K_x \sum_z B_{x-z} K_z \pi_z
\leq \sum_z \pi_z = 1$, and therefore
\begin{eqnarray*}
        \mu \left(\frac{1}{\hat\pi^\mu_x} - (2L+1)\right) \leq 1,
\end{eqnarray*}
which means that
\begin{eqnarray}
        \hat\pi^\mu_x \geq \frac{1}{2L+1+\frac{1}{\mu}} \geq \mu.
		\label{ineq:min_pi}
\end{eqnarray}
For $x\in [l,L]$, we have (recall the first three steps of~(\ref{eq:bd_fit3})),
\begin{eqnarray}
        m_x(\pi) & = & b K_x \sum_z K_z \pi_z
                + (1-b) K_x \sum_{z=-L}^{x-M} K_z \pi_z \nnb \\
        & \geq & b K_L \sum_z K_z \pi_z
                + (1-b) K_L \sum_{z=-L}^{x-M} K_z \pi_z \ \ \ \ \
		\mbox{ since $K_x$ is decreasing in $[0,L]$} \nnb \\
        & \geq & b K_L^2 + (1-b) K_L^2 \nnb \\
	& = & K_L^2. \label{eq:min_fit}
\end{eqnarray}
Since $K$ and $B$ are both symmetric, the same bound applies if
$x\in [-L,-l]$. Therefore~(\ref{ineq:min_pi}) and~(\ref{eq:min_fit})
imply the following estimate of the mean fitness $\overline m_\pi$:
\begin{eqnarray}
        \overline m_{\hat\pi^\mu} & = & \sum_x m_x(\hat\pi^\mu) \hat\pi^\mu_x
	\geq \sum_{x=l}^L m_x(\hat\pi^\mu) \hat\pi^\mu_x
		+ \sum_{x=-L}^{-l} m_x(\hat\pi^\mu) \hat\pi^\mu_x \nnb \\
        & \geq & 2 \sum_{x=l}^L K_L^2 \mu
        = 2 K_L^2 \mu (L-l). \label{ineq:mean_fit}
\end{eqnarray}
This crude lower bound on $\overline m_{\hat\pi^\mu}$ depends on $\mu$,
but in Lemma~\ref{lem:largemu} we will improve it such that
it does not depend on $\mu$ for sufficiently small $b$. For now,
we use it to establish an estimate on $\hat\pi^\mu([-l+1,l-1])$ that
will be needed for the proof of Lemma~\ref{lem:largemu} below.
Condition~(\ref{eq:cond3}) implies that for $x\in [-l+1,l-1]$,
\begin{eqnarray}
        \mu \left(\frac{1}{\hat\pi^\mu_x} - (2L+1)\right)
        = \overline m_{\hat\pi^\mu} - \hat m^\mu_x. \label{eq:cond5}
\end{eqnarray}
For $x\in [-l+1,l-1]$, $\sum_z B_{x-z} K_z \pi_z = b \sum_z K_z \pi_z$ is
constant, thus the maximum fitness is attained at $x=0$, with
\begin{eqnarray}
        m_0(\pi) = b K_0 \sum_z K_z \pi_z \leq b, \label{eq:max_fit}
\end{eqnarray}
since $K_x$ is increasing in $[-L,0]$ and decreasing in $[0,L]$.
If $b \leq 2 K_L^2 \mu (L-l)$, then we can apply the
estimate~(\ref{ineq:mean_fit}) on $\overline m_\pi$,
and bound the right hand side of~(\ref{eq:cond5}): 
\begin{eqnarray*}
 	\overline m_{\hat\pi^\mu} - \hat m^\mu_x
	\geq \overline m_{\hat\pi^\mu} - \hat m^\mu_0
	\geq 2 K_L^2 \mu (L-l) - b,
\end{eqnarray*}
which is $\geq 0$ if $b \leq 2 K_L^2 \mu (L-l)$.
Thus for $b \leq 2 K_L^2 \mu (L-l)$,
(\ref{eq:cond5}) implies that $\frac{1}{\hat\pi^\mu_x} - (2L+1) \geq 0$, i.e.
\begin{eqnarray*}
        \hat\pi^\mu_x \leq \frac{1}{2L+1}.
\end{eqnarray*}
Hence we can bound the mass in $[-l+1,l-1]$:
\begin{eqnarray}
	\hat\pi^\mu([-l+1,l-1]) \leq \frac{2l+1}{2L+1}. \label{ineq:mass_mid}
\end{eqnarray}
Before we state the theorem of this section, we establish
the following lemma, which improves upon the
bound~(\ref{ineq:mean_fit}) such that it does not depend on $\mu$:
                                                                                
\begin{LEM} For $\mu$ and $b$ in the region
$R_1 = \{(\mu,b): 0 \leq b \leq \min( 4 \mu K_L^2 (L-l),
\frac{K_L^2 (L-l)}{4 (2L+1)^3}) \}$, there is a positive constant
$c_1$ that depends on $L$, $l$, and $K$ but not on $\mu$,
such that $\overline m_{\hat\pi^\mu} > c_1$ for any stationary distribution
$\hat\pi^\mu$.
\label{lem:largemu}
\end{LEM}
\proof
We notationally suppress the dependence of $\hat\pi^\mu$,
$\hat m^\mu_x$, and $\overline m_{\hat\pi^\mu}$ on $\mu$. Suppose that
\begin{eqnarray}
        \hat\pi_x \hat\pi_y < \delta \
	\mbox{ for all $x\in [-L,-l]$  and $y\in [x+M,L]$}, \label{eq:cond4}
\end{eqnarray}
where
\begin{eqnarray}
	\delta = \frac{K_L^2 (L-l)}{2 (2L+1)^3} - b \label{def2:delta}
\end{eqnarray}
is a positive constant if $(\mu,b) \in R_1$.
The pairs $(x,y)$, with $x\in [-L,-l]$ and $y\in [x+M,L]$, are exactly those
that contribute weight $1$ to
the calculation of the mean fitness $\overline m_{\hat\pi}$ as defined
in~(\ref{def:mbar}), while the pairs $(x',y')$ with $x'\in [-L,-l]$
and $y'\in [-L,x+M-1]$ contribute only weight $b$.
Condition~(\ref{eq:cond4}) implies that
\begin{eqnarray}
        \overline m_{\hat\pi}
	& = & \sum_{|x-z|\geq M} K_x K_z{\hat\pi}_x {\hat\pi}_z
                + b \sum_{|x-z|< M} K_x K_z{\hat\pi}_x {\hat\pi}_z \nnb \\
        & \leq & K_0^2 \sum_{|x-z|\geq M} \delta + b \sum_{|x-z|< M} K_0^2
                \nonumber \\
        & \leq & 2 (L-l+1)^2 \delta + 2 b (2L+1)^2 \nnb \\
        & \leq & 2 (2L+1)^2 (\delta+b). \label{upper:mpi}
\end{eqnarray}
Recall from~(\ref{ineq:mass_mid}) that
$\hat\pi([-l+1,l-1]) \leq \frac{2l+1}{2L+1}$, which is equivalent to
$\hat\pi([-L,-l]\cup [l,L]) \geq \frac{2(L-l)}{2L+1}$. Therefore
either $\hat\pi([-L,-l]) \geq \frac{L-l}{2L+1}$ or
$\hat\pi([l,L]) \geq \frac{L-l}{2L+1}$.
Suppose $\hat\pi([l,L]) \geq \frac{L-l}{2L+1}$, then we can bound
the fitness of site $-L$:
\begin{eqnarray}
        \hat m_{-L} & = & K_{-L} \sum_{z=-L}^L B_{-L-z} K_z \hat\pi_z
	\geq K_{-L} \sum_{z=l}^{L} K_z \hat\pi_z \nnb \\
        & \geq & K_L^2 \sum_{z=l}^{L} \hat\pi_z
        = K_L^2 \hat\pi([l,L])
        \geq K_L^2 \frac{L-l}{2L+1}. \label{ineq:fit_-L}
\end{eqnarray}
Since $\delta+b = \frac{K_L^2 (L-l)}{2 (2L+1)^3}$, (\ref{upper:mpi})
and~(\ref{ineq:fit_-L}) together imply that
\begin{eqnarray*}
        \hat m_{-L} - \overline m_{\hat\pi}
        \geq K_L^2 \frac{L-l}{2L+1} - 2 (2L+1)^2 (\delta+b)
        = 0.
\end{eqnarray*}
As a result of the inequality above, condition~(\ref{eq:cond3}) implies that
\begin{eqnarray*}
        \frac{1}{\hat\pi_{-L}} - (2L+1) \leq 0,
\end{eqnarray*}
hence
\begin{eqnarray}
	\hat\pi_{-L} \geq \frac{1}{2L+1}. \label{ineq:pi_-L}
\end{eqnarray}
Combining the two bounds~(\ref{ineq:fit_-L}) and~(\ref{ineq:pi_-L})
on $\hat\pi_{-L}$ and $\hat m_{-L}$, we have
\begin{eqnarray}
        \overline m_{\hat\pi} = \sum_x \hat\pi_x \hat m_x
        \geq \hat\pi_{-L} \hat m_{-L}
        \geq K_L^2 \frac{L-l}{(2L+1)^2}. \label{ineq:mpi1}
\end{eqnarray}
By similar reasoning, if $\hat\pi([-L,-l]) \geq \frac{L-l}{2L+1}$,
then we have
$\hat m_L \geq K_L^2 \frac{L-l}{2L+1}$ and $\hat\pi_L > \frac{1}{2L+1}$,
which also implies that
$\overline m_{\hat\pi} \geq K_L^2 \frac{L-l}{(2L+1)^2}$. Therefore
condition~(\ref{eq:cond4}) implies that
$\overline m_{\hat\pi} \geq K_L^2 \frac{L-l}{(2L+1)^2}$.
This would actually contradict~(\ref{upper:mpi}) for sufficiently small
$\delta$ and $b$, but it does not matter since in that case,
it just says that condition~(\ref{eq:cond4}) is impossible and the
analysis below applies.
                                                                                
On the other hand, if condition~(\ref{eq:cond4}) is not satisfied,
i.e. there is at least one pair of phenotypes, say $(\tilde x,\tilde y)$,
with $\tilde x \in [-L,-l]$ and $\tilde y \in [\tilde x+M,L]$, such that
$\hat\pi_{\tilde x} \hat\pi_{\tilde y} > \delta$. Then it is easy to see
that
\begin{eqnarray}
        \overline m_{\hat\pi} \geq \sum_{|x-z|\geq M} K_x K_z\pi_x \pi_z
        \geq K_L^2 \hat\pi_{\tilde x} \hat\pi_{\tilde y}
        \geq K_L^2 \delta
        \geq \frac{K_L^4 (L-l)}{4(2L+1)^3}, \label{ineq:mpi2}
\end{eqnarray}
the last inequality due to~(\ref{def2:delta}) and the requirement
$b\leq \frac{K_L^2(L-l)}{4(2L+1)^3}$.
Therefore combining~(\ref{ineq:mpi1}) and~(\ref{ineq:mpi2}), we conclude that
\begin{eqnarray*}
        \overline m_{\hat\pi} \geq \min\left(
                \frac{K_L^2 (L-l)}{(2L+1)^2},
                \frac{K_L^4 (L-l)}{4 (2L+1)^3}\right) > 0.
\end{eqnarray*}
This bound is uniform for any positive $(\mu,b)$ satisfying the condition
$b\leq 2 K_L^2 \mu (L-l)$ and $b\leq \frac{K_L^2 (L-l)}{4 (2L+1)^3}$
and the conclusion of the lemma follows.
\qed

\vspace{.3cm}

\begin{THM}
\label{thm:largemu}
If $K$ is symmetric and unimodal with $K_0=1$, and
$B_x = b + (1-b) 1_{\{|x| \geq M\}}$ with $L < M \leq 2L$,
then with the constant $c_1 = c_1(K,L,l)$ defined in Lemma~\ref{lem:largemu},
we have $\hat\pi^\mu_x \leq \frac{2\mu}{c_1}$
for $x\in [-l+1,l-1]$ and $(\mu,b)$ lying in
\begin{eqnarray*}
        R = \left\{(\mu,b):
        0 \leq b \leq \min\left(4 \mu K_L^2 (L-l),
                \frac{c_1}{2},
                \frac{K_L^2 (L-l)}{4 (2L+1)^3}\right) \right\},
\end{eqnarray*}
\end{THM}
\proof
Recall from~(\ref{eq:max_fit}) that
for $x\in [-l+1,l-1]$, the maximum fitness is attained at $x=0$, and
$m_0(\pi) \leq b$.
If $b \leq \frac{c_1}{2}$, with $c_1$ from Lemma~\ref{lem:largemu},
then $m_0(\pi) \leq \frac{c_1}{2}$.
Hence for $\mu$ and $b$ lying in $R \subset R_1$, where $R_1$ is defined in
Lemma~\ref{lem:largemu}, Lemma~\ref{lem:largemu} implies that
\begin{eqnarray}
	\overline m_{\hat\pi^\mu} - \hat m^\mu_0 \geq \frac{c_1}{2}.
	\label{ineq:fit_diff}
\end{eqnarray}
Condition~(\ref{eq:cond3})
and~(\ref{ineq:fit_diff}) together imply that for $x\in [-l+1,l-1]$
and $(\mu,b) \in R$,
\begin{eqnarray*}
        \mu \left(\frac{1}{\hat\pi^\mu_x} - (2L+1)\right)
        = \overline m_{\hat\pi^\mu} - \hat m^\mu_x \geq \frac{c_1}{2}.
\end{eqnarray*}
Therefore for $x\in [-l+1,l-1]$ and $(\mu,b) \in R$,
\begin{eqnarray*}
        \hat\pi^\mu_x \leq \frac{1}{\frac{c_1}{2\mu} + 2L+1}
        \leq \frac{2\mu}{c_1},
\end{eqnarray*}
and the proof is complete.
\qed

\vspace{.3cm}

\begin{REM}
If $b=0$, then Theorem~\ref{thm:largemu} works for any $\mu$, no matter how
small.
\end{REM}

\begin{REM}
The constant $c_1$ in Lemma~\ref{lem:largemu} is small.
\end{REM}

\subsection{Small $\mu$: Existence of $\delta$-like Stationary Measure}
\label{sec:sm_mu1}
Results from Chapter~\ref{sec:one-dim} for one-dimensional systems
indicate that when $b$ and $\mu$ are sufficiently small, there
should be at least two stationary distributions, one resembling the
$\delta$-measure, the other bimodal and having little mass in the middle.
In this section, we show the existence of $\delta$-like stationary
distributions. More specifically, we establish the following:

\begin{PRP}
\label{prop:sm_mu1}
Define $k=\min_x |K_x - K_{x-1}|$. If $\mu < \frac{b k}{8} \epsilon_1$,
then the set $A_1$ is an invariant set for the
dynamical system~(\ref{eq:ODE3}), where we define
\begin{eqnarray*}
	A_1=\{\pi \in \Delta:
		\pi_x \leq \epsilon_1 \mbox{ for all } x \neq 0\},
\end{eqnarray*}
with $\epsilon_1 \leq \min(\frac{1}{4L},\frac{b k}{2 K_l^2 (L-l+1)},
\frac{k}{16 L K_1})$.
\end{PRP}
\proof Let $\mu < \frac{b k}{8} \epsilon_1$ be fixed.
For $\pi \in A_1$, since $\pi_x < \epsilon_1 < \frac{1}{4L}$, we have
$\pi_0 \geq 1-2L\epsilon_1 > \frac{1}{2}$.
Recall that $l=M-L$ and that for $x\in [-l+1,l-1]$,
\begin{eqnarray*}
        m_x = b K_x \sum_z K_z \pi_z.
\end{eqnarray*}
Since $K_x$ is increasing in $[-L,0]$ and decreasing in $[0,L]$ with $K_0=1$,
sites in $[0,l-1]$ have decreasing fitness,
and sites in $[-l+1,0]$ have increasing fitness, and
we also have the following estimates:
\begin{eqnarray}
        m_0 - m_1 & = & b (K_0 - K_1) \sum_z K_z \pi_z
        \geq b k K_0 \pi_0 \geq \frac{b k}{2} \ \ \ \
		\mbox{ since $\pi_0 \geq \frac{1}{2}$}, \label{ineq:2est1}\\
        m_1 & \leq & b K_1 K_0 \leq b K_1. \label{ineq:2est2}
\end{eqnarray}
For $\pi \in A_1$ and $x \in [l,L]$,
\begin{eqnarray}
        m_l - m_x & = & K_l b \sum_z K_z \pi_z
                - K_x \left[ b \sum_z K_z \pi_z
                + (1-b) \sum_{z=-L}^{x-M} K_z \pi_z \right] \nnb \\
        & = & (K_l - K_x) b \sum_z K_z \pi_z
                - K_x (1-b) \sum_{z=-L}^{x-M} K_z \pi_z \nnb \\
        & \geq & k b K_0 \pi_0 - K_l \sum_{z=-L}^{-l} K_z \pi_z \nnb \\
        & \geq & \frac{k b}{2} - K_l^2 (L-l+1) \epsilon_1 \ \ \ \ \ \mbox{
		by~(\ref{ineq:2est1})} \nnb \\
	& \geq & 0, \label{ineq:fit_l_x}
\end{eqnarray}
since we assume $\epsilon_1 \leq \frac{k b}{2 K_l^2 (L-l+1)}$.
By the same calculation as in~(\ref{ineq:fit_l_x}), $m_{-l} - m_x \geq 0$ for
$x\in [-L,-l]$ as well.
Therefore sites in $[-L,-l]$ are less fit than site $-l$, and sites in
$[l,L]$ are less fit than site $l$; furthermore, sites in $[-l,0]$ have
increasing fitness, while sites in $[0,l]$ have decreasing fitness.
In particular, among sites in $[-L,-1]\cup [1,L]$, sites $1$ and $-1$ have
maximum fitness.
Any measure $\pi \in A_1$ looks like the $\delta$-measure.
For such measures, the mean fitness $\overline m_\pi$ is close to but less
than $m_0$. We estimate the difference between $\overline m_\pi$ and $m_1$:
\begin{eqnarray*}
        \overline m_\pi - m_1 = \sum_x m_x \pi_x - m_1
	\geq m_0 \pi_0 - m_1
	\geq m_0 (1-2L\epsilon_1) - m_1,
\end{eqnarray*}
since $\pi_0 \geq 1 - 2L\epsilon_1$ from the beginning of the proof.
Now using estimates~(\ref{ineq:2est1}) and~(\ref{ineq:2est2}), we
continue estimating $\overline m_\pi - m_1$: 
\begin{eqnarray*}
        m_0 (1-2L\epsilon_1) - m_1
	= (1-2L\epsilon_1) (m_0 - m_1) - 2L\epsilon_1 m_1
        \geq (1-2L\epsilon_1)\frac{b k}{2} - 2L\epsilon_1 b K_1.
\end{eqnarray*}
Recalling from the beginning of the proof that $1-2L\epsilon_1 > \frac{1}{2}$,
we use the assumption $\epsilon_1 \leq \frac{k}{16 L K_1}$ to estimate
the right hand side of the above inequality:
\begin{eqnarray*}
        (1-2L\epsilon_1)\frac{b k}{2} - 2L\epsilon_1 b K_1
	\geq \frac{b k}{4} - 2L b K_1 \frac{k}{16 L K_1}
	\geq \frac{b k}{8}.
\end{eqnarray*}
Therefore
\begin{eqnarray}
	\overline m_\pi - m_1 \geq \frac{b k}{8}. \label{ineq:diff_fit}
\end{eqnarray}
Since among sites in $[-L,-1]\cup [1,L]$, sites $1$ and $-1$ have
maximum fitness,~(\ref{ineq:diff_fit}) implies that for all $x\neq 0$,
\begin{eqnarray*}
        m_x - \overline m_\pi \leq - \frac{b k}{8}.
\end{eqnarray*}
                                                                                
We write $\partial A_1 = B_1  \cup B_2 \cup B_3$, where
\begin{eqnarray*}
	B_1 & = & \{\pi \in \partial A_1: \pi_x = 0 \ \mbox{ for some $x$
		and $\pi_x \neq \epsilon_1$ for all $x$} \}, \\
	B_2 & = & \{\pi \in \partial A_1: \pi_x = \epsilon_1 \
		\mbox{ for some $x$ and $\pi_x \neq 0$ for all $x$} \}, \\
	B_3 & = & \{\pi \in \partial A_1: \pi_x = \epsilon_1 \
		\mbox{ for some $x$ and $\pi_y = 0$ for some $y$} \}.
\end{eqnarray*}
For $\pi \in B_1$, we have shown following~(\ref{ineq:muneq0})
that $\partial_t \pi_x > 0$ at $x$ where $\pi_x = 0$. Therefore
$\partial_t \pi$ points toward the interior of $\Delta$.
For $\pi \in B_2$, we have for $x$ where $\pi_x = \epsilon_1$,
\begin{eqnarray}
        \partial_t \pi_x & = & \pi_x (m_x - \overline m_\pi)
                + \mu(1-(2L+1)\pi_x) \nnb \\
        & \leq & - \epsilon_1 \frac{b k}{8} + \mu \nnb \\
        & < & 0, \label{ineq:muneqe}
\end{eqnarray}
since $\mu < \frac{b k}{8} \epsilon_1$. Thus $\partial_t \pi$
also points toward the interior of $\Delta$. For $\pi \in B_3$,
we can apply~(\ref{ineq:muneq0}) to sites $x$ where $\pi_x = 0$
and~(\ref{ineq:muneqe}) to sites $y$ where $\pi_y = \epsilon_1$,
and conclude that $\partial_t \pi$ also points
toward the interior of $\Delta$. Therefore the set $A_1$ is invariant
for the dynamical system~(\ref{eq:ODE3}), as required.
\qed

\subsection{Small $\mu$: Existence of Bimodal Stationary Measure}
\label{sec:sm_mu2}
Results from Chapter~\ref{sec:sm_mu1} show that the set $A_1$, members of
which resemble the $\delta$-measure, is invariant for the dynamical
system~(\ref{eq:ODE3}). In this section, we assume that
$\frac{b}{1-b} \leq \frac{K_{M/2}^2}{8}$ and 
show that there is another
invariant set
\begin{eqnarray}
	A_2 & = & \left\{\pi: \pi_x \leq \epsilon_2 
		\mbox{ for all } x \notin \{p,-p\},
		\mbox{ and $| \log \frac{\pi_p}{\pi_{-p}} |
		\leq \epsilon_3$} \right\}, \label{def:A2}
\end{eqnarray}
where $p=M/2$ and $L < M \leq 2L$,
\begin{eqnarray}
        c_2 & = & \frac{1}{2} \min\left( (1-b) \frac{K_p (K_p-K_{p+1})}{8},
                b,\frac{K_p^2}{16}\right), \label{def2:c2} \\
	\epsilon_2 & = & \min\left(\frac{1}{8(2L+1)},
                \frac{K_p (K_p-K_{p+1})}{8 K_{p+1}^2 (L-l+1)},
                \frac{K_p^2}{16 K_l (L-l+1)}, \right. \nnb \\
	& & \ \ \ \ \ \ \ \ \ \ \left.
                \frac{c_2}{16 (2L+1) K_p}\right), \label{def:e2} \\
	\epsilon_3 & = & \min\left(\log 2,
		\log\left(1+\frac{c_2}{4 K_p^2}\right)\right). \label{def:e3}
\end{eqnarray}
Apparently, $M$ must be even; this is such that
$B_{p-(-p)} = 1$ but $B_{(p-1)-(-p)}=b$. Notice that since $M \leq 2L$,
$l = M-L \leq M - \frac{M}{2} = \frac{M}{2} = p$,
and also $p = \frac{M}{2} \leq L$. Thus $l \leq p \leq L$.
Members of $A_2$ are bimodal distributions with very sharp peaks at sites
$p$ and $-p$.

For $\pi \in A_2$, $\pi_p + \pi_{-p} = 1 - \sum_{x \neq p,-p} \pi_x
\geq 1-(2L-1) \frac{1}{8(2L+1)}$ by the condition
$\pi_x \leq \frac{1}{8(2L+1)}$ for $x \notin \{p,-p\}$, therefore
$\pi_p + \pi_{-p} \geq \frac{7}{8}$.
Now since $| \log \frac{\pi_p}{\pi_{-p}} | \leq \log 2$, we have
$\frac{1}{2} \leq \frac{\pi_p}{\pi_{-p}} \leq 2$.
This means that
\begin{eqnarray}
	\min(\pi_p,\pi_{-p}) > \frac{1}{4}, \label{ineq:pi_p}
\end{eqnarray}
for otherwise, say
$\pi_p \leq \frac{1}{4}$, then $\pi_{-p} > \frac{5}{8}$, which means
that $\frac{\pi_p}{\pi_{-p}} < \frac{2}{5} < \frac{1}{2}$.

We use the same idea that we used to establish the invariance property of
$A_1$ to show that $A_2$ is also an invariant set for the dynamical system,
for $\mu < \frac{3 c_2}{8} \epsilon_2$.

\begin{LEM}
For any $\pi \in A_2$, the following estimate holds:
\begin{eqnarray}
        \min(m_p,m_{-p}) - m_x > \frac{c_2}{2}
        \mbox{ for $x\in [-L,L]\backslash \{p,-p\}$},
                \label{ineq:diff_fit4}
\end{eqnarray}
where $A_2$ and $c_2$ are defined in~(\ref{def:A2}) and ~(\ref{def2:c2}),
respectively.
\end{LEM}
\proof
As in the proof of Proposition~\ref{prop:sm_mu1}, we first establish
a few bounds on fitness of various sites for $\pi\in A_2$.
For $x=p$,
\begin{eqnarray}
        m_p & = & K_p \left[ b \sum_z K_z \pi_z
                + (1-b) \sum_{z=-L}^{-p} K_z \pi_z \right] \nnb \\
        & \geq & K_p (1-b) K_{-p} \pi_{-p} \geq \frac{(1-b) K_p^2}{4},
	\label{ineq:fit_p}
\end{eqnarray}
where we use~(\ref{ineq:pi_p}) in the last inequality.
For $x\in [p+1,L]$,
\begin{eqnarray*}
        \lefteqn{ m_p - m_x } \nnb \\
	& = & K_p \left[ b \sum_z K_z \pi_z
                + (1-b) \sum_{z=-L}^{-p} K_z \pi_z \right]
                - K_x \left[ b \sum_z K_z \pi_z
                + (1-b) \sum_{z=-L}^{x-M} K_z \pi_z \right] \nonumber \\
        & = & (K_p - K_x) b \sum_z K_z \pi_z
                + (K_p - K_x) (1-b) \sum_{z=-L}^{-p} K_z \pi_z
                - K_x (1-b) \sum_{z=-p+1}^{x-M} K_z \pi_z.
\end{eqnarray*}
{
\psfrag{big}{$> \frac{1}{4}$}
\psfrag{small}{$\leq \epsilon_2$}
\psfrag{L}{$L$}
\psfrag{-L}{$-L$}
\psfrag{p}{$p$}                                     
\psfrag{-p}{$-p$}        
\psfrag{l}{$l$}
\psfrag{-l}{$-l$}
\begin{figure}[h!]
\centering
\includegraphics[height = 1.5in, width = 5.5in]{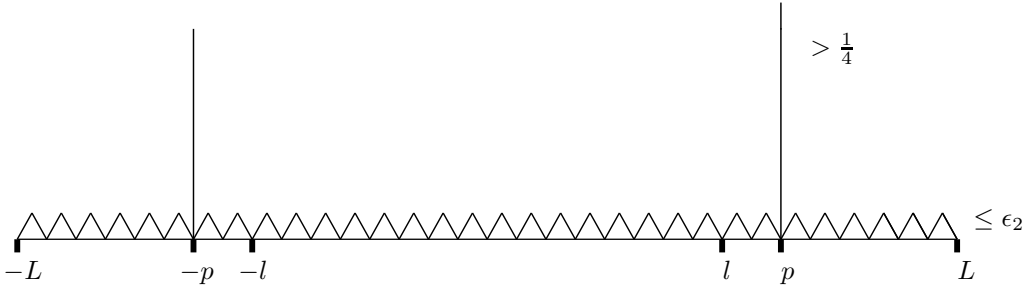}
\caption{Illustration of $A_2$}
\end{figure}
}
Since $x\in [p+1,L]$ and $K_x$ is decreasing in $[0,L]$, we have
$K_p - K_x \geq K_p - K_{p+1} > 0$. We apply these facts, along
with~(\ref{ineq:pi_p}) and
the requirements on $\pi_x$ for $\pi \in A_2$, to the right hand
side and obtain
\begin{eqnarray}
        m_p - m_x & \geq & 0 + (K_p-K_{p+1}) (1-b) K_{-p} \frac{1}{4}
                - K_{p+1} (1-b) (L-l+1) K_{-p+1} \epsilon_2 \nonumber \\
        & = & (1-b) \left[ \frac{K_p (K_p-K_{p+1})}{4}
                - K_{p+1} K_{p-1} (L-l+1) \epsilon_2 \right] \nonumber \\
        & \geq & (1-b) \frac{K_p (K_p-K_{p+1})}{8}
        \label{ineq:tinye1}
\end{eqnarray}
for $x\in [p+1,L]$ since
$\epsilon_2 \leq \frac{K_p (K_p-K_{p+1})}{8 K_{p+1} K_{p-1} (L-l+1)}$.
For $x\in [0,l-1]$,
\begin{eqnarray*}
        m_x = b K_x \sum_z K_z \pi_z \leq b K_x K_0 \leq b,
\end{eqnarray*}
therefore using~(\ref{ineq:fit_p}) and the inequality above, we have
\begin{eqnarray}
        m_p - m_x \geq \frac{(1-b) K_p^2}{4} - b \geq b
        \label{ineq:tinye2}
\end{eqnarray}
since $\frac{b}{1-b} \leq \frac{K_p^2}{8}$.
For $x\in [l,p-1]$,
\begin{eqnarray*}
        m_x & = & K_x \left[ b \sum_z K_z \pi_z
                + (1-b) \sum_{z=-L}^{x-M} K_z \pi_z \right] \\
	& \leq & K_l (b K_0 \sum_z \pi_z + K_0 \sum_{z=-L}^{-l} \epsilon_2) \\
        & \leq & K_l (b + (L-l+1) \epsilon_2),
\end{eqnarray*}
therefore using~(\ref{ineq:fit_p}) and the inequality above, we have
\begin{eqnarray*}
        m_p - m_x & \geq & \frac{(1-b)K_p^2}{4} - K_l(b+(L-l+1)\epsilon_2) \\
	& = & \frac{K_p^2}{4} - b \left(\frac{K_p^2}{4}+K_l\right)
                - K_l (L-l+1) \epsilon_2.
\end{eqnarray*}
Since $b \leq \frac{K_p^2}{2K_p^2 + K_l K_0}$, we have
$\frac{K_p^2}{4} - b \left(\frac{K_p^2}{4}+K_l\right) \geq \frac{K_p^2}{8}$.
Furthermore, since $\epsilon_2 \leq \frac{K_p^2}{16 K_l (L-l+1)}$, we have
\begin{eqnarray}
	m_p - m_x \geq \frac{K_p^2}{16}. \label{ineq:tinye3}
\end{eqnarray}
The estimates~(\ref{ineq:tinye1}),~(\ref{ineq:tinye2}),
and~(\ref{ineq:tinye3}) compare $m_p$ with $m_x$ for $x\in [0,L]$.
Similar calculations comparing $m_{-p}$ with $m_x$ for $x\in [-L,0]$
yield similar results. Then recalling the definition of
$c_2$ in~(\ref{def2:c2}), we have
\begin{eqnarray}
	& & m_p - m_x > c_2 \mbox{ for $x\in [0,L]\backslash \{p\}$}, \nnb \\
	& \mbox{ and } & m_{-p} - m_x > c_2
		\mbox{ for $x\in [-L,0]\backslash \{-p\}$}.
		\label{ineq:diff_fit2}
\end{eqnarray}

To establish the lemma, it suffices to compare $m_p$ and $m_{-p}$:
\begin{eqnarray*}
	| m_p - m_{-p} | & = & \left| K_p \left[ b \sum_z K_z \pi_z
                + (1-b) \sum_{z=-L}^{-p} K_z \pi_z \right] \right. \nnb \\
	& & \ \ \ \ \	- \left. K_{-p} \left[ b \sum_z K_z \pi_z
                + (1-b) \sum_{z=p}^{L} K_z \pi_z \right] \right| \\
	& = & K_p (1-b) \left| \sum_{z=p}^{L} K_z (\pi_{-z} - \pi_z) \right| \\
	& \leq & K_p (1-b) \sum_{z=p}^{L} K_z | \pi_{-z} - \pi_z |.
\end{eqnarray*}
Since $K_z \leq K_p$ for $z \in [p,L]$, $K_z \leq K_{-p} = K_p$ for
$z \in [-L,-p]$, and $\pi_z \leq \epsilon_2$
for $z\neq -p,p$, we have
\begin{eqnarray}
	| m_p - m_{-p} | & \leq & K_p^2 (1-b)
		( | \pi_{-p} - \pi_p | + (L-p) 2\epsilon_2) \nnb \\
	& \leq & K_p^2 | \pi_{-p} - \pi_p | + 2 K_p^2 (L-p) \epsilon_2.
	\label{ineq:m_p-m_-p}
\end{eqnarray}
We treat the two terms in the above sum separately. We first deal with
the second term $2 K_p^2 (L-p) \epsilon_2$. Since
$L-p < 2L+1$ and $K_p \geq K_p^2$, we have
$\frac{c_2}{16(2L+1)K_p} \leq \frac{c_2}{8(L-p)K_p^2}$. Then the definition
of $\epsilon_2$ in~(\ref{def:e2}) means that
$\epsilon_2 \leq \frac{c_2}{16(2L+1)K_p} \leq \frac{c_2}{8(L-p)K_p^2}$,
which implies
\begin{eqnarray}
	2 K_p^2 (L-p) \epsilon_2 \leq \frac{c_2}{4}. \label{ineq:m_p-m_-p1}
\end{eqnarray}
For the first term in~(\ref{ineq:m_p-m_-p}), we can divide into two
cases: $\pi_p \geq \pi_{-p}$ and $\pi_p \leq \pi_{-p}$.
If $\pi_p \geq \pi_{-p}$, then
\begin{eqnarray*}
	0 \leq \pi_p - \pi_{-p}
	= \pi_p - e^{\epsilon_3} \pi_{-p} + e^{\epsilon_3} \pi_{-p} - \pi_{-p}
	= (\pi_p - e^{\epsilon_3} \pi_{-p}) + (e^{\epsilon_3}-1) \pi_{-p}.
\end{eqnarray*}
The definition of $\epsilon_3$ in~(\ref{def:e3}) implies that
$\pi_p \leq e^{\epsilon_3} \pi_{-p}$ and
$e^{\epsilon_3}-1 \leq \frac{c_2}{4 K_p^2}$, so continuing the calculation
in the line above, we obtain
\begin{eqnarray}
	|\pi_p - \pi_{-p}| = \pi_p - \pi_{-p}
	\leq 0 + \frac{c_2}{4 K_p^2} \pi_{-p}
	\leq \frac{c_2}{4 K_p^2}. \label{ineq:m_p-m_-p2}
\end{eqnarray}
If $\pi_p \leq \pi_{-p}$, we get the same bound. Therefore
applying~(\ref{ineq:m_p-m_-p1}) and~(\ref{ineq:m_p-m_-p2})
to~(\ref{ineq:m_p-m_-p}), we have
\begin{eqnarray*}
        | m_p - m_{-p} | \leq \frac{c_2}{4} + \frac{c_2}{4} = \frac{c_2}{2}.
\end{eqnarray*}
This result means that the estimate in~(\ref{ineq:diff_fit2}) can
be generalized to
\begin{eqnarray*}
        \min(m_p,m_{-p}) - m_x > \frac{c_2}{2}
	\mbox{ for $x\in [-L,L]\backslash \{p,-p\}$},
\end{eqnarray*}
as required by the lemma. \qed

\vspace{.3cm}

\begin{PRP}
\label{prop:sm_mu2}
If $\mu < \frac{3 c_2}{8} \epsilon_2$,
then the set $A_2$ defined in~(\ref{def:A2}) is an invariant set for the
dynamical system~(\ref{eq:ODE3}).
\end{PRP}
\proof As in the proof of Proposition~\ref{prop:sm_mu1}, it suffices to
observe that $\partial_t \pi_x >0$ where $\pi_x = 0$ as shown by the argument
following~(\ref{ineq:muneqe}), and
check that the following inequalities hold:
\begin{eqnarray}
	& & \left.\left(\partial_t \pi_x\right)\right|_{
		\pi_x = \epsilon_2,\pi \in A_2} < 0 \mbox{ for $x\neq p,-p$},
		\label{ineq:smmu1} \\
	& & \left.\left(\partial_t \log \frac{\pi_p}{\pi_{-p}}\right)\right|_{
		\log \frac{\pi_p}{\pi_{-p}} = \epsilon_3,\pi \in A_2} < 0,
		\label{ineq:smmu2} \\
	& \mbox{and} & \left.\left(\partial_t
		\log \frac{\pi_{-p}}{\pi_p}\right)\right|_{
                \log \frac{\pi_{-p}}{\pi_p} = \epsilon_3,\pi \in A_2} < 0.
		\label{ineq:smmu3}
\end{eqnarray}
The estimate~(\ref{ineq:diff_fit4}) means that for $\pi \in A_2$
and $x \notin \{p,-p\}$, $m_x$ is significantly smaller than
$\min(m_p,m_{-p})$. But the mean fitness $\overline m_\pi$ cannot be
much smaller than $\min(m_p,m_{-p})$:
\begin{eqnarray*}
        \min(m_p,m_{-p}) - \overline m_\pi
	& = & \min(m_p,m_{-p}) - \sum_x m_x \pi_x \\
        & \leq & \min(m_p,m_{-p}) - m_p \pi_p - m_{-p} \pi_{-p} \\
	& \leq & \min(m_p,m_{-p}) - \min(m_p,m_{-p}) (\pi_p + \pi_{-p}) \\
        & = & (1-\pi_p-\pi_{-p}) \min(m_p,m_{-p}) \\
	& \leq & (2L-1)\epsilon_2 \min(m_p,m_{-p})
\end{eqnarray*}
by the requirement that $\pi_x \leq \epsilon_2$ for $x\neq p,-p$
and $\pi\in A_2$. Since \\ $m_p = K_p \sum_z B_{p-z} K_z \pi_z \leq K_p K_0$
and similarly $m_{-p} \leq K_p K_0$, we have
\begin{eqnarray*}
        \min(m_p,m_{-p}) - \overline m_\pi \leq (2L-1) \epsilon_2 K_p K_0
	= (2L-1) \epsilon_2 K_p.
\end{eqnarray*}
We use the definition of $\epsilon_2$ in~(\ref{def:e2}) to obtain
$\epsilon_2 \leq \frac{c_2}{16(2L+1) K_p} < \frac{c_2}{8 (2L-1) K_p}$,
which, when applied to the estimate in the line above, implies that
\begin{eqnarray}
        \min(m_p,m_{-p}) - \overline m_\pi < \frac{c_2}{8}.
		\label{ineq:diff_fit3}
\end{eqnarray}
Now~(\ref{ineq:diff_fit4}) and~(\ref{ineq:diff_fit3}) imply that for
$\pi \in A_2$ and $x \notin \{p,-p\}$,
\begin{eqnarray*}
        m_x - \overline m_\pi < m_x + \frac{c_2}{8} - \min(m_p,m_{-p})
	< - \frac{c_2}{2} + \frac{c_2}{8} = - \frac{3 c_2}{8}.
\end{eqnarray*}
Thus for $\pi\in \partial A_2 \cap \{\pi_z = \epsilon_2 \ \mbox{ for some } \ 
z \neq \{p,-p\}\}$, we have for $x$ where $\pi_x = \epsilon_2$,
\begin{eqnarray*}
        \left.\left(\partial_t \pi_x\right)\right|_{\pi_x = \epsilon_2} & = &
		\pi_x (m_x - \overline m_\pi)
                + \mu(1-(2L+1)\pi_x)  |_{\pi_x = \epsilon_2}\\
        & \leq & - \epsilon_2 \frac{3 c_2}{8} + \mu \\
        & < & 0,
\end{eqnarray*}
since $\mu < \frac{3 c_2}{8} \epsilon_2$. This verifies~(\ref{ineq:smmu1}).

Now we deal with $\pi\in \partial A_2 \cap \{\log \frac{\pi_p}{\pi_{-p}}
= \epsilon_3 \}$. By~(\ref{eq:ODE3}), we have
\begin{eqnarray}
	\partial_t \log \frac{\pi_p}{\pi_{-p}} & = &
		(m_p - \overline m_\pi + \frac{\mu}{\pi_p} - \mu (2L+1)) -
		(m_{-p} - \overline m_\pi + \frac{\mu}{\pi_{-p}} - \mu (2L+1))
		\nnb \\
	& = & m_p - m_{-p} + \mu \left( \frac{1}{\pi_p} - \frac{1}{\pi_{-p}}
		\right). \label{eq:pi_p:pi_-p}
\end{eqnarray}
If $\frac{\pi_p}{\pi_{-p}} = e^{\epsilon_3} > 1$, then
$\pi_p > \pi_{-p}$, which means that
$\frac{1}{\pi_p} - \frac{1}{\pi_{-p}} < 0$. Therefore it remains
to check the sign of $m_p - m_{-p}$:
\begin{eqnarray}
        \lefteqn{ m_p - m_{-p} } \nnb \\
	& = & K_p \left[ b \sum_z K_z \pi_z
                + (1-b) \sum_{z=-L}^{-p} K_z \pi_z \right]
                - K_{-p} \left[ b \sum_z K_z \pi_z
                + (1-b) \sum_{z=p}^{L} K_z \pi_z \right] \nnb \\
        & = & K_p (1-b) \left( K_p (\pi_{-p} - \pi_p)
		+ \sum_{z=p+1}^{L} K_z (\pi_{-z} - \pi_z) \right).
		\label{ineq:m_p-m_-p5}
\end{eqnarray}
If $\frac{\pi_p}{\pi_{-p}} = e^{\epsilon_3}$, then
\begin{eqnarray}
	\pi_{-p} - \pi_p = \pi_{-p} (1-e^{\epsilon_3})
	= \pi_{-p} \max(-1,-\frac{c_2}{4 K_p^2})
	\leq -\frac{c_2}{16 K_p^2}
	\label{ineq:m_p-m_-p3}
\end{eqnarray} 
by the definition of $\epsilon_3$ in~(\ref{def:e3}) and the fact
$\pi_{-p} > \frac{1}{4}$ established in~(\ref{ineq:pi_p}).
For $\pi \in A_2$ and $z\geq p+1$,
$\pi_{-z}-\pi_z \leq 2 \epsilon_2$, therefore
\begin{eqnarray}
	\sum_{z=p+1}^{L} K_z (\pi_{-z} - \pi_z) \leq (L-p) K_p 2 \epsilon_2.
	\label{ineq:m_p-m_-p4}
\end{eqnarray}
Applying~(\ref{ineq:m_p-m_-p3}) and~(\ref{ineq:m_p-m_-p4})
to~(\ref{ineq:m_p-m_-p5}), and using
the requirement $\epsilon_2 \leq \frac{c_2}{16 (2L+1) K_p}$ in~(\ref{def:e2}),
which implies $\epsilon_2 \leq \frac{c_2}{32 K_p^2 (L-p)}$, we conclude that
$m_p - m_{-p} < 0$ if $\frac{\pi_p}{\pi_{-p}} = e^{\epsilon_3}$.
Therefore~(\ref{eq:pi_p:pi_-p}) implies that
\begin{eqnarray*}
        \left.\left(\partial_t \log \frac{\pi_p}{\pi_{-p}} \right)\right|_{
		\frac{\pi_p}{\pi_{-p}} = e^{\epsilon_3}} < 0,
\end{eqnarray*}
which verifies~(\ref{ineq:smmu2}). The verification of~(\ref{ineq:smmu3})
is similar, and the proof is complete. \qed

\vspace{.3cm}

Propositions~\ref{prop:sm_mu1} and~\ref{prop:sm_mu2} imply that if
$b$ is sufficiently small and $\mu < c_3 b^2$, where $c_3$ is a small
constant dependent on $K$ and $L$, then~(\ref{eq:ODE3}) has
at least two stationary distributions, one resembling a $\delta$-measure
and the other bimodal and having little mass in the middle. But
Theorem~\ref{thm:largemu} imply that if $b$ is sufficiently small
and $\mu > c_4 b$, where $c_4$ is a large constant dependent on $K$ and $L$,
then all stationary distributions are bimodal and have very little mass
in the middle. This phenomenon is illustrated in figure~\ref{fig:transition}.
We conjecture that there is a phase transition between a unique stationary
distribution and two stationary distributions in the behaviour
of~(\ref{eq:ODE3}) for small $\mu$ and $b$. But since we cannot come up
with a straightforward comparison argument in either $\mu$ or $b$, this
remains a conjecture.

{
\psfrag{t1}{2 stationary distributions}
\psfrag{t2}{(one with speciation,}
\psfrag{t3}{one without speciation)}
\psfrag{t4}{Unique stationary distribution}
\psfrag{t5}{(speciation)}
\psfrag{ee}{$\mu$}
\psfrag{bb}{$b$}
\begin{figure}[h!]
\centering
\includegraphics[height = 2.5in, width = 3in]{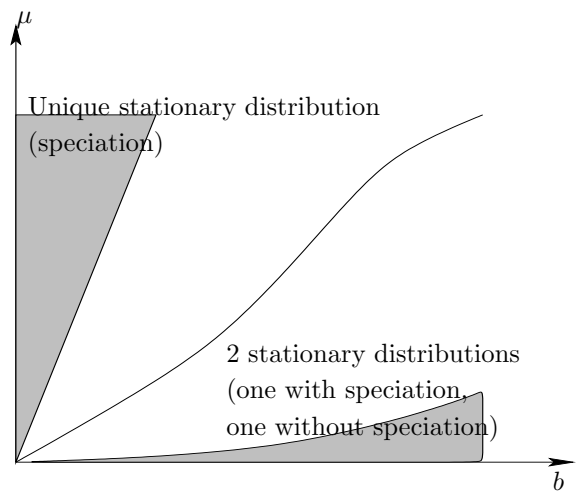}
\caption{Conjecture of phase transition in~(\ref{eq:ODE3}) for small
$\mu$ and $b$ (Shaded region is where Theorem~\ref{thm:largemu},
Propositions~\ref{prop:sm_mu1} and~\ref{prop:sm_mu2} work)}
\label{fig:transition}
\end{figure}
}

\chapter{Stationary Distributions}
In Chapter~\ref{ch:selmut}, we examined the large-time behaviour
of the deterministic dynamical system~(\ref{eq:ODE2}), i.e.
$\lim_{t\rightarrow\infty} \lim_{N\rightarrow\infty} \pi^N(t)$.
In this chapter, however, we will take the limit $t\rightarrow\infty$
first and examine the behaviour of $\nu^{\mu,N}$, the stationary distribution
of $\pi^N$; more specifically, we will do this by taking the limits
$N\rightarrow\infty$ then $\mu\rightarrow 0$ and examine
$\lim_{\mu\rightarrow 0}\lim_{N\rightarrow\infty} \nu^{\mu,N}$.
For this, we consider the case of symmetric, strictly positive,
and unimodal $K$, with $K_x$ strictly
decreasing for $x\in [0,L]$ and $K_0 = 1$, and $B_x = 1_{\{|x| \geq M\} }$
with $L < M \leq 2L$. Define $l=-L+M$, then $m_x = 0$ for $x\in [-l+1,l-1]$.
Define $c_3 = \sup_{\pi} \overline m_\pi$. For the strong selection model
with $N$ particles and mutation rate $\mu$ described in
Chapter~\ref{ch:st-model}, we observe that
this continuous-time finite-state Markov process has the property that
all states communicate, and therefore it has a unique stationary
distribution $\nu^{\mu,N}$ [Durrett 1991]. Let $(\Gcal^N,\Dcal(\Gcal^N))$
denote the generator associated with this Markov process, then for all
$F\in C^\infty(\Pcal(E)) \subset \Dcal(\Gcal^N)$, we have
\begin{eqnarray}
	\int \Gcal^N F(\pi) \nu^{\mu,N}(d\pi) = 0. \label{cond:GN:stat}
\end{eqnarray}

Let $(\Gcal,\Dcal(\Gcal))$ denote the generator associated with the
deterministic process described by the ODE~(\ref{eq:ODE}):
\begin{eqnarray*}
        \partial_t \pi_x = \pi_x ( m(x,\pi) - \overline m_{\pi} )
                + \mu ( 1 - (2L+1) \pi_x).
\end{eqnarray*}
We calculate the effect of the generator ${\cal G}^N$ on a $C^\infty$-function
\\
$F(\pi^N) = F(\pi^N_{-L},\ldots,\pi^N_0,\ldots,\pi^N_L)$:
\begin{eqnarray*}
        \Gcal^N F(\pi^N) & = & \sum_x \sum_{y\neq x} \left[
                F\left(\pi^N - \frac{\delta_x}{N} + \frac{\delta_y}{N}\right)
                - F(\pi^N) \right] \left(N \pi^N_x m(y,\pi^N) \pi^N_y
                + N \mu \pi^N_x\right) \nnb \\
        & = & N \sum_x \sum_{y\neq x} \left[
                F\left( \pi^N_{-L},\ldots,\pi^N_x-\frac{1}{N},\ldots,
                        \pi^N_y+\frac{1}{N},\ldots,\pi^N_L \right)
                - F(\pi^N) \right] \nnb \\
	& & \ \ \ \ \ \times \pi^N_x \left(m(y,\pi^N) \pi^N_y
                + \mu \right). \nnb
\end{eqnarray*}
Performing a Taylor expansion on $F$, we continue the above computation:
\begin{eqnarray}
	\lefteqn{ \Gcal^N F(\pi^N) } \nnb \\
        & = & N \sum_x \sum_{y\neq x} \left[
                - \frac{1}{N} \frac{\partial F}{\partial \pi^N_x}(\pi^N)
                + \frac{1}{N} \frac{\partial F}{\partial \pi^N_y}(\pi^N)
                + O(N^{-2}) \right] \pi^N_x \left( m(y,\pi^N) \pi^N_y
                + \mu \right) \nnb \\
        & = & \sum_x \sum_{y\neq x} \left[
                \frac{\partial F}{\partial \pi^N_y}(\pi^N)
                - \frac{\partial F}{\partial \pi^N_x}(\pi^N) \right]
                \pi^N_x \left( m(y,\pi^N) \pi^N_y + \mu \right)
                + R_1(N,F) (\pi^N). \label{eq:GN}
\end{eqnarray}
Since $\sup_{x,y,\pi} \pi^N_x (m(y,\pi^N)+\mu) O(N^{-2}) = O(N^{-2})$, we have
\begin{eqnarray*}
	R_1(N,F)(\pi) \leq  \frac{C \|F\|}{N},
\end{eqnarray*}
where we suppress the dependence on $L$. Therefore $R_1(N,F)(\pi)$ is
$O(N^{-1})$ uniformly for all $\pi \in {\cal P}(E)$.
Here we use the norm
\[\|F\| = \int \left( 1+\sum_{k=-L}^L \xi_k^2 \right) |\hat F(\xi)|^2 d\xi \]
associated with the Sobolev space $H^2(\Rbold^E)$ for $F$, where
$\hat F$ denotes the Fourier transform of $F$.
The first term in~(\ref{eq:GN}) is in fact equal to $\Gcal F(\pi^N)$ by
the following computation:
\begin{eqnarray*}
        \lefteqn{\sum_x \sum_{y\neq x} \left[
                \frac{\partial F}{\partial \pi_y}(\pi)
                - \frac{\partial F}{\partial \pi_x}(\pi) \right]
                \pi_x ( m(y,\pi) \pi_y + \mu )} \nnb \\
        & = & \sum_x \sum_y \left[
                \frac{\partial F}{\partial \pi_y}(\pi)
                - \frac{\partial F}{\partial \pi_x}(\pi) \right]
                \pi_x ( m(y,\pi) \pi_y + \mu ) \nnb \\
        & = & \sum_x \pi_x \sum_y \frac{\partial F}{\partial \pi_y}(\pi)
                        ( m(y,\pi) \pi_y + \mu )
                - \sum_y ( m(y,\pi) \pi_y + \mu ) \sum_x
                        \frac{\partial F}{\partial \pi_x}(\pi) \pi_x
                        \nnb \\                
        & = & \sum_x \frac{\partial F}{\partial \pi_x}(\pi)
                        ( m(x,\pi) \pi_x + \mu )
                - ( \overline m_\pi + (2L+1) \mu ) \sum_x \pi_x
                        \frac{\partial F}{\partial \pi_x}(\pi) \nonumber \\
        & = & \sum_x [\pi_x ( m(x,\pi) - \overline m_\pi)
                + \mu ( 1 - (2L+1) \pi_x)]     
                \frac{\partial F}{\partial \pi_x}(\pi) \\
	& = & \Gcal F(\pi). \nnb
\end{eqnarray*}

Therefore~(\ref{eq:GN}) can be written in the following much-simplified
form
\begin{eqnarray*}
	\Gcal^N F(\pi^N) = \Gcal F(\pi^N) + R_1(N,F) (\pi^N).	
\end{eqnarray*}
Then~(\ref{cond:GN:stat}) implies
\begin{eqnarray*}
	\left|\int {\cal G} F(\pi) \nu^{\mu,N}(d\pi)\right|
	= \left|\int R_1(N,F) (\pi) \nu^{\mu,N}(d\pi)\right|
	\leq \frac{C \|F\|}{N} \left|\int \nu^{\mu,N}(d\pi)\right|
	= \frac{C \|F\|}{N},
\end{eqnarray*}
i.e.,
\begin{eqnarray}
	\int {\cal G} F(\pi) \nu^{\mu,N}(d\pi) = O(N^{-1}). \label{eq:nuNe}
\end{eqnarray}
Since $E$ is compact, so is ${\cal P}(E)$ and $\Pcal(\Pcal(E))$,
therefore for each $\mu$, we can take a sequence
$N_k(\mu)$ such that $\nu^{\mu,N_k(\mu)}$
converges weakly to some $\nu^\mu \in \Pcal(\Pcal(E))$. By~(\ref{eq:nuNe}),
$\nu^\mu$ satisfies: for all $F\in C^\infty(\Pcal(E))$,
\begin{eqnarray}
	\int {\cal G} F(\pi) \nu^{\mu}(d\pi)
	& = & \int\sum_x [\pi_x ( m(x,\pi) - \overline m_\pi)
                + \mu ( 1 - (2L+1) \pi_x)]
                \frac{\partial F}{\partial \pi_x}(\pi) \nu^{\mu}(d\pi) \nnb \\
	& = & 0. \label{eq:nue}
\end{eqnarray}
Therefore $\nu^\mu$ is an stationary distribution for the
deterministic flow $\Gcal$.
Now we take a sequence $\mu_k \rightarrow 0$, such that $\nu^{\mu_k}$
converges weakly to some $\nu^0 \in \Pcal(\Pcal(E))$, and by~(\ref{eq:nue}) and
the following estimate: as $\mu\rightarrow 0$,
\begin{eqnarray*}
        \left| \int \sum_x \mu ( 1 - (2L+1) \pi_x)
        \frac{\partial F}{\partial \pi_x}(\pi) \nu^\mu(d\pi) \right|
        \leq \mu C(F,L) \left| \int \nu^\mu(d\pi) \right|
        = \mu C(F,L) \rightarrow 0,
\end{eqnarray*}
$\nu^0$ satisfies:
\begin{eqnarray}
        \int \sum_x \pi_x ( m(x,\pi) - \overline m_\pi)
                \frac{\partial F}{\partial \pi_x}(\pi) \nu^0(d\pi) = 0.
		\label{eq:nu0}
\end{eqnarray}
After establishing several lemmas, we will use
the above characterization of $\nu^\mu$ and $\nu^0$ to prove the following:

\begin{THM} \label{thm:stat_meas}
Suppose $K$ is symmetric, unimodal, and strictly
decreasing for $x\in [0,L]$, with $K_0 = 1$, and $B_x = 1_{\{|x| \geq M\} }$
with $L < M \leq 2L$. Define $l=-L+M$.
If $\nu^\mu$ is a weak limit point of $\nu^{\mu,N}$, then 
for any $z \in [-l+1,l-1]$, we have
\begin{eqnarray*}
        \nu^\mu \{\pi: \pi_z \geq \delta \} \leq
		\frac{1}{\delta} \frac{\mu}{\frac{\delta_2}{2}+\mu(2L+1)},
\end{eqnarray*}
where $\delta_2 = \min\left(\frac{K_L^2}{2(2L+1)},\frac{K_L^4}{4 (2L+1)K_l^2},
\frac{K_L^2}{2(2L+1)^2}\right)$.
\end{THM}

\begin{COR}
Under the same assumption on $K$ and $B$ as in Theorem~\ref{thm:stat_meas},
we have
\begin{eqnarray*}
        \nu^0 \{\pi: \pi_x = 0 \ \forall x \in [-l+1,l-1]\} = 1
\end{eqnarray*}
if $\nu^0$ is a weak limit point of $\nu^\mu$, and consequently,
for any $\delta>0$,
\begin{eqnarray*}
        \nu^{\mu_i,N_j(\mu_i)}
        \{\pi: \pi_x < \delta \ \forall x \in [-l+1,l-1]\}
        \geq 1-\delta
\end{eqnarray*}
for some sufficiently large $i$ and $j=j_i$.
\label{cor:stat_meas}
\end{COR}

\noindent {\bf Proof of Corollary~\ref{cor:stat_meas}}. \hspace{2mm}
We recall that if a sequence of random variables
$X_n \Rightarrow X_\infty$, then
$\liminf_{n\rightarrow \infty} P(X_n \in A) \geq P(X_\infty \in A)$ for
any open set $A$. Thus Theorem~\ref{thm:stat_meas} implies that for any
$z \in [-l+1,l-1]$ and sequence $\mu_i$ such that
$\nu^{\mu_i} \Rightarrow \nu^0$,
\begin{eqnarray*}
	\nu^0 \{\pi: \pi_z > \delta' \}
	\leq \liminf_{i\rightarrow \infty} \nu^{\mu_i} \{\pi: \pi_z > \delta'\}
	\leq \liminf_{i\rightarrow \infty}
		\frac{1}{\delta'} \frac{\mu_i}{\frac{\delta_2}{2}+\mu_i(2L+1)}
	= 0.
\end{eqnarray*}
This holds for any positive $\delta'$, so $\nu^0 \{\pi: \pi_z > 0\} = 0$
and therefore,
\begin{eqnarray*}
        \nu^0 \{\pi: \pi_x = 0 \ \forall x \in [-l+1,l-1]\} = 1.
\end{eqnarray*}
\qed

\vspace{.3cm}

\begin{LEM}
\label{lem:lb:pi}
If $\mu < 1$, then
\begin{eqnarray*}
	\nu^\mu \left\{\pi: \pi_x > \frac{\mu}{4(2L+2)}
		\ \forall x\in [-L,L] \right\} = 1.
\end{eqnarray*}
\end{LEM}

\proof Let $z\in [-L,L]$ be an arbitrary site.
Let $\delta < \frac{\mu}{2(2L+2)}$ be a small positive constant and
$f\in C^\infty(\mathbb{R})$ be a function that satisfies
the following requirements:
\begin{enumerate}
\item[(a)] $f'(x) = 1$ for $x\leq \frac{\delta}{2}$;
\item[(b)] $f'(x) = 0$ for $x\geq \delta$;
\item[(c)] $f'(x) \in [0,1]$ for $x \in [\frac{\delta}{2},\delta]$; and
\item[(d)] $f(1)=0$.
\end{enumerate}
Define $F(\pi_{-L},\ldots,\pi_L) = f(\pi_z)$. Then~(\ref{eq:nue}) implies:
\begin{eqnarray}
	\int [\pi_z ( m_z - \overline m_\pi) + \mu
                ( 1 - (2L+1) \pi_z)] f'(\pi_z)
                \nu^\mu(d\pi) = 0. \label{eq:lemma8.1_1}
\end{eqnarray}
Since $m_z \geq 0$ and $\overline m_\pi$ is bounded above by $1$
uniformly in $\pi$, we have $m_z - \overline m_\pi > -1$.
The integrand in the above integral is nonzero only for
$\pi_z\in [0,\delta]$, therefore~(\ref{eq:lemma8.1_1}) can be rewritten as:
\begin{eqnarray}
        \int_{\{ \pi: \pi_z \leq \delta \} }
		[\pi_z ( m_z - \overline m_\pi) + \mu
                ( 1 - (2L+1) \pi_z)] f'(\pi_z)
                \nu^\mu(d\pi) = 0. \label{eq:lemma8.1_4}
\end{eqnarray}
Furthermore, the integrand is bounded below:
\begin{eqnarray}
	\pi_z ( m_z - \overline m_\pi) + \mu ( 1 - (2L+1) \pi_z)
	& \geq & - \pi_z + \mu ( 1 - (2L+1) \pi_z) \nnb \\
	& = & - (1+\mu(2L+1))\pi_z + \mu.
	\label{eq:lemma8.1_2}
\end{eqnarray}
Since $\mu<1$, we have $\delta < \frac{\mu}{2(1+2L+1)}
< \frac{\mu}{2(1+\mu(2L+1))}$, and therefore if
$\pi_z \in [0,\delta] \subset [0,\frac{\mu}{2(1+\mu(2L+1))}]$, then
$(1+\mu(2L+1))\pi_z < \frac{\mu}{2}$. Thus~(\ref{eq:lemma8.1_2}) implies that
\begin{eqnarray}
	\pi_z ( m_z - \overline m_\pi) + \mu ( 1 - (2L+1) \pi_z)
	> \frac{\mu}{2}. \label{eq:lemma8.1_3}
\end{eqnarray}
Applying the estimate~(\ref{eq:lemma8.1_3}) to~(\ref{eq:lemma8.1_4}),
we obtain
\begin{eqnarray*}
	0 & \geq & \frac{\mu}{2} \int_{\{ \pi: \pi_z \leq \delta \} }
		f'(\pi_z) \nu^\mu(d\pi) \\
	& \geq & \frac{\mu}{2} \int_{\{ \pi: \pi_z \leq \delta/2 \} }
                f'(\pi_z) \nu^\mu(d\pi) \ \ \ \ \ \mbox{since $f'(x)\geq 0$
			for $x\in [\delta/2,\delta]$} \\
	& = & \frac{\mu}{2} \int_{\{ \pi: \pi_z \leq \delta/2 \} }
		\nu^\mu(d\pi) \ \ \ \ \ \ \ \ \ \mbox{since $f'(x) = 1$
                        for $x\in [0,\delta/2]$} \\
	& = & \frac{\mu}{2} \nu^\mu
		\left\{\pi: \pi_z \leq \frac{\delta}{2} \right\}.
\end{eqnarray*}
Therefore
\begin{eqnarray*}
	\nu^\mu \left\{\pi: \pi_z \leq \frac{\delta}{2} \right\} = 0
\end{eqnarray*}
which implies that
\begin{eqnarray*}
	\nu^\mu \left\{\pi: \pi_x \leq \frac{\delta}{2} \mbox{ for some }
		x \in [-L,L] \right\}
	\leq \sum_{x=-L}^L \nu^\mu \left\{\pi:
		\pi_x \leq \frac{\delta}{2} \right\} = 0.
\end{eqnarray*}     
Thus the lemma follows. \qed

\vspace{.3cm}

Lemma~\ref{lem:lb:pi} implies that $\nu^\mu$-a.s.,
\begin{eqnarray}
	\overline m_\pi \geq \pi_{-L} K_{-L} K_L \pi_L
		\geq \frac{\mu^2 K_L^2}{16(2L+2)^2}. \label{ineq:mbar_pi:crude}
\end{eqnarray}
Recall that $\delta_2 = \min\left(\frac{K_L^2}{2(2L+1)},
\frac{K_L^4}{4 (2L+1)K_l^2}, \frac{K_L^2}{2(2L+1)^2}\right)$
is independent of $\mu$. We define
$\delta_1 = \frac{\mu^2 K_L^2}{16(2L+2)^2}$ and
\begin{eqnarray}
	A = \{\pi: \overline m_\pi < \delta_2\}, \label{eq:defA}
\end{eqnarray}
and the function $\psi:\Pcal(E) \rightarrow \mathbb{R}$,
\begin{eqnarray}
        \psi(\pi) = \sum_x \pi_x ( m_x - \overline m_\pi)^2
                + \mu (2L+1) \sum_x m_x
                        \left( \frac{1}{2L+1} - \pi_x \right).
                \label{eq:defpsi}
\end{eqnarray}
We recall from a formula following~(\ref{def:V}) that
$\frac{\partial \overline m_\pi}{\partial \pi_x} = 2 m_x$, and observe that
\begin{eqnarray*}
	\Gcal \overline m_\pi & = & \sum_x [\pi_x (m_x-\overline m_\pi)
		+ \mu ( 1 - (2L+1) \pi_x)]
                \frac{\partial \overline m_\pi}{\partial \pi_x} \\
        & = & \sum_x [\pi_x ( m_x - \overline m_\pi)
		+ \mu ( 1 - (2L+1) \pi_x)] 2 m_x\\
	& = & 2 \sum_x \pi_x m_x ( m_x - \overline m_\pi)
                + \mu (2L+1) \sum_x m_x
                        \left( \frac{1}{2L+1} - \pi_x \right).
\end{eqnarray*}
Adding $- 2 \sum_x \pi_x \overline m_\pi ( m_x - \overline m_\pi) = 0$
to the right hand side, we continue as follows:
\begin{eqnarray*}
        \Gcal \overline m_\pi & = & 2 \left[
		\sum_x \pi_x m_x ( m_x - \overline m_\pi)
                - \sum_x \pi_x \overline m_\pi ( m_x - \overline m_\pi)
			\right. \nnb \\
        & &	\ \ \ \ \ + \left. \mu (2L+1) \sum_x m_x
                        \left( \frac{1}{2L+1} - \pi_x \right) \right] \\
        & = & 2 \left[ \sum_x \pi_x ( m_x - \overline m_\pi)^2
                + \mu (2L+1) \sum_x m_x
                        \left( \frac{1}{2L+1} - \pi_x \right) \right].
\end{eqnarray*}
Therefore
\begin{eqnarray}
	\Gcal \overline m_\pi = 2 \psi (\pi). \label{eq:Gmbar_pi}
\end{eqnarray}
We will establish in the following two lemmas that on the set $A$ defined
in~(\ref{eq:defA}),
$\psi(\pi)$ is bounded below by $\frac{1}{4(2L+1)} \overline m_\pi^2$.
We write $A$ as a disjoint union of two sets $A_1$ and $A_2$,
where
\begin{eqnarray*}
        A_1 & = & \left\{\pi\in A: \pi_x \leq \frac{1}{2L+1}
		\ \mbox{ for all } x \mbox{ with } m_x \neq 0 \right\}, \\
        A_2 & = & \left\{\pi\in A: \pi_x > \frac{1}{2L+1}
		\ \mbox{ for some } x \mbox{ with } m_x \neq 0 \right\},
\end{eqnarray*}
and prove a lemma for each case.

\begin{LEM} \label{lem:A1}
For any $\pi\in A_1$, we have
$\psi(\pi)\geq \frac{1}{2L+1}\overline m_\pi^2$.
\end{LEM}

\proof For any $\pi\in A_1$, there are two cases:
\begin{enumerate}
\item[Case 1] $\pi_x \leq \frac{1}{2L+1}$ for all $x\in [-L,-l]\cup [l,L]$;
\item[Case 2] $\pi_x > \frac{1}{2L+1}$ and $m_x=0$ for some
	(possibly more than one)
	$x\in [-L,-l]\cup [l,L]$, but for any $x$ with $m_x\neq 0$,
	$\pi_x$ is still $\leq \frac{1}{2L+1}$.
\end{enumerate}
For Case 1,~(\ref{eq:defpsi}) implies
\begin{eqnarray*}
        \psi(\pi) \geq \sum_{x=-l+1}^{l-1} \pi_x (m_x - \overline m_\pi)^2
                + \mu (2L+1) \sum_{x:m_x\neq 0} m_x
                        \left(\frac{1}{2L+1}-\pi_x\right). \\
\end{eqnarray*}
Since $m_x = 0$ for $x\in [-l+1,l-1]$ and $\pi_x \leq \frac{1}{2L+1}$
for all $x \notin [-l,l]$, the second sum on the right hand side is
nonnegative, and therefore
\begin{eqnarray}
        \psi(\pi) \geq \sum_{x=-l+1}^{l-1} \pi_x (m_x-\overline m_\pi)^2
	= \sum_{x=-l+1}^{l-1} \pi_x \overline m_\pi^2, \label{ineq:psi1}
\end{eqnarray}
Using
\begin{eqnarray*}
	\sum_{x=-l+1}^{l-1} \pi_x
	= 1-\sum_{x=-L}^{-l} \pi_x - \sum_{x=l}^L \pi_x
	\geq \frac{2l-1}{2L+1},
\end{eqnarray*}
(\ref{ineq:psi1}) implies that
\begin{eqnarray}
	\psi(\pi) \geq \frac{2l-1}{2L+1} \overline m_\pi^2
	\geq \frac{1}{2L+1} \overline m_\pi^2. \label{ineq:psi1a}
\end{eqnarray}

We now turn to Case 2. Suppose $y$ is a site with $\pi_y > \frac{1}{2L+1}$
and $m_y = 0$. We bound $\psi(\pi)$ in~(\ref{eq:defpsi}) a bit differently
from Case 1. Since $\pi_x \leq \frac{1}{2L+1}$ for any $x$ with $m_x \neq 0$,
the second sum in~(\ref{eq:defpsi}) is nonnegative, and therefore
\begin{eqnarray}
	\psi(\pi) \geq \sum_x \pi_x (m_x - \overline m_\pi)^2
	\geq \pi_y (m_y - \overline m_\pi)^2
	\geq \frac{1}{2L+1} \overline m_\pi^2. \label{ineq:psi2}
\end{eqnarray}
Therefore~(\ref{ineq:psi1a}) and~(\ref{ineq:psi2}) imply the lemma. \qed

\vspace{.3cm}

\begin{LEM} \label{lem:A2} If $\mu < \frac{K_L^4}{8(2L+1)^4 K_l^2}$,
then for any $\pi\in A_2$, we have
$\psi(\pi)\geq \frac{\overline m_\pi^2}{4(2L+1)}$.
\end{LEM}

\proof For $\pi\in A_2$, the main inconvenience
is that the second term in~(\ref{eq:defpsi}), i.e.  the term involving
$\frac{1}{2L+1}-\pi_x$, can be negative.
We divide into three cases and show that in each case, we have
\begin{eqnarray*}
        \psi(\pi) \geq \frac{\overline m_\pi^2}{4(2L+1)}.
\end{eqnarray*}
\begin{enumerate}
\item[Case 1.] For all $x \in [-L,-l]$, $\pi_x \leq \frac{1}{2L+1}$, and
	$i \in [l,L]$ is the rightmost site with
	$\pi_\cdot > \frac{1}{2L+1}$, i.e.  there is no $x$ to the
	right of $i$ with $\pi_x > \frac{1}{2L+1}$;
\item[Case 2.] For all $x \in [l,L]$, $\pi_x \leq \frac{1}{2L+1}$, and
	$i \in [-L,-l]$ is the leftmost site with
	$\pi_\cdot > \frac{1}{2L+1}$, i.e.  there is no $x$ to the
	left of $i$ with $\pi_x > \frac{1}{2L+1}$;
\item[Case 3.] $i_1 \in [-L,-l]$ is the leftmost site with
	$\pi_\cdot > \frac{1}{2L+1}$, and
	$i_2 \in [l,L]$ is the rightmost site with
        $\pi_\cdot > \frac{1}{2L+1}$.
\end{enumerate}

{
\psfrag{big}{$> \frac{1}{2L+1}$}
\psfrag{small}{$\leq \frac{1}{2L+1}$}
\psfrag{L}{$L$}
\psfrag{-L}{$-L$}
\psfrag{i}{$i$}
\psfrag{i-M}{$i-M$}
\psfrag{l}{$l$}
\psfrag{-l}{$-l$}
\begin{figure}[h!]
\centering
\includegraphics[height = 2in, width = 5.5in]{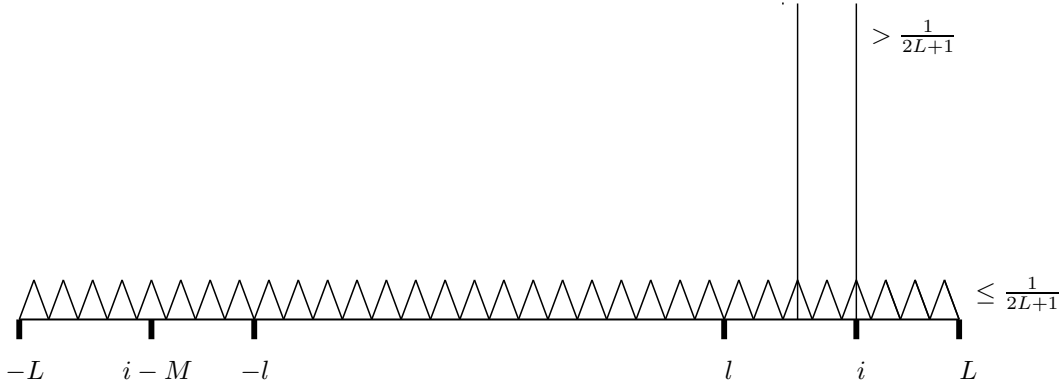}
\caption{Illustration of Case 1}
\end{figure}
}

We first deal with Case 1. The terms in the first sum of the
definition of $\psi$ in~(\ref{eq:defpsi}) are squares and therefore
nonnegative, hence we can throw some of them away and obtain the following:
\begin{eqnarray}
        \psi(\pi) & \geq & \sum_{x=-L}^{i-M} \pi_x ( m_x - \overline m_\pi)^2
                + \mu (2L+1) \sum_{x=-L}^L m_x
                        \left( \frac{1}{2L+1} - \pi_x \right) \nnb \\
	& & \ \ \ \ \	+ \pi_i ( m_i - \overline m_\pi)^2. \label{ineq:psi3}
\end{eqnarray}
For $x\in [-L,i-M]$, $B_{x-i} = 1$, therefore
\begin{eqnarray*}
        m_x = K_x \sum_{z=x+M}^L B_{x-z} K_z \pi_z \geq K_L^2 \pi_i
	> \frac{K_L^2}{2L+1}.
\end{eqnarray*}
Applying the above and the requirement
$\overline m_\pi < \frac{K_L^2}{2(2L+1)}$ for $\pi\in A$
to~(\ref{ineq:psi3}), we obtain:
\begin{eqnarray}
	\sum_{x=-L}^{i-M} \pi_x ( m_x - \overline m_\pi)^2
	& \geq & \sum_{x=-L}^{i-M} \pi_x \left(\frac{K_L^2}{2L+1}
		- \frac{K_L^2}{2(2L+1)}\right)^2	 \nnb \\
	& \geq & \frac{K_L^4}{4(2L+1)^2} \sum_{x=-L}^{i-M} \pi_x.
		\label{ineq:psi4}
\end{eqnarray}
This deals with the first term in~(\ref{ineq:psi3}). For the second term
in~(\ref{ineq:psi3}), we observe that $m_x = 0$ for $x \in [-l,l]$, and
$\frac{1}{2L+1}-\pi_x \geq 0$ for $x \in [-L,-l] \cup [i+1,L]$, therefore
\begin{eqnarray*}
	\sum_{x=-L}^L m_x \left( \frac{1}{2L+1} - \pi_x \right)
	\geq \sum_{x=l}^i m_x \left( \frac{1}{2L+1} - \pi_x \right).
\end{eqnarray*}
Applying the universal bound $\frac{1}{2L+1}-\pi_x \geq -1$ to right hand side
above, we obtain
\begin{eqnarray}
        \sum_{x=-L}^L m_x \left( \frac{1}{2L+1} - \pi_x \right)  
        \geq -\sum_{x=l}^i m_x. \label{ineq:psi9}
\end{eqnarray}
Now applying~(\ref{ineq:psi4}),~(\ref{ineq:psi9}), and the requirement
$\pi_i > \frac{1}{2L+1}$ to~(\ref{ineq:psi3}), we obtain
\begin{eqnarray}
	\psi(\pi) \geq \frac{K_L^4}{4(2L+1)^2} \sum_{x=-L}^{i-M} \pi_x
		- \mu (2L+1) \sum_{x=l}^i m_x
		+ \frac{1}{2L+1} ( m_i - \overline m_\pi)^2. \label{ineq:psi5}
\end{eqnarray}
We observe that
\begin{eqnarray}
	\sum_{x=l}^i m_x & = & \sum_{x=l}^i K_x \sum_{y=-L}^{x-M} K_y \pi_y
	\leq K_l^2 \sum_{x=l}^i \sum_{y=-L}^{x-M} \pi_y
	\leq K_l^2 \sum_{x=l}^i \sum_{y=-L}^{i-M} \pi_y \nnb \\
	& \leq & K_l^2 (L-l+1) \sum_{x=-L}^{i-M} \pi_x
	\leq K_l^2 (2L+1) \pi([-L,i-M]). \label{ineq:psi10}
\end{eqnarray}
Therefore the computation in~(\ref{ineq:psi5}) can be continued as follows:
\begin{eqnarray}
	\psi(\pi) & \geq & \frac{K_L^4}{4(2L+1)^2} \pi([-L,i-M])
		- \mu K_l^2 (2L+1)^2 \pi([-L,i-M]) \nnb \\
	& & \ \ \ \ \	+ \frac{1}{2L+1} ( m_i - \overline m_\pi)^2 \nnb \\
	& \geq & \frac{K_L^4}{8(2L+1)^2} \pi([-L,i-M])
		+ \frac{1}{2L+1} ( m_i - \overline m_\pi)^2 \label{ineq:psi8}
\end{eqnarray}
since $\mu < \frac{K_L^4}{8(2L+1)^4 K_l^2}$.
The right hand side of the above is a sum of two nonnegative terms,
both of which cannot be small at the same time. Indeed,
notice that
\begin{eqnarray*}
	m_i = K_i \sum_{x=-L}^{i-M} K_x \pi_x \leq K_l^2 \pi([-L,i-M]).
\end{eqnarray*}
Therefore, if $\pi([-L,i-M]) < \frac{1}{2K_l^2} \overline m_\pi$,
then $m_i < \frac{1}{2} \overline m_\pi$, and we have
$( m_i - \overline m_\pi)^2 > \frac{1}{4} \overline m_\pi^2$,
hence the second term in~(\ref{ineq:psi8}) alone implies
$\psi(\pi) \geq \frac{\overline m_\pi^2}{4(2L+1)}$.
Otherwise, $\pi([-L,i-M]) \geq \frac{1}{2K_l^2} \overline m_\pi$,
then the first term in~(\ref{ineq:psi8}) implies
$\psi(\pi) \geq \frac{K_L^4 \overline m_\pi}{16 (2L+1)^2 K_l^2}$.
So we have the following estimate on $\psi(\pi)$:
\begin{eqnarray}
	\psi(\pi) \geq \min\left(\frac{K_L^4 \overline m_\pi}{
		16 (2L+1)^2 K_l^2}, \frac{\overline m_\pi^2}{4(2L+1)} \right).
	\label{ineq:psi6}
\end{eqnarray}
Since $\delta_2 \leq \frac{K_L^4}{4 (2L+1)K_l^2}$, on the set
$A = \{\pi: \overline m_\pi < \delta_2\}$, we have
$\frac{K_L^4 \overline m_\pi}{16 (2L+1)^2 K_l^2} \geq 
\frac{\overline m_\pi^2}{4(2L+1)}$, therefore~(\ref{ineq:psi6}) implies
\begin{eqnarray*}
	\psi(\pi) \geq \frac{\overline m_\pi^2}{4(2L+1)}.
\end{eqnarray*}

{
\psfrag{big}{$> \frac{1}{2L+1}$}
\psfrag{small}{$\leq \frac{1}{2L+1}$}
\psfrag{L}{$L$}
\psfrag{-L}{$-L$}
\psfrag{i1}{$i_2$}
\psfrag{i1-M}{$i_2-M$}
\psfrag{i2}{$i_1$} 
\psfrag{i2-M}{$i_1+M$}
\psfrag{l}{$l$}
\psfrag{-l}{$-l$}
\begin{figure}[h!]
\centering
\includegraphics[height = 2in, width = 5.5in]{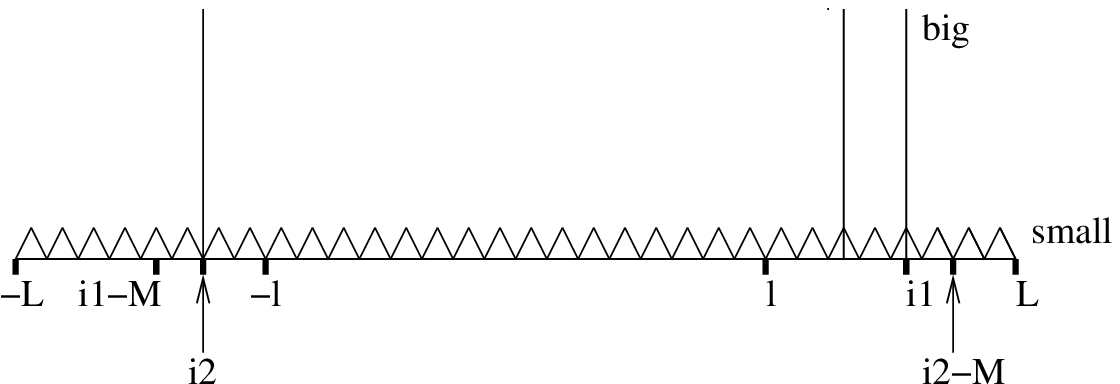}
\caption{Illustration of Case 1}
\end{figure}
}

Case 2 follows by exactly the same argument. For Case 3, we first
observe that if $i_2-i_1 \geq M$ then $B_{i_2 - i_1} = 1$ and on
$A=\{\pi: \overline m_\pi < \delta_2\}$, where by definition
$\delta_2 \leq \frac{K_L^2}{2(2L+1)^2}$, we have
$\delta_2 > \overline m_\pi \geq \pi_{i_1} K_{i_1} K_{i_2} \pi_{i_2}
\geq \frac{K_L^2}{(2L+1)^2}$, which is impossible. Thus $i_2-i_1 < M$, and
\begin{eqnarray}
        \psi(\pi) & \geq & \sum_{x=-L}^{i_2-M} \pi_x ( m_x - \overline m_\pi)^2
                + \pi_{i_2} ( m_{i_2} - \overline m_\pi)^2
		+ \sum_{x=i_1+M}^{L} \pi_x ( m_x - \overline m_\pi)^2 \nnb \\
	& & 	+ \pi_{i_1} ( m_{i_1} - \overline m_\pi)^2
		+ \mu (2L+1) \sum_{x=-L}^L m_x
                        \left( \frac{1}{2L+1} - \pi_x \right).
		\label{ineq:psi7}
\end{eqnarray}
We can use techniques similar to those used for Case 1 (leading
to~(\ref{ineq:psi4}) and~(\ref{ineq:psi10})) to obtain the following
bounds:
\begin{eqnarray*}
	& & \sum_{x=-L}^{i_2-M} \pi_x ( m_x - \overline m_\pi)^2
	+ \sum_{x=i_1+M}^{L} \pi_x ( m_x - \overline m_\pi)^2 \\
	& & \ \ \ \ \
	\geq \frac{K_L^4}{4(2L+1)^2} (\pi([-L,i_2-M]) + \pi([i_1+M,L]))
\end{eqnarray*}
and
\begin{eqnarray*}
	\sum_{x=-L}^L m_x \left( \frac{1}{2L+1} - \pi_x \right)
	& \geq & \sum_{x=l}^{i_2} m_x \left( \frac{1}{2L+1} - \pi_x \right)
	+ \sum_{x=i_1}^{-l} m_x \left( \frac{1}{2L+1} - \pi_x \right) \\
	& \geq & - \sum_{x=l}^{i_2} m_x - \sum_{x=i_1}^{-l} m_x \\
	& \geq & - K_l^2 (2L+1) \pi([-L,i_2-M] \cup [i_1+M,L]).
\end{eqnarray*}
Since $\mu < \frac{K_L^4}{8(2L+1)^4 K_l^2}$, the above two estimates
applied to~(\ref{ineq:psi7}) implies that
\begin{eqnarray*}
        \psi(\pi) & \geq &
		\frac{K_L^4}{4(2L+1)^2} \pi([-L,i_2-M] \cup [i_1+M,L]) \\
	& &	- \mu (2L+1) K_l^2 (2L+1) \pi([-L,i_2-M] \cup [i_1+M,L]) \\
	& &	+ \frac{1}{2L+1} ( m_{i_2} - \overline m_\pi)^2
		+ \frac{1}{2L+1} ( m_{i_1} - \overline m_\pi)^2 \nnb \\
	& \geq &	\frac{K_L^4}{8(2L+1)^2} \pi([-L,i_2-M])
		+ \frac{1}{2L+1} ( m_{i_2} - \overline m_\pi)^2 \nnb \\
	& &	+ \frac{K_L^4}{8(2L+1)^2} \pi([i_1+M,L])
                + \frac{1}{2L+1} ( m_{i_1} - \overline m_\pi)^2.
\end{eqnarray*}
We can now apply the technique leading to~(\ref{ineq:psi6}) to the sum
of the first two terms above, and then to the sum of the last two terms,
to obtain:
\begin{eqnarray*}
        \psi(\pi) \geq 2 \min\left(\frac{K_L^4 \overline m_\pi}{
                16 (2L+1)^2 K_l^2}, \frac{\overline m_\pi^2}{4(2L+1)} \right),
\end{eqnarray*}
which by virtue of
$\overline m_\pi < \delta_2 \leq \frac{K_L^4}{4(2L+1)K_l^2}$ on $A$ implies
\begin{eqnarray*}
        \psi(\pi) \geq \frac{\overline m_\pi^2}{2(2L+1)}.
\end{eqnarray*}
Thus we have established the necessary bound on $\psi(\pi)$ for $\pi \in A_2$
in all three cases. \qed

\vspace{.3cm}

\begin{LEM} \label{lem:ineq:mbar_pi}
If $\delta_1 = \frac{\mu^2 K_L^2}{16(2L+2)^2}$ and
$\delta_2 = \min(\frac{K_L^2}{2(2L+1)},\frac{K_L^4}{4 (2L+1)K_l^2},
\frac{K_L^2}{2(2L+1)^2})$, then
\begin{eqnarray*}
        \nu^\mu \left\{\pi:
                \overline m_\pi \leq \frac{\delta_2}{2} \right\} = 0.
\end{eqnarray*}
\end{LEM}
\proof Recall from~(\ref{ineq:mbar_pi:crude}) that $\nu^\mu$-a.s.,
\begin{eqnarray}
        \overline m_\pi \geq \delta_1. \label{ineq:delta1}
\end{eqnarray}
Let $f\in C^\infty(\mathbb{R})$ to be a function that satisfies the
following requirements:
\begin{enumerate}
\item[(a)]
	$f'(x)$ increases from 0 to $\frac{1}{\delta_1^2}$
        for $x \in [0,\delta_1]$;
\item[(b)] $f'(x) = \frac{1}{x^2}$ for
	$\delta_1 \leq x \leq \frac{\delta_2}{2}$;
\item[(c)] $f'(x)$ decreases from $\frac{4}{\delta_2^2}$ to $0$
        for $x \in [\frac{\delta_2}{2},\frac{2\delta_2}{3}]$;
\item[(d)] $f'(x) = 0$ for $x\geq \frac{2\delta_2}{3}$.
\end{enumerate}
Define $F(\pi_{-L},\ldots,\pi_L) = f(\overline m_\pi)$.
Then~(\ref{eq:nue}) implies:
\begin{eqnarray*}
        \int f'(\overline m_\pi) \Gcal \overline m_\pi \nu^\mu(d\pi) = 0.
\end{eqnarray*}
Substituting~(\ref{eq:Gmbar_pi}) into the above equation, we obtain
\begin{eqnarray}
	\int f'(\overline m_\pi) \psi(\pi) \nu^\mu(d\pi) = 0.
	\label{eq:lemma8.5_1}
\end{eqnarray}
Lemmas~\ref{lem:A1} and~\ref{lem:A2} imply that
$\psi(\pi)$ is bounded below by $\frac{1}{4(2L+1)} \overline m_\pi^2$
on the set $A = \{\pi: \overline m_\pi < \delta_2\}$,
defined in~(\ref{eq:defA}).
Applying this fact and~(\ref{ineq:delta1}) to~(\ref{eq:lemma8.5_1}),
we obtain
\begin{eqnarray*}
        0 & = & 2 \int_{
                \{\pi: 0 \leq \overline m_\pi \leq \frac{2\delta_2}{3} \} }
                \psi(\pi) f'(\overline m_\pi) \nu^\mu(d\pi)
                \ \ \ \ \ \mbox{since $f'(x)$ only nonzero for $x$ in
                        $[0,\frac{2\delta_2}{3}]$} \\
        & = & 2 \int_{ \{\pi: \delta_1 \leq \overline m_\pi
                        \leq \frac{2 \delta_2}{3} \} }
                \psi(\pi) f'(\overline m_\pi) \nu^\mu(d\pi)
                \ \ \ \ \ \mbox{since $\nu^\mu\{\pi: \overline m_\pi
                        < \delta_1\} = 0$ by~(\ref{ineq:mbar_pi:crude})} \\
	& \geq & 2 \int_{ \{\pi: \delta_1 \leq \overline m_\pi
                        \leq \frac{\delta_2}{2} \} }
                \psi(\pi) f'(\overline m_\pi) \nu^\mu(d\pi)
                \ \ \ \ \ \mbox{since $\psi(\pi) f'(\overline m_\pi)>0$
                        if $\overline m_\pi \in
                        [\frac{\delta_2}{2},\frac{2\delta_2}{3}]$} \\
        & \geq & 2 \int_{ \{\pi: \delta_1 \leq \overline m_\pi
                        \leq \frac{\delta_2}{2} \} }
                \frac{\overline m_\pi^2}{4(2L+1)} \frac{1}{\overline m_\pi^2}
                        \nu^\mu(d\pi) \ \ \ \ \ \mbox{by the bound
                        on $\psi(\overline m_\pi)$ for $\pi\in A$} \\
        & = & \frac{1}{2(2L+1)} \int_{ \{\pi: \delta_1 \leq \overline m_\pi
                        \leq \frac{\delta_2}{2} \} } \nu^\mu(d\pi).
\end{eqnarray*}
Therefore
\begin{eqnarray*}
        \nu^\mu \left\{\pi:
                \overline m_\pi \leq \frac{\delta_2}{2} \right\}
        = \nu^\mu \left\{\pi:
                0 \leq \overline m_\pi < \delta_1 \right\}
        + \nu^\mu \left\{\pi: \delta_1 \leq \overline m_\pi
                        \leq \frac{\delta_2}{2} \right\}
        = 0+0 = 0,
\end{eqnarray*}
as required. \qed

\vspace{.3cm}

\noindent {\bf Proof of Theorem~\ref{thm:stat_meas}}. \hspace{2mm}
For an arbitrary site $z\in [-l+1,l-1]$, we have $m_z = 0$. If we take
$F(\pi) = \pi_z$, then by~(\ref{eq:nue}), we have
\begin{eqnarray*}
        0 & = & \int \pi_z ( m(z,\pi) - \overline m_\pi) + \mu(1-(2L+1)\pi_z)
		\nu^\mu(d\pi) \\
        & = & \int \mu - (\overline m_\pi + \mu(2L+1)) \pi_z \nu^\mu(d\pi),
\end{eqnarray*}
so
\begin{eqnarray*}
	\mu & = & \int (\overline m_\pi + \mu(2L+1)) \pi_z \nu^\mu(d\pi) \\
	& = & \int_{\{\pi:\overline m_\pi > \delta_2/2\}}
		(\overline m_\pi + \mu(2L+1)) \pi_z \nu^\mu(d\pi) \\
	& & \ \ \ \ \ + \int_{\{\pi:\overline m_\pi \leq \delta_2/2\}}
		(\overline m_\pi + \mu(2L+1)) \pi_z \nu^\mu(d\pi),
\end{eqnarray*}
where $\delta_2$ is as defined in Lemma~\ref{lem:ineq:mbar_pi}. The same
lemma shows that                     
$\{\pi: \overline m_\pi \leq \delta_2/2\}$ has $\nu^\mu$-measure 0, so
the second integral in the above equation is 0, thus
\begin{eqnarray*}
	\mu = \int_{\{\pi:\overline m_\pi > \delta_2/2\}}
                (\overline m_\pi + \mu(2L+1)) \pi_z \nu^\mu(d\pi)
	\geq \int \left(\frac{\delta_2}{2} + \mu(2L+1)\right) \pi_z
		\nu^\mu(d\pi),
\end{eqnarray*}
i.e.
\begin{eqnarray*}
	\int \pi_z \nu^\mu(d\pi)
	\leq \frac{\mu}{\frac{\delta_2}{2} + \mu(2L+1)}.
\end{eqnarray*}
The observation $\nu^\mu \{\pi: \pi_z \geq \delta \}
\leq \frac{1}{\delta} \int \pi_z \nu^\mu(d\pi)$ completes the proof.
\qed

\appendix
\chapter{A Result on the Conditioned Dieckmann-Doebeli Model}
\label{ch:appen}
In this section, we deal with a special case of the conditioned
Dieckmann-Doebeli model described in~(\ref{eq:badmodel:V2})
and~(\ref{eq:badmodel}), and show that in this special case, there
exist symmetric bimodal stationary distributions. Let
$E = [-L,L]\cap \Zbold$, $A$ be a Markov transition matrix associated with
mutation, $\pi(t) \in \Pcal(E)$ for all $t \in \Zbold^+$,
and $M$ be an \emph{even} constant such that $L-1 \leq M < 2(L-1)$, then the
equation of the discrete-time dynamical system is as follows:
\begin{eqnarray}
	\pi_x(t+1) & = & \sum_y A(y,x) \frac{\pi_y(t) V_y(\pi(t))}{
                \sum_z \pi_z(t) V_z(\pi(t))}, \nnb \\
        \mbox{where } V_x(\pi) & = & V^{(2)}_x(\pi)
		= \frac{K_x}{\sum_z C_{x-z} \pi_z}, \nnb \\
	K_x & = & 1_{\{|x| \leq L-1\}}, \nnb \\
	\mbox{and } C_x & = & 1_{\{|x| \leq M\}}. \label{eq:badmodel2}
\end{eqnarray}
Every step of~(\ref{eq:badmodel2}) can be divided into three sub-steps:
\begin{eqnarray}
	& \mbox{Resampling} &	\pi'_x (t) = \pi_x(t) V_x(\pi(t)); \nnb \\
	& \mbox{Mutation} &	\pi''_x (t) = \sum_y A(y,x) \pi'_x(t); \nnb \\
	& \mbox{Normalization}&	\pi_x(t+1) = \frac{\pi''_x (t)}{
		\sum_y \pi''_y (t)}. \label{eq:3steps}
\end{eqnarray}
Notice that performing the normalization step before the mutation step does
not change the model, but for this section, we will use the step order
in~(\ref{eq:3steps}).

If $K$, $C$, and $A(y,x)=A(y-x)$ are symmetric about $0$, then the map
$\pi(t) \mapsto \pi(t+1)$ maps the set of symmetric probability measures
on $[-L,L]$ to itself.
Therefore by Brouwer's fixed point theorem, there exists \emph{symmetric}
stationary distribution(s). We first derive a few simple facts about symmetric
stationary distributions in the no-mutation case,
i.e. when $A = I$. In this case, if $\nu^0$ is a stationary
distribution of~(\ref{eq:badmodel2}), we must have
\begin{eqnarray}
        & & \mbox{ if $\nu^0_x \neq 0$, then }
        V_x(\nu^0) = \frac{K_x}{\sum_z C_{x-z} \nu^0_z}
	= \frac{K_x}{\sum_{z=x-M}^{x+M} \nu^0_z}
                \mbox{ is a constant}. \label{eq:nu0cond2}
\end{eqnarray}
This is the same condition as~(\ref{cond:Vconst}).
Since $K_x = 0$ outside the interval $[-L+1,L-1]$, $V_x = 0$ outside
that interval, i.e. at $x=\pm L$,
therefore the support of any stationary distribution
$\nu^0$ must lie in $[-L+1,L-1]$. We restrict our attention to sites
in $[-L+1,L-1]$. If the competition intensity function $C$ is rectangular, 
as in~(\ref{eq:badmodel2}), two sites either compete
at intensity 1 or they do not compete against each other at all. If $M = L-1$,
then $0$ is the only site that competes with every other site in $[-L+1,L-1]$;
but if $M=2(L-1)-1$, then for any site $x \in [-L+2,L-2]$, $x$ competes with
every site in $[-L+1,L-1]$. Define $l=-(L-1-M)$, then $[-l,l]$ contains the
sites that compete with every site in $[-L+1,L-1]$, therefore
\begin{eqnarray}
        V_x(\nu^0) = \frac{K_x}{\sum_{z=x-M}^{x+M} \nu^0_z} = \left\{
        \begin{array}{ll}       1 & x\in [-l,l] \\
                                0 & x \in (\infty,-L]\cup [L,\infty)
        \end{array} \right. .
	\label{cond:V1}
\end{eqnarray}
We also observe that $\frac{M}{2}>l$ because $M-2l = M+2(L-1-M) = 2(L-1)-M>0$.
If there is some mass in $[-L+1,-\frac{M}{2}-1] \cup [\frac{M}{2}+1,L-1]$,
then since no site in $[\frac{M}{2}+1,L-1]$ competes with any site in
$[-L+1,-\frac{M}{2}-1]$, we have
\begin{eqnarray*}                 
	\sum_{z=x-M}^{x+M} \nu^0_z < 1
\end{eqnarray*}                   
for $x\in [-L+1,-\frac{M}{2}-1] \cup [\frac{M}{2}+1,L-1]$,
which by the definition of $V_x(\pi)$
in~(\ref{eq:badmodel2}) implies the following:
\begin{eqnarray}
        \ V_x(\nu^0) > 1 \ \mbox{ for }
		\ x\in \left[-L+1,-\frac{M}{2}-1\right]
		\cup \left[\frac{M}{2}+1,L-1\right].
	\label{cond:V2}
\end{eqnarray}
Combining the results on fitness $V_x(\pi)$ in~(\ref{cond:V1})
and~(\ref{cond:V2}) and condition~(\ref{eq:nu0cond2}), we conclude
that stationary distributions with rectangular $K$ and $C$ as defined
in~(\ref{eq:badmodel2}) must have all the mass falling in
either $\left[-L+1,-\frac{M}{2}-1\right] \cup \left[\frac{M}{2}+1,L-1\right]$
or $\left[-\frac{M}{2},\frac{M}{2}\right]$.

{
\psfrag{L}{$L$}
\psfrag{-L}{$-L$}
\psfrag{i1}{$\frac{M}{2}$}
\psfrag{i2}{$-\frac{M}{2}$}
\psfrag{l}{$l$}
\psfrag{-l}{$-l$}
\begin{figure}[t!]
\centering
\includegraphics[height = 0.5in, width = 5.5in]{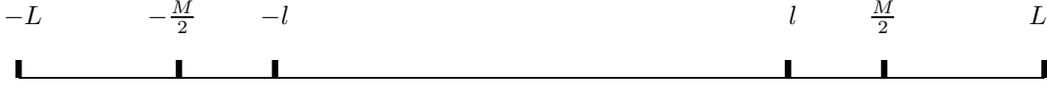}
\caption{Relative locations of various sites of interest}
\end{figure}
}

Now we turn to the case with small mutations.
We take $A^\mu$ to be an operator that corresponds to a small
1-step mutation, i.e. convolution with
$\mu \delta_{-1} + (1-2\mu) \delta_0 + \mu \delta_1$, and define
$\nu^\mu$ to be a stationary distribution of~(\ref{eq:badmodel2})
with mutation kernel
$A^\mu$. Any stationary distribution $\nu^\mu$ of~(\ref{eq:badmodel2})
satisfies the following condition:
\begin{eqnarray}
        \forall x \in [-L,L], \ \mu \nu^\mu_{x-1} V_{x-1}
        	+ (1-2\mu) \nu^\mu_x V_x
        	+ \mu \nu^\mu_{x+1} V_{x+1}
        = \overline V \nu^\mu_x, \label{eq:nuecond}
\end{eqnarray}
where $V_y=V_y(\nu^\mu)$ and
\begin{eqnarray*}
        \overline V = \overline V(\nu^\mu) = \sum_{z: \nu^\mu_z\neq 0}
                \nu^\mu_z V_z(\nu^\mu)
\end{eqnarray*}
is the normalization constant. Condition~(\ref{eq:nuecond}) implies that
$\nu^\mu$ is nowhere zero in $[-L,L]$; otherwise, say $\nu^\mu_{z} = 0$, then
then $\nu^\mu_{z-1} = \nu^\mu_{z+1} = 0$ as well, which by induction
means that $\nu^\mu_x = 0$ for all $x$, a clear contradiction.
Since the support of $\nu^\mu$ has expanded on both sides each by 1 site
compared to $\nu^0$, the sites where $V_\cdot(\nu^\mu)$ is constant 1
should correspondingly contract by 1 site on each side:
\begin{eqnarray}
        V_x(\nu^\mu) = \frac{K_x}{\sum_{z=x-M}^{x+M} \nu^\mu_z}
        = \left\{
        \begin{array}{ll}       1 & x\in [-l+1,l-1] \\
                                >1 & x\in [-L+1,-l] \cup [l,L-1] \\
                                0 & x \in (\infty,-L]\cup [L,\infty)
        \end{array} \right. . \label{eq:facts:V}
\end{eqnarray}
Therefore, for sites in the middle, i.e. $x\in [-l+2,l-2]$,
$\nu^\mu$ satisfies:
\begin{eqnarray*}
        \mu \nu^\mu_{x-1} + (1-2\mu) \nu^\mu_x + \mu \nu^\mu_{x+1}
        = \overline V(\nu^\mu) \nu^\mu_x.
\end{eqnarray*}

Since every site $x$ in $[-L+1,L-1]$ competes with all sites
lying on the same side (with respect to the origin) as $x$,
and the stationary distributions we consider are symmetric,
$\sum_{z=x-M}^{x+M} \nu^\mu_x$ is at least
$\frac{1}{2}$ for $x \in [-L+1,L-1]$,
therefore
\begin{eqnarray}
	V_x(\nu^\mu) \leq 2 \mbox{ for } x \in [-L+1,L-1],
	\label{ineq:Vupper1}
\end{eqnarray}
i.e. the normalization constant $\overline V(\nu^\mu)$ is bounded above by $2$.
We also need a nontrivial lower bound of $\overline V(\nu^\mu)$ for
\emph{symmetric} $\nu^\mu$ that is uniform for small $\mu$ for the proof
of the upcoming theorem.
We rewrite condition~(\ref{eq:nuecond}) for $x=L$ and $x=L-1$,
taking into account the fact $V_L(\nu^\mu) = 0$ from~(\ref{eq:facts:V}):
\begin{eqnarray*}
        \mu \nu^\mu_{L-1} V_{L-1}(\nu^\mu)
	& = & \overline V(\nu^\mu) \nu^\mu_L, \\
	\mu \nu^\mu_{L-2} V_{L-2}(\nu^\mu)
		+ (1-2\mu) \nu^\mu_{L-1} V_{L-1}(\nu^\mu)
	& = & \overline V(\nu^\mu) \nu^\mu_{L-1}.
\end{eqnarray*}
Dividing both sides of the above two equations, we get
\begin{eqnarray*}
	\frac{\nu^\mu_L}{\nu^\mu_{L-1}}
	& = & \frac{\mu \nu^\mu_{L-1} V_{L-1}(\nu^\mu)}{
		\mu \nu^\mu_{L-2} V_{L-2}(\nu^\mu)
		+ (1-2\mu) \nu^\mu_{L-1} V_{L-1}(\nu^\mu)} \\
	& \leq & \frac{\mu \nu^\mu_{L-1} V_{L-1}(\nu^\mu)}{
		(1-2\mu) \nu^\mu_{L-1} V_{L-1}(\nu^\mu)}
	= \frac{\mu}{1-2\mu},
\end{eqnarray*}
which is $<1$ if $\mu < \frac{1}{3}$. This means that
$\nu^\mu_L = \nu^\mu_{-L} \leq \frac{1}{4}$, for otherwise, 
$\nu^\mu_L = \nu^\mu_{-L} > \frac{1}{4}$ implies that
$\nu^\mu_{L-1} = \nu^\mu_{-L+1} \leq \frac{1}{2} - \nu^\mu_L < \frac{1}{4}$,
hence $\frac{\nu^\mu_L}{\nu^\mu_{L-1}} > 1$, a contradiction.
Since $V_x(\nu^\mu) \geq 1$ for $x\in [-L+1,L-1]$ by~(\ref{eq:facts:V}),
we have
\begin{eqnarray*}
        \overline V(\nu^\mu) \geq \sum_{z=-L+1}^{L-1} \nu^\mu_z V_z(\nu^\mu)
	\geq \sum_{z=-L+1}^{L-1} \nu^\mu_x.
\end{eqnarray*}
The fact $\nu^\mu_L = \nu^\mu_{-L} \leq \frac{1}{4}$ then implies
\begin{eqnarray}
	\overline V(\nu^\mu) \geq 1 - \nu^\mu_{-L} - \nu^\mu_L
	\geq \frac{1}{2}. \label{ineq:Mnu:crude}
\end{eqnarray} 
In particular, $\overline V(\nu^\mu)$ is bounded between $\frac{1}{2}$ and $2$.
We take a sequence $\mu_n \rightarrow 0$, such that
$\nu^n = \nu^{\mu_n}$ converges to some $\nu^0$,
then since $\overline V(\nu)$ is a continuous function of $\nu$,
$\overline V_n = \overline V(\nu^n)$ also converges to a positive constant
$\overline V$. We will prove the following:

\begin{THM}
	If $\nu^n = \nu^{\mu_n}$ is
	a convergent sequence of \emph{symmetric} stationary distributions
	for the conditioned Dieckmann-Doebeli model in~(\ref{eq:badmodel2}),
	then $\nu^n_{\left[-\frac{M}{2},\frac{M}{2}\right]} \rightarrow 0$
	as $n \rightarrow \infty$.
\label{thm:append}
\end{THM}

\proof 
Since $\nu^n \rightarrow \nu^0$ and
$\overline V(\nu^n) \rightarrow \overline V$
as $n \rightarrow \infty$, condition~(\ref{eq:nuecond}) converges to
the following:
\begin{eqnarray*}
        \forall x \in [-L,L], \ \nu^0_x V_x(\nu^0) = \overline V \nu^0_x.
\end{eqnarray*}
Therefore $\nu^0$ is in fact an stationary distribution for the no-mutation
case, i.e.~(\ref{eq:badmodel2}) with $A=I$.
If some mass of $\nu^0$ lies in
$\left[-L+1,-\frac{M}{2}-1\right] \cup \left[\frac{M}{2}+1,L-1\right]$,
then $\nu^0_{\left[-\frac{M}{2},\frac{M}{2}\right]}$ must be zero,
which means that
$\nu^n_{\left[-\frac{M}{2},\frac{M}{2}\right]} \rightarrow 0$
when $n \rightarrow \infty$ as required.
Therefore it suffices to show that
$\nu^0_{\left[-L+1,-\frac{M}{2}-1\right] \cup \left[\frac{M}{2}+1,L-1\right]}
> 0$.

We assume, toward a contradiction, that
$\nu^0_{\left[-L+1,-\frac{M}{2}-1\right]}
	= \nu^0_{\left[\frac{M}{2}+1,L-1\right]} = 0$.
Then for any positive $\delta$, we have
$\nu^n_{\left[-L+1,-\frac{M}{2}-1\right]}
	= \nu^n_{\left[\frac{M}{2}+1,L-1\right]} < \delta$
for sufficiently large $n$.
We first derive more refined (than~(\ref{ineq:Mnu:crude})
and~(\ref{ineq:Vupper1}) respectively)
lower and upper bounds for $\overline V_n$. Because of the supposition
$\nu^n_{\left[\frac{M}{2}+1,L-1\right]}<\delta$,~(\ref{eq:nuecond})
and~(\ref{eq:facts:V}) with $x=L$ and the bound $\nu^n_{L-1}<\delta$
imply that
\begin{eqnarray}
        \nu^n_L<\frac{\delta \mu_n V_{L-1}(\nu^n)}{\overline V_n}
	\leq \frac{2 \delta \mu_n}{\overline V_n}
	\label{ineq:nun_L2}
\end{eqnarray}
by~(\ref{ineq:Vupper1}).
Applying the estimate~(\ref{ineq:Mnu:crude}) to the right hand side,
we have
\begin{eqnarray}
        \nu^n_L<\frac{2 \delta \mu_n}{\overline V_n} \leq 4 \delta \mu_n.
                \label{ineq:nun_L}
\end{eqnarray}
Therefore
\begin{eqnarray}
        \overline V_n & = & \sum_{z: \nu^n_z\neq 0} \nu^n_z V_z(\nu^n) \nnb \\
        & = & \sum_{z=-L+1}^{L-1} \nu^n_z V_z(\nu^n) \ \ \ \ \ \mbox{ since
                $V_{-L}(\nu^n) = V_L(\nu^n) = 0$ by~(\ref{eq:facts:V})} \nnb \\
        & \geq & \sum_{z=-L+1}^{L-1} \nu^n_z \ \ \ \ \ \mbox{ since
                $V_x(\nu^n) \geq 1$ for $x \in [-L+1,L-1]$} \nnb \\
        & = & 1 - \nu^n_{-L} - \nu^n_L \nnb \\
        & \geq & 1-8\delta \mu_n, \label{ineq:Mnu:fine}
\end{eqnarray}
by~(\ref{ineq:nun_L}).

For $x\in \left[-\frac{M}{2},\frac{M}{2}\right]$,
\begin{eqnarray}
	V_x(\nu^n) = \frac{K_x}{\sum_{z=x-M}^{x+M} \nu^n_z}
	\leq \frac{1}{\sum_{z=-\frac{M}{2}}^{\frac{M}{2}} \nu^n_z}
	= \frac{1}{1-2\nu^n_{\left[\frac{M}{2}+1,L\right]}}
	\leq \frac{1}{1-2(\delta+4\delta\mu_n)} \label{ineq:Vupper2}
\end{eqnarray}
by~(\ref{ineq:nun_L}) and the supposition
$\nu^n_{\left[\frac{M}{2}+1,L-1\right]} < \delta$.
Therefore
\begin{eqnarray}
	\overline V(\nu^n) & = & \sum_{x=-L}^L \nu^n_x V_x(\nu^n)
	= \sum_{x=-\frac{M}{2}}^{\frac{M}{2}} \nu^n_x V_x(\nu^n)
		+ 2 \sum_{x=\frac{M}{2}+1}^{L-1} \nu^n_x V_x(\nu^n) \nnb \\
	& \leq & \frac{1}{1-2(\delta+4\delta\mu_n)} + 4\delta
	\label{ineq:Vupper3}
\end{eqnarray}
using~(\ref{ineq:Vupper2}) for the first sum and~(\ref{ineq:Vupper1})
for the second.

Let $r=\nu^n_{\left[-\frac{M}{2},-l-1\right]}
=\nu^n_{\left[l+1,\frac{M}{2}\right]}$, then
\begin{eqnarray}
	V_{L-1}(\nu^n) & = & \frac{K_{L-1}}{\sum_{z=L-1-M}^{L-1+M} \nu^n_z}
	= \frac{1}{\sum_{z=-l}^{L} \nu^n_z}
	= \frac{1}{1-\nu^n_{ \left[-\frac{M}{2},-l-1\right] }
		- \nu^n_{\left[-L+1,-\frac{M}{2} \right]} - \nu^n_{-L}} \nnb \\
	& \geq & \frac{1}{1-r}. \label{ineq:VL-1}
\end{eqnarray}
Condition~(\ref{eq:nuecond}) applied with $x=L-1$ implies
\begin{eqnarray*}
        (1-2\mu_n) \nu^n_{L-1} V_{L-1}(\nu^n)
	\leq \overline V(\nu^n) \nu^n_{L-1},
\end{eqnarray*}
hence
\begin{eqnarray*}
	V_{L-1}(\nu^n) \leq \frac{\overline V(\nu^n)}{1-2\mu_n}.
\end{eqnarray*}
Then~(\ref{ineq:Vupper3}) and~(\ref{ineq:VL-1}) imply
\begin{eqnarray*}
	\frac{1}{1-r}
		\leq \frac{1}{1-2\mu_n} \left(
			\frac{1}{1-2(\delta+4\delta\mu_n)} + 4\delta \right)
\end{eqnarray*}
i.e.
\begin{eqnarray*}
	1-r \geq \frac{1-2\mu_n}{
			\frac{1}{1-2(\delta+4\delta\mu_n)} + 4\delta}
		= \frac{(1-2\mu_n)(1-2\delta-8\delta\mu_n)}{
			1+4\delta-8\delta^2-32\delta^2 \mu_n}
		\geq \frac{1-3\delta}{1+4\delta} = 1-\frac{7\delta}{1+4\delta}
\end{eqnarray*}
for sufficiently small $\mu_n$. Therefore
\begin{eqnarray*}
	r = \nu^n_{\left[-\frac{M}{2},-l-1\right]}
	=\nu^n_{\left[l+1,\frac{M}{2}\right]}
	\leq \frac{7\delta}{1+4\delta} \leq 7\delta.
\end{eqnarray*} 
The above inequality and the supposition 
$\nu^n_{\left[-L+1,-\frac{M}{2}-1\right]}
        = \nu^n_{\left[\frac{M}{2}+1,L-1\right]} < \delta$
imply that
\begin{eqnarray}
	\nu^n_{[-L+1,-l-1]} = \nu^n_{[l+1,L-1]} < 8\delta. \label{ineq:nu_sides}
\end{eqnarray}

Let $p=\nu^n_l$. We bound $V_l(\nu^n)$. Since site $-L$ is the only site
in the support of $\nu^n$ that does not compete with site $l$, we have
\begin{eqnarray}
        1 < V_l(\nu^n) = \frac{K_l}{\sum_{z=-L+1}^L \nu^n_z}
        = \frac{K_l}{1 - \nu^n_{-L}}
        \leq \frac{1}{1 - 4 \delta \mu_n} \label{ineq:V_l}
\end{eqnarray}
by~(\ref{ineq:nun_L}). We will
establish the following lemma after the proof of the present theorem:
\begin{LEM}
	Let $\nu^n = \nu^{\mu_n}$ be as in Theorem~\ref{thm:append},
	and suppose $\nu^n_{[-L+1,-l-1]}=\nu^n_{[l+1,L-1]}<8\delta$,
	then
\begin{enumerate}
\item $\nu^n_{-l}$ and $\nu^n_l$ is bounded away from 0 as
        $n \rightarrow \infty$;
\item $\nu^n_{l-1} \leq \nu^n_{l}+\delta$ for sufficiently
        large $n$.
\end{enumerate}
\label{lem:append}
\end{LEM}
Since $\nu^n_{[-L+1,-l-1]} < 8\delta$ by~(\ref{ineq:nu_sides}), we have
from~(\ref{ineq:nun_L}),
\begin{eqnarray}
	V_{l+1}(\nu^n) \leq \frac{1}{1-8\delta-4\delta\mu_n}
	\label{ineq:V_l+1}
\end{eqnarray}
because sites $-L$ and $-L+1$ are the only sites in $[-L,L]$ that
do not compete with site $l+1$.
We estimate $\sum_y A(y,l) \nu^n_y V_y(\nu^n) - \nu^n_l$:
\begin{eqnarray*}
        \lefteqn{ \sum_y A(y,l) \nu^n_y V_y(\nu^n) - \nu^n_l } \nnb \\
	& = & \mu_n \nu^n_{l-1} V_{l-1}(\nu^n) + (1-2\mu_n) \nu^n_l V_l(\nu^n)
		+ \mu_n \nu^n_{l+1} V_{l+1}(\nu^n) - \nu^n_l \\
        & \leq & \mu_n(p+\delta) + \frac{(1-2\mu_n)p}{1-4\delta\mu_n}
                + \frac{\mu_n \delta}{1-8\delta-4\delta\mu_n} - p,
\end{eqnarray*}
using $\nu^n_{l-1} \leq p+\delta$, $V_{l-1}(\nu^n) = 1$,~(\ref{ineq:V_l}),
$\nu^n_{l+1} \leq \delta$, and~(\ref{ineq:V_l+1}).
Simplifying the right hand side of the above, we get
\begin{eqnarray*}
	\sum_y A(y,l) \nu^n_y V_y(\nu^n) - \nu^n_l
        & \leq & \frac{-4\mu_n^2 \delta p + p - \mu_n p}{1-4\delta \mu_n}
                + \mu_n \delta + \frac{\mu_n \delta}{
			1-8\delta-4\delta\mu_n} - p \\
        & = & \frac{- 4\mu_n^2 \delta p - \mu_n p + 4\mu_n \delta p}{
                1-4\delta \mu_n} + \mu_n \delta
                + \frac{\mu_n \delta}{1-8\delta-4\delta\mu_n} \\
        & \leq & -\frac{\mu_n p (1-4\delta)}{1-4\delta \mu_n}
                + \mu_n \delta + \frac{\mu_n \delta}{1-8\delta-4\delta\mu_n} \\
        & \leq & -\frac{\mu_n p}{2}
\end{eqnarray*}
for sufficiently small $\delta$ and $\mu_n$. This estimate means that
after resampling and mutation (and before normalization), $\nu^n_l$
decreases by at least $\mu_n p/2$. On the other hand, $\nu^n_{l+1}$ can
only increase: let $q=\nu^n_{l+1}$, then
$q \leq \delta$ and
\begin{eqnarray*}           
        \lefteqn{ \sum_y A(y,l+1) \nu^n_y V_y(\nu^n) - \nu^n_{l+1} } \nnb \\
	& = & \mu_n \nu^n_{l} V_{l}(\nu^n)
		+ (1-2\mu_n) \nu^n_{l+1} V_{l+1}(\nu^n)
                + \mu_n \nu^n_{l+2} V_{l+2}(\nu^n) - \nu^n_{l+1} \\
        & \geq & \mu_n p  + (1-2\mu_n) q - q,        
\end{eqnarray*}             
since $V_{l}(\nu^n) \geq 1$ and $V_{l+1}(\nu^n) \geq 1$.
Therefore
\begin{eqnarray*}
	\sum_y A(y,l+1) \nu^n_y V_y(\nu^n) - \nu^n_{l+1} \geq \mu_n (p - 2q)
	\geq \mu_n(p-2\delta) > 0
\end{eqnarray*}
if $\delta$ is small enough.
After normalization, i.e. dividing by $\overline V(\nu^n)$,
$\nu^n_l$ and $\nu^n_{l+1}$
cannot possibly return to their original values. This contradicts the
assumption that $\nu^n$ is an stationary distribution for~(\ref{eq:badmodel2})
with mutation kernel $A^{\mu_n}$, and the proof is complete. \qed

\vspace{.3cm}

\noindent {\bf Proof of Lemma~\ref{lem:append}}. \hspace{2mm}
Define
\begin{eqnarray}
        \zeta_x = \nu^n_x V_x(\nu^n). \label{def:zeta}
\end{eqnarray}
Notice that $\zeta_x$ depends on $n$, but notationally we suppress this
dependence. For $x \in [-l+1,l-1]$, $\zeta_{x} = \nu^n_{x}$ since
$V_x(\nu^n) = 1$ from~(\ref{eq:facts:V}).
For $x \in [-l+1,l-1]$, we rewrite condition~(\ref{eq:nuecond}) as follows:
\begin{eqnarray}
        \mu_n \zeta_{x-1} + (1-2\mu_n-\overline V_n) \zeta_x
                + \mu_n \zeta_{x+1} = 0. \label{eq:nuecond2}
\end{eqnarray}
This is a recurrence relation with general solution
$\zeta_x = A\beta_1^x + B \beta_2^x$, where $\beta_1$ and $\beta_2$ are
the two roots of the quadratic polynomial
$\mu_n + (1-2\mu_n-\overline V_n) r + \mu_n r^2$; or
$\zeta_x = (A+ Bx) \beta_1^x$, where $\beta_1$ is the double root of the
polynomial. Elementary calculation shows that for the solution
$\zeta_x = (A+ Bx) \beta_1^x$ to satisfy the symmetry requirement for
$L\geq 1$, either $B=0$ or $\beta_1=0$; $\beta_1=0$ leads to the solution
of $\zeta_x = 0$, and $B = 0$ leads to the conclusion $\beta_1 = 1$ and
$\zeta_x = A$; both these two scenarios will be included in Case 2 below.
For the solution
$\zeta_x = A\beta_1^x + B \beta_2^x$, simple calculation leads to:
\begin{eqnarray*}
        \beta_1,\beta_2 = \frac{1}{2\mu_n} \left(
        2\mu_n + \overline V_n -1 \pm \sqrt{(\overline V_n-1)^2
                + 4\mu_n(\overline V_n-1)}\right).
\end{eqnarray*}
We divide into three cases:
\begin{enumerate}
\item If $\beta_1$ and $\beta_2$ are two real roots,
        then since $\zeta$ is symmetric,
        we must have $\beta_1 = 1/\beta_2$ with $\beta_1>0$, and the solution
        is $\zeta_x = A (\beta_1^x + \beta_1^{-x})$ for $x \in [-l,l]$.
\item If $\beta_1 = \beta_2$, then the solution is $\zeta_x = A$
        for $x \in [-l,l]$.
\item If $\beta_1$ and $\beta_2$ are complex roots, then we write
        $\beta_1=\gamma e^{i\theta}$ and $\beta_2 = \gamma e^{-i \theta}$, and
        the solution is $\zeta_x = A \gamma \cos (x \theta)$
        for $x \in [-l,l]$. Define
\begin{eqnarray}
        \alpha_n = \overline V_n - 1, \label{def:alpha_n}
\end{eqnarray}
        then for $\beta_1$ and $\beta_2$ to be complex,
        $\alpha_n^2 + 4 \mu_n \alpha_n = \alpha_n (\alpha_n+4\mu_n)<0$,
        which means that
\begin{eqnarray}
        & \mbox{either} & \mbox{$\alpha_n < 0$ and $\alpha_n+4\mu_n > 0$}
                \label{ineq:roots1} \\
        & \mbox{or} & \mbox{$\alpha_n > 0$ and $\alpha_n+4\mu_n < 0$}.
                \label{ineq:roots2}
\end{eqnarray}
        Now (\ref{ineq:roots2}) is clearly impossible since $\mu_n \geq 0$,
        and~(\ref{ineq:roots1}) implies that $\alpha_n < 0$. Furthermore,
\begin{eqnarray}
        \tan \theta & = & \sqrt{ \frac{-\alpha_n^2-4\mu_n \alpha_n}{
                (\alpha_n+2\mu_n)^2}}
        = \left( -\frac{\alpha_n^2+4\mu_n\alpha_n+4\mu_n^2}{
                \alpha_n^2+4\mu_n\alpha_n} \right)^{-1/2} \nnb \\
        & = & \left( -1 - \frac{4\mu_n^2}{
                \alpha_n^2+4\mu_n\alpha_n} \right)^{-1/2} \nnb \\
        & = & \left( -1 + \frac{1}{-\frac{\alpha_n}{\mu_n}
                + (-\frac{\alpha_n}{2\mu_n})^2} \right)^{-1/2}. \label{eq:tan}
\end{eqnarray}
        Hence ~(\ref{ineq:Mnu:fine}) and~(\ref{def:alpha_n}) imply that
        $\alpha_n > -8\delta\mu_n$, and since $\alpha_n < 0$,we have
\begin{eqnarray*}
        0 < -\alpha_n / \mu_n < 8\delta.
\end{eqnarray*}
        We conclude from~(\ref{eq:tan}) that for sufficiently
        small $\delta$, $\tan \theta$ is also very small,
\end{enumerate}

Note that $-l$ and $l$ are the boundary sites for~(\ref{eq:nuecond2}),
therefore statements about $\zeta_x$ in the three cases above all hold for
$x\in [-l,l]$, even though~(\ref{eq:nuecond2}) holds for only
$x\in [-l+1,l-1]$.

In case 1, $\zeta$ is a linear combination of two convex functions, therefore
$\zeta$ is convex for $x \in [-l,l]$. In case 2, $\zeta$ is flat
for $x \in [-l,l]$. And in case 3, $\zeta$ is concave
for $x \in [-l,l]$, but since $\theta$ is small for small $\delta$, it
is close to being flat for small $\delta$.
Therefore recalling the definition of $\zeta_x$
in~(\ref{def:zeta}) and using~(\ref{ineq:V_l}), we have
\begin{eqnarray*}
        (1 - 4 \delta \mu_n) \zeta_l \leq \nu^n_l < \zeta_l.
\end{eqnarray*}
In summary, for $x\in [-l,l]$, $\zeta$ is convex, or flat, or nearly flat
for small enough $\delta$; $\nu^n_x = \zeta_x$ for $x\in [-l+1,l-1]$
and $\nu^n_l = \nu^n_{-l}$ is smaller than but very close to
$\zeta_l = \zeta_{-l}$; furthermore, by~(\ref{ineq:nun_L})
and~(\ref{ineq:nu_sides}), we have
\begin{eqnarray*}
        \nu^n_{[-l,l]} = 1 - \nu^n_{-L} - \nu^n_{[-L+1,-l-1]}
                - \nu^n_{[l+1,L-1]} - \nu^n_L
        \geq  1 - 8 \delta \mu_n - 16 \delta,
\end{eqnarray*}
i.e. $\nu^n$ has almost all its mass on $[-l,l]$. We can then
use the symmetry assumption on $\nu^n$ to arrive at the conclusion of the
lemma.
\qed

%%%%%%%%%%%%%%%%%%%%%%%%%%%%%%%%%%%%%%%%%%%%%%%%%%%%%%%%%%%%%%%%%%%%%%%%
%%%
%%% The end matter of your thesis goes here.  The bibliography is
%%% generally mandatory, while any appendices and an index are optional.

%% Bibliography
%% Most any citation style is acceptable, as long as it is consistent.
%% Simple numbering (plain style) is generally sufficient.
% \bibliographystyle{plain}
% \thesisbibliography{...}

\newpage

\newcommand{\newpar}           {\vspace{0.1in}\hspace{-0.47in}}

\newpar
{\bf\huge Bibliography}

\vspace{0.6in}

\newpar
[Burkholder 1973]
 Distribution function inequalities for martingales.
D. L. Burkholder. 
Ann. Probability  1  (1973), 19--42.

\newpar
[B\"urger 2000]
 The mathematical theory of selection, recombination, and mutation.
Reinhard B\"urger. John Wiley, 2000.

\newpar
[Dawson 1993]
Measure-valued Markov Processes,
Lecture Notes in Mathematics, 1541.
D. A. Dawson.  Springer-Verlag, Berlin, 1993.

\newpar
[Dawson and Perkins 1998]
Long-time behavior and coexistence in a mutually catalytic branching model.
D. A. Dawson and E. A. Perkins.
Ann. Probab.  26  (1998),  no. 3, 1088--1138.

\newpar
[Del Moral 1998]
Measure-valued processes and interacting particle systems: Application to nonlinear filtering problems.
P. Del Moral.
Ann. Appl. Probab.  8  (1998),  no. 2, 438--495.

\newpar
[Dieckmann and Doebeli 1999]
On the origin of species by sympatric speciation.
U. Dieckmann and M. Doebeli.
Nature v. 400 no. 6742 (July 22 1999) p. 354-7.

\newpar
[Durrett 1991]
Probability : theory and examples.
R. Durrett. Wadsworth \& Brooks, 1991.

\newpar
[Durrett 1995]
Ten lectures on particle systems,
Lecture Notes in Mathematics 1608
R. Durrett.
Springer, Berlin, 1995.

\newpar
[Durrett and Neuhauser 1994]
Particle systems and reaction-diffusion equations.
R. Durrett and C. Neuhauser.
Ann. Probab.  22  (1994),  no. 1, 289--333.

\newpar
[Ethier and Kurtz 1981]
The infinitely-many-neutral-alleles diffusion model.
S. N. Ethier and T. G. Kurtz.
Adv. in Appl. Probab.  13  (1981), no. 3, 429--452.

\newpar
[Ethier and Kurtz 1994]
Convergence to Fleming-Viot processes in the weak atomic topology.
S. N. Ethier and T. G. Kurtz.
Stochastic Process. Appl.  54  (1994),  no. 1, 1--27.

\newpar
[Liggett 1985]
Interacting particle systems.
Thomas M. Liggett.
Springer-Verlag, 1985.

\newpar
[Liggett 1999]
Stochastic interacting systems : contact, voter, and exclusion processes.
Thomas M. Liggett.
Springer, 1999.

\newpar
[Mayr 1963]
Animal species and evolution.
Enst Mayr.
Harvard University Press, 1963.

\newpar
[Perkins 2002]
Dawson-Watanabe superprocesses and measure-valued diffusions,
Lectures on probability theory and statistics (Saint-Flour, 1999),  125--324,
Lecture Notes in Mathematics 1781
E. A. Perkins.
Springer, Berlin, 2002.

\newpar
[Taylor 1996]
Partial differential equations.
M. E. Taylor.
Springer, 1996.

\newpar
[Wiggins 1988]
Global bifurcations and chaos : analytical methods.
Stephen Wiggins.
Springer-Verlag, 1988.

%% Appendices
%% Index

%%%%%%%%%%%%%%%%%%%%%%%%%%%%%%%%%%%%%%%%%%%%%%%%%%%%%%%%%%%%%%%%%%%%%%%%

\end{document}